\documentclass[12pt,a4paper]{article}
\usepackage{amsfonts}

\usepackage{amssymb}

\usepackage{amsmath}

\newtheorem{theorem}{Theorem}[section]
\newtheorem{definition}[theorem]{Definition}
\newtheorem{lemma}[theorem]{Lemma}

\begin{document}
\title{A note
 on Artin-Markov normal form theorem for braid groups\footnote{Supported by the
NNSF of China (No.10771077) and the NSF of Guangdong Province
(No.06025062).} }
\author{
Yuqun Chen and Qiuhui Mo  \\
\\
{\small \ School of Mathematical Sciences}\\
{\small \ South China Normal University}\\
{\small \ Guangzhou 510631}\\
{\small \ P. R. China}\\
{\small \ yqchen@scnu.edu.cn} \\
{\small \ scnuhuashimomo@126.com}}

\date{}

\maketitle \noindent\textbf{Abstract:} In a recent paper by L. A.
Bokut, V. V. Chaynikov and K. P. Shum  in 2007, Braid group $B_n$ is
represented by Artin-Burau's relations. For such a representation,
it is told that all other compositions can be checked in the same
way. In this note, we support this claim and check all compositions.


\noindent {\bf AMS} Mathematics Subject Classification(2000): 20E06,
20F05, 16S15, 13P10

\section{Introduction}

In two classical papers by A. A. Markov \cite{ma} and E. Artin
\cite{ar}, the normal form theorem for braid groups was established.
Recently, a new proof of this result had been published by L. A.
Bokut, V. V. Chaynikov and K. P. Shum  \cite{bchas}. The last  proof
is based on Gr\"{o}bner-Shirshov bases theory for non-commutative
(and non-associative) algebras invented by A. I. Shirshov
\cite{s62}. In what follows, we will use the presentation of this
theory from L. A. Bokut and  K. P. Shum  \cite{bs2}, and L. A. Bokut
and Y. Q. Chen  \cite{BoCh1}. To be more precise, some compositions
($S$-polynomials) of Artin-Burau's relations were checked in
\cite{bchas}, and it was claimed that all other compositions can be
checked in the same way. In this paper, we support the last claim
and show that all compositions of Artin-Markov's relations are
trivial.

\section{Basic notation and results}

We first cite some concepts and results from the literature
\cite{s62, B72, B76} which are related to the Gr\"{o}bner-Shirshov
bases for associative algebras.

Let $k$ be a field, $k\langle X\rangle$ the free associative algebra
over $k$ generated by $X$ and $ X^{*}$ the free monoid generated by
$X$, where the empty word is the identity which is denoted by 1. For
a word $w\in X^*$, we denote the length of $w$ by $|w|$. Let $X^*$
be a well ordered set. Let $f\in k\langle X\rangle$ with the leading
word $\bar{f}$. Then we call $f$  monic if $\bar{f}$ has coefficient
1.

A well order $>$ on $X^*$ is called monomial if it is compatible
with the multiplication of words, that is, for $u, v\in X^*$, we
have
$$
u > v \Rightarrow w_{1}uw_{2} > w_{1}vw_{2},  \ for \  all \
 w_{1}, \ w_{2}\in  X^*.
$$
A standard example of monomial order on $X^*$ is the deg-lex order
to compare two words first by degree and then lexicographically,
where $X$ is a linearly ordered set.

Let $f$ and $g$ be two monic polynomials in \textmd{k}$\langle
X\rangle$ and $<$ a monomial order on $X^*$. Then, there are two
kinds of compositions:

$ (i)$ If \ $w$ is a word such that $w=\bar{f}b=a\bar{g}$ for some
$a,b\in X^*$ with $|\bar{f}|+|\bar{g}|>|w|$, then the polynomial
 $ (f,g)_w=fb-ag$ is called the intersection composition of $f$ and
$g$ with respect to $w$.

$ (ii)$ If  $w=\bar{f}=a\bar{g}b$ for some $a,b\in X^*$, then the
polynomial $ (f,g)_w=f - agb$ is called the inclusion composition of
$f$ and $g$ with respect to $w$.

Let $S\subset k\langle X\rangle$ such that every $s\in S$ is monic.
Then the composition $ (f,g)_w$ is called trivial modulo $ (S,w)$ if
$ (f,g)_w=\sum\alpha_i a_i s_i b_i$, where each $\alpha_i\in k$,
$a_i,b_i\in X^{*}, \ s_i\in S$ and $a_i \overline{s_i} b_i<w$. If
this is the case, then we write
$$
 (f,g)_w\equiv0\quad mod (S,w).
$$
In general, for $p,q\in k\langle X\rangle$, we write $ p\equiv
q\quad mod (S,w) $ which means that $p-q=\sum\alpha_i a_i s_i b_i $,
where each $\alpha_i\in k,a_i,b_i\in X^{*}, \ s_i\in S$ and $a_i
\overline{s_i} b_i<w$.

A set $S\subset k\langle X\rangle$ is called a Gr\"{o}bner-Shirshov
basis with respect to the monomial order $<$ if any composition of
polynomials in $S$ is trivial modulo $S$.

The following lemma was first proved by Shirshov \cite{s62} for free
Lie algebras  (with deg-lex order) (see also Bokut \cite{B72}).
Bokut \cite{B76} specialized the approach of Shirshov to associative
algebras  (see also Bergman \cite{Be}). For the case of commutative
polynomials, this lemma is known as the Buchberger's Theorem
\cite{Bu70}.

\begin{lemma}\label{CDL}
 (Composition-Diamond Lemma) \ Let $k$ be a field, $A=k \langle
X|S\rangle=k\langle X\rangle/Id (S)$ and $<$ a monomial order on
$X^*$, where $Id (S)$ is the ideal of $k \langle X\rangle$ generated
by $S$. Then  $S $ is a Gr\"{o}bner-Shirshov basis if and only if
$Irr (S) = \{ u \in X^* |  u \neq a\bar{s}b ,s\in S,a ,b \in X^*\}$
is a linear basis of the algebra $A=k\langle X | S \rangle.$
\end{lemma}

Let $A=sgp\langle X|S\rangle$ be a semigroup presentation. By abuse
of notation, $S$ is also a subset of $k\langle X \rangle$. Suppose
that $S$ is a Gr\"{o}bner-Shirshov basis of $k\langle X|S\rangle$.
Then the $S$-irreducible set $Irr(S)=\{u\in X^*|u\neq
a\overline{f}b,\ a ,b \in X^*,\ f\in S\}$ is a linear basis of
$k\langle X|S\rangle$ which is also a normal form of $A$.

\begin{definition}(\cite{bchas})
Let $\Sigma=\{\sigma_{1},\cdots,\sigma_{n-1}\}$ be a finite
alphabet. Then, the following group presentation defines the
n-strand braid group:
$$
B_{n}=gp\langle \Sigma \ | \ \sigma_{t+1}\sigma_{t}\sigma_{t+1}=
\sigma_{t}\sigma_{t+1}\sigma_{t} \ \ (1\leq t\leq n-2), \ \
\sigma_{i}\sigma_{j}=\sigma_{j}\sigma_{i} \ \ (2\leq 1+j<i\leq
n-1)\rangle.
$$
\end{definition}

We now introduce the concept of {\it inverse tower order} of words.
\begin{definition}\label{d2.3}(\cite{bchas})
Let $X=Y\dot{\cup}Z$, words $Y^*$ and the letters in Z  be well
ordered. Suppose that the order on $Y^*$ is monomial. Then, any word
in $X$ has the form $u=u_0z_1u_1\cdots z_{k}u_{k}$, where $k\geq 0,\
u_i\in Y^*,\ z_{i}\in Z$. Define the inverse weight of the word
$u\in
 X^*$ by
$$
inwt(u)=( k, u_{k}, z_{k}, \cdots,u_1, z_{1}, u_{0} ).
$$
Now we order the inverse weights lexicographically as follows
$$
u>v\Leftrightarrow inwt(u)>inwt(v).
$$
Then we call the above order the inverse tower order. Clearly, this
order is  a monomial order on $X^*$.
\end{definition}

In case $Y=T\dot{\cup} U$ and $Y^*$ is endowed with the inverse
tower order, we call the order of words in $X$ the inverse tower
order of words relative to the presentation
$$
X=(T\dot{\cup} U)\dot{\cup} Z.
$$
In general, we can define the inverse tower order of $X$-words
relative to the presentation
$$
X=(\cdots
(X^{(n)}\dot{\cup}X^{(n-1)})\dot{\cup}\cdots)\dot{\cup}X^{(0)},
$$
where $X^{(n)}$-words are endowed by a monomial order.

In the braid group $B_{n}$, we now introduce a new set of generators
which are called  the Artin-Burau generators. We set
\begin{eqnarray*}
&&s_{i,i+1}=\sigma_{i}^{2},\ \ \
s_{i,j+1}=\sigma_{j}\cdots\sigma_{i+1}\sigma_{i}^{2}\sigma_{i+1}^{-1}\cdots\sigma_{j}^{-1},
\ \ \ \  1\leq i<j\leq n-1;\\
&&\sigma_{i,j+1}=\sigma_{i}^{-1}\cdots\sigma_{j}^{-1},\ \ \ 1\leq
i\leq j\leq n-1;\ \ \ \sigma_{ii}=1,\ \ \ \{a,b\}=b^{-1}ab.
\end{eqnarray*}
Form the set
$$
S_{j}=\{s_{i,j},s_{i,j}^{-1}, \ 1\leq i,j<n\} \ \mbox{ and } \
\Sigma^{-1}=\{\sigma_{1}^{-1},\cdots\sigma_{n-1}^{-1}\}.
$$
Then the set
$$
S=S_{n}\cup S_{n-1}\cup \cdots \cup S_{2}\cup\Sigma^{-1}
$$
generates $B_{n}$ as a semigroup.

Now we order the set S in the following way:
$$
S_{n}<S_{n-1}<\cdots < S_{2}<\Sigma^{-1},
$$
and
$$
s_{1,j}^{-1}< s_{1,j}< s_{2,j}^{-1}<\cdots<s_{j-1,j} \ , \ \ \
\sigma_{1}^{-1}<\sigma_{2}^{-1}<\cdots\sigma_{n-1}^{-1}.
$$
With the above notation, we now order the $S$-words by using the
inverse tower order, according to the fixed presentation of $S$ as
the union of $S_{j}$ and $\Sigma^{-1}$. We order the $S_{n}$-words
by the $deg-inlex$ order, i.e., we first compare the words by length
and then by inverse lexicographical order, starting from their last
letters.

\begin{lemma}(\cite{bchas})
The following Artin-Markov relations hold in the braid group
$B_{n}$. For $\delta=\pm1$,
\begin{eqnarray}
&&\sigma_{k}^{-1}s_{i,j}^{\delta}=s_{i,j}^{\delta}\sigma_{k}^{-1},
k\neq i-1,i,j-1,j\\
&&\sigma_{i}^{-1}s_{i,i+1}^{\delta}=s_{i,i+1}^{\delta}\sigma_{1}^{-1}\\
&&\sigma_{i-1}^{-1}s_{i,j}^{\delta}=s_{i-1,j}^{\delta}\sigma_{i-1}^{-1}\\
&&\sigma_{i}^{-1}s_{i,j}^{\delta}=\{s_{i+1,j}^{\delta},s_{i,i+1}\}\sigma_{i}^{-1}\\
&&\sigma_{j-1}^{-1}s_{i,j}^{\delta}=s_{i,j-1}^{\delta}\sigma_{j-1}^{-1}\\
&&\sigma_{j}^{-1}s_{i,j}^{\delta}=\{s_{i,j+1}^{\delta},s_{j,j+1}\}\sigma_{j}^{-1}
\end{eqnarray}
for $\ i<j<k<l,\ \varepsilon = \pm1 $,
\begin{eqnarray}
&&s_{j,k}^{-1}s_{k,l}^{\varepsilon}=\{s_{k,l}^{\varepsilon},s_{j,l}^{-1}\}s_{j,k}^{-1}\\
&&s_{j,k}s_{k,l}^{\varepsilon}=\{s_{k,l}^{\varepsilon},s_{j,l}s_{k,l}\}s_{j,k}\\
&&s_{j,k}^{-1}s_{j,l}^{\varepsilon}=\{s_{j,l}^{\varepsilon},s_{k,l}^{-1}s_{j,l}^{-1}\}s_{j,k}^{-1}\\
&&s_{j,k}s_{j,l}^{\varepsilon}=\{s_{j,l}^{\varepsilon},s_{k,l}\}s_{j,k}\\
&&s_{i,k}^{-1}s_{j,l}^{\varepsilon}=\{s_{j,l}^{\varepsilon},s_{k,l}s_{i,l}s_{k,l}^{-1}s_{i,l}^{-1}\}s_{i,k}^{-1}\\
&&s_{i,k}s_{j,l}^{\varepsilon}=\{s_{j,l}^{\varepsilon},s_{i,l}^{-1}s_{k,l}^{-1}s_{i,l}s_{k,l}\}s_{i,k}
\end{eqnarray}
for $j<i<k<l$ or $ i<k<j<l, \ and  \ \varepsilon, \ \delta=\pm 1 $,
\begin{eqnarray}
&&s_{i,k}^{\delta}s_{j,l}^{\varepsilon}=s_{j,l}^{\varepsilon}s_{i,k}^{\delta}
\end{eqnarray}
and
\begin{eqnarray}
&&\sigma_{j}^{-1}\sigma_{k}^{-1}=\sigma_{k}^{-1}\sigma_{j}^{-1}, \  \ j<k-1\\
&&\sigma_{j,j+1}\sigma_{k,j+1}=\sigma_{k,j+1}\sigma_{j-1,j} \ ,\ \ j<k\\
&&\sigma_{i}^{-2}=s_{i,i+l}^{-1}\\
&&s_{i,j}^{\pm1}s_{i,j}^{\mp1}=1
\end{eqnarray}
\end{lemma}

 \ \

 The following theorem is from L. A.
Bokut, V. V. Chaynikov and K. P. Shum  \cite{bchas} in which the
authors checked some compositions in (1)-(17) and claimed that all
compositions can be checked to be trivial in the same way. Here, we
support the claim and check all compositions.

\begin{theorem}(\cite{bchas})
The Artin-Markov relations (1)-(17) form a minimal
Gr\"{o}bner-Shirshov basis of the braid group $B_{n}$ in terms of
the Artin-Burau generators relative to the inverse tower order of
words.
\end{theorem}

Proof: We check all compositions step by step.

 Denote by $(i)\wedge (j)$ the composition of the type
$(i)$ and type $(j)$, and $s_{i,j}=(i,j)=(ij)$.

$(1)\wedge(7)$ \\
Letn $f=\sigma_{q}^{-1}(jk)^{-1}-(jk)^{-1}\sigma_{q}^{-1},\ \ q\neq
j,j-1,k,k-1,\
g=(jk)^{-1}(kl)^{\varepsilon}-\{{(kl)}^{\varepsilon},{(jl)}^{-1}\}{(jk)}^{-1},\
\ j<k<l$. Then $w=\sigma_{q}^{-1}(jk)^{-1}(kl)^{\varepsilon}$ and
$$
(f,g)_{w}=\sigma_{q}^{-1}\{{(kl)}^{\varepsilon},{(jl)}^{-1}\}{(jk)}^{-1}-
(jk)^{-1}\sigma_{q}^{-1}{(kl)}^{\varepsilon}.
$$
There are three cases to consider.
\begin{enumerate}
\item[1)]\ $q\neq l,l-1$. This case is trivial.
\item[2)]\ $q=l$.
\begin{eqnarray*}
&&\sigma_{q}^{-1}\{{(kl)}^{\varepsilon},{(jl)}^{-1}\}{(jk)}^{-1}\\
&=&\sigma_{l}^{-1}\{{(kl)}^{\varepsilon},{(jl)}^{-1}\}{(jk)}^{-1}\\
&\equiv&\{\{{(k,l+1)}^{\varepsilon},(l,l+1)\},\{{(j,l+1)}^{-1},(l,l+1)\}\}\sigma_{l}^{-1}{(jk)}^{-1}\\
&\equiv&\{{(k,l+1)}^{\varepsilon},{(j,l+1)}^{-1}(l,l+1)\}{(jk)}^{-1}\sigma_{l}^{-1}\
\ \ \ \ \ and \\
&&\\
&&(jk)^{-1}\sigma_{q}^{-1}(kl)^{\varepsilon}=(jk)^{-1}\sigma_{l}^{-1}(kl)^{\varepsilon}\\
&\equiv&(jk)^{-1}\{(k,l+1)^{\varepsilon},(l,l+1)\}\sigma_{l}^{-1}\\
&\equiv&\{\{(k,l+1)^{\varepsilon},(j,l+1)^{-1}\},(l,l+1)\}(jk)^{-1}\sigma_{l}^{-1}\\
&\equiv&\{(k,l+1)^{\varepsilon},(j,l+1)^{-1}(l,l+1)\}(jk)^{-1}\sigma_{l}^{-1}.
\end{eqnarray*}
Thus, $(f,g)_{w}\equiv0$.
\item[3)]\ $q=l-1.$
\begin{eqnarray*}
(f,g)_{w}&=&\sigma_{l-1}^{-1}\{(kl)^{\varepsilon},(jl)^{-1}\}(jk)^{-1}-
(jk)^{-1}\sigma_{l-1}^{-1}(kl)^{\varepsilon}\\
&\equiv&\{(k,l-l)^{\varepsilon},
(j,l-1)^{-1}\}(jk)^{-1}\sigma_{l-1}^{-1}- (jk)^{-1}
(k,l-l)^{\varepsilon}\sigma_{l-1}^{-1}\\
&\equiv&\{(k,l-l)^{\varepsilon},
(j,l-1)^{-1}\}(jk)^{-1}\sigma_{l-1}^{-1}-\{(k,l-l)^{\varepsilon},(j,l-1)^{-1}\}(jk)^{-1}\sigma_{l-1}^{-1}\\
&\equiv&0.
\end{eqnarray*}
\end{enumerate}

$(1)\wedge(8)$\\
Letn $f=\sigma_{q}^{-1}(jk)-(jk)\sigma_{q}^{-1},q\neq j,j-1,k,k-1,\
g=(jk)(kl)^{\varepsilon}-\{(kl)^{\varepsilon},(jl)(kl)\}(jk),j<k<l$.
Then $w=\sigma_{q}^{-1}(jk)(kl)^{\varepsilon}$ and
$$
(f,g)_w=\sigma_{q}^{-1}\{(kl)^{\varepsilon},(jl)(kl)\}(jk)-(jk)\sigma_{q}^{-1}(kl)^{\varepsilon}.
$$
There are three cases to consider.
\begin{enumerate}
\item[1)]\ The case $q\neq l,l-1$ is trivial.
\item[2)]\ $q=l$.
\begin{eqnarray*}
&&\sigma_{l}^{-1}\{(kl)^{\varepsilon},(jl)(kl)\}(jk)\\
&\equiv&\{(k,l+1)^{\varepsilon},(j,l+1)(k,l+1)(l,l+1)\}(jk)\sigma_{l}^{-1} \ \ \ \ \  and \\
&&\\
&&(jk)\sigma_{l}^{-1}(kl)^{\varepsilon}\\
&\equiv&(jk)\{(k,l+1)^{\varepsilon},(l,l+1)\}\sigma_{l}^{-1}\\
&\equiv&\{(k,l+l)^{\varepsilon},(j,l+1)(k,l+1)(l,l+1)\}(jk)\sigma_{l}^{-1}.\\
\end{eqnarray*}
\item[3)]\ $q=l-1$.
\begin{eqnarray*}
&&\sigma_{l-1}^{-1}\{(kl)^{\varepsilon},(jk)(kl)\}(jk)\\
&\equiv&\{(k,l-l)^{\varepsilon},(j,l-1)(k,l-1)\}(jk)\sigma_{l-1}^{-1}\
\ \ \ \ \ and\\
&&\\
&&(jk)\sigma_{l-1}^{-1}(kl)^{\varepsilon}\\
&\equiv&(jk)(k,l-1)^{\varepsilon}\sigma_{l-1}^{-1}\\
&\equiv&\{(k,l-l)^{\varepsilon},(j,l-1)(k,l-1)\}(jk)\sigma_{l-1}^{-1}.\\
\end{eqnarray*}
\end{enumerate}

$(1)\wedge(9)$\\
Let $f=\sigma_{q}^{-1}(jk)^{-1}-(jk)^{-1}\sigma_{q}^{-1}, q\neq
j,j-1,k,k-1,\
g=(jk)^{-1}(jl)^{\varepsilon}-\{(jl)^{\varepsilon},(kl)^{-1}(jl)^{-1}\}(jk)^{-1},j<k<l$.
Then $w=\sigma_{q}^{-1}(jk)^{-1}(jl)^{\varepsilon}$ and
$$
(f,g)_w\equiv\sigma_{q}^{-1}\{(jl)^{\varepsilon},(kl)^{-1}(jl)^{-1}\}(jk)^{-1}-(jk)^{-1}\sigma_{q}^{-1}(jl)^{\varepsilon}.
$$
There are three cases to consider.
\begin{enumerate}
\item[1)]\ The case $q\neq l,l-1$ is trivial.
\item[2)]\ $q=l$.
\begin{eqnarray*}
&&\sigma_{l}^{-1}\{(jl)^{\varepsilon},(kl)^{-1}(jl)^{-1}\}(jk)^{-1}\\
&\equiv&\{(j,l+1)^{\varepsilon},(k,l+1)^{-1}(j,l+l)^{-1}(l,l+1)\}(jk)^{-1}\sigma_{l}^{-1} \ \ \ \ and \\
&&\\
&&(jk)^{-1}\sigma_{l}^{-1}(j,l)^{\varepsilon}\\
&\equiv&(jk)^{-1}\{(j,l+1)^{\varepsilon},(l,l+1)\}\sigma_{l}^{-1}\\
&\equiv&\{(j,l+1)^{\varepsilon},(k,l+l)^{-1}(j,l+1)^{-1}(l,l+1)\}(jk)^{-1}\sigma_{l}^{-1}.\\
\end{eqnarray*}
\item[3)]\ $q=l-1$.
\begin{eqnarray*}
&&\sigma_{l-1}^{-1}\{(jl)^{\varepsilon},(kl)^{-1}(jl)^{-1}\}(jk)^{-1}\\
&\equiv&\{(j,l-1)^{\varepsilon},(k,l-1)^{-1}(j,l-1)^{-1}\}(jk)^{-1}\sigma_{l-1}^{-1} \ \ \ \ \ \ and\\
&&\\
&&(jk)^{-1}\sigma_{l-1}^{-1}(jl)^{\varepsilon}\\
&\equiv&(jk)^{-1}(j,l-1)^{\varepsilon}\sigma_{l-1}^{-1}\\
&\equiv&\{(j,l-1)^{\varepsilon},(k,l-1)^{-1}(j,l-1)^{-1}\}(jk)^{-1}\sigma_{l-1}^{-1}.\\
\end{eqnarray*}
\end{enumerate}

$(1)\wedge(10)$\\
Let $f=\sigma_{q}^{-1}(jk)-(jk)\sigma_{q}^{-1},q \neq j,j-1,k,k-1,\
g=(jk)(jl)^{\varepsilon}-\{(jl)^{\varepsilon},(k,l)\}(jk),j<k<l$.
Then $w=\sigma_{q}^{-1}(jk)(jl)^{\varepsilon}$ and
$$
(f,g)_w=\sigma_{q}^{-1}\{(jl)^{\varepsilon},(kl)\}(jk)-(jk)\sigma_{q}^{-1}(jl)^{\varepsilon}.
$$
There are three cases to consider.
\begin{enumerate}
\item[1)]\ The case $q\neq l,l-1$ is trivial.
\item[2)]\ $q=l$.
$$
(f,g)_w\equiv\{(j,l+1)^{\varepsilon},(k,l+1)(l,l+1)\}(jk)\sigma_{q}^{-1}-
\{(j,l+1)^{\varepsilon},(k,l+1)(l,l+1)\}(jk)\sigma_{q}^{-1}\equiv0.
$$
\item[3)] \ $q=l-1$.
$$
(f,g)_w\equiv\{(j,l-1)^{\varepsilon},(k,l-1)\}(jk)\sigma_{l-1}^{-1}-
\{(j,l-1)^{\varepsilon},(k,l-1)\}(jk)\sigma_{l-1}^{-1}\equiv0.
$$
\end{enumerate}

$(1)\wedge(11)$ \\
Let $f=\sigma_{q}^{-1}(ik)^{-1}-(ik)^{-1}\sigma_{q}^{-1},q\neq
i,i-1,k,k-1,\
g=(ik)^{-1}(jl)^{\varepsilon}-\{(jl)^{\varepsilon},(kl)(il)(kl)^{-1}(il)^{-1}\},i<j<k<l$.
Then $w=\sigma_{q}^{-1}(ik)^{-1}(jl)^{\varepsilon}$ and
$$
(f,g)_w=\sigma_{q}^{-1}\{(jl)^{\varepsilon},(kl)(il)(kl)^{-1}(il)^{-1}\}(ik)^{-1}-
(ik)^{-1}\sigma_{q}^{-1}(jl)^{\varepsilon}.
$$
There are five cases to consider.
\begin{enumerate}
\item[1)]\ The case $q\neq l,l-1,j,j-1$ is trivial.
\item[2)] \ $q=j$.
\begin{eqnarray*}
&&\sigma_{j}^{-1}\{(jl)^{\varepsilon},(kl)(il)(kl)^{-1}(il)^{-1}\}(ik)^{-1}\\
&\equiv&\{(j+1,l)^{\varepsilon},(j,j+1)(kl)(il)(kl)^{-1}(il)^{-1}\}(ik)^{-1}\sigma_{j}^{-1}\\
&\equiv&\{(j+1,l)^{\varepsilon},(kl)(il)(kl)^{-1}(il)^{-1}(j,j+1)\}(ik)^{-1}\sigma_{j}^{-1} \ \ \ \ \ and\\
&&\\
&&(ik)^{-1}\sigma_{j}^{-1}(jl)^{\varepsilon}\\
&\equiv&(ik)^{-1}\{(j+1,l)^{\varepsilon},(j,j+1)\}\sigma_{j}^{-1}\\
&\equiv&\{(j+1,l)^{\varepsilon},(kl)(il)(kl)^{-1}(il)^{-1}(j,j+1)\}(ik)^{-1}\sigma_{j}^{-1}.\\
\end{eqnarray*}
\item[3)]\ $q=j-1$.
\begin{eqnarray*}
(f,g)_w&\equiv&\{(j-1,l)^{\varepsilon},(kl)(il)(kl)^{-1}(il)^{-1}\}(ik)^{-1}\sigma_{j-1}^{-1}\\
&&-\{(j-1,l)^{\varepsilon},(kl)(il)(kl)^{-1}(il)^{-1}\}(ik)^{-1}\sigma_{j-1}^{-1}\\
&\equiv&0.
\end{eqnarray*}
\item[4)]\ $q=l$.
\begin{eqnarray*}
(f,g)_w&\equiv&\{(j,l+1)^{\varepsilon},(k,l+1)(i,l+1)(k,l+1)^{-1}(i,l+1)^{-1}(l,l+1)\}(ik)^{-1}\sigma_{l}^{-1}\\
&&-\{(j,l+1)^{\varepsilon},(k,l+1)(i,l+1)(k,l+1)^{-1}(i,l+1)^{-1}(l,l+1)\}(ik)^{-1}\sigma_{l}^{-1}\\
&\equiv&0.
\end{eqnarray*}
\item[5)]\ $q=l-1$.
\begin{eqnarray*}
(f,g)_w&\equiv&\{(j,l-1)^{\varepsilon},(k,l-1)(i,l-1)(k,l-1)^{-1}(i,l-1)^{-1}\}(ik)^{-1}\sigma_{l-1}^{-1}\\
&&-\{(j,l-1)^{\varepsilon},(k,l-1)(i,l-1)(k,l-1)^{-1}(i,l-1)^{-1}\}(ik)^{-1}\sigma_{l-1}^{-1}\\
&\equiv&0.
\end{eqnarray*}
\end{enumerate}

$(1)\wedge(12)$\\
Let $f=\sigma_{q}^{-1}(ik)-(ik)\sigma_{q}^{-1}, q\neq i,i-1,k,k-1,\
g=(ik)(jl)^{\varepsilon}-\{(jl)^{\varepsilon},(il)^{-1}(kl)^{-1}(il)(kl)\}(ik),i<j<k<l$.
Then $w=\sigma_{q}^{-1}(ik)(jl)^{\varepsilon}$ and
$$
(f,g)_w=\sigma_{q}^{-1}\{(jl)^{\varepsilon},(il)^{-1}(kl)^{-1}(il)(kl)\}(ik)-(ik)\sigma_{q}^{-1}\{(jl)^{\varepsilon}.
$$
There are five cases to consider.
\begin{enumerate}
\item[1)]\ The case $q\neq l,l-1,j,j-1$ is trivial.
\item[2)]\ $q=j$.
\begin{eqnarray*}
&&\sigma_{j}^{-1}\{(jl)^{\varepsilon},(il)^{-1}(kl)^{-1}(il)(kl)\}(ik)\\
&\equiv&\{(j+1,l)^{\varepsilon},(j,j+1)(il)^{-1}(kl)^{-1}(il)(kl)\}(ik)\sigma_{j}^{-1}\\
&\equiv&\{(j+1,l)^{\varepsilon},(il)^{-1}(kl)^{-1}(il)(kl)(j,j+1)\}(ik)\sigma_{j}^{-1} \ \ \ \ \ \ and \\
&&\\
&&(ik)\sigma_{j}^{-1}(jl)^{\varepsilon}\\
&=&(ik)\{(j+1,l)^{\varepsilon},(j,j+1)\}\sigma_{j}^{-1}\\
&\equiv&\{(j+1,l)^{\varepsilon},(il)^{-1}(kl)^{-1}(il)(kl)(j,j+1)\}(ik)\sigma_{j}^{-1}.\\
\end{eqnarray*}
\item[3)]\ $q=j-1$.
\begin{eqnarray*}
(f,g)_w&\equiv&\{(j-1,l)^{\varepsilon},(il)^{-1}(kl)^{-1}(il)(kl)\}(ik)\sigma_{j-1}^{-1}-\\
&&\{(j-1,l)^{\varepsilon},(il)^{-1}(kl)^{-1}(il)(kl)\}(ik)\sigma_{j-1}^{-1}\\
&\equiv&0.
\end{eqnarray*}
\item[4)]\ $q=l$.
\begin{eqnarray*}
(f,g)_w&\equiv&\{(j,l+1)^{\varepsilon},(i,l+1)^{-1}(k,l+1)^{-1}(i,l+1)(k,l+1)(l,l+1)\}(ik)\sigma_{l}^{-1}\\
&&-\{(j,l+1)^{\varepsilon},(i,l+1)^{-1}(k,l+1)^{-1}(i,l+1)(k,l+1)(l,l+1)\}(ik)\sigma_{l}^{-1}\\
&\equiv&0.
\end{eqnarray*}
\item[5)]\ $q=l-1$.
\begin{eqnarray*}
(f,g)_w&\equiv&\{(j,l-1)^{\varepsilon},(i,l-1)^{-1}(k,l-1)^{-1}(i,l-1)(k,l-1)\}(ik)\sigma_{l-1}^{-1}\\
&&-\{(j,l-1)^{\varepsilon},(i,l-1)^{-1}(k,l-1)^{-1}(i,l-1)(k,l-1)\}(ik)\sigma_{l-1}^{-1}\\
&\equiv&0.
\end{eqnarray*}
\end{enumerate}

$(1)\wedge(13)$\\
Let$f=\sigma_{q}^{-1}(ik)\delta-(ik)^\delta\sigma_{q}^{-1},q\neq
i,i-1,k,k-1,\
g=(ik)^\delta(jl)^{\varepsilon}-(jl)^{\varepsilon}(ik)^\delta,j<i<k<l
\ \ or\ \  i<k<j<l$. Then
$w=\sigma_{q}^{-1}(ik)^{\delta}(jl)^{\varepsilon}$ and
$$
(f,g)_w=\sigma_{q}^{-1}(jl)^{\varepsilon}(ik)^{\delta}-(ik)^{\delta}\sigma_{q}^{-1}(jl)^{\varepsilon}.
$$
There are five cases to consider.
\begin{enumerate}
\item[1)]\ The case $q\neq l,l-1,j,j-1$ is trivial.
\item[2)]\ $q=j$.
\begin{eqnarray*}
&&\sigma_{j}^{-1}(jl)^{\varepsilon}(ik)^{\delta}\\
&\equiv&\{(j+1,l)^{\varepsilon},(j,j+1)\}\sigma_{j}^{-1}(ik)^{\delta}\\
&\equiv&\{(j+1,l)^{\varepsilon},(j,j+1)\}(ik)^{\delta}\sigma_{j}^{-1} \ \ \ \ \ \ and\\
&&\\
&&(ik)^\delta\sigma_{j}^{-1}(jl)^{\varepsilon}\\
&\equiv&(ik)^\delta\{(j+1,l)^{\varepsilon},(j,j+1)\}\sigma_{j}^{-1}\\
&\equiv&\{(j+1,l)^{\varepsilon},(j,j+1)\}(ik)^{\delta}\sigma_{j}^{-1}.
\end{eqnarray*}
\item[3)]\ $q=j-1$.
$$
(f,g)_w\equiv(j-1,l)^{\varepsilon}(ik)^\delta\sigma_{j-1}^{-1}-
(j-1,l)^{\varepsilon}(ik)^\delta\sigma_{j-1}^{-1}\equiv0.
$$
\item[4)]\ $q=l$.
$$
(f,g)_w\equiv\{(j,l+1)^{\varepsilon},(l,l+1)\}(ik)^\delta\sigma_{l}^{-1}-
\{(j,l+1)^{\varepsilon},(l,l+1)\}(ik)^\delta\sigma_{l}^{-1}\equiv0.
$$
\item[5)]\ $q=l-1$.
$$
(f,g)_w\equiv(j,l-1)^{\varepsilon}(ik)^\delta\sigma_{l-1}^{-1}-
(j,l-1)^{\varepsilon}(ik)^\delta\sigma_{l-1}^{-1}\equiv0.
$$
\end{enumerate}

$(2)\wedge(7)$\\
Let $f=\sigma_{j}^{-1}(j,j+1)^{-1}-(j,j+1)^{-1}\sigma_{j}^{-1},\
g=(j,j+1)^{-1}(j+1,l)^{\varepsilon}-\{(j+1,l)^{\varepsilon},(jl)^{-1}\}
(j,j+1)^{-1},\ j+1<l$. Then
$w=\sigma_{j}^{-1}(j,j+1)^{-1}(j+1,l)^{\varepsilon}$ and
\begin{eqnarray*}
(f,g)_w&=&\sigma_{j}^{-1}\{(j+1,l)^{\varepsilon},(jl)^{-1}\}(j,j+1)^{-1}-(j,j+1)^{-1}\sigma_{j}^{-1}(j+1,l)^{\varepsilon}\\
&\equiv&\{(jl)^{\varepsilon},\{(j+1,l)^{-1},(j,j+1)\}\}(j,j+1)^{-1}\sigma_{j}^{-1}-\\
&&(j,j+1)^{-1}(jl)^{\varepsilon}\sigma_{j}^{-1}\\
&\equiv&\{(jl)^{\varepsilon},(j+1,l)^{-1}(jl)^{-1}\}(j,j+1)^{-1}\sigma_{j}^{-1}-
\{(jl)^{\varepsilon},(j+1,l)^{-1}(jl)^{-1}\}(j,j+1)^{-1}\sigma_{j}^{-1}\\
&\equiv&0.
\end{eqnarray*}

$(2)\wedge(8)$\\
Let $f=\sigma_{j}^{-1}(j,j+1)-(j,j+1)\sigma_{j}^{-1},\
g=(j,j+1)(j+1,l)^{\varepsilon}-\{(j+1,l)^{\varepsilon},(jl)(j+1,l)\}(j,j+1),j+1<l$.
Then $w=\sigma_{j}^{-1}(j,j+1)(j+1,l)^{\varepsilon}$ and
\begin{eqnarray*}
(f,g)_w&=&\sigma_{j}^{-1}\{(j+1,l)^{\varepsilon},(jl)(j+1,l)\}(j,j+1)-(j,j+1)\sigma_{j}^{-1}(j+1,l)^{\varepsilon}\\
&\equiv&\{(jl)^{\varepsilon},\{(j+1,l),(j,j+1)\}(jl)\}(j,j+1)\sigma_{j}^{-1}-(j,j+1)(jl)^{\varepsilon}\sigma_{j}^{-1}\\
&\equiv&(jl)^{-1}(j,j+1)^{-1}(j+1,l)^{-1}(j,j+1)(jl)^{\varepsilon}(j,j+1)^{-1}(j+1,l)(j,j+1)\\
&&(jl)(j,j+1)\sigma_{j}^{-1}-\{(jl)^{\varepsilon},(j+1,l)\}(j,j+1)\sigma_{j}^{-1}\\
&\equiv&(jl)^{-1}(jl)(j+1,l)^{-1}(jl)^{-1}(j,j+1)^{-1}(j,j+1)(jl)^{\varepsilon}(jl)(j+1,l)(jl)^{-1}(j,j+1)^{-1}\\
&&(j,j+1)(jl)(j,j+1)\sigma_{j}^{-1}-\{(jl)^{\varepsilon},(j+1,l)\}(j,j+1)\sigma_{j}^{-1}\\
&\equiv&0.
\end{eqnarray*}

$(2)\wedge(9)$\\
Let $f=\sigma_{j}^{-1}(j,j+1)^{-1}-(j,j+1)^{-1}\sigma_{j}^{-1},\
g=(j,j+1)^{-1}(jl)^{\varepsilon}-\{(jl)^{\varepsilon},(j+1,l)^{-1}(jl)^{-1}\}(j,j+1)^{-1},j+1<l$.
Then $w=\sigma_{j}^{-1}(j,j+1)^{-1}(jl)^{\varepsilon}$ and
\begin{eqnarray*}
(f,g)_w&=&\sigma_{j}^{-1}\{(jl)^{\varepsilon},(j+1,l)^{-1}(jl)^{-1}\}(j,j+1)^{-1}-
(j,j+1)^{-1}\sigma_{j}^{-1}(jl)^{\varepsilon}\\
&\equiv&\{\{(j+1,l)^{\varepsilon},(j,j+1)\},(jl)^{-1}\{(j+1,l)^{-1},(j,j+1)\}\}(j,j+1)^{-1}\sigma_{j}^{-1}\\
&&-(j,j+1)^{-1}\{(j+1,l)^{\varepsilon},(j,j+1)\}\sigma_{j}^{-1}\\
&\equiv&(j,j+1)^{-1}(j+1,l)(j,j+1)(jl)(j,j+1)^{-1}(j+1,l)^{\varepsilon}(j,j+1)(jl)^{-1}(j,j+1)^{-1}\sigma_{j}^{-1}\\
&& -\{(j+1,l)^{\varepsilon},(jl)^{-1}(j,j+1)\}(j,j+1)^{-1}\sigma_{j}^{-1}\\
&\equiv&(jl)(j+1,l)(jl)^{-1}(j,j+1)^{-1}(j,j+1)(jl)(j,j+1)^{-1}(j+1,l)^{\varepsilon}(j,j+1)(jl)^{-1}(jl)\\
&&(j+1,l)^{-1}(jl)^{-1}(j,j+1)^{-1}\sigma_{j}^{-1}-(j,j+1)^{-1}\{(j+1,l)^{\varepsilon},(jl)^{-1}\}\sigma_{j}^{-1}\\
&\equiv&(jl)(j+1,l)(jl)(j+1,l)^{\varepsilon}(jl)^{-1}(j,j+1)^{-1}(j,j+1)(j+1,l)^{-1}(jl)^{-1}(j,j+1)^{-1}\sigma_{j}^{-1}\\
&&-\{\{(j+1,l)^{\varepsilon},(jl)^{-1}\},\{(jl)^{-1},(j+1,l)^{-1}(jl)^{-1}\}\}(j,j+1)^{-1}\sigma_{j}^{-1}\\
&\equiv&0.
\end{eqnarray*}

$(2)\wedge(10)$\\
Let $f=\sigma_{j}^{-1}(j,j+1)-(j,j+1)\sigma_{j}^{-1},\
g=(j,j+1)(jl)^{\varepsilon}-\{(jl)^{\varepsilon},(j+1,l)\}(j,j+1),j+1<l$.
Then $w=\sigma_{j}^{-1}(j,j+1)(jl)^{\varepsilon}$ and
\begin{eqnarray*}
(f,g)_w&=&\sigma_{j}^{-1}\{(jl)^{\varepsilon},(j+1,l)\}(j,j+1)-(j,j+1)\sigma_{j}^{-1}(jl)^{\varepsilon}\\
&\equiv&\sigma_{j}^{-1}\{(jl)^{\varepsilon},(j+1,l)\}(j,j+1)-(j,j+1)\{(j+1,l)^{\varepsilon},(j,j+1)\}\sigma_{j}^{-1}\\
&\equiv&\{\{(j+1,l)^{\varepsilon},(j,j+1)\},(jl)\}(j,j+1)\sigma_{j}^{-1}-(j+1,l)^{\varepsilon}(j,j+1)\sigma_{j}^{-1}\\
&\equiv&(jl)^{-1}(j,j+1)^{-1}(j+1,l)^{\varepsilon}(j,j+1)(jl)(j,j+1)\sigma_{j}^{-1}-
(j+1,l)^{\varepsilon}(j,j+1)\sigma_{j}^{-1}\\
&\equiv&(jl)^{-1}(jl)(j+1,l)^{\varepsilon}(jl)^{-1}(j,j+1)^{-1}(j,j+1)(jl)(j,j+1)\sigma_{j}^{-1}-\\
&&(j+1,l)^{\varepsilon}(j,j+1)\sigma_{j}^{-1}\\
&\equiv&0.
\end{eqnarray*}

$(2)\wedge(13)$\\
Let $f=\sigma_{i}^{-1}(i,i+1)^\delta-(i,i+1)^\delta\sigma_{i}^{-1},\
g=(i,i+1)^\delta(jl)^{\varepsilon}-(jl)^{\varepsilon}(i,i+1)^\delta,j<i<i+1<l
\ \ \ or\ \ \  i<i+1<j<l$. Then
$w=\sigma_{i}^{-1}(i,i+1)^\delta(jl)^{\varepsilon}$ and
\begin{eqnarray*}
(f,g)_w&=&\sigma_{i}^{-1}(jl)^{\varepsilon}(i,i+1)^\delta-(i,i+1)^\delta\sigma_{i}^{-1}(jl)^{\varepsilon}\\
&\equiv&(jl)^{\varepsilon}\sigma_{i}^{-1}(i,i+1)^\delta-(i,i+1)^\delta(jl)^{\varepsilon}\sigma_{i}^{-1}\\
&\equiv&(jl)^{\varepsilon}(i,i+1)^\delta\sigma_{i}^{-1}-(jl)^{\varepsilon}(i,i+1)^\delta\sigma_{i}^{-1}\\
&\equiv&0.
\end{eqnarray*}

$(3)\wedge(7)$\\
Let $f=\sigma_{j-1}^{-1}(jk)^{-1}-(j-1,k)^{-1}\sigma_{j-1}^{-1},\
g=(jk)^{-1}(kl)^{\varepsilon}-\{(kl)^{\varepsilon},(jl)^{-1}\}(jk)^{-1},j<k<l$.
Then $w=\sigma_{j-1}^{-1}(jk)^{-1}(kl)^{\varepsilon}$ and
\begin{eqnarray*}
(f,g)_w&=&\sigma_{j-1}^{-1}\{(kl)^{\varepsilon},(jl)^{-1}\}(jk)^{-1}-(j-1,k)^{-1}\sigma_{j-1}^{-1}(kl)^{\varepsilon}\\
&\equiv&\{(kl)^{\varepsilon},(j-1,l)^{-1}\}(j-1,k)^{-1}\sigma_{j-1}^{-1}-
(j-1,k)^{-1}(kl)^{\varepsilon}\sigma_{j-1}^{-1}\\
&\equiv&\{(kl)^{\varepsilon},(j-1,l)^{-1}\}(j-1,k)^{-1}\sigma_{j-1}^{-1}-
\{(kl)^{\varepsilon},(j-1,l)^{-1}\}(j-1,k)^{-1}\sigma_{j-1}^{-1}\\
&\equiv&0.
\end{eqnarray*}

$(3)\wedge(8)$\\
Let $f=\sigma_{j-1}^{-1}(jk)-(j-1,k)\sigma_{j-1}^{-1},\
g=(jk)(kl)^{\varepsilon}-\{(kl)^{\varepsilon},(jl)(kl)\}(jk),j<k<l$.
Then $w=\sigma_{j-1}^{-1}(jk)(kl)^{\varepsilon}$ and
\begin{eqnarray*}
(f,g)_w&=&\sigma_{j-1}^{-1}\{(kl)^{\varepsilon},(jl)(kl)\}(jk)-(j-1,k)\sigma_{j-1}^{-1}(kl)^{\varepsilon}\\
&\equiv&\{(kl)^{\varepsilon},(j-1,l)(kl)\}(j-1,k)\sigma_{j-1}^{-1}-(j-1,k)(kl)^{\varepsilon}\sigma_{j-1}^{-1}\\
&\equiv&\{(kl)^{\varepsilon},(j-1,l)(kl)\}(j-1,k)\sigma_{j-1}^{-1}-
\{(kl)^{\varepsilon},(j-1,l)(kl)\}(j-1,k)\sigma_{j-1}^{-1}\\
&\equiv&0.
\end{eqnarray*}

$(3)\wedge(9)$\\
Let $f=\sigma_{j-1}^{-1}(jk)^{-1}(j-1,k)^{-1}\sigma_{j-1}^{-1},\
g=(jk)^{-1}(jl)^{\varepsilon}-\{(jl)^{\varepsilon},(kl)^{-1}(jl)^{-1}\}(jk)^{-1},j<k<l$.
Then $w=\sigma_{j-1}^{-1}(jk)^{-1}(jl)^{\varepsilon}$ and
\begin{eqnarray*}
(f,g)_w&=&\sigma_{j-1}^{-1}\{(jl)^{\varepsilon},(kl)^{-1}(jl)^{-1}\}(jk)^{-1}-
(j-1,k)^{-1}\sigma_{j-1}^{-1}(jl)^{\varepsilon}\\
&\equiv&\{(j-1,l)^{\varepsilon},(kl)^{-1}(j-1,l)^{-1}\}(j-1,k)^{-1}\sigma_{j-1}^{-1}-
(j-1,k)^{-1}(j-1,l)^{\varepsilon}\sigma_{j-1}^{-1}\\
&\equiv&\{(j-1,l)^{\varepsilon},(kl)^{-1}(j-1,l)^{-1}\}(j-1,k)^{-1}\sigma_{j-1}^{-1}-\\
&&\{(j-1,l)^{\varepsilon},(kl)^{-1}(j-1,l)^{-1}\}(j-1,k)^{-1}\sigma_{j-1}^{-1}\\
&\equiv&0.
\end{eqnarray*}

$(3)\wedge(10)$\\
Let $f=\sigma_{j-1}^{-1}(jk)-(j-1,k)\sigma_{j-1}^{-1},\
g=(jk)(jl)^{\varepsilon}-\{(jl)^{\varepsilon},(kl)\}(jk),j<k<l$.
Then $w=\sigma_{j-1}^{-1}(jk)(jl)^{\varepsilon}$ and
\begin{eqnarray*}
(f,g)_w&=&\sigma_{j-1}^{-1}\{(jl)^{\varepsilon},(kl)\}(jk)-(j-1,k)\sigma_{j-1}^{-1}(jl)^{\varepsilon}\\
&\equiv&\{(j-1,l)^{\varepsilon},(kl)\}(j-1,k)\sigma_{j-1}^{-1}-(j-1,k)(j-1,l)^{\varepsilon}\sigma_{j-1}^{-1}\\
&\equiv&0.
\end{eqnarray*}

$(3)\wedge(11)$\\
Let $f=\sigma_{i-1}^{-1}(ik)^{-1}-(i-1,k)^{-1}\sigma_{i-1}^{-1},\
g=(ik)^{-1}(jl)^{\varepsilon}-\{(jl)^{\varepsilon},(kl)(il)(kl)^{-1}(il)^{-1}\}(ik)^{-1},i<j<k<l$.
Then $w=\sigma_{i-1}^{-1}(ik)^{-1}(jl)^{\varepsilon}$ and
\begin{eqnarray*}
(f,g)_w&=&\sigma_{i-1}^{-1}\{(jl)^{\varepsilon},(kl)(il)(kl)^{-1}(il)^{-1}\}(ik)^{-1}-
(i-1,k)^{-1}\sigma_{i-1}^{-1}(jl)^{\varepsilon}\\
&\equiv&\{(jl)^{\varepsilon},(kl)(i-1,l)(kl)^{-1}(i-1,l)^{-1}\}(i-1,k)^{-1}\sigma_{i-1}^{-1}-(
i-1,k)^{-1}(jl)^{\varepsilon}\sigma_{i-1}^{-1}\\
&\equiv&0.
\end{eqnarray*}

$(3)\wedge(12)$\\
Let $f=\sigma_{i-1}^{-1}(ik)-(i-1,k)\sigma_{i-1}^{-1},\
g=(ik)(jl)^{\varepsilon}-\{(jl)^{\varepsilon},(il)^{-1}(kl)^{-1}(il)(kl)\}(ik),i<j<k<l$.
Then $w=\sigma_{i-1}^{-1}(ik)(jl)^{\varepsilon}$ and
\begin{eqnarray*}
(f,g)_w&=&\sigma_{i-1}^{-1}\{(jl)^{\varepsilon},(il)^{-1}(kl)^{-1}(il)(kl)\}(ik)-
(i-1,k)\sigma_{i-1}^{-1}(jl)^{\varepsilon}\\
&\equiv&\{(jl)^{\varepsilon},(i-1,l)^{-1}(kl)^{-1}(i-1,l)(kl)\}(i-1,k)\sigma_{i-1}^{-1}-
(i-1,k)(jl)^{\varepsilon}\sigma_{i-1}^{-1}\\
&\equiv&0.
\end{eqnarray*}

$(3)\wedge(13)$\\
$f=\sigma_{i-1}^{-1}(ik)^\delta-(i-1,k)^\delta\sigma_{i-1}^{-1}, \
g=(ik)^\delta(jl)^{\varepsilon}-(jl)^{\varepsilon}(ik)^\delta$. Then
$w=\sigma_{i-1}^{-1}(ik)^\delta(jl)^{\varepsilon}$ and
$$
(f,g)_w=\sigma_{i-1}^{-1}(jl)^{\varepsilon}
(ik)^\delta-(i-1,k)^\delta\sigma_{i-1}^{-1}(jl)^{\varepsilon}.
$$
There are two cases to consider.
\begin{enumerate}
\item[1)]\ $j<i<k<l$. In this case, there are two subcases to consider.
\begin{enumerate}
\item[a)]\ $i-1=j$.
\begin{eqnarray*}
&&(f,g)_w\\
&=&\sigma_{i-1}^{-1}(i-1,l)^{\varepsilon}(ik)^\delta-(i-1,k)^\delta
\sigma_{i-1}^{-1}(i-1,l)^{\varepsilon}\\
&\equiv&\{(il)^{\varepsilon},(i-1,i)\}(i-1,k)^\delta\sigma_{i-1}^{-1}-
(i-1,k)^\delta\{(il)^{\varepsilon},(i-1,i)\}\sigma_{i-1}^{-1}\\
&\equiv&(i-1,i)^{-1}(il)^{\varepsilon}(i-1,i)(i-1,k)^\delta\sigma_{i-1}^{-1}
- (i-1,k)^\delta(i-1,l)(il)^{\varepsilon}
(i-1,l)^{-1}\sigma_{i-1}^{-1}.
\end{eqnarray*}
If $\delta=1$, then
\begin{eqnarray*}
(f,g)_w&\equiv&(i-1,l)(il)^{\varepsilon}(i-1,l)^{-1}(i-1,k)^\delta
\sigma_{i-1}^{-1}-\{(i-1,l),(kl)\}\{(il)^{\varepsilon},\\
&&(i-1,l)^{-1}(kl)^{-1}(i-1,l)(kl)\}\{(i-1,l)^{-1},(kl)\}(i-1,k)^\delta\sigma_{i-1}^{-1}\\
&\equiv&0.
\end{eqnarray*}
If $\delta=-1$, then
\begin{eqnarray*}
(f,g)_w&\equiv&(i-1,l)(il)^{\varepsilon}(i-1,l)^{-1}(i-1,k)^\delta
\sigma_{i-1}^{-1}-\{(i-1,l),(kl)^{-1}(i-1,l)^{-1}\}\\
&&\{(il)^{\varepsilon},(kl)(i-1,l)(kl)^{-1}(i-1,l)^{-1}\}\{(i-1,l)^{-1}(kl)^{-1}(i-1,l)^{-1}\}\\
&&(i-1,k)^\delta\sigma_{i-1}^{-1}\\
&\equiv&0.
\end{eqnarray*}
\item[b)]\ $j<i-1$.
\begin{eqnarray*}
(f,g)_w&=&\sigma_{i-1}^{-1}(jl)^{\varepsilon}(ik)^\delta-(i-1,k)^\delta\sigma_{i-1}^{-1}(jl)^{\varepsilon}\\
&\equiv&(jl)^{\varepsilon}(i-1,k)^\delta\sigma_{i-1}^{-1}-(i-1,k)^\delta(jl)^{\varepsilon}\sigma_{i-1}^{-1}\\
&\equiv&(jl)^{\varepsilon}(i-1,k)^\delta\sigma_{i-1}^{-1}-(jl)^{\varepsilon}(i-1,k)^\delta\sigma_{i-1}^{-1}\\
&\equiv&0.
\end{eqnarray*}
\end{enumerate}
\item[2)]\ $i<k<j<l$.
\begin{eqnarray*}
(f,g)_w&=&\sigma_{i-1}^{-1}(jl)^{\varepsilon}(ik)^\delta-(i-1,k)^\delta\sigma_{i-1}^{-1}(jl)^{\varepsilon}\\
&\equiv&(jl)^{\varepsilon}(i-1,k)^\delta\sigma_{i-1}^{-1}-(i-1,k)^\delta(jl)^{\varepsilon}\sigma_{i-1}^{-1}\\
&\equiv&(jl)^{\varepsilon}(i-1,k)^\delta\sigma_{i-1}^{-1}-(jl)^{\varepsilon}(i-1,k)^\delta\sigma_{i-1}^{-1}\\
&\equiv&0.
\end{eqnarray*}
\end{enumerate}

$(4)\wedge(7)$\\
Let
$f=\sigma_{j}^{-1}(jk)^{-1}-\{(j+1,k)^{-1}(j,j+1)\}\sigma_{j}^{-1},
g=(jk)^{-1}(kl)^{\varepsilon}-\{(kl)^{\varepsilon},(jl)^{-1}\}(jk)^{-1},j+1<k<l$.\
Then $w=\sigma_{j}^{-1}(jk)^{-1}(kl)^{\varepsilon}$ \ and
\begin{eqnarray*}
(f,g)_w&=&\sigma_{j}^{-1}\{(kl)^{\varepsilon},(jl)^{-1}\}(jk)^{-1}-\{(j+1,k)^{-1},(j,j+1)\}
\sigma_{j}^{-1}(kl)^{\varepsilon}\\
&\equiv&\{(kl)^{\varepsilon},\{(j+1,l)^{-1},(j,j+1)\}\}\{(j+1,k)^{-1},(j,j+1)\}\sigma_{j}^{-1}-\{(j+1,k)^{-1}\\
&&(j,j+1)\}(kl)^{\varepsilon}\sigma_{j}^{-1}\\
&\equiv&(j,j+1)^{-1}(j+1,l)(j,j+1)(kl)^{\varepsilon}(j,j+1)^{-1}(j+1,l)^{-1}(j+1,k)^{-1}\\
&&(j,j+1)\sigma_{j}^{-1}
-(j,j+1)^{-1}(j+1,k)^{-1}(j,j+1)(kl)^{\varepsilon}\sigma_{j}^{-1}\\
&\equiv&(jl)(j+1,l)(jl)^{-1}(kl)^{\varepsilon}(jl)(j+1,l)^{-1}(jl)^{-1}(j,j+1)^{-1}(j+1,k)^{-1}\\
&&(j,j+1)\sigma_{j}^{-1}-(j,j+1)^{-1}(j+1,k)^{-1}(kl)^{\varepsilon}(j,j+1)\sigma_{j}^{-1}\\
&\equiv&(jl)(j+1,l)(jl)^{-1}(kl)^{\varepsilon}(jl)(j+1,l)^{-1}(jl)^{-1}(jk)(j+1,k)^{-1}(jk)^{-1}(j,j+1)^{-1}\\
&&(j,j+1)\sigma_{j}^{-1}-(j,j+1)^{-1}(j+1,l)(kl)^{\varepsilon}(j+1,l)^{-1}(j+1,k)^{-1}(j,j+1)\sigma_{j}^{-1}\\
&\equiv&(jl)(j+1,l)(jl)^{-1}(kl)^{\varepsilon}(jl)(j+1,l)^{-1}(jl)^{-1}(jk)(j+1,k)^{-1}(jk)^{-1}\sigma_{j}^{-1}-(jl)\\
&&(j+1,l)(jl)^{-1}(kl)^{\varepsilon}(jl)(j+1,l)^{-1}(jl)^{-1}(j,j+1)^{-1}(j+1,k)^{-1}(j,j+1)\sigma_{j}^{-1}\\
&\equiv&(jl)(j+1,l)(jl)^{-1}(kl)^{\varepsilon}(jl)(j+1,l)^{-1}(jl)^{-1}(jk)(j+1,k)^{-1}(jk)^{-1}\sigma_{j}^{-1}\\
&&-(jl)(j+1,l)(jl)^{-1}(kl)^{\varepsilon}(jl)(j+1,l)^{-1}(jl)^{-1}(jk)(j+1,k)^{-1}(jk)^{-1}\\
&&(j,j+1)^{-1}(j,j+1)\sigma_{j}^{-1}\\
&\equiv&0.
\end{eqnarray*}

$(4)\wedge(8)$\\
Let $f=\sigma_{j}^{-1}(jk)-\{(j+1,k),(j,j+1)\}\sigma_{j}^{-1},\
g=(jk)(kl)^{\varepsilon}-\{(kl)^{\varepsilon},(jl)(kl)\}(jk),j+1<k<l$.
Then $w=\sigma_{j}^{-1}(jk)(kl)^{\varepsilon}$ and
$$
(f,g)_w=\sigma_{j}^{-1}\{(kl)^{\varepsilon},(jl)(kl)\}(jk)-\{(j+1,k),(j,j+1)\}\sigma_{j}^{-1}(kl)^{\varepsilon}\equiv0
$$
since
\begin{eqnarray*}
&&\sigma_{j}^{-1}\{(kl)^{\varepsilon},(jl)(kl)\}(jk)\\
&\equiv&\{(kl)^{\varepsilon},\{(j+1,l),(j,j+1)\}(kl)\}\{(j+1,k),(j,j+1)\}\sigma_{j}^{-1}\\
&\equiv&(kl)^{-1}(j,j+1)^{-1}(j+1,l)^{-1}(j,j+1)^{-1}(kl)^{\varepsilon}(j,j+1)^{-1}(j+1,l)(j,j+1)(kl)(j,j+1)^{-1}\\
&&(j+1,k)(j,j+1)\sigma_{j}^{-1}\\
&\equiv&(kl)^{-1}(jl)(j+1,l)^{-1}(jl)^{-1}(kl)^{\varepsilon}(jl)(j+1,l)(jl)^{-1}(kl)(jk)(j+1,k)(jk)^{-1}
\sigma_{j}^{-1} \ \ \ \ \ and\\
&&\\
&&\{(j+1,k),(j,j+1)\}\sigma_{j}^{-1}(kl)^{\varepsilon}\\
&\equiv&(j,j+1)^{-1}(j+1,k)(j,j+1)(kl)^{\varepsilon}\sigma_{j}^{-1}\\
&\equiv&(j,j+1)^{-1}(j+1,k)(kl)^{\varepsilon}(j,j+1)\sigma_{j}^{-1}\\
&\equiv&(j,j+1)^{-1}\{(kl)^{\varepsilon},(j+1,l)(kl)\}(j+1,k)(j,j+1)\sigma_{j}^{-1}\\
&\equiv&\{(kl)^{\varepsilon},\{(j+1,l),(jl)^{-1}\}(kl)\}.\{(j+1,k)(jk)^{-1}\}(j,j+1)^{-1}(j,j+1)\sigma_{j}^{-1}\\
&\equiv&(kl)^{-1}(jl)(j+1,l)^{-1}(jl)^{-1}(kl)^{\varepsilon}(jl)(j+1,l)(jl)^{-1}(kl)(jk)(j+1,k)(jk)^{-1}
\sigma_{j}^{-1}.
\end{eqnarray*}

$(4)\wedge(9)$\\
Let
$f=\sigma_{j}^{-1}(jk)^{-1}-\{(j+1,k)^{-1},(j,j+1)\}\sigma_{j}^{-1},\
g=(jk)^{-1}(jl)^{\varepsilon}-\{(jl)^{\varepsilon},(kl)^{-1}(jl)^{-1}\}(jk)^{-1},j+1<k<l$.
Then $w=\sigma_{j}^{-1}(jk)^{-1}(jl)^{\varepsilon}$ and
$$
(f,g)_w=\sigma_{j}^{-1}\{(jl)^{\varepsilon},(kl)^{-1}(jl)^{-1}\}(jk)^{-1}-
\{(j+1,k)^{-1}(j,j+1)\}\sigma_{j}^{-1}(jl)^{\varepsilon}\equiv0
$$
since
\begin{eqnarray*}
&&\sigma_{j}^{-1}\{(jl)^{\varepsilon},(kl)^{-1}(jl)^{-1}\}(jk)^{-1}\\
&\equiv&\{\{(j+1,l)^{\varepsilon},(j,j+1)\},(kl)^{-1}\{(j+1,l)^{-1},(j,j+1)\}\}\{(j+1,k)^{-1},(j,j+1)\sigma_{j}^{-1}\\
&\equiv&(j,j+1)^{-1}(j+1,l)(j,j+1)(kl)(j,j+1)^{-1}(j+1,l)^{\varepsilon}(j,j+1)(kl)^{-1}(j,j+1)^{-1}\\
&&(j+1,l)^{-1}(j+1,k)^{-1}(j,j+1)\sigma_{j}^{-1}\\
&\equiv&(jl)(j+1,l)(jl)^{-1}(kl)(jl)(j+1,l)^{\varepsilon}(jl)^{-1}(kl)^{-1}(jl)(j+1,l)^{-1}(jl)^{-1}(jk)(j+1,k)^{-1}\\
&&(jk)^{-1}\sigma_{j}^{-1}  \ \ \ \ \ \ \ \ and \\
&&\\
&&\{(j+1,k)^{-1},(j,j+1)\}\sigma_{j}^{-1}(jl)^{\varepsilon}\\
&\equiv&(j,j+1)^{-1}(j+1,k)^{-1}(j,j+1)(j,j+1)^{-1}(j+1,l)^{\varepsilon}(j,j+1)\sigma_{j}^{-1}\\
&\equiv&(j,j+1)^{-1}\{(j+1,l)^{\varepsilon},(kl)^{-1}(j+1,l)^{-1}\}(j+1,k)^{-1}(j,j+1)\sigma_{j}^{-1}\\
&\equiv&\{\{(j+1,l)^{\varepsilon},(jl)^{-1}\},(kl)^{-1}\{(j+1,l)^{-1},(jl)^{-1}\}\}\{(j+1,k)^{-1},(jk)^{-1}\}
\sigma_{j}^{-1}\\
&\equiv&(jl)(j+1,l)(jl)^{-1}(kl)(jl)(j+1,l)^{\varepsilon}(jl)^{-1}(kl)^{-1}(jl)(j+1,l)^{-1}(jl)^{-1}(jk)(j+1,k)^{-1}\\
&&(jk)^{-1}\sigma_{j}^{-1}.
\end{eqnarray*}

$(4)\wedge(10)$\\
Let $f=\sigma_{j}^{-1}(jk)-\{(j+1,k),(j,j+1)\}\sigma_{j}^{-1},\
g=(jk)(jl)^{\varepsilon}-\{(jl)^{\varepsilon},(kl)\}(jk),j+1<k<l$.
Then $w=\sigma_{j}^{-1}(jk)(jl)^{\varepsilon}$ and
\begin{eqnarray*}
(f,g)_w&=&\sigma_{j}^{-1}\{(jl)^{\varepsilon},(kl)\}(jk)-\{(j+1,k),(j,j+1)\}\sigma_{j}^{-1}(jl)^{\varepsilon}\\
&\equiv&\{(j+1,l)^{\varepsilon},(j,j+1)(kl)\}\{(j+1,k),(j,j+1)\}\sigma_{j}^{-1}\\
&&-(j,j+1)^{-1}(j+1,k)(j+1,l)^{\varepsilon}(j,j+1)\sigma_{j}^{-1}\\
&\equiv&\{(j+1,l)^{\varepsilon},(j,j+1)(kl)\}\{(j+1,k),(j,j+1)\}\sigma_{j}^{-1}\\
&&-(j,j+1)^{-1}\{(j+1,l)^{\varepsilon},(kl)\}(j+1,k)(j,j+1)\sigma_{j}^{-1}\\
&\equiv&(kl)^{-1}(jl)(j+1,l)^{\varepsilon}(jl)^{-1}(kl)(jk)(j+1,k)(jk)^{-1}\sigma_{j}^{-1}\\
&&-(kl)^{-1}(jl)(j+1,l)^{\varepsilon}(jl)^{-1}(kl)(jk)(j+1,k)(jk)^{-1}\sigma_{j}^{-1}\\
&\equiv&0.
\end{eqnarray*}

$(4)\wedge(11)$\\
Let
$f=\sigma_{i}^{-1}(ik)^{-1}-\{(i+1,k)^{-1},(i,i+1)\}\sigma_{i}^{-1},\
g=(ik)^{-1}(jl)^{\varepsilon}-\{(jl)^{\varepsilon},(kl)(il)(kl)^{-1}(il)^{-1}\}(ik)^{-1},i<j<k<l$.
Then $w=\sigma_{i}^{-1}(ik)^{-1}(jl)^{\varepsilon}$ and
$$
(f,g)_w=\sigma_{i}^{-1}\{(jl)^{\varepsilon},(kl)(il)(kl)^{-1}(il)^{-1}\}(ik)^{-1}-
\{(i+1,k)^{-1},(i,i+1)\}\sigma_{i}^{-1}(jl)^{\varepsilon}.
$$
There are two cases to consider.
\begin{enumerate}
\item[1)]\ $i=j-1$.
\begin{eqnarray*}
&&\sigma_{j-1}^{-1}\{(jl)^{\varepsilon},(kl)(j-1,l)(kl)^{-1}(j-1,l)^{-1}\}(j-1,k)^{-1}\\
&\equiv&\{(j-1,l)^{\varepsilon},(kl)\{(jl),(j-1,j)\}(kl)^{-1}\{(jl)^{-1},(j-1,j)\}\}\{(jk)^{-1},(j-1,j)\}\sigma_{j-1}^{-1}\\
&\equiv&(j-1,j)^{-1}(jl)(j-1,j)(kl)(j-1,j)^{-1}(jl)^{-1}(j-1,j)(kl)^{-1}(j-1,l)^{\varepsilon}(kl)\\
&&(j-1,j)^{-1}(jl)(j-1,j)(kl)^{-1}(j-1,j)^{-1}(jl)^{-1}(jk)^{-1}(j-1,j)\sigma_{j-1}^{-1}\\
&\equiv&\{(j-1,l)^{\varepsilon},(kl)(j-1,l)(jl)(j-1,l)^{-1}(kl)^{-1}(j-1,l)(jl)^{-1}(j-1,l)^{-1}\}\\
&&(j-1,k)(jk)^{-1}(j-1,k)^{-1}\sigma_{j-1}^{-1} \ \ \ \ \ \ \ and \\
&&\\
&&\{(jk)^{-1},(j-1,j)\}\sigma_{j-1}^{-1}(jl)^{\varepsilon}\\
&\equiv&(j-1,j)^{-1}(jk)^{-1}(j-1,j)(j-1,l)^{\varepsilon}\sigma_{j-1}^{-1}\\
&\equiv&(j-1,j)^{-1}(jk)^{-1}\{(j-1,l)^{\varepsilon},(jl)\}(j-1,j)\sigma_{j-1}^{-1}\\
&\equiv&(j-1,j)^{-1}\{(j-1,l)^{\varepsilon},\{(jl),(kl)^{-1}(jl)^{-1}\}\}(jk)^{-1}(j-1,j)\sigma_{j-1}^{-1}\\
&\equiv&\{\{(j-1,l)^{\varepsilon},(jl)^{-1}(j-1,l)^{-1}\},\{\{(jl),(j-1,l)^{-1}\},(kl)^{-1}\{(jl)^{-1},(j-1,l)^{-1}\}\}\}\\
&&\{(jk)^{-1},(j-1,k)^{-1}\}\sigma_{j-1}^{-1}\\
&\equiv&\{(j-1,l)^{\varepsilon},(kl)(j-1,l)(jl)(j-1,l)^{-1}(kl)^{-1}(j-1,l)(jl)^{-1}(j-1,l)^{-1}\}\\
&&(j-1,k)(jk)^{-1}(j-1,k)^{-1}\sigma_{j-1}^{-1}.
\end{eqnarray*}
\item[2)]\ $i<j-1$.
\begin{eqnarray*}
&&\sigma_{i}^{-1}\{(jl)^{\varepsilon},(kl)(il)(kl)^{-1}(il)^{-1}\}(ik)^{-1}\\
&\equiv&\{(jl)^{\varepsilon},(kl)\{(i+1,l),(i,i+1)\}(kl)^{-1}\{(i+1,l)^{-1},(i,i+1)\}\}\\
&&\{(i+1,k)^{-1},(i,i+1)\}\sigma_{i}^{-1}\\
&\equiv&(i,i+1)^{-1}(i+1,l)(i,i+1)(kl)(i,i+1)^{-1}(i+1,l)^{-1}(i,i+1)(kl)^{-1}(jl)^{\varepsilon}(kl)\\
&&(i,i+1)^{-1}(i+1,l)(i,i+1)(kl)^{-1}(i,i+1)^{-1}(i+1,l)^{-1}(i+1,k)^{-1}(i,i+1)\sigma_{i}^{-1}\\
&\equiv&\{(jl)^{\varepsilon},(kl)(il)(i+1,l)(il)^{-1}(kl)^{-1}(il)(i+1,l)^{-1}(il)^{-1}\}\\
&&(ik)(i+1,k)^{-1}(ik)^{-1}\sigma_{i}^{-1}\ \ \ \ \ and \\
&&\\
&&\{(i+1,k)^{-1},(i,i+1)\}\sigma_{i}^{-1}(jl)^{\varepsilon}\\
&\equiv&(i,i+1)^{-1}(i+1,k)^{-1}(i,i+1)(jl)^{\varepsilon}\sigma_{i}^{-1}\\
&\equiv&(i,i+1)^{-1}(i+1,k)^{-1}(jl)^{\varepsilon}(i,i+1)\sigma_{i}^{-1}\\
&\equiv&(i,i+1)^{-1}(i+1,l)(kl)(i+1,l)^{-1}(kl)^{-1}(jl)^{\varepsilon}(kl)(i+1,l)(kl)^{-1}\\
&&(i+1,l)^{-1}(i+1,k)^{-1}(i,i+1)\sigma_{i}^{-1}\\
&\equiv&\{(jl)^{\varepsilon},(kl)(il)(i+1,l)(il)^{-1}(kl)^{-1}(il)(i+1,l)^{-1}(il)^{-1}\}\\
&&(ik)(i+1,k)^{-1}(ik)^{-1}\sigma_{i}^{-1}.
\end{eqnarray*}
\end{enumerate}

$(4)\wedge(12)$\\
Let $f=\sigma_{i}^{-1}(ik)-\{(i+1,k),(i,i+1)\}\sigma_{i}^{-1},\
g=(ik)(jl)^{\varepsilon}-\{(jl)^{\varepsilon},(il)^{-1}(kl)^{-1}(il)(kl)\}(ik),i<j<k<l$.
Then $w=\sigma_{i}^{-1}(ik)(jl)^{\varepsilon}$ and
$$
(f,g)_w=\sigma_{i}^{-1}\{(jl)^{\varepsilon},(il)^{-1}(kl)^{-1}(il)(kl)\}(ik)-
\{(i+1,k),(i,i+1)\}\sigma_{i}^{-1}(jl)^{\varepsilon}.
$$
There are two cases to consider.
\begin{enumerate}
\item[1)]\ $i=j-1$.
\begin{eqnarray*}
&&\sigma_{j-1}^{-1}\{(jl)^{\varepsilon},(j-1,l)^{-1}(kl)^{-1}(j-1,l)(kl)\}(j-1,k)\\
&\equiv&\{(j-1,l)^{\varepsilon},\{(jl)^{-1},(j-1,j)\}(kl)^{-1}\{(jl),(j-1,j)\}(kl)\}\{(jk),(j-1,j)\}\sigma_{j-1}^{-1}\\
&\equiv&(kl)^{-1}(j-1,j)^{-1}(jl)^{-1}(j-1,j)(kl)(j-1,j)^{-1}(jl)(j-1,j)(j-1,l)^{\varepsilon}(j-1,j)^{-1}\\
&&(jl)^{-1}(j-1,j)(kl)^{-1}(j-1,j)^{-1}(jl)(j-1,j)(kl)(j-1,j)^{-1}(jk)(j-1,j)\sigma_{j-1}^{-1}\\
&\equiv&\{(j-1,l)^{\varepsilon},(jl)^{-1}(j-1,l)^{-1}(kl)^{-1}(j-1,l)(jl)(j-1,l)^{-1}(kl)\}\\
&&(j-1,k)(jk)(j-1,k)^{-1}\sigma_{j-1}^{-1}\ \ \ \ \ and \\
&&\\
&&\{(jk),(j-1,j)\}\sigma_{j-1}^{-1}(jl)^{\varepsilon}\\
&\equiv&(j-1,j)^{-1}(jk)(j-1,j)(j-1,l)^{\varepsilon}\sigma_{j-1}^{-1}\\
&\equiv&(j-1,j)^{-1}(jk)\{(j-1,l)^{\varepsilon},(jl)\}(j-1,j)\sigma_{j-1}^{-1}\\
&\equiv&(j-1,j)^{-1}\{(j-1,l)^{\varepsilon},\{(jl),(kl)\}\}(jk)(j-1,j)\sigma_{j-1}^{-1}\\
&\equiv&\{\{(j-1,l)^{\varepsilon},(jl)^{-1}(j-1,l)^{-1}\},\{\{(jl),(j-1,l)^{-1}\},(kl)\}\}
\{(jk),(j-1,k)^{-1}\}\sigma_{j-1}^{-1}\\
&\equiv&\{(j-1,l)^{\varepsilon},(jl)^{-1}(j-1,l)^{-1}(kl)^{-1}(j-1,l)(jl)(j-1,l)^{-1}(kl)\}\\
&&(j-1,k)(jk)(j-1,k)^{-1}\sigma_{j-1}^{-1}.
\end{eqnarray*}
\item[2)]\ $i<j-1$.
\begin{eqnarray*}
&&\sigma_{i}^{-1}\{(jl)^{\varepsilon},(il)^{-1}(kl)^{-1}(il)(kl)\}(ik)\\
&\equiv&\{(jl)^{\varepsilon},\{(i+1,l)^{-1},(i,i+1)\}(kl)^{-1}\{(i+1,l)(i,i+1)\}\\
&&(kl)\}\{(i+1,k)(i.i+1)\}\sigma_{i}^{-1}\\
&\equiv&(kl)^{-1}(i,i+1)^{-1}(i+1,l)^{-1}(i,i+1)(kl)(i,i+1)^{-1}(i+1,l)(i,i+1)(jl)^{\varepsilon}\\
&&(i,i+1)^{-1}(i+1,l)^{-1}(i,i+1)(kl)^{-1}(i,i+1)^{-1}(i+1,l)(i,i+1)(kl)\\
&&(i,i+1)^{-1}(i+1,k)(i,i+1)\sigma_{i}^{-1}\\
&\equiv&\{(jl)^{\varepsilon},(il)(i+1,l)^{-1}(il)^{-1}(kl)^{-1}(il)(i+1,l)(il)^{-1}(kl)\}(ik)(i+1,k)\\
&&(ik)^{-1}\sigma_{i}^{-1} \ \ \ \  and \\
&&\\
&&\{(i+1,k),(i,i+1)\}\sigma_{i}^{-1}(jl)^{\varepsilon}\\
&\equiv&(i,i+1)^{-1}(i+1,k)(i,i+1)(jl)^{\varepsilon}\sigma_{i}^{-1}\\
&\equiv&(i,i+1)^{-1}(i+1,k)(jl)^{\varepsilon}(i,i+1)\sigma_{i}^{-1}\\
&\equiv&(i,i+1)^{-1}(kl)^{-1}(i+1,l)^{-1}(kl)(i+1,l)(jl)^{\varepsilon}(i+1,l)^{-1}(kl)^{-1}\\
&&(i+1,l)(kl)(i+1,k)(i,i+1)\sigma_{i}^{-1}\\
&\equiv&\{(jl)^{\varepsilon},(il)(i+1,l)^{-1}(il)^{-1}(kl)^{-1}(il)(i+1,l)(il)^{-1}(kl)\}(ik)(i+1,k)\\
&&(ik)^{-1}\sigma_{i}^{-1}.
\end{eqnarray*}
\end{enumerate}

$(4)\wedge(13)$\\
Let
$f=\sigma_{i}^{-1}(ik)^\delta-\{(i+1,k)^\delta,(i,i+1)\}\sigma_{i}^{-1},\
g=(ik)^\delta(jl)^{\varepsilon}-(jl)^{\varepsilon}(ik)^\delta,j<i<i+1<k<l,or,i<i+1<k<l$.
Then $w=\sigma_{i}^{-1}(ik)^\delta(jl)^{\varepsilon}$ and
\begin{eqnarray*}
(f,g)_w&=&\sigma_{i}^{-1}(jl)^{\varepsilon}(ik)^\delta-\{(i+1,k)^\delta,(i,i+1)\}\sigma_{i}^{-1}(jl)^{\varepsilon}\\
&\equiv&(jl)^{\varepsilon}\{(i+1,k)^\delta,(i,i+1)\}\sigma_{i}^{-1}-(jl)^{\varepsilon}\{(i+1,k)^\delta,(i,i+1)\}
\sigma_{i}^{-1}\\
&\equiv&0.
\end{eqnarray*}

$(5)\wedge(7)$\\
Let $f=\sigma_{k-1}^{-1}(jk)^{-1}-(j,k-1)^{-1}\sigma_{k-1}^{-1},\
g=(jk)^{-1}(kl)^{\varepsilon}-\{(kl)^{\varepsilon},(jl)^{-1}\}(jk)^{-1},j<k-1<k<l$.
Then $w=\sigma_{k-1}^{-1}(jk)^{-1}(kl)^{\varepsilon}$ and
\begin{eqnarray*}
(f,g)_w&=&\sigma_{k-1}^{-1}\{(kl)^{\varepsilon},(jl)^{-1}\}(jk)^{-1}-(j,k-1)^{-1}\sigma_{k-1}^{-1}(kl)^{\varepsilon}\\
&\equiv&\{(k-1,l)^{\varepsilon},(jl)^{-1}\}(j,k-1)^{-1}\sigma_{k-1}^{-1}-
(j,k-1)^{-1}(k-1,l)^{\varepsilon}\sigma_{k-1}^{-1}\\
&\equiv&0.
\end{eqnarray*}

$(5)\wedge(8)$\\
Let $f=\sigma_{k-1}^{-1}(jk)-(j,k-1)\sigma_{k-1}^{-1},\
g=(jk)(kl)^{\varepsilon}-\{(kl)^{\varepsilon},(jl)(kl)\}(jk),j<k-1<k<l$.
Then $w=\sigma_{k-1}^{-1}(jk)(kl)^{\varepsilon}$ and
\begin{eqnarray*}
(f,g)_w&=&\sigma_{k-1}^{-1}\{(kl)^{\varepsilon},(jl)(kl)\}(jk)-(j,k-1)\sigma_{k-1}^{-1}(kl)^{\varepsilon}\\
&\equiv&\{(k-1,l)^{\varepsilon},(jl)(k-1,l)\}(j,k-1)\sigma_{k-1}^{-1}-
(j,k-1)(k-1,l)^{\varepsilon}\sigma_{k-1}^{-1}\\
&\equiv&0.
\end{eqnarray*}

$(5)\wedge(9)$\\
Let $f=\sigma_{k-1}^{-1}(jk)^{-1}-(j,k-1)^{-1}\sigma_{k-1}^{-1},\
g=(jk)^{-1}(jl)^{\varepsilon}-\{(jl)^{\varepsilon},(kl)^{-1}(jl)^{-1}\}(jk)^{-1},j<k-1<k<l$.
Then $w=\sigma_{k-1}^{-1}(jk)^{-1}(jl)^{\varepsilon}$ and
\begin{eqnarray*}
(f,g)_w&=&\sigma_{k-1}^{-1}\{(jl)^{\varepsilon},(kl)^{-1}(jl)^{-1}\}(jk)^{-1}-
(j,k-1)^{-1}\sigma_{k-1}^{-1}(jl)^{\varepsilon}\\
&\equiv&\{(jl)^{\varepsilon},(k-1,l)^{-1}(jl)^{-1}\}(j,k-1)^{-1}\sigma_{k-1}^{-1}-
(j,k-1)^{-1}(jl)^{\varepsilon}\sigma_{k-1}^{-1}\\
&\equiv&0.
\end{eqnarray*}

$(5)\wedge(10)$\\
Let $f=\sigma_{k-1}^{-1}(jk)-(j,k-1)\sigma_{k-1}^{-1},\
g=(jk)(jl)^{\varepsilon}-\{(jl)^{\varepsilon},(kl)\}(jk),j<k-1<k<l$.
Then $w=\sigma_{k-1}^{-1}(jk)(jl)^{\varepsilon}$ and
\begin{eqnarray*}
(f,g)_w&=&\sigma_{k-1}^{-1}\{(jl)^{\varepsilon},(kl)\}(jk)-(j,k-1)\sigma_{k-1}^{-1}(jl)^{\varepsilon}\\
&\equiv&\{(jl)^{\varepsilon},(k-1,l)\}(j,k-1)\sigma_{k-1}^{-1}-(j,k-1)(jl)^{\varepsilon}\sigma_{k-1}^{-1}\\
&\equiv&0.
\end{eqnarray*}

$(5)\wedge(11)$\\
Let $f=\sigma_{k-1}^{-1}(ik)^{-1}-(i,k-1)^{-1}\sigma_{k-1}^{-1},
g=(ik)^{-1}(jl)^{\varepsilon}-\{(jl)^{\varepsilon},(kl)(il)(kl)^{-1}(il)^{-1}\}(ik)^{-1},i<j<k<l.$\
Then
$w=\sigma_{k-1}^{-1}(ik)^{-1}(jl)^{\varepsilon}$ \ and\\
$$
(f,g)_w=\sigma_{k-1}^{-1}\{(jl)^{\varepsilon},(kl)(il)(kl)^{-1}(il)^{-1}\}(ik)^{-1}
-(i,k-1)^{-1}\sigma_{k-1}^{-1}(jl)^{\varepsilon}.
$$
 There are two cases to consider.
\begin{enumerate}
\item[1)]\ $j=k-1$.
\begin{eqnarray*}
&& \sigma_{k-1}^{-1}\{(k-1,l)^{\varepsilon},(kl)(il)(kl)^{-1}(il)^{-1}\}(ik)^{-1}\\
&\equiv&\{\{(kl)^{\varepsilon},(k-1,k)\},(k-1,l)(il)(k-1,l)^{-1}(il)^{-1}\}(i,k-1)^{-1}\sigma_{k-1}^{-1}\\
&\equiv&(il)(k-1,l)(il)^{-1}(k-1,l)^{-1}(k-1,k)^{-1}(kl)^{\varepsilon}(k-1,k)(k-1,l)(il)\\
&&(k-1,l)^{-1}
(il)^{-1}(i,k-1)^{-1}\sigma_{k-1}^{-1}\\
&\equiv&(il)(k-1,l)(il)^{-1}(kl)^{\varepsilon}(il)(k-1,l)^{-1}(il)^{-1}(i,k-1)^{-1}\sigma_{k-1}^{-1}\ \ \ \ \ and \\
&&\\
&&(i,k-1)^{-1}\sigma_{k-1}^{-1}(k-1,l)^{\varepsilon}\\
&\equiv&(i,k-1)^{-1}\{(kl)^{\varepsilon},(k-1,k)\}\sigma_{k-1}^{-1}\\
&\equiv&\{(kl)^{\varepsilon},\{(k-1,k),(ik)^{-1}\}(i,k-1)^{-1}\sigma_{k-1}^{-1}\\
&\equiv&(ik)(k-1,k)^{-1}(ik)^{-1}(kl)^{\varepsilon}(ik)(k-1,k)(ik)^{-1}(i,k-1)^{-1}\sigma_{k-1}^{-1}\\
&\equiv&(ik)(k-1,k)^{-1}\{(kl)^{\varepsilon},(il)^{-1}\}(k-1,k)(ik)^{-1}(i,k-1)^{-1}\sigma_{k-1}^{-1}\\
&\equiv&(ik)\{\{(kl)^{\varepsilon},(k-1,l)^{-1}\},(il)^{-1}\}(ik)^{-1}(i,k-1)^{-1}\sigma_{k-1}^{-1}\\
&\equiv&\{\{\{(kl)^{\varepsilon},(il)(kl)\},\{(k-1,l)^{-1},(il)^{-1}(kl)^{-1}(il)(kl)\}\},\{(il)^{-1},(kl)\}\}\\
&&(i,k-1)^{-1}\sigma_{k-1}^{-1}\\
&\equiv&(il)(k-1,l)(il)^{-1}(kl)^{\varepsilon}(il)(k-1,l)^{-1}(il)^{-1}(i,k-1)^{-1}\sigma_{k-1}^{-1}.
\end{eqnarray*}
\item[2)]\ $j<k-1$.
\begin{eqnarray*}
(f,g)_w&=&\sigma_{k-1}^{-1}\{(jl)^{\varepsilon},(kl)(il)(kl)^{-1}(il)^{-1}\}(ik)^{-1}
-(i,k-1)^{-1}\sigma_{k-1}^{-1}(jl)^{\varepsilon}\\
&\equiv&\{(jl)^{\varepsilon},(k-1,l)(il)(k-1,l)^{-1}(il)^{-1}\}(i,k-1)^{-1}\sigma_{k-1}^{-1}
-(i,k-1)^{-1}\\
&&(jl)^{\varepsilon}\sigma_{k-1}^{-1}\\
&\equiv&0.
\end{eqnarray*}
\end{enumerate}

$(5)\wedge(12)$\\
Let $f=\sigma_{k-1}^{-1}(ik)-(i,k-1)\sigma_{k-1}^{-1},\
g=(ik)(jl)^{\varepsilon}-\{(jl)^{\varepsilon},(il)^{-1}(kl)^{-1}(il)(kl)\}(ik),i<j<k<l$.
Then $w=\sigma_{k-1}^{-1}(ik)(jl)^{\varepsilon}$ and
$$
(f,g)_w=\sigma_{k-1}^{-1}\{(jl)^{\varepsilon},(il)^{-1}(kl)^{-1}(il)(kl)\}(ik)-
(i,k-1)\sigma_{k-1}^{-1}(jl)^{\varepsilon}.
$$
There are two cases to consider.
\begin{enumerate}
\item[1)]\ $j=k-1$.
\begin{eqnarray*}
&&\sigma_{k-1}^{-1}\{(k-1,l)^{\varepsilon},(il)^{-1}(kl)^{-1}(il)(kl)\}(ik)\\
&\equiv&\{\{(kl)^{\varepsilon},(k-1,k)\}(il)^{-1}(k-1,l)^{-1}(il)(k-1,l)\}(i,k-1)\sigma_{k-1}^{-1}\\
&\equiv&\{(kl)^{\varepsilon},(k-1,l)^{-1}(il)^{-1}(k-1,l)^{-1}(il)(k-1,l)\}(i,k-1)\sigma_{k-1}^{-1} \ \  \ \ and \\
&&\\
&&(i,k-1)\sigma_{k-1}^{-1}(k-1,l)^{\varepsilon}\\
&\equiv&(i,k-1)\{(kl)^{\varepsilon},(k-1,k)\}\sigma_{k-1}^{-1}\\
&\equiv&\{(kl)^{\varepsilon},\{(k-1,k)(ik)(k-1,k)\}\}(i,k-1)\sigma_{k-1}^{-1}\\
&\equiv&(k-1,k)^{-1}(ik)^{-1}(k-1,k)^{-1}(ik)(k-1,k)(kl)^{\varepsilon}(k-1,k)^{-1}(ik)^{-1}(k-1,k)(ik)\\
&&(k-1,k)(i,k-1)\sigma_{k-1}^{-1}\\
&\equiv&(k-1,k)^{-1}(ik)^{-1}(k-1,k)^{-1}(ik)\{(kl)^{\varepsilon},(k-1,l)(kl)\}(ik)^{-1}(k-1,k)(ik)(k-1,k)\\
&&(i,k-1)\sigma_{k-1}^{-1}\\
&\equiv&(k-1,k)^{-1}(ik)^{-1}(k-1,k)^{-1}\{\{(kl)^{\varepsilon},(il)(kl)\},\{(k-1,l),(il)^{-1}(kl)^{-1}(il)(kl)\}\\
&&\{(kl),(il)(kl)\}\}(k-1,k)(ik)(k-1,k)(i,k-1)\sigma_{k-1}^{-1}\\
&\equiv&(k-1,k)^{-1}(ik)^{-1}\{\{(kl)^{\varepsilon},(k-1,l)^{-1}\},(il)\{(k-1,l),(kl)^{-1}(k-1,l)^{-1}\}\\
&&\{(kl),(k-1,l)^{-1}\}\}(ik)(k-1,k)(i,k-1)\sigma_{k-1}^{-1}\\
&\equiv&(k-1,k)^{-1}\{\{\{(kl)^{\varepsilon},(il)^{-1}\},\{(k-1,l)^{-1},(kl)(il)(kl)^{-1}(il)^{-1}\}\},\{(il),(kl)^{-1}\\
&&(il)^{-1}\}\{(k-1,l)(kl)(il)(kl)^{-1}(il)^{-1}\}\{(kl)(il)^{-1}\}\}(k-1,k)(i,k-1)\sigma_{k-1}^{-1}\\
&\equiv&\{\{(kl)^{\varepsilon},(k-1,l)^{-1}\},\{\{(k-1,l)^{-1},(kl)^{-1}(k-1,l)^{-1}\},\{(kl),(k-1,l)^{-1}\}(il)\}\\
&&\{\{(k-1,l),(kl)^{-1}(k-1,l)^{-1}\},\{(kl),(k-1,l)^{-1}\}\}(i,k-1)\sigma_{k-1}^{-1}\\
&\equiv&\{(kl)^{\varepsilon},(k-1,l)^{-1}(il)^{-1}(k-1,l)^{-1}(il)(k-1,l)\}(i,k-1)\sigma_{k-1}^{-1}.
\end{eqnarray*}
\item[2)]\ $j<k-1$.
$$
(f,g)_w=\{(jl)^{\varepsilon},(il)^{-1}(k-1,l)^{-1}(il)(k-1,l)\}(i,k-1)\sigma_{k-1}^{-1}-
(i,k-1)(jl)^{\varepsilon}\sigma_{k-1}^{-1}\equiv0.
$$
\end{enumerate}

$(5)\wedge(13)$ \\
Let
$f=\sigma_{k-1}^{-1}(ik)^{\delta}-(i,k-1)^{\delta}\sigma_{k-1}^{-1},\
g=(ik)^\delta(jl)^{\varepsilon}-(jl)^{\varepsilon}(ik)^{\delta},j<i<k-1<k<l\
\ or \ \  i<k-1<k<j<l$. Then
$w=\sigma_{k-1}^{-1}(ik)^{\delta}(jl)^{\varepsilon}$ and

\begin{eqnarray*}
(f,g)_w&=&\sigma_{k-1}^{-1}(jl)^{\varepsilon}(ik)^{\delta}-(i,k-1)^{\delta}\sigma_{k-1}^{-1}(jl)^{\varepsilon}\\
&\equiv&(jl)^{\varepsilon}(i,k-1)^{\delta}\sigma_{k-1}^{-1}-(i,k-1)^{\delta}(jl)^{\varepsilon}\sigma_{k-1}^{-1}\\
&\equiv&(jl)^{\varepsilon}(i,k-1)^{\delta}\sigma_{k-1}^{-1}-(jl)^{\varepsilon}(i,k-1)^{\delta}\sigma_{k-1}^{-1}\\
&\equiv&0.
\end{eqnarray*}

$(6)\wedge(7)$\\
Let
$f=\sigma_{k}^{-1}(jk)^{-1}-\{(j,k+1)^{-1},(k,k+1)\}\sigma_{k}^{-1},\
g=(jk)^{-1}(kl)^{\varepsilon}-\{(kl)^{\varepsilon},(jl)^{-1}\}(jk)^{-1},j<k<l$.
Then $w=\sigma_{k}^{-1}(jk)^{-1}(kl)^{\varepsilon}$ and
$$
(f,g)_w=\sigma_{k}^{-1}\{(kl)^{\varepsilon},(jl)^{-1}\}(jk)^{-1}-
\{(j,k+1)^{-1},(k,k+1)\}\sigma_{k}^{-1}(kl)^{\varepsilon}.
$$
There are two cases to consider.
\begin{enumerate}
\item[1)]\ $k=l-1$.
\begin{eqnarray*}
&&\sigma_{k}^{-1}\{(k,k+1)^{\varepsilon},(j,k+1)^{-1}\}(jk)^{-1}\\
&\equiv&\{(k,k+1)^{\varepsilon},(jk)^{-1}\}\{(j,k+1)^{-1},(k,k+1)\}\sigma_{k}^{-1}\\
&\equiv&(jk)(k,k+1)^{\varepsilon}(jk)^{-1}(k,k+1)^{-1},(j,k+1)^{-1}(k,k+1)\sigma_{k}^{-1}\\
&\equiv&(k,k+1)^{-1}(j,k+1)^{-1}(k,k+1)^{\varepsilon}(k,k+1)\sigma_{k}^{-1} \ \ \ and\\
&&\\
&&\{(j,k+1)^{-1},(k,k+1)\}\sigma_{k}^{-1}(kl)^{\varepsilon}\\
&\equiv&(k,k+1)^{-1}(j,k+1)^{-1}(k,k+1)(k,k+1)^{\varepsilon}\sigma_{k}^{-1}\\
&\equiv&(k,k+1)^{-1}(j,k+1)^{-1}(k,k+1)^{\varepsilon}(k,k+1)\sigma_{k}^{-1}.
\end{eqnarray*}
\item[2)]\  $k<l-1$, then   $j<k<k+1<l$.
\begin{eqnarray*}
&&\sigma_{k}^{-1}\{(kl)^{\varepsilon},(jl)^{-1}\}(jk)^{-1}\\
&\equiv&\{\{(k+1,l)^{\varepsilon},(k,k+1)\},(jl)^{-1}\}\{(j,k+1)^{-1},(k,k+1)\}\sigma_{k}^{-1}\\
&\equiv&(jl)(k,k+1)^{-1}(k+1,l)^{\varepsilon}(k,k+1)(jl)^{-1}(k,k+1)^{-1}(j,k+1)^{-1}(k,k+1)\}\sigma_{k}^{-1}\\
&\equiv&(jl)(kl)(k+1,l)^{\varepsilon}(kl)^{-1}(jl)^{-1}(k,k+1)^{-1}(j,k+1)^{-1}(k,k+1)\}\sigma_{k}^{-1}  \ \ \ \  and \\
&&\\
&&\{(j,k+1)^{-1},(k,k+1)\}\sigma_{k}^{-1}(kl)^{\varepsilon}\\
&\equiv&(k,k+1)^{-1}(j,k+1)^{-1}(k,k+1)(k,k+1)^{-1}(k+1,l)^{\varepsilon}(k,k+1)\sigma_{k}^{-1}\\
&\equiv&(k,k+1)^{-1}\{(k+1,l)^{\varepsilon},(jl)^{-1}\}(j,k+1)^{-1}(k,k+1)\}\sigma_{k}^{-1}\\
&\equiv&\{\{(k+1,l)^{\varepsilon},(kl)^{-1}\},(jl)^{-1}\}(k,k+1)^{-1}(j,k+1)^{-1}(k,k+1)\sigma_{k}^{-1}\\
&\equiv&(jl)(kl)(k+1,l)^{\varepsilon}(kl)^{-1}(jl)^{-1}(k,k+1)^{-1}(j,k+1)^{-1}(k,k+1)\sigma_{k}^{-1}.
\end{eqnarray*}
\end{enumerate}

$(6)\wedge(8)$\\
Let $f=\sigma_{k}^{-1}(jk)-\{(j,k+1),(k,k+1)\}\sigma_{k}^{-1},
g=(jk)(kl)^{\varepsilon}-\{(kl)^{\varepsilon},(jl)(kl)\}(jk),j<k<l$.\
Then $w=\sigma_{k}^{-1}(jk)(kl)^{\varepsilon}$\  and
$$
(f,g)_w=\sigma_{k}^{-1}\{(kl)^{\varepsilon},(jl)(kl)\}(jk)-\{(j,k+1),(k,k+1)\}\sigma_{k}^{-1}(kl)^{\varepsilon}.
$$
There are two cases to consider.
\begin{enumerate}
\item[1)]\ $l=k+1$.
\begin{eqnarray*}
(f,g)_w&=&\sigma_{k}^{-1}\{(k,k+1)^{\varepsilon},(j,k+1)(k,k+1)\}(jk)-
\{(j,k+1),(k,k+1)\}\sigma_{k}^{-1}\\
&&(k,k+1)^{\varepsilon}\\
&\equiv&\{(k,k+1)^{\varepsilon},(jk)(k,k+1)\}\{(j,k+1)(k,k+1)\}\sigma_{k}^{-1}\\
&&-(k,k+1)^{-1}(j,k+1),(k,k+1)(k,k+1)^{\varepsilon}\sigma_{k}^{-1}\\
&\equiv&(k,k+1)^{-1}(jk)^{-1}(k,k+1)^{\varepsilon}(jk)(j,k+1)(k,k+1)\sigma_{k}^{-1}-(k,k+1)^{-1}\\
&&(j,k+1),(k,k+1)^{\varepsilon+1}\sigma_{k}^{-1}\\
&\equiv&(k,k+1)^{-1}(j,k+1)(k,k+1)^{\varepsilon
+1}\sigma_{k}^{-1}-(k,k+1)^{-1}(j,k+1)\\
&&(k,k+1)^{\varepsilon+1}\sigma_{k}^{-1}\\
&\equiv&0.
\end{eqnarray*}
\item[2)]\ $l<k+1,j<k<k+1<l$.
\begin{eqnarray*}
(f,g)_w&=&\{\{(k+1,l)^{\varepsilon},(k,k+1)\},(jl)\{(k+1,l),(k,k+1)\}\}\{(j,k+1),\\
&&(k,k+1)\}\sigma_{k}^{-1}
-(k,k+1)^{-1}(j,k+1)(k+1,l)^{\varepsilon}(k,k+1)\sigma_{k}^{-1}\\
&\equiv&(k,k+1)^{-1}(k+1,l)^{-1}(k,k+1)(jl)^{-1}(jl)(k,k+1)^{-1}(k+1,l)\\
&&(j,k+1)(k,k+1)\sigma_{k}^{-1}-(k,k+1)^{-1}\{(k+1,l)^{\varepsilon},(jl)\\
&&(k+1,l)\}(j,k+1)(k,k+1)\sigma_{k}^{-1}\\
&\equiv&\{\{(k+1,l)^{\varepsilon},(kl)^{-1}\},(jl)\{(k+1,l)(kl)^{-1}\}\}(k,k+1)^{-1}(j,k+1)\\
&&(k,k+1)\sigma_{k}^{-1}-\{\{(k+1,l)^{\varepsilon},(kl)^{-1}\},(jl)\{(k+1,l)(kl)^{-1}\}\}\\
&&(k,k+1)^{-1}(j,k+1)(k,k+1)\sigma_{k}^{-1}\\
&\equiv&0.
\end{eqnarray*}
\end{enumerate}

$(6)\wedge(9)$\\
Let
$f=\sigma_{k}^{-1}(jk)^{-1}-\{(j,k+1)^{-1},(k,k+1)\}\sigma_{k}^{-1},\
g=(jk)^{-1}(jl)^{\varepsilon}-\{(jl)^{\varepsilon},(kl)^{-1}(jl)^{-1}\}(jk)^{-1},j<k<l$.
Then $w=\sigma_{k}^{-1}(jk)^{-1}(jl)^{\varepsilon}$ and
$$
(f,g)_w=\sigma_{k}^{-1}\{(jl)^{\varepsilon},(kl)^{-1}(jl)^{-1}\}(jk)^{-1}-
\{(j,k-1)^{-1},(k,k+1)\}\sigma_{k}^{-1}(jl)^{\varepsilon}.
$$
There are two cases to consider.
\begin{enumerate}
\item[1)]\ $k=l-1$.
\begin{eqnarray*}
&&\sigma_{k}^{-1}\{(j,k+1)^{\varepsilon},(k,k+1)^{-1}(j,k+1)^{-1}\}(jk)^{-1}\\
&\equiv&\{(jk)^{\varepsilon},(k,k+1)^{-1}(jk)^{-1}\}\{(j,k+1)^{-1}(k,k+1)\}\sigma_{k}^{-1}\\
&\equiv&(jk)(k,k+1)(jk)^{\varepsilon}(k,k+1)^{-1}(jk)^{-1}(k,k+1)^{-1}(j,k+1)^{-1}(k,k+1)\sigma_{k}^{-1}.
\end{eqnarray*}
If $\varepsilon=1$ , then
\begin{eqnarray*}
&&\sigma_{k}^{-1}\{(j,k+1)^{\varepsilon},(k,k+1)^{-1}(j,k+1)^{-1}\}(jk)^{-1}\\
&\equiv&(jk)(j,k+1)^{-1}\sigma_{k}^{-1}\equiv(k,k+1)^{-1}(j,k+1)^{-1}(k,k+1)(jk)\sigma_{k}^{-1}.
\end{eqnarray*}
If $\varepsilon=-1$, then
\begin{eqnarray*}
&&\sigma_{k}^{-1}\{(j,k+1)^{\varepsilon},(k,k+1)^{-1}(j,k+1)^{-1}\}(jk)^{-1}\\
&\equiv&(k,k+1)^{-1}(j,k+1)^{-1}(k,k+1)(j,k+1)(jk)^{-1}(k,k+1)^{-1}(j,k+1)^{-1}(k,k+1)\sigma_{k}^{-1}\\
&\equiv&(k,k+1)^{-1}(j,k+1)^{-1}(k,k+1)(j,k+1)(j,k+1)(k,k+1)^{-1}(j,k+1)^{-1}(jk)^{-1}\\
&&(j,k+1)^{-1}(k,k+1)\sigma_{k}^{-1}\\
&\equiv&(k,k+1)^{-1}(j,k+1)^{-1}(k,k+1)(jk)^{-1}\sigma_{k}^{-1}.
\end{eqnarray*}
Hence,
\begin{eqnarray*}
&&\sigma_{k}^{-1}\{(j,k+1)^{\varepsilon},(k,k+1)^{-1}(j,k+1)^{-1}\}(jk)^{-1}\\
&\equiv&(k,k+1)^{-1}(j,k+1)^{-1}(k,k+1)(jk)^{\varepsilon}\sigma_{k}^{-1} \ \ \ \  and \\
&&\\
&&\{(j,k+1)^{-1},(k,k+1)\}\sigma_{k}^{-1}(j,k+1)^{\varepsilon}\\
&\equiv&(k,k+1)^{-1}(j,k+1)^{-1}(k,k+1)(jk)^{\varepsilon}\sigma_{k}^{-1}.
\end{eqnarray*}
\item[2)]\ $k<l-1,j<k<k+1<l$.
\begin{eqnarray*}
(f,g)_w&=&\{(jl)^{\varepsilon},\{(k+1,l)^{-1},(k,k+1)\}(jl)^{-1}\}\{(j,k+1)^{-1},(k,k+1)\}\sigma_{k}^{-1}\\
&&-(k,k+1)^{-1}(j,k+1)^{-1}(k,k+1)(jl)^{\varepsilon}\sigma_{k}^{-1}\\
&\equiv&(jl)(k,k+1)^{-1}(k+1,l)(k,k+1)(jl)^{\varepsilon}(k,k+1)^{-1}(k+1,l)^{-1}\\
&&(k,k+1)(jl)^{-1}(k,k+1)^{-1}(j,k+1)^{-1}(k,k+1)\sigma_{k}^{-1}-(k,k+1)^{-1}\\
&&(j,k+1)^{-1}(jl)^{\varepsilon}(k,k+1)\sigma_{k}^{-1}\\
&\equiv&(jl)(kl)(k+1,l)(kl)^{-1}(jl)^{\varepsilon}(kl)(k+1,l)^{-1}(kl)^{-1}(jl)^{-1}(j,k+1)^{-1}\sigma_{k}^{-1}\\
&&-(k,k+1)^{-1}\{(jl)^{\varepsilon},(k+1,l)^{-1}(jl)^{-1}\}(j,k+1)^{-1}(k,k+1)\sigma_{k}^{-1}\\
&\equiv&(jl)(kl)(k+1,l)(kl)^{-1}(jl)^{\varepsilon}(kl)(k+1,l)^{-1}(kl)^{-1}(jl)^{-1}(j,k+1)^{-1}\sigma_{k}^{-1}\\
&&-\{(jl)^{\varepsilon},\{(k+1,l)^{-1},(kl)^{-1}\}(jl)^{-1}\}(k,k+1)^{-1}(j,k+1)^{-1}(k,k+1)\sigma_{k}^{-1}\\
&\equiv&0.
\end{eqnarray*}
\end{enumerate}

$(6)\wedge(10)$\\
Let $f=\sigma_{k}^{-1}(jk)-\{(j,k+1),(k,k+1)\}\sigma_{k}^{-1},\
g=(jk)(jl)^{\varepsilon}-\{(jl)^{\varepsilon},(kl)\}(jk),j<k<l$.
Then $w=\sigma_{k}^{-1}(jk)(jl)^{\varepsilon}$ and
$$
(f,g)_w=\sigma_{k}^{-1}\{(jl)^{\varepsilon},(kl)\}(jk)-\{(j,k+1),(k,k+1)\}\sigma_{k}^{-1}(jl)^{\varepsilon}.
$$
There are two cases to consider.
\begin{enumerate}
\item[1)]\ $l=k+1$.
\begin{eqnarray*}
&&\sigma_{k}^{-1}\{(j,k+1)^{\varepsilon},(k,k+1)\}(jk)\\
&\equiv&\{(jk)^{\varepsilon},(k,k+1)\}\{(j,k+1),(k,k+1)\}\sigma_{k}^{-1}\\
&\equiv&(k,k+1)^{-1}(jk)^{\varepsilon}(j,k+1)(k,k+1)\sigma_{k}^{-1}.
\end{eqnarray*}
If \ $\varepsilon=1$, then
\begin{eqnarray*}
&&\sigma_{k}^{-1}\{(j,k+1)^{\varepsilon},(k,k+1)\}(jk)\\
&\equiv&(k,k+1)^{-1}(k,k+1)^{-1}(j,k+1)(k,k+1)(jk)(k,k+1)\sigma_{k}^{-1}\\
&\equiv&(k,k+1)^{-1}(j,k+1)(k,k+1)(jk)\sigma_{k}^{-1}.
\end{eqnarray*}
If \ $\varepsilon=-1$, then
\begin{eqnarray*}
&&\sigma_{k}^{-1}\{(j,k+1)^{\varepsilon},(k,k+1)\}(jk)\\
&\equiv&(k,k+1)^{-1}(j,k+1)(k,k+1)(j,k+1)(k,k+1)^{-1}(j,k+1)^{-1}(jk)^{-1}(k,k+1)\sigma_{k}^{-1}\\
&\equiv&(k,k+1)^{-1}(j,k+1)(k,k+1)(jk)^{-1}\sigma_{k}^{-1}.
\end{eqnarray*}
Hence,
$$
\sigma_{k}^{-1}\{(j,k+1)^{\varepsilon},(k,k+1)\}(jk)\equiv(k,k+1)^{-1}(j,k+1)(k,k+1)(jk)^{\varepsilon}\sigma_{k}^{-1}.
$$
Also,
$$
\{(j,k+1),(k,k+1)\}\sigma_{k}^{-1}(j,k+1)^{\varepsilon}\equiv(k,k+1)^{-1}(j,k+1)(k,k+1)(jk)^{\varepsilon}\sigma_{k}^{-1}.
$$
\item[2)]\ $l<k+1$.
\begin{eqnarray*}
(f,g)_w&\equiv&\{(jl)^{\varepsilon},\{(k+1,l),(k,k+1)\}\}\{(j,k+1),(k,k+1)\}\\
&&-(k,k+1)^{-1}(j,k+1)(k,k+1)(jl)^{\varepsilon}\sigma_{k}^{-1}\\
&\equiv&(k,k+1)^{-1}(k+1,l)^{-1}(k,k+1)(jl)^{\varepsilon}(k,k+1)^{-1}(k+1,l)(j,k+1)(k,k+1)\\
&&-(k,k+1)^{-1}(j,k+1)(jl)^{\varepsilon}(k,k+1)\sigma_{k}^{-1}\\
&\equiv&(kl)(k+1,l)^{-1}(kl)^{-1}(jl)^{\varepsilon}(kl)(k+1,l)(kl)^{-1}(k,k+1)^{-1}(j,k+1)(k,k+1)\\
&&-(k,k+1)^{-1}\{(jl)^{\varepsilon},(k+1,l)\}(j,k+1)(k,k+1)\sigma_{k}^{-1}\\
&\equiv&(kl)(k+1,l)^{-1}(kl)^{-1}(jl)^{\varepsilon}(kl)(k+1,l)(kl)^{-1}(k,k+1)^{-1}(j,k+1)(k,k+1)\\
&&-\{(jl)^{\varepsilon},\{(k+1,l),(kl)^{-1}\}\}(k,k+1)^{-1}(j,k+1)(k,k+1)\sigma_{k}^{-1}\\
&\equiv&0.
\end{eqnarray*}
\end{enumerate}

$(6)\wedge(11)$\\
Let
$f=\sigma_{k}^{-1}(ik)^{-1}-\{(i,k+1)^{-1},(k,k+1)\}\sigma_{k}^{-1},\
g=(ik)^{-1}(jl)^{\varepsilon}-\{(jl)^{\varepsilon},(kl)(il)(kl)^{-1}(il)^{-1}\}(ik)^{-1},i<j<k<l$.
Then $w=\sigma_{k}^{-1}(ik)^{-1}(jl)^{\varepsilon}$ and
$$
(f,g)_w=\sigma_{k}^{-1}\{(jl)^{\varepsilon},(kl)(il)(kl)^{-1}(il)^{-1}\}(ik)^{-1}-
\{(i,k+1)^{-1},(k,k+1)\}\sigma_{k}^{-1}(jl)^{\varepsilon}.
$$
There are two cases to consider.
\begin{enumerate}
\item[1)]\ $k=l-1$.
\begin{eqnarray*}
&&\sigma_{k}^{-1}\{(j,k+1)^{\varepsilon},(k,k+1)(i,k+1)(k,k+1)^{-1}(i,k+1)^{-1}\}(ik)^{-1}\\
&\equiv&\{(jk)^{\varepsilon},(k,k+1)(ik)(k,k+1)^{-1}(ik)^{-1}\}\{(i,k+1)^{-1},(k,k+1)\}\sigma_{k}^{-1}\\
&\equiv&(ik)(k,k+1)(ik)^{-1}(k,k+1)^{-1}(jk)^{\varepsilon}(k,k+1)(ik)(k,k+1)^{-1}(ik)^{-1}(k,k+1)^{-1}\\
&&(i,k+1)^{-1}(k,k+1)\sigma_{k}^{-1}\\
&\equiv&(k,k+1)^{-1}(i,k+1)^{-1}(k,k+1)(i,k+1)(jk)^{\varepsilon}(k,k+1)(k,k+1)^{-1}(i,k+1)^{-1}\\
&&(k,k+1)^{-1}(i,k+1)(k,k+1)(k,k+1)^{-1}(i,k+1)^{-1}(k,k+1)\sigma_{k}^{-1}\\
&\equiv&(k,k+1)^{-1}(i,k+1)^{-1}(k,k+1)(jk)^{\varepsilon}\sigma_{k}^{-1} \ \ \ \ \ and\\
&&\\
&&\{(i,k+1)^{-1},(k,k+1)\}\sigma_{k}^{-1}(j,k+1)^{\varepsilon}\\
&\equiv&(k,k+1)^{-1}(i,k+1)^{-1}(k,k+1)(jk)^{\varepsilon}\sigma_{k}^{-1}.
\end{eqnarray*}
\item[2)]\ $k<l+1, \ i<j<k<k+1<l$.
\begin{eqnarray*}
(f,g)_w&\equiv&\{(jl)^{\varepsilon},\{(k+1,l),(k,k+1)\}(il)\{(k+1,l)^{-1},(k,k+1)\}(il)^{-1}\}\\
&&\{(i,k+1)^{-1},(k,k+1)\}\sigma_{k}^{-1}-(k,k+1)^{-1}(i,k+1)^{-1}(k,k+1)(jl)^{\varepsilon}\sigma_{k}^{-1}\\
&\equiv&\{(jl)^{\varepsilon},\{(k+1,l),(kl)^{-1}\}(il)\{(k+1,l)^{-1},(kl)^{-1}\}(il)^{-1}\}(k,k+1)^{-1}\\
&&(i,k+1)^{-1}(k,k+1)\sigma_{k}^{-1}\\
&&-\{(jl)^{\varepsilon},\{(k+1,l),(kl)^{-1}\}(il)\{(k+1,l)^{-1},(kl)^{-1}\}(il)^{-1}\}(k,k+1)^{-1}\\
&&(i,k+1)^{-1}(k,k+1)\sigma_{k}^{-1}\\
&\equiv&0.
\end{eqnarray*}
\end{enumerate}

$(6)\wedge(12)$\\
Let $f=\sigma_{k}^{-1}(ik)-\{(i,k+1),(k,k+1)\}\sigma_{k}^{-1},\
g=(ik)(jl)^{\varepsilon}-\{(jl)^{\varepsilon},(il)^{-1}(kl)^{-1}(il)(kl)\}(ik),i<j<k<l$.
Then $w=\sigma_{k}^{-1}(ik)(jl)^{\varepsilon}$ and
$$
(f,g)_w=\sigma_{k}^{-1}\{(jl)^{\varepsilon},(il)^{-1}(kl)^{-1}(il)(kl)\}(ik)-
\{(i,k+1),(k,k+1)\}\sigma_{k}^{-1}(jl)^{\varepsilon}.
$$
There are two cases to consider.
\begin{enumerate}
\item[1)]\ $k=l-1,i<j<k<l=k+1$.
\begin{eqnarray*}
&&\sigma_{k}^{-1}\{(j,k+1)^{\varepsilon},(i,k+1)^{-1}(k,k+1)^{-1}(i,k+1)(k,k+1)\}(ik)\\
&\equiv&\{(jk)^{\varepsilon},(ik)^{-1}(k,k+1)^{-1}(ik)(k,k+1)\}\{(i,k+1),(k,k+1)\}\sigma_{k}^{-1}\\
&\equiv&(k,k+1)^{-1}(ik)^{-1}(k,k+1)(ik)(jk)^{\varepsilon}(ik)^{-1}(k,k+1)^{-1}(ik)(i,k+1)(k,k+1)\sigma_{k}^{-1}\\
&\equiv&(k,k+1)^{-1}(i,k+1)(k,k+1)(i,k+1)^{-1}(jk)^{\varepsilon}(i,k+1)\sigma_{k}^{-1}\\
&\equiv&(k,k+1)^{-1}(i,k+1)(k,k+1)(jk)^{\varepsilon}\sigma_{k}^{-1} \ \ \ \ \  and \\
&&\\
&&\{(i,k+1),(k,k+1)\}\sigma_{k}^{-1}(j,k+1)^{\varepsilon}\\
&&\equiv(k,k+1)^{-1}(i,k+1)(k,k+1)(jk)^{\varepsilon}\sigma_{k}^{-1}.
\end{eqnarray*}
\item[2)]\ $k<l-1,i<j<k<k+1<l$.
\begin{eqnarray*}
&&\sigma_{k}^{-1}\{(jl)^{\varepsilon},(il)^{-1}(kl)^{-1}(il)(kl)\}(ik)\\
&\equiv&\{(jl)^{\varepsilon},(il)^{-1}\{(k+1,l)^{-1},(k,k+1)\}(il)\{(k+1,l),(k,k+1)\}\}\\
&&\{(i,k+1).(k,k+1)\}\sigma_{k}^{-1} \ \ \ \ \ \ and \\
&&\\
&&\{(i,k+1),(k,k+1)\}\sigma_{k}^{-1}(jl)^{\varepsilon}\\
&\equiv&(k,k+1)^{-1}(i,k+1)(k,k+1)(jl)^{\varepsilon}\sigma_{k}^{-1}\\
&\equiv&(k,k+1)^{-1}\{(jl)^{\varepsilon},(il)^{-1}(k+1,l)^{-1}(il)(k+1,l)\}(i,k+1)(k,k+1)\sigma_{k}^{-1}\\
&\equiv&\{(jl)^{\varepsilon},(il)^{-1}\{(k+1,l)^{-1},(k,k+1)\}(il)\{(k+1,l),(k,k+1)\}\}\\
&&\{(i,k+1).(k,k+1)\}\sigma_{k}^{-1}.
\end{eqnarray*}
\end{enumerate}

$(6)\wedge(13)$\\
Let
$f=\sigma_{k}^{-1}(ik)^\delta-\{(i,k+1)^\delta,(k,k+1)\}\sigma_{k}^{-1},\
g=(ik)^\delta(jl)^{\varepsilon}-(jl)^{\varepsilon}(ik)^\delta,j<i<k<l$.
Then $w=\sigma_{k}^{-1}(ik)^\delta(jl)^{\varepsilon}$ and
$$
(f,g)_w=\sigma_{k}^{-1}(jl)^{\varepsilon}(ik)^\delta-\{(i,k+1)^\delta,(k,k+1)\}\sigma_{k}^{-1}(jl)^{\varepsilon}.
$$
There are two cases to consider.
\begin{enumerate}
\item[1)]\ $k=l-1 \ ,\ \ j<i<k<l=k+1$.
$$
\sigma_{k}^{-1}(j,k+1)^{\varepsilon}(ik)^\delta\equiv(jk)^{\varepsilon}\{(i,k+1)^\delta,(k,k+1)\}\sigma_{k}^{-1}.
$$
If $\varepsilon=1$, then
\begin{eqnarray*}
&&\sigma_{k}^{-1}(j,k+1)^{\varepsilon}(ik)^\delta\\
&\equiv&\{\{(i,k+1)^\delta,(j,k+1)^{-1}(k,k+1)^{-1}(j,k+1)(k,k+1)\},\\
&&\{(k,k+1),(j,k+1)(k,k+1)\}\}(jk)^{\varepsilon}\sigma_{k}^{-1}\\
&\equiv&(k,k+1)^{-1}(i,k+1)^\delta(k,k+1)(jk)^{\varepsilon}\sigma_{k}^{-1}.
\end{eqnarray*}
If $\varepsilon=-1$, then
\begin{eqnarray*}
&&\sigma_{k}^{-1}(j,k+1)^{\varepsilon}(ik)^\delta\\
&\equiv&\{\{(i,k+1)^\delta,(k,k+1)(j,k+1)(k,k+1)^{-1}(j,k+1)^{-1}\},\{(k,k+1),(j,k+1)^{-1}\}\}\\
&&(jk)^{\varepsilon}\sigma_{k}^{-1}\\
&\equiv&(k,k+1)^{-1}(i,k+1)^\delta(k,k+1)(jk)^{\varepsilon}\sigma_{k}^{-1}.\\
\end{eqnarray*}
On the other hand,
$$
\{(i,k+1)^\delta,(k,k+1)\}\sigma_{k}^{-1}(j,k+1)^{\varepsilon}\equiv(k,k+1)^{-1}(i,k+1)^\delta
(k,k+1)(jk)^{\varepsilon}\sigma_{k}^{-1}.
$$
\item[2)]\ $k<l-1,j<i<k<k+1<l$.
\begin{eqnarray*}
(f,g)_w&\equiv&(jl)^{\varepsilon}\{(i,k+1)^\delta,(k,k+1)\}\sigma_{k}^{-1}-\\
&&(k,k+1)^{-1}(i,k+1)^\delta(k,k+1)(jl)^{\varepsilon}\sigma_{k}^{-1}\\
&\equiv&(jl)^{\varepsilon}(k,k+1)^{-1}(i,k+1)^\delta(k,k+1)\sigma_{k}^{-1}-\\
&&(jl)^{\varepsilon}(k,k+1)^{-1}(i,k+1)^\delta(k,k+1)\sigma_{k}^{-1}\\
&\equiv&0.
\end{eqnarray*}
\end{enumerate}

$(6)\wedge(13)$\\
Let
$f=\sigma_{k}^{-1}(ik)^\delta-\{(i,k+1)^\delta,(k,k+1)\}\sigma_{k}^{-1},\
g=(ik)^\delta(jl)^{\varepsilon}-(jl)^{\varepsilon}(ik)^\delta,i<k<j<l$.
Then $w=\sigma_{k}^{-1}(ik)^\delta(jl)^{\varepsilon}$ and
$$
(f,g)_w=\sigma_{k}^{-1}(jl)^{\varepsilon}(ik)^\delta-\{(i,k+1)^\delta,(k,k+1)\}\sigma_{k}^{-1}(jl)^{\varepsilon}.
$$
There are two cases to consider.
\begin{enumerate}
\item[1)]\ $k=j-1,i<k<k+1=j<l$.
\begin{eqnarray*}
&&\sigma_{k}^{-1}(k+1,l)^{\varepsilon}(ik)^\delta\\
&\equiv&(kl)^{\varepsilon}\{(i,k+1)^\delta,(k,k+1)\}\sigma_{k}^{-1}\\
&\equiv&(kl)^{\varepsilon}(k,k+1)^{-1}(i,k+1)^\delta(k,k+1)\sigma_{k}^{-1}
\ \ \ \  and \\
&&\\
&&\{(i,k+1)^\delta,(k,k+1)\}\sigma_{k}^{-1}(k+1,l)^{\varepsilon}\\
&\equiv&(k,k+1)^{-1}(i,k+1)^\delta(k,k+1)(kl)^{\varepsilon}\sigma_{k}^{-1}\\
&\equiv&(k,k+1)^{-1}(i,k+1)^\delta\{(kl)^{\varepsilon},(k+1,l)\}(k,k+1)\sigma_{k}^{-1}.
\end{eqnarray*}
If $\varepsilon=1$, then
\begin{eqnarray*}
&&\{(i,k+1)^\delta,(k,k+1)\}\sigma_{k}^{-1}(k+1,l)^{\varepsilon}\\
&\equiv&(k,k+1)^{-1}\{\{(kl)^{\varepsilon},(il)^{-1}(k+1,l)^{-1}(il)(k+1,l)\},\{(k+1,l),(il)(k+1,l)\}\}\\
&&(i,k+1)^\delta(k,k+1)\sigma_{k}^{-1}\\
&\equiv&(k,k+1)^{-1}\{(kl)^{\varepsilon},(k+1,l)\}(i,k+1)^\delta(k,k+1)\sigma_{k}^{-1}\\
&\equiv&\{\{(kl)^{\varepsilon},(k+1,l)^{-1}(kl)^{-1}\},\{(k+1,l),(kl)^{-1}\}\}(k,k+1)^{-1}(i,k+1)^\delta(k,k+1)\sigma_{k}^{-1}\\
&\equiv&(kl)^{\varepsilon}(k,k+1)^{-1}(i,k+1)^\delta(k,k+1)\sigma_{k}^{-1}.
\end{eqnarray*}
If $ \varepsilon=-1$, then
\begin{eqnarray*}
&&\{(i,k+1)^\delta,(k,k+1)\}\sigma_{k}^{-1}(k+1,l)^{\varepsilon}\\
&\equiv&(k,k+1)^{-1}\{\{(kl)^{\varepsilon},(k+1,l)(il)(k+1,l)^{-1}(il)^{-1}\},\{(k+1,l),(il)^{-1}\}\}\\
&&(i,k+1)^\delta(k,k+1)\sigma_{k}^{-1}\\
&\equiv&(k,k+1)^{-1}\{(kl)^{\varepsilon},(k+1,l)\}(i,k+1)^\delta(k,k+1)\sigma_{k}^{-1}\\
&\equiv&(kl)^{\varepsilon}(k,k+1)^{-1}(i,k+1)^\delta(k,k+1)\sigma_{k}^{-1}.
\end{eqnarray*}
\item[2)]\ $k<j-1,i<k<k+1<j<l$.
$$
(f,g)_{w}\equiv(jl)^{\varepsilon}\{(i,k+1)^\delta,(k,k+1)\}\sigma_{k}^{-1}-
(jl)^{\varepsilon}\{(i,k+1)^\delta,(k,k+1)\}\sigma_{k}^{-1}\equiv0.
$$
\end{enumerate}

$(7)\wedge(7)$\\
Let $f=(ij)^{-1}(jk)^{-1}-\{(jk)^{-1},(ik)^{-1}\}(ij)^{-1},\
g=(jk)^{-1}(kl)^{\varepsilon}-\{(kl)^{\varepsilon},(jl)^{-1}\}(jk)^{-1},i<j<k<l$.
Then $w=(ij)^{-1}(jk)^{-1}(kl)^{\varepsilon}$ and
$$
(f,g)_w=(ij)^{-1}\{(kl)^{\varepsilon},(jl)^{-1}\}(jk)^{-1}-\{(jk)^{-1},(ik)^{-1}\}(ij)^{-1}(kl)^{\varepsilon}\equiv
0
$$
since
\begin{eqnarray*}
&&(ij)^{-1}\{(kl)^{\varepsilon},(jl)^{-1}\}(jk)^{-1}\\
&\equiv&\{(kl)^{\varepsilon},\{(jl)^{-1},(il)^{-1}\}\}\{(jk)^{-1},(ik)^{-1}\}(ij)^{-1} \ \ \ \ \ \ \ and \\
&&\\
&&\{(jk)^{-1},(ik)^{-1}\}(ij)^{-1}(kl)^{\varepsilon}\\
&\equiv&(ik)(jk)^{-1}(ik)^{-1}(ij)^{-1}(kl)^{\varepsilon}\\
&\equiv&(ik)(jk)^{-1}\{(kl)^{\varepsilon},(il)^{-1}\}(ik)^{-1}(ij)^{-1}\\
&\equiv&(ik)\{\{(kl)^{\varepsilon},(jl)^{-1}\},(il)^{-1}\}(jk)^{-1}(ik)^{-1}(ij)^{-1}\\
&\equiv&\{\{\{(kl)^{\varepsilon},(il)(kl)\},\{(jl)^{-1},(il)^{-1}(kl)^{-1}(il)(kl)\}\},\{(il)^{-1}(kl)\}\}\\
&&(ik)(jk)^{-1}(ik)^{-1}(ij)^{-1}\\
&\equiv&\{(kl)^{\varepsilon},\{(jl)^{-1},(il)^{-1}\}\{(jk)^{-1},(ik)^{-1}\}(ij)^{-1}.
\end{eqnarray*}

$(7)\wedge(8)$\\
Let $f=(ij)^{-1}(jk)-\{(jk),(ik)^{-1}\}(ij)^{-1},\
g=(jk)(kl)^{\varepsilon}-\{(kl)^{\varepsilon},(jl)(kl)\}(jk),
i<j<k<l$. Then $w=(ij)^{-1}(jk)(kl)^{\varepsilon}$ and
$$
(f,g)_w=(ij)^{-1}\{(kl)^{\varepsilon},(jl)(kl)\}(jk)-\{(jk),(ik)^{-1}\}(ij)^{-1}(kl)^{\varepsilon}\equiv
0
$$
since
\begin{eqnarray*}
&&(ij)^{-1}\{(kl)^{\varepsilon},(jl)(kl)\}(jk)\\
&\equiv&\{(kl)^{\varepsilon},\{(jl),(il)^{-1}\}(kl)\}\{(jk),(ik)^{-1}\}(ij)^{-1} \ \ \ \ \ \ and\\
&&\\
&&\{(jk),(ik)^{-1}\}(ij)^{-1}(kl)^{\varepsilon}\\
&\equiv&(ik)(jk)\{(kl)^{\varepsilon},(il)^{-1}\}(ik)^{-1}(ij)^{-1}\\
&\equiv&(ik)\{\{(kl)^{\varepsilon},(jl)(kl)\},(il)^{-1}\}(jk)(ik)^{-1}(ij)^{-1}\\
&\equiv&\{\{\{(kl)^{\varepsilon},(il)(kl)\},\{(jl),(il)^{-1}(kl)^{-1}(il)(kl)\}\\
&&\{(kl),(il)(kl)\}\},\{(il)^{-1},(kl)^{-1}(il)^{-1}\}\}(ik)(jk)(ik)^{-1}(ij)^{-1}\\
&\equiv&\{(kl)^{\varepsilon},\{(jl),(il)^{-1}\}(kl)\}\{(jk),(ik)^{-1}\}(ij)^{-1}.
\end{eqnarray*}

$(7)\wedge(9)$\\
Let $f=(ij)^{-1}(jk)^{-1}-\{(jk)^{-1},(ik)^{-1}\}(ij)^{-1},\
g=(jk)^{-1}(jl)^{\varepsilon}-\{(jl)^{\varepsilon},(kl)^{-1}(jl)^{-1}\}(jk)^{-1},i<j<k<l$.
Then $w=(ij)^{-1}(jk)^{-1}(jl)^{\varepsilon}$ and
$$
(f,g)_w=(ij)^{-1}\{(jl)^{\varepsilon},(kl)^{-1}(jl)^{-1}\}(jk)^{-1}-\{(jk)^{-1},(ik)^{-1}\}(ij)^{-1}(jl)^{\varepsilon}\equiv
0
$$
since
\begin{eqnarray*}
&&(ij)^{-1}\{(jl)^{\varepsilon},(kl)^{-1}(jl)^{-1}\}(jk)^{-1}\\
&\equiv&\{(jl)^{\varepsilon},(il)^{-1}(kl)^{-1}\{(jl)^{-1},(il)^{-1}\}\}\{(jk)^{-1},(ik)^{-1}\}(ij)^{-1} \ \ \ \ \ \ and\\
&&\\
&&\{(jk)^{-1},(ik)^{-1}\}(ij)^{-1}(jl)^{\varepsilon}\\
&\equiv&(ik)(jk)(ik)^{-1}\{(jl)^{\varepsilon},(il)^{-1}\}(ij)^{-1}\\
&\equiv&(ik)(jk)^{-1}\{\{(jl)^{\varepsilon},(kl)(il)(kl)^{-1}(il)^{-1}\},\{(il)^{-1},(kl)^{-1}(il)^{-1}\}\}
(ik)^{-1}(ij)^{-1}\\
&\equiv&(ik)(jk)^{-1}\{(jl)^{\varepsilon},(il)^{-1}\}(ik)^{-1}(ij)^{-1}\\
&\equiv&(ik)\{\{(jl)^{\varepsilon},(kl)^{-1}(jl)^{-1}\},(il)^{-1}\}(jk)^{-1}(ik)^{-1}(ij)^{-1}\\
&\equiv&\{\{\{(jl)^{\varepsilon},(il)^{-1}(kl)^{-1}(il)(kl)\},\{(kl)^{-1},(il)(kl)\}\\
&&\{(jl)^{-1},(il)^{-1}(kl)^{-1}(il)(kl)\}\},\{(il)^{-1}(kl)\}\}(ik)(jk)^{-1}(ik)^{-1}(ij)^{-1}\\
&\equiv&\{(jl)^{\varepsilon},(il)^{-1}(kl)^{-1}\{(jl)^{-1},(il)^{-1}\}\}\{(jk)^{-1},(ik)^{-1}\}(ij)^{-1}.
\end{eqnarray*}

$(7)\wedge(10)$\\
Let $f=(ij)^{-1}(jk)-\{(jk),(ik)^{-1}\}(ij)^{-1},\
g=(jk)(jl)^{\varepsilon}-\{(jl)^{\varepsilon},(kl)\}(jk),i<j<k<l$.
Then $ w=(ij)^{-1}(jk)(jl)^{\varepsilon}$ and
$$
(f,g)_w=(ij)^{-1}\{(jl)^{\varepsilon},(kl)\}(jk)-\{(jk),(ik)^{-1}\}(ij)^{-1}(jl)^{\varepsilon}\equiv
0
$$
since
\begin{eqnarray*}
&&(ij)^{-1}\{(jl)^{\varepsilon},(kl)\}(jk)\\
&\equiv&\{\{(jl)^{\varepsilon},(il)^{-1}\},(kl)\}\{(jk),(ik)^{-1}\}(ij)^{-1}
\ \ \ \ \ \ \ and\\
&&\\
&&\{(jk),(ik)^{-1}\}(ij)^{-1}(jl)^{\varepsilon}\\
&\equiv&(ik)(jk)\{(jl)^{\varepsilon},(il)^{-1}\}(ik)^{-1}(ij)^{-1}\\
&\equiv&(ik)\{\{(jl)^{\varepsilon},(kl)\},(il)^{-1}\}(jk)(ik)^{-1}(ij)^{-1}\\
&\equiv&\{\{\{(jl)^{\varepsilon},(il)^{-1}(kl)^{-1}(il)(kl)\},\{(kl),(il)(kl)\}\},\{(il)^{-1}(kl)\}\}\\
&&(ik)(jk)(ik)^{-1}(ij)^{-1}\\
&\equiv&\{\{(jl)^{\varepsilon},(il)^{-1}\},(kl)\}\{(jk),(ik)^{-1}\}(ij)^{-1}.
\end{eqnarray*}

$(7)\wedge(11)$\\
Let $f=(pi)^{-1}(ik)^{-1}-\{(ik)^{-1},(pk)^{-1}\}(pi)^{-1},\
g=(ik)^{-1}(jl)^{\varepsilon}-\{(jl)^{\varepsilon},(kl)(il)(kl)^{-1}(il)^{-1}\}(ik)^{-1},
p<i<j<k<l$. Then $w=(pi)^{-1}(ik)^{-1}(jl)^{\varepsilon}$ and
$$
(f,g)_w=(pi)^{-1}\{(jl)^{\varepsilon},(kl)(il)(kl)^{-1}(il)^{-1}\}(ik)^{-1}-
\{(ik)^{-1},(pk)^{-1}\}(pi)^{-1}(jl)^{\varepsilon}\equiv 0
$$
since
\begin{eqnarray*}
&&(pi)^{-1}\{(jl)^{\varepsilon},(kl)(il)(kl)^{-1}(il)^{-1}\}(ik)^{-1}\\
&\equiv&\{(jl)^{\varepsilon},(kl)\{(il),(pl)^{-1}\}(kl)^{-1}\{(il)^{-1},(pl)^{-1}\}\}\{(ik)^{-1},(pk)^{-1}\}(pi)^{-1} \ \ \ \ \ \ \ and \\
&&\\
&&\{(ik)^{-1},(pk)^{-1}\}(pi)^{-1}(jl)^{\varepsilon}\\
&\equiv&(pk)(ik)^{-1}(pk)^{-1}(jl)^{\varepsilon}(pi)^{-1}\\
&\equiv&(pk)(ik)^{-1}\{(jl)^{\varepsilon},(kl)(pl)(kl)^{-1}(pl)^{-1}\}(pk)^{-1}(pi)^{-1}\\
&\equiv&(pk)\{\{(jl)^{\varepsilon},(kl)(il)(kl)^{-1}(il)^{-1}\},\{(kl)(il)^{-1}\}(pl)\\
&&\{(kl)^{-1},(il)^{-1}\}(pl)^{-1}\}(ik)^{-1}(pk)^{-1}(pi)^{-1}\\
&\equiv&(pk)\{(jl)^{\varepsilon},(kl)\{(kl)^{-1}(il)^{-1}(pl)^{-1}\}\}(ik)^{-1}(pk)^{-1}(pi)^{-1}\\
&\equiv&\{\{(jl)^{\varepsilon},(pl)^{-1}(kl)^{-1}(pl)(kl)\},\{(kl),(pl)(kl)\}\\
&&\{\{(kl)^{-1},(pl)(kl)\},\{(il)^{-1},(pl)^{-1}(kl)^{-1}(pl)(kl)\}\{(pl)^{-1},(kl)\}\}(pk)(ik)^{-1}(pk)^{-1}(pi)^{-1}\\
&\equiv&\{(jl)^{\varepsilon},(kl)\{(il),(pl)^{-1}\}(kl)^{-1}\{(il)^{-1},(pl)^{-1}\}\}\{(ik)^{-1},(pk)^{-1}\}(pi)^{-1}.
\end{eqnarray*}

$(7)\wedge(12)$\\
Let $f=(pi)^{-1}(ik)-\{(ik),(pk)^{-1}\}(pi)^{-1},\
g=(ik)(jl)^{\varepsilon}-\{(jl)^{\varepsilon},(il)^{-1}(kl)^{-1}(il)(kl)\}(ik),p<i<j<k<l$.
Then $w=(pi)^{-1}(ik)(jl)^{\varepsilon}$ and
$$
(f,g)_w=(pi)^{-1}\{(jl)^{\varepsilon},(il)^{-1}(kl)^{-1}(il)(kl)\}(ik)-\{(ik),(pk)^{-1}\}(pi)^{-1}(jl)^{\varepsilon}\equiv
0
$$
since
\begin{eqnarray*}
&&(pi)^{-1}\{(jl)^{\varepsilon},(il)^{-1}(kl)^{-1}(il)(kl)\}(ik)\\
&\equiv&\{(jl)^{\varepsilon},\{(il)^{-1},(pl)^{-1}\}(kl)^{-1}\{(il),(pl)^{-1}\}(kl)\}\{(ik),(pk)^{-1}\}(pi)^{-1} \ \ \ \ and\\
&&\\
&&\{(ik),(pk)^{-1}\}(pi)^{-1}(jl)^{\varepsilon}\\
&\equiv&(pk)(ik)\{(jl)^{\varepsilon},(kl)(pl)(kl)^{-1}(pl)^{-1}\}(pk)^{-1}(pi)^{-1}\\
&\equiv&(pk)\{\{(jl)^{\varepsilon},(il)^{-1}(kl)^{-1}(il)(kl)\},\{(kl),(il)(kl)\}(pl)\\
&&\{(kl)^{-1},(il)(kl)\}(pl)^{-1}\}(ik)(pk)^{-1}(pi)^{-1}\\
&\equiv&(pk)\{(jl)^{\varepsilon},(kl)\{(kl)^{-1},(il)(kl)(pl)^{-1}\}\}(ik)(pk)^{-1}(pi)^{-1}\\
&\equiv&\{\{(jl)^{\varepsilon},(pl)^{-1}(kl)^{-1}(pl)(kl)\},\{(kl),(pl)(kl)\}\{\{(kl)^{-1},(pl)(kl)\},\\
&&\{(il),(pl)^{-1}(kl)^{-1}(pl)(kl)\}\{(kl),(pl)(kl)\}\{(pl)^{-1},(kl)\}\}\}(pk)(ik)(pk)^{-1}(pi)^{-1}\\
&\equiv&\{(jl)^{\varepsilon},\{(il)^{-1},(pl)^{-1}\}(kl)^{-1}\{(il),(pl)^{-1}\}(kl)\}\{(ik),(pk)^{-1}\}(pi)^{-1}.
\end{eqnarray*}

$(7)\wedge(13)$\\
$f=(pi)^{-1}(ik)^\delta-\{(ik)^\delta,(pk)^{-1}\}(pi)^{-1},
g=(ik)^\delta(jl)^{\varepsilon}-(jl)^{\varepsilon}(ik)^\delta.$\
Then $w=(pi)^{-1}(ik)^\delta(jl)^{\varepsilon}$\ and
$$
(f,g)_w=(pi)^{-1}(jl)^{\varepsilon}(ik)^\delta-\{(ik)^\delta,(pk)^{-1}\}(pi)^{-1}(jl)^{\varepsilon}.
$$
There are two cases to consider.
\begin{enumerate}
\item[1)]\ $p<i<k,j<i<k<l.$\ In this case, there are three subcases to consider.
\begin{enumerate}
\item[a)]\ $p<j,p<j<i<k<l$.
\begin{eqnarray*}
&&(pi)^{-1}(jl)^{\varepsilon}(ik)^\delta\\
&\equiv&\{(jl)^{\varepsilon},(il)(pl)(il)^{-1}(pl)^{-1}\}\{(ik)^\delta,(pk)^{-1}\}(pi)^{-1}\ \ \ and\\
&&\\
&&\{(ik)^\delta,(pk)^{-1}\}(pi)^{-1}(jl)^{\varepsilon}\\
&\equiv&(pk)(ik)^\delta(pk)^{-1}\{(jl)^{\varepsilon},(il)(pl)(il)^{-1}(pl)^{-1}\}(pi)^{-1}\\
&\equiv&(pk)(ik)^\delta\{\{(jl)^{\varepsilon},(kl)(pl)(kl)^{-1}(pl)^{-1}\},\{(il),(kl)(pl)(kl)^{-1}(pl)^{-1}\}\\
&&\{(pl),(kl)^{-1}(pl)^{-1}\}\{(il)^{-1},(kl)(pl)(kl)^{-1}(pl)^{-1}\}\{(pl)^{-1},(kl)^{-1}\\
&&(pl)^{-1}\}(pk)^{-1}(pi)^{-1}\\
&\equiv&(pk)(ik)^\delta\{(jl)^{\varepsilon},(il)(kl)(pl)(kl)^{-1}(il)^{-1}(pl)^{-1}\}(pk)^{-1}(pi)^{-1}.
\end{eqnarray*}
If $\delta=1$, then
\begin{eqnarray*}
&&\{(ik)^\delta,(pk)^{-1}\}(pi)^{-1}(jl)^{\varepsilon}\\
&\equiv&(pk)\{(jl)^{\varepsilon},\{(il),(kl)\}\{(kl),(il)(kl)\}(pl)\{(kl)^{-1},(il)(kl)\}\}\\
&&(il)^{-1},(kl)\}(pl)^{-1}\}(ik)^\delta(pk)^{-1}(pi)^{-1}\\
&\equiv&(pk)\{(jl)^{\varepsilon},(il)(kl)(pl)(kl)^{-1}(il)^{-1}(pl)^{-1}\}(ik)^\delta(pk)^{-1}(pi)^{-1}\\
&\equiv&\{\{(jl)^{\varepsilon},(pl)^{-1}(kl)^{-1}(pl)(kl)\},\{\{(pl),(kl)\},\{(kl)^{-1},(pl)(kl)\}\\
&&\{(il)^{-1},(pl)^{-1}(kl)^{-1}(pl)(kl)\}\}\{(pl),(kl)\}\}\{(ik)^\delta,(pk)^{-1}\}(pi)^{-1}\\
&\equiv&\{(jl)^{\varepsilon},(il)(pl)(il)^{-1}(pl)^{-1}\}\{(ik)^\delta,(pk)^{-1}\}(pi)^{-1}.
\end{eqnarray*}
If $\delta=-1$, then
\begin{eqnarray*}
&&\{(ik)^\delta,(pk)^{-1}\}(pi)^{-1}(jl)^{\varepsilon}\\
&\equiv&(pk)\{(jl)^{\varepsilon},\{(il),(kl)^{-1}(il)^{-1}\}\{(kl),(il)^{-1}\}(pl)\{(kl)^{-1},(il)^{-1}\}\\
&&\{(il)^{-1},(kl)^{-1}(il)^{-1}\}(pl)^{-1}\}(ik)^\delta(pk)^{-1}(pi)^{-1}\\
&\equiv&(pk)\{(jl)^{\varepsilon},(il)(kl)(pl)(kl)^{-1}(il)^{-1}(pl)^{-1}(ik)^\delta(pk)^{-1}(pi)^{-1}\\
&\equiv&\{(jl)^{\varepsilon},(il)(pl)(il)^{-1}(pl)^{-1}\}\{(ik)^\delta,(pk)^{-1}\}(pi)^{-1}.
\end{eqnarray*}
\item[b)]\ $p=j,p=j<i<k<l$.
\begin{eqnarray*}
&&(pi)^{-1}(jl)^{\varepsilon}(ik)^\delta\\
&=&(ji)^{-1}(jl)^{\varepsilon}(ik)^\delta\\
&\equiv&\{(jl)^{\varepsilon},(il)^{-1}(jl)^{-1}\}\{(ik)^\delta,(jk)^{-1}\}(ji)^{-1}\ \ \ and \\
&&\\
&&\{(ik)^\delta,(pk)^{-1}\}(pi)^{-1}(jl)^{\varepsilon}\\
&\equiv&(jk)(ik)^\delta(jk)^{-1}(ji)^{-1}(jl)^{\varepsilon}\\
&\equiv&(jk)(ik)^\delta(jk)^{-1}\{(jl)^{\varepsilon},(il)^{-1}(jl)^{-1}\}(ji)^{-1}\\
&\equiv&(jk)(ik)^\delta\{\{(jl)^{\varepsilon},(kl)^{-1}(jl)^{-1}\},\{(il)^{-1},(kl)(jl)(kl)^{-1}jl)^{-1}\}\\
&&\{(jl)^{-1},(kl)^{-1}(jl)^{-1}\}(jk)^{-1}(jl)^{-1}\\
&\equiv&(jk)(ik)^\delta\{(jl)^{\varepsilon},(kl)^{-1}(il)^{-1}(jl)^{-1}\}(jk)^{-1}(ji)^{-1}.
\end{eqnarray*}
If $\delta=1$, then
\begin{eqnarray*}
&&\{(ik)^\delta,(pk)^{-1}\}(pi)^{-1}(jl)^{\varepsilon}\\
&\equiv&(jk)\{(jl)^{\varepsilon},\{(kl)^{-1},(il)(kl)\}\{(il)^{-1},(kl)\}(jl)^{-1}\}(ik)^\delta(jk)^{-1}(ji)^{-1}\\
&\equiv&(jk)\{(jl)^{\varepsilon},(kl)^{-1}(il)^{-1}(jl)^{-1}\}(ik)^\delta(jk)^{-1}(ji)^{-1}\\
&\equiv&\{\{(jl)^{\varepsilon},(kl)\}\{(kl)^{-1},(jl)(kl)\}\{(il)^{-1},(jl)^{-1}(kl)^{-1}(ji)(kl)\}\\
&&\{(jl)^{-1},(kl)\}\}(jk)(ik)^\delta(jk)^{-1}(ji)^{-1}\\
&\equiv&\{(jl)^{\varepsilon},(il)^{-1}(jl)^{-1}\}\{(ik)^\delta,(jk)^{-1}\}(ji)^{-1}.
\end{eqnarray*}
If $\delta=-1$, then
\begin{eqnarray*}
&&\{(ik)^\delta,(pk)^{-1}\}(pi)^{-1}(jl)^{\varepsilon}\\
&\equiv&(jk)\{(jl)^{\varepsilon},\{(kl)^{-1},(il)^{-1}\}\{(il)^{-1},(kl)^{-1}(il)^{-1}\}(jl)^{-1}\}
(ik)^\delta(jk)^{-1}(ji)^{-1}\\
&\equiv&(jk)\{(jl)^{\varepsilon},(kl)^{-1}(il)^{-1}(jl)^{-1}\}(ik)^\delta(jk)^{-1}(ji)^{-1}\\
&\equiv&\{(jl)^{\varepsilon},(il)^{-1}(jl)^{-1}\}\{(ik)^\delta,(jk)^{-1}\}(ji)^{-1}.
\end{eqnarray*}
\item[c)]\ $p>j,j<p<i<k<l$.
\begin{eqnarray*}
(f,g)_w&\equiv&(jl)^{\varepsilon}\{(ik)^\delta,(pk)^{-1}\}(pi)^{-1}-(pk)(ik)^\delta(pk)^{-1}(pi)^{-1}(jl)^{\varepsilon}\\
&\equiv&(jl)^{\varepsilon}(pk)(ik)^\delta(pk)^{-1}(pi)^{-1}-(jl)^{\varepsilon}(pk)(ik)^\delta(pk)^{-1}(pi)^{-1}\\
&\equiv&0.
\end{eqnarray*}
\end{enumerate}
\item[2)]\ $p<i<k<j<l$.
\begin{eqnarray*}
(f,g)_w&=&(pi)^{-1}(jl)^{\varepsilon}(ik)^\delta-\{(ik)^\delta,(pk)^{-1}\}(pi)^{-1}(jl)^{\varepsilon}\\
&\equiv&(jl)^{\varepsilon}(pi)^{-1}(ik)^\delta-(jl)^{\varepsilon}\{(ik)^\delta,(pk)^{-1}\}(pi)^{-1}\\
&\equiv&0.
\end{eqnarray*}
\end{enumerate}

$(8)\wedge(7)$\\
Let $f=(ij)(jk)^{-1}-\{(jk)^{-1},(ik)(jk)\}(ij),\
g=(jk)^{-1}(kl)^{\varepsilon}-\{(kl)^{\varepsilon},(jl)^{-1}\}(jk)^{-1},i<j<k<l$.
Then $w=(ij)(jk)^{-1}(kl)^{\varepsilon}$ and
$$
(f,g)_w=(ij)\{(kl)^{\varepsilon},(jl)^{-1}\}(jk)^{-1}-\{(jk)^{-1},(ik)(jk)\}(ij)(kl)^{\varepsilon}\equiv
0
$$
since
\begin{eqnarray*}
&&(ij)\{(kl)^{\varepsilon},(jl)^{-1}\}(jk)^{-1}\\
&\equiv&\{(kl)^{\varepsilon},\{(jl)^{-1},(il)(jl)\}\}\{(jk)^{-1},(ik)(jk)\}(ij)\ \ \ \ \ and \\
&&\\
&&\{(jk)^{-1},(ik)(jk)\}(ij)(kl)^{\varepsilon}\\
&\equiv&(jk)^{-1}(ik)^{-1}(jk)^{-1}(ik)(jk)(kl)^{\varepsilon}(ij)\\
&\equiv&(jk)^{-1}(ik)^{-1}(jk)^{-1}(ik)\{(kl)^{\varepsilon},(jl)(kl)\}(jk)(ij)\\
&\equiv&(jk)^{-1}(ik)^{-1}(jk)^{-1}\{\{(kl)^{\varepsilon},(il)(kl)\},\{(jl),(il)^{-1}(kl)^{-1}(il)(kl)\}\\
&&\{(kl),(il)(kl)\}\}(ik)(jk)(ij)\\
&\equiv&(jk)^{-1}(ik)^{-1}\{\{(kl)^{\varepsilon},(jl)^{-1}\},(il)\{(jl),(kl)^{-1}(jl)^{-1}\}\{(kl),(jl)^{-1}\}\}\\
&&(jk)^{-1}(ik)(jk)(ij)\\
&\equiv&(jk)^{-1}\{\{\{(kl)^{\varepsilon},(il)^{-1}\},\{(jl)^{-1},(kl)(il)(kl)^{-1}(il)^{-1}\}\},\\
&&\{(il),(kl)^{-1}(il)^{-1}\}\{(jl),(kl)(il)(kl)^{-1}(il)^{-1}\}\{(kl),(il)^{-1}\}\}(ik)^{-1}(jk)^{-1}(ik)(jk)(ij)\\
&\equiv&\{\{(kl)^{\varepsilon},(jl)^{-1}\},(il)^{-1}\{(kl)^{-1},(jl)^{-1}\}
\{(jl)^{-1},(kl)^{-1}(jl)^{-1}\}\{(kl),(jl)^{-1}\}(il)\\
&&\{(kl)^{-1},(jl)^{-1}\}\{(jl),(kl)^{-1}(jl)^{-1}\}\{(kl),(jl)^{-1}\}\}(jk)^{-1}(ik)^{-1}(jk)^{-1}(ik)(jk)(ij)\\
&\equiv&\{(kl)^{\varepsilon},\{(jl)^{-1},(il)(jl)\}\}\{(jk)^{-1},(ik)(jk)\}(ij).
\end{eqnarray*}

$(8)\wedge (8)$\\
Let $f=(ij)(jk)-\{(jk),(ik)(jk)\}(ij),\
g=(jk)(kl)^{\varepsilon}-\{(kl)^{\varepsilon},(jl)(kl)\}(jk),\ \
i<j<k<l$. Then $w=(ij)(jk)(kl)^{\varepsilon}$ and
$$
(f,g)_{w}=(ij)\{(kl)^{\varepsilon},(jl)(kl)\}(jk)-\{(jk),(ik)(jk)\}(ij)(kl)\equiv
0
$$
since
\begin{eqnarray*}
&&(ij)\{(kl)^{\varepsilon},(jl)(kl)\}(jk)\\
&\equiv& \{(kl)^{\varepsilon},\{(jl),(il)(jl)\}(kl)\}\{jk,(ik)(jk)\}(ij) \ \ \ \ \ \ and\\
&&\\
&&\{(jk),(ik)(jk)\}(ij)(kl)\\
&\equiv&(jk)^{-1}(ik)^{-1}(jk)(ik)(jk)(kl)(ij)\\
&\equiv&(jk)^{-1}(ik)^{-1}(jk)(ik)\{(kl),(jl)(kl)\}(jk)(ij)\\
&\equiv&(jk)^{-1}(ik)^{-1}(jk)\{\{(kl),(il)(kl)\},\{(jl),(il)^{-1}(kl)^{-1}(il)(kl)\}\{(kl),(il)(kl)\}\}(ik)(jk)(ij)\\
&\equiv&(jk)^{-1}(ik)^{-1}\{\{(kl),(jl)(kl)\},(il)\{(jl),(kl)\}\{(kl),(jl)(kl)\}\}(jk)(ik)(jk)(ij)\\
&\equiv&(jk)^{-1}\{\{(kl),(il)^{-1}\},\{(jl),(kl)(il)(kl)^{-1}(il)^{-1}\}\{(kl),(il)^{-1}\}\{(il),(kl)^{-1}(il)^{-1}\}\\
&&\{(jl),(kl)(il)(kl)^{-1}(il)^{-1}\}\{(kl),(il)^{-1}\}\}(ik)^{-1}(jk)(ik)(jk)(ij)\\
&\equiv&\{\{(kl),(jl)^{-1}\},\{\{(jl),(kl)^{-1}(jl)^{-1}\},\{(kl),(jl)^{-1}\}(il)\}\{(jl),(kl)^{-1}(jl)^{-1}\}\\
&&\{(kl),(jl)^{-1}\}\}\{(jk),(ik)(jk)\}(ij)\\
&\equiv&\{(kl)^{\varepsilon},\{(jl),(il)(jl)\}(kl)\}\{jk,(ik)(jk)\}(ij).
\end{eqnarray*}

$(8)\wedge(9)$\\
Let $f=(ij)(jk)^{-1}-\{(jk)^{-1},(ik)(jk)\}(ij),\
g=(jk)^{-1}(jl)^{\varepsilon}-\{(jl)^{\varepsilon},(kl)^{-1}(jl)^{-1}\}(jk)^{-1},\
i<j<k<l$. Then $w=(ij)(jk)^{-1}(jl)^{\varepsilon}$ and
$$
(f,g)_w=(ij)\{(jl)^{\varepsilon},(kl)^{-1}(jl)^{-1}\}(jk)^{-1}-\{(jk)^{-1},(ik)(jk)\}(ij)(jl)^{\varepsilon}\equiv
0
$$
since
\begin{eqnarray*}
&&(ij)\{(jl)^{\varepsilon},(kl)^{-1}(jl)^{-1}\}(jk)^{-1}\\
&\equiv&\{(jl)^{\varepsilon},(il)(jl)(kl)^{-1}\{(jl)^{-1},(il)(jl)\}\}\{(jk)^{-1},(ik)(jk)\}(ij)\ \ \ \ \ and \\
&&\\
&&\{(jk)^{-1},(ik)(jk)\}(ij)(jl)^{\varepsilon}\\
&\equiv&(jk)^{-1}(ik)^{-1}(jk)^{-1}(ik)(jk)\{(jl)^{\varepsilon},(il)(jl)\}(ij)\\
&\equiv&(jk)^{-1}(ik)^{-1}(jk)^{-1}(ik)\{\{(jl)^{\varepsilon},(kl)\},(il)\{(jl),(kl)\}\}(jk)(ij)\\
&\equiv&(jk)^{-1}(ik)^{-1}(jk)^{-1}\{\{\{(jl)^{\varepsilon},(il)^{-1}(kl)^{-1}(il)(kl)\},
\{(kl),(il)(kl)\}\},\{(il),(kl)\}\\
&&\{\{(jl),(il)^{-1}(kl)^{-1}(il)(kl)\},\{(kl),(il)(kl)\}\}\}(ik)(jk)(ij)\\
&\equiv&(jk)^{-1}(ik)^{-1}\{\{(jl)^{\varepsilon},(kl)^{-1}(jl)^{-1}\},(il)\{(jl),(kl)^{-1}(jl)^{-1}\}\\
&&\{(kl),(jl)^{-1}\}\}(jk)^{-1}(ik)(jk)(ij)\\
&\equiv&(jk)^{-1}\{\{\{(jl)^{\varepsilon},(kl)(il)(kl)^{-1}(il)^{-1}\},\{(kl)^{-1},(il)^{-1}\}\\
&&\{(jl)^{-1},(kl)(il)(kl)^{-1}(il)^{-1}\}\},\{(il),(kl)^{-1}(il)^{-1}\}\\
&&\{(jl)^{-1},(kl)(il)(kl)^{-1}(il)^{-1}\}\{(kl),(il)^{-1}\}\}(ik)^{-1}(jk)^{-1}(ik)(jk)(ij)\\
&\equiv&(jk)^{-1}\{(jl)^{\varepsilon},(kl)(il)(kl)^{-1}(il)^{-1}(kl)^{-1}(jl)^{-1}(kl)(il)(kl)^{-1}(jl)(kl)\}\\
&&(ik)^{-1}(jk)^{-1}(ik)(jk)(ij)\\
&\equiv&\{\{(jl)^{\varepsilon},(kl)^{-1}(jl)^{-1}\},\{(kl),(jl)^{-1}\}(il)\{(kl)^{-1},(jl)^{-1}\}\\
&&(il)^{-1}\{(kl)^{-1},(jl)^{-1}\}\{(jl)^{-1},(kl)^{-1}(jl)^{-1}\}\{(kl),(jl)^{-1}\}(il)\{(kl)^{-1},(jl)^{-1}\}\\
&&\{(jl),(kl)^{-1}(jl)^{-1}\}\{(kl),(jl)^{-1}\}\}(jk)^{-1}(ik)^{-1}(jk)^{-1}(ik)(jk)(ij)\\
&\equiv&\{(jl)^{\varepsilon},(il)(jl)(kl)^{-1}\{(jl)^{-1},(il)(jl)\}\}\{(jk)^{-1},(ik)(jk)\}(ij).
\end{eqnarray*}

$(8)\wedge(10)$\\
Let $f=(ij)(jk)-\{(jk),(ik)(jk)\}(ij),\
g=(jk)(jl)^{\varepsilon}-\{(jl)^{\varepsilon},(kl)\}(jk),\ i<j<k<l$
Then $w=(ij)(jk)(jl)^{\varepsilon}$ and
$$
(f,g)_w=(ij)\{(jl)^{\varepsilon},(kl)\}(jk)-\{(jk),(ik)(jk)\}(ij)(jl)^{\varepsilon}\equiv
0
$$
since
\begin{eqnarray*}
&&(ij)\{(jl)^{\varepsilon},(kl)\}(jk)\equiv\{(jl)^{\varepsilon},(il)(jl)(kl)\}\{(jk),(ik)(jk)\}(ij)\ \ \ \ \ and \\
&&\\
&&\{(jk),(ik)(jk)\}(ij)(jl)^{\varepsilon}\\
&\equiv&(jk)^{-1}(ik)^{-1}(jk)\{(jl)^{\varepsilon},(il)(jl)(kl)\}(ik)(jk)(ij)\\
&\equiv&(jk)^{-1}(ik)^{-1}\{\{(jl)^{\varepsilon},(kl)\}(il)\{(jl),(kl)\}\{(kl),(jl)(kl)\}\}(jk)(ik)(jk)(ij)\\
&\equiv&(jk)^{-1}\{\{\{(jl)^{\varepsilon},(kl)(il)(kl)^{-1}(il)^{-1}\},\{(kl),(il)^{-1}\}\},\{(il),(kl)^{-1}(il)^{-1}\}\\
&&\{(jl),(kl)(il)(kl)^{-1}(il)^{-1}\}\{(kl)(il)^{-1}\}\}(ik)^{-1}(jk)(ik)(jk)(ij)\\
&\equiv&\{\{(jl)^{\varepsilon},(kl)^{-1}(jl)^{-1}\},\{(kl),(jl)^{-1}\}(il)\{(jl),(kl)^{-1}(jl)^{-1}\}\\
&&\{(kl),(jl)^{-1}\}\}\{(jk),(ik)(jk)\}(ij)\\
&\equiv&\{(jl)^{\varepsilon},(il)(jl)(kl)\}\{(jk),(ik)(jk)\}(ij).
\end{eqnarray*}

$(8)\wedge(11)$\\
Let $f=(pi)(ik)^{-1}-\{(ik)^{-1},(pk)(ik)\}(pi),\
g=(ik)^{-1}(jl)^{\varepsilon}-\{(jl)^{\varepsilon},(kl)(il)(kl)^{-1}(il)^{-1}\}(ik)^{-1},\
p<i<j<k<l$. Then $w=(pi)(ik)^{-1}(jl)^{\varepsilon}$ and
$$
(f,g)_w=(pi)\{(jl)^{\varepsilon},(kl)(il)(kl)^{-1}(il)^{-1}\}(ik)^{-1}-\{(ik)^{-1},(pk)(ik)\}(pi)(jl)^{\varepsilon}\equiv
0
$$
since
\begin{eqnarray*}
&&(pi)\{(jl)^{\varepsilon},(kl)(il)(kl)^{-1}(il)^{-1}\}(ik)^{-1}\\
&\equiv&\{(jl)^{\varepsilon},(kl)\{(il),(pl)(il)\}(kl)^{-1}\{(il)^{-1},(pl)(il)\}\}\{(ik)^{-1},(pk)(ik)\}(pi) \ \ \ and \\
&&\\
&&\{(ik)^{-1},(pk)(ik)\}(pi)(jl)^{\varepsilon}\\
&\equiv&(ik)^{-1}(pk)^{-1}(ik)^{-1}(pk)\{(jl)^{\varepsilon},(il)^{-1}(kl)^{-1}(il)(kl)\}(ik)(pi)\\
&\equiv&(ik)^{-1}(pk)^{-1}(ik)^{-1}\{\{(jl)^{\varepsilon},(pl)^{-1}(kl)^{-1}(pl)(kl)\},\{(il)^{-1},(pl)^{-1}\\
&&(kl)^{-1}(pl)(kl)\}\{(kl)^{-1},(pl)(kl)\}\{(il),(pl)^{-1}(kl)^{-1}(pl)(kl)\}\{(kl),(pl)(kl)\}\}(pk)(ik)(pi)\\
&\equiv&(ik)^{-1}(pk)^{-1}\{\{(jl)^{\varepsilon},(kl)(il)(kl)^{-1}(il)^{-1}\},\{(il)^{-1},(kl)^{-1}(il)^{-1}\}\\
&&(pl)^{-1}\{(kl)^{-1},(il)^{-1}\}(pl)\{(il),(kl)^{-1}(il)^{-1}\}\{(kl),(il)^{-1}\}(ik)^{-1}(pk)(ik)(pi)\\
&\equiv&(ik)^{-1}\{\{(jl)^{\varepsilon},(kl)(pl)(kl)^{-1}(pl)^{-1}\},\{\{(pl)^{-1},(kl)^{-1}(pl)^{-1}\},\\
&&\{(il),(kl)(pl)(kl)^{-1}(pl)^{-1}\}\}\{\{(pl),(kl)^{-1}(pl)^{-1}\},\{(il),(kl)(pl)(kl)^{-1}(pl)^{-1}\}\\
&&\{(kl),(pl)^{-1}\}\}\}(pk)^{-1}(ik)^{-1}(pk)(ik)(pi)\\
&\equiv&\{\{(jl)^{\varepsilon},(kl)(il)(kl)^{-1}(il)^{-1}\},\{(il)^{-1},(kl)^{-1}(il)^{-1}\}\\
&&\{(kl),(il)^{-1}\}(pl)^{-1}\{(kl)^{-1},(il)^{-1}\}\{(il),(kl)^{-1}(il)^{-1}\}\{(kl),(il)^{-1}\}(pl)\\
&&\{(kl)^{-1},(il)^{-1}\}(pl)^{-1}\{(kl)^{-1},(il)^{-1}\}\{(il)^{-1},(kl)^{-1}(il)^{-1}\}\{(kl),(il)^{-1}\}\\
&&(pl)\{(kl)^{-1},(il)^{-1}\}\{(il),(kl)^{-1}(il)^{-1}\}\{(kl),(il)^{-1}\}\}\{(ik)^{-1}(pk)(ik)\}(pi)\\
&\equiv&\{(jl)^{\varepsilon},(kl)\{(il),(pl)(il)\}(kl)^{-1}\{(il)^{-1},(pl)(il)\}\}\{(ik)^{-1},(pk)(ik)\}(pi).
\end{eqnarray*}

$(8)\wedge(12)$\\
Letn $f=(pi)(ik)-\{(ik),(pk)(ik)\}(pi),\
g=(ik)(jl)^{\varepsilon}-\{(jl)^{\varepsilon},(il)^{-1}(kl)^{-1}(il)(kl)\}(ik),\
p<i<j<k<l$. Then $w=(pi)(ik)(jl)^{\varepsilon}$ and
$$
(f,g)_w=(pi)\{(jl)^{\varepsilon},(il)^{-1}(kl)^{-1}(il)(kl)\}(ik)-\{(ik),(pk)(ik)\}(pi)(jl)^{\varepsilon}\equiv
0
$$
since
\begin{eqnarray*}
&&(pi)\{(jl)^{\varepsilon},(il)^{-1}(kl)^{-1}(il)(kl)\}(ik)\\
&\equiv&\{(jl)^{\varepsilon},\{(il)^{-1},(pl)(il)\}(kl)^{-1}\{(il),(pl)(il)\}(kl)\}\{(ik),(pk)(ik)\}(pi) \ \ \ \ and\\
&&\\
&&\{(ik),(pk)(jk)\}(pi)(jl)^{\varepsilon}\\
&\equiv&(ik)^{-1}(pk)^{-1}(ik)\{(jl)^{\varepsilon},(il)^{-1}(pl)^{-1}(kl)^{-1}(pl)(il)(kl)\}(pk)(ik)(pi)\\
&\equiv&(ik)^{-1}(pk)^{-1}\{\{(jl)^{\varepsilon},(il)^{-1}(kl)^{-1}(il)(kl)\},\{(il)^{-1},(kl)\}(pl)^{-1}\\
&&\{(kl)^{-1},(il)(kl)\}(pl)\{(il),(kl)\}\{(kl),(il)(kl)\}\}(ik)(pk)(ik)(pi)\\
&\equiv&(ik)^{-1}\{\{(jl)^{\varepsilon},(kl)(pl)(kl)^{-1}(pl)^{-1}\},\{(il)^{-1},(kl)(pl)(kl)^{-1}(pl)^{-1}\}\\
&&\{(pl)^{-1},(kl)^{-1}(pl)^{-1}\}\{(kl)^{-1},(pl)^{-1}\}\{(il)^{-1},(kl)(pl)(kl)^{-1}(pl)^{-1}\}\}\\
&&\{(kl)^{-1},(pl)^{-1}\}\{(il),(kl)(pl)(kl)^{-1}(pl)^{-1}\}\{(kl),(pl)^{-1}\}\{(pl),(kl)^{-1}(pl)^{-1}\}\\
&&\{(il),(kl)(pl)(kl)^{-1}(pl)^{-1}\}\{(kl),(pl)^{-1}\}\}(pk)^{-1}(ik)(pk)(ik)(pi)\\
&\equiv&\{\{(jl)^{\varepsilon},(kl)(il)(kl)^{-1}(il)^{-1}\},\{(il)^{-1},(kl)^{-1}(il)^{-1}\}\\
&&\{\{(kl)^{-1},(il)^{-1}\},\{\{(il),(kl)^{-1}(il)^{-1}\},\{(kl),(il)^{-1}\}(pl)\}\}\\
&&\{(il),(kl)^{-1}(il)^{-1}\}\{(kl),(il)^{-1}\}\}(ik)^{-1}(pk)^{-1}(ik)(pk)(ik)(pi)\\
&\equiv&\{(jl)^{\varepsilon},\{(il)^{-1},(pl)(il)\}(kl)^{-1}\{(il),(pl)(il)\}(kl)\}\{(ik),(pk)(ik)\}(pi).
\end{eqnarray*}

$(8)\wedge(13)$\\
Let $f=(pi)(ik)^\delta-\{(ik)^\delta,(pk)(ik)\}(pi),
g=(ik)^\delta(jl)^{\varepsilon}-(jl)^{\varepsilon}(ik)^\delta$.\
Then $w=(pi)(ik)^\delta(jl)^{\varepsilon}$\ and
$$
(f,g)_w=(pi)(jl)^{\varepsilon}(ik)^\delta-\{(ik)^\delta,(pk)(ik)\}(pi)(jl)^{\varepsilon}.
$$
There are two cases to consider.
\begin{enumerate}
\item[1)]\ $p<i<k,j<i<k<l$.\ In this case,there are three subcases to
consider.
\begin{enumerate}
\item[(a)]\  $p<j,p<j<i<k<l$.
\begin{eqnarray*}
&&(pi)(jl)^{\varepsilon}(ik)^\delta\\
&\equiv&\{(jl)^{\varepsilon},(pl)^{-1}(il)^{-1}(pl)(il)\}\{(ik)^\delta,(pk)(ik)\}(pi)\ \ \ and \\
&&\\
&&\{(ik)^\delta,(pk)(ik)\}(pi)(jl)^{\varepsilon}\\
&\equiv&(ik)^{-1}(pk)^{-1}(ik)^\delta(pk)(ik)\{(jl)^{\varepsilon},(pl)^{-1}(il)^{-1}(pl)(il)\}(pi)\\
&\equiv&(ik)^{-1}(pk)^{-1}(ik)^\delta(pk)\{(jl)^{\varepsilon},(pl)^{-1}\{(il)^{-1},(kl)\}(pl)\{(il),(kl)\}\}(ik)(pi)\\
&\equiv&(ik)^{-1}(pk)^{-1}(ik)^\delta\{\{(jl)^{\varepsilon},(pl)^{-1}(kl)^{-1}(pl)(kl)\}\{(pl)^{-1},(kl)^{-1}\}\\
&&\{\{(il)^{-1},(pl)^{-1}(kl)^{-1}(pl)(kl)\},\{(kl),(pl)(kl)\}\}\\
&&\{(pl),(kl)\}\{\{(il),(pl)^{-1}(kl)^{-1}(pl)(kl)\},\{(kl),(pl)(kl)\}\}(pk)(ik)(pi)\\
&\equiv&(ik)^{-1}(pk)^{-1}(ik)^\delta\{(jl)^{\varepsilon},(pl)^{-1}(kl)^{-1}(il)^{-1}(pl)(il)(kl)\}(pk)(ik)(pi).
\end{eqnarray*}
If $\delta=1$, then
\begin{eqnarray*}
&&\{(ik)^\delta,(pk)(ik)\}(pi)(jl)^{\varepsilon}\\
&\equiv&(ik)^{-1}(pk)^{-1}\{(jl)^{\varepsilon},(pl)^{-1}\{(kl)^{-1},(il)(kl)\}\{(il)^{-1},(kl)\}(pl)\\
&&\{(il),(kl)\}\{(kl),(il)(kl)\}\}(ik)^\delta(pk)(ik)(pi)\\
&\equiv&(ik)^{-1}(pk)^{-1}\{(jl)^{\varepsilon},(pl)^{-1}\{(pl),(il)(kl)\}\}(ik)^\delta(pk)(ik)(pi)\\
&\equiv&(ik)^{-1}\{\{(jl)^{\varepsilon},(kl)(pl)(kl)^{-1}(pl)^{-1}\},\{(pl)^{-1},(kl)^{-1}(pl)^{-1}\}\\
&&\{\{(pl),(kl)^{-1}(pl)^{-1}\},\{(il),(kl)(pl)(kl)^{-1}(pl)^{-1}\}\{(kl),(pl)^{-1}\}\}\}\\
&&(pk)^{-1}(ik)^\delta(pk)(ik)(pi)\\
&\equiv&\{(jl)^{\varepsilon},(pl)^{-1}\{(pl),\{(kl)^{-1},(il)^{-1}\}\{(il),(kl)^{-1}(il)^{-1}\}\\
&&\{(kl),(il)^{-1}\}\}\}(ik)^{-1}(pk)^{-1}(ik)^\delta(pk)(ik)(pi)\\
&\equiv&\{(jl)^{\varepsilon},(pl)^{-1}(il)^{-1}(pl)(il)\}\{(ik)^\delta,(pk)(ik)\}(pi).
\end{eqnarray*}
If $ \delta=-1$, then
\begin{eqnarray*}
&&\{(ik)^\delta,(pk)(ik)\}(pi)(jl)^{\varepsilon}\\
&\equiv&(ik)^{-1}(pk)^{-1}\{(jl)^{\varepsilon},(pl)^{-1}\{(kl)^{-1},(il)^{-1}\}\{(il)^{-1},(kl)^{-1}\\
&&(il)^{-1}\}(pl)\{(il),(kl)^{-1}(il)^{-1}\}\{(kl),(il)^{-1}\}(ik)^\delta(pk)(ik)(pi)\\
&\equiv&\{(jl)^{\varepsilon},(pl)^{-1}(il)^{-1}(pl)(il)\}\{(ik)^\delta,(pk)(ik)\}(pi).
\end{eqnarray*}
\item[(b)]\ $p=j,p<j<i<k<l$.
\begin{eqnarray*}
&&(pi)(jl)^{\varepsilon}(ik)^\delta\\
&\equiv&(ji)(jl)^{\varepsilon}(ik)^\delta\equiv\{(jl)^{\varepsilon},(il)\}\{(ik)^\delta,(jk)(ik)\}(ji)\ \ \ and \\
&&\\
&&\{(ik)^\delta,(pk)(ik)\}(pi)(jl)^{\varepsilon}=(ik)^{-1}(jk)^{-1}(ik)^\delta(jk)(ik)(ji)(jl)^{\varepsilon}\\
&\equiv&(ik)^{-1}(jk)^{-1}(ik)^\delta(jk)(ik)\{(jl)^{\varepsilon},(il)\}(ji)\\
&\equiv&(ik)^{-1}(jk)^{-1}(ik)^\delta(jk)\{(jl)^{\varepsilon},\{(il),(kl)\}\}(ik)(ji)\\
&\equiv&(ik)^{-1}(jk)^{-1}(ik)^\delta\{\{(jl)^{\varepsilon},(kl)\}\{\{(il),(jl)^{-1}(kl)^{-1}(jl)(kl)\},\\
&&\{(kl),(jl)(kl)\}\}\}(jk)(ik)(ji)\\
&\equiv&(ik)^{-1}(jk)^{-1}(ik)^\delta\{(jl)^{\varepsilon},(il)(kl)\}(jk)(ik)(ji).
\end{eqnarray*}
If $\delta=1$, then
\begin{eqnarray*}
&&\{(ik)^\delta,(pk)(ik)\}(pi)(jl)^{\varepsilon}\\
&\equiv&(ik)^{-1}(jk)^{-1}\{(jl)^{\varepsilon},\{(il),(kl)\}\{(kl),(il)(kl)\}\}(ik)^\delta(jk)(ik)(ji)\\
&\equiv&(ik)^{-1}(jk)^{-1}\{(jl)^{\varepsilon},(il)(kl)\}(ik)^\delta(jk)(ik)(ji)\\
&\equiv&(ik)^{-1}\{\{(jl)^{\varepsilon},(kl)^{-1}(jl)^{-1}\},\{(il),(kl)(jl)(kl)^{-1}(jl)^{-1}\}\\
&&\{(kl),(jl)^{-1}\}\}(jk)^{-1}(ik)^\delta(jk)(ik)(ji)\\
&\equiv&\{(jl)^{\varepsilon},(jl)^{-1}\{(kl)^{-1},(il)^{-1}\}\{(il),(kl)^{-1}(il)^{-1}\}\\
&&\{(kl),(il)^{-1}\}\}(ik)^{-1}(jk)^{-1}(ik)^\delta(jk)(ik)(ji)\\
&\equiv&\{(jl)^{\varepsilon},(il)\}\{(ik)^\delta,(jk)(ik)\}(ji).
\end{eqnarray*}
If $\delta=-1$, then
\begin{eqnarray*}
&&\{(ik)^\delta,(pk)(ik)\}(pi)(jl)^{\varepsilon}\\
&\equiv&(ik)^{-1}(jk)^{-1}\{(jl)^{\varepsilon},\{(il),(kl)^{-1}(il)^{-1}\}\{(kl),(il)^{-1}\}\}(ik)^\delta(jk)(ik)(ji)\\
&\equiv&(ik)^{-1}(jk)^{-1}\{(jl)^{\varepsilon},(il)(kl)\}(ik)^\delta(jk)(ik)(ji)\\
&\equiv&\{(jl)^{\varepsilon},(il)\}\{(ik)^\delta,(jk)(ik)\}(ji).
\end{eqnarray*}
\item[(c)]\ $p>j,j<p<i<k<l$.
\begin{eqnarray*}
(f,g)_w&\equiv&(jl)^{\varepsilon}(pi)(ik)^\delta-(jl)^{\varepsilon}\{(ik)^\delta,(pk)(ik)\}(pi)\\
&\equiv&0.
\end{eqnarray*}
\end{enumerate}
\item[2)]\ $p<i<k<j<l$.
\begin{eqnarray*}
(f,g)_w&=&(pi)(jl)^{\varepsilon}(ik)^\delta-\{(ik)^\delta,(pk)(ik)\}(pi)(jl)^{\varepsilon}\\
&\equiv&(jl)^{\varepsilon}(pi)(ik)^\delta-(jl)^{\varepsilon}\{(ik)^\delta,(pk)(ik)\}(pi)\\
&\equiv&0.
\end{eqnarray*}
\end{enumerate}

$(9)\wedge(7)$\\
Let
$f=(ji)^{-1}(jk)^{-1}-\{(jk)^{-1},(ik)^{-1}(jk)^{-1}\}(ji)^{-1},\
g=(jk)^{-1}(kl)^{\varepsilon}-\{(kl)^{\varepsilon},(jl)^{-1}\}(jk)^{-1},\
j<i<k<l$. Then $w=(ji)^{-1}(jk)^{-1}(kl)^{\varepsilon}$ and
$$
(f,g)_w=(ji)^{-1}\{(kl)^{\varepsilon},(jl)^{-1}\}(jk)^{-1}-\{(jk)^{-1},(ik)^{-1}(jk)^{-1}\}(ji)^{-1}(kl)^{\varepsilon}\equiv
0
$$
since
\begin{eqnarray*}
&&(ji)^{-1}\{(kl)^{\varepsilon},(jl)^{-1}\}(jk)^{-1}\\
&\equiv&\{(kl)^{\varepsilon},\{(jl)^{-1},(il)^{-1}(jl)^{-1}\}\}\{(jk)^{-1},(ik)^{-1}(jk)^{-1}\}(ji)^{-1} \ \ \ \ \ and\\
&&\\
&&\{(jk)^{-1},(ik)^{-1}(jk)^{-1}\}(ji)^{-1}(kl)^{\varepsilon}\\
&\equiv&(jk)(ik)(jk)^{-1}(ik)^{-1}(jk)^{-1}(kl)^{\varepsilon}(ji)^{-1}\\
&\equiv&(jk)(ik)(jk)^{-1}(ik)^{-1}\{(kl)^{\varepsilon},(jl)^{-1}\}(jk)^{-1}(ji)^{-1}\\
&\equiv&(jk)(ik)(jk)^{-1}\{\{(kl)^{\varepsilon},(il)^{-1}\},(jl)^{-1}\}(ik)^{-1}(jk)^{-1}(ji)^{-1}\\
&\equiv&(jk)(ik)\{\{(kl)^{\varepsilon},(jl)^{-1}\},\{(il)^{-1},(kl)(jl)(kl)^{-1}(jl)^{-1}\}\\
&&\{(jl)^{-1},(kl)^{-1}(jl)^{-1}\}\}(jk)^{-1}(ik)^{-1}(jk)^{-1}(ji)^{-1}\\
&\equiv&(jk)\{\{(kl)^{\varepsilon},(il)(kl)\},(jl)^{-1}\{(kl)^{-1},(il)(kl)\}\{(il)^{-1},(kl)\}(jl)^{-1}\}\\
&&(ik)(jk)^{-1}(ik)^{-1}(jk)^{-1}(ji)^{-1}\\
&\equiv&\{\{(kl)^{\varepsilon},(jl)(kl)\},\{(il),(jl)^{-1}(kl)^{-1}(jl)(kl)\}
\{(kl),(jl)(kl)\}\{(jl)^{-1},(kl)\}\{(kl)^{-1},\\
&&(jl)(kl)\}\{(il)^{-1},(jl)^{-1}(kl)^{-1}(jl)(kl)\}\{(jl)^{-1},(kl)\}\}(jk)(ik)(jk)^{-1}(ik)^{-1}(jk)^{-1}(ji)^{-1}\\
&\equiv&\{(kl)^{\varepsilon},\{(jl)^{-1},(il)^{-1}(jl)^{-1}\}\}\{(jk)^{-1},(ik)^{-1}(jk)^{-1}\}(ji)^{-1}.
\end{eqnarray*}

$(9)\wedge(8)$\\
let $f=(ji)^{-1}(jk)-\{(jk),(ik)^{-1}(jk)^{-1}\}(ji)^{-1},\
g=(jk)(kl)^{\varepsilon}-\{(kl)^{\varepsilon},(jl)(kl)\}(jk),\
j<i<k<l$. Then $w=(ji)^{-1}(jk)(kl)^{\varepsilon}$ and
$$
(f,g)_w=(ji)^{-1}\{(kl)^{\varepsilon},(jl)(kl)\}(jk)-\{(jk),(ik)^{-1}(jk)^{-1}\}(ji)^{-1}(kl)^{\varepsilon}\equiv
0
$$
since
\begin{eqnarray*}
&&(ji)^{-1}\{(kl)^{\varepsilon},(jl)(kl)\}(jk)\\
&\equiv&\{(kl)^{\varepsilon},\{(jl),(il)^{-1}(jl)^{-1}\}(kl)\}\{(jk),(ik)^{-1}(jk)^{-1}\}(ji)^{-1}\\
&\equiv&\{(kl)^{\varepsilon},(jl)(il)(jl)(il)^{-1}(jl)^{-1}(kl)\}(jk)(ik)(jk)(ik)^{-1}(jk)^{-1}(ji)^{-1} \ \ \ \ and\\
&&\\
&&\{(jk),(ik)^{-1}(jk)^{-1}\}(ji)^{-1}(kl)^{\varepsilon}\\
&\equiv&(jk)(ik)(jk)\{(kl)^{\varepsilon},(il)^{-1}(jl)^{-1}\}(ik)^{-1}(jk)^{-1}(ji)^{-1}\\
&\equiv&(jk)(ik)\{\{(kl)^{\varepsilon},(jl)(kl)\},\{(il)^{-1},(jl)^{-1}(kl)^{-1}(jl)(kl)\}\\
&&\{(jl)^{-1},(kl)\}\}(jk)(ik)^{-1}(jk)^{-1}(ji)^{-1}\\
&\equiv&(jk)\{\{(kl)^{\varepsilon},(il)(kl)\},(jl)\{(il)^{-1},(kl)\}(jl)^{-1}\}(ik)(jk)(ik)^{-1}(jk)^{-1}(ji)^{-1}\\
&\equiv&\{\{(kl)^{\varepsilon},(jl)(kl)\},\{\{(jl),(kl)\},\{(kl)^{-1},(jl)(kl)\}\{(il)^{-1},(jl)^{-1}(kl)^{-1}\\
&&(jl)(kl)\}\}\{(kl),(jl)(kl)\}\{(jl)^{-1},(kl)\}\}(jk)(ik)(jk)(ik)^{-1}(jk)^{-1}(ji)^{-1}\\
&\equiv&\{(kl)^{\varepsilon},(jl)(il)(jl)(il)^{-1}(jl)^{-1}(kl)\}(jk)(ik)(jk)(ik)^{-1}(jk)^{-1}(ji)^{-1}.
\end{eqnarray*}

$(9)\wedge(9)$\\
Let
$f=(ji)^{-1}(jk)^{-1}-\{(jk)^{-1},(ik)^{-1}(jk)^{-1}\}(ji)^{-1},\
g=(jk)^{-1}(jl)^{\varepsilon}-\{(jl)^{\varepsilon},(kl)^{-1}(jl)^{-1}\}(jk)^{-1},\
j<i<k<l$.
Then $w=(ji)^{-1}(jk)^{-1}(jl)^{\varepsilon}$ and
$$
(f,g)_w=(ji)^{-1}\{(jl)^{\varepsilon},(kl)^{-1}(jl)^{-1}\}(jk)^{-1}-
\{(jk)^{-1},(ik)^{-1}(jk)^{-1}\}(ji)^{-1}(jl)^{\varepsilon}\equiv 0
$$
since
\begin{eqnarray*}
&&(ji)^{-1}\{(jl)^{\varepsilon},(kl)^{-1}(jl)^{-1}\}(jk)^{-1}\\
&\equiv&\{\{(jl)^{\varepsilon},(il)^{-1}(jl)^{-1}\},(kl)^{-1}\{(jl)^{-1},(il)^{-1}(jl)^{-1}\}\}
\{(jk)^{-1},(ik)^{-1}(jk)^{-1}\}(ji)^{-1}\\
&\equiv&\{(jl)^{\varepsilon},(il)^{-1}(jl)^{-1}(kl)^{-1}(jl)(il)(jl)^{-1}(il)^{-1}(jl)^{-1}\}\\
&&\{(jk)^{-1},(ik)^{-1}(jk)^{-1}\}(ji)^{-1} \ \ \ \ \ \ \ \ \ and\\
&&\\
&&\{(jk)^{-1},(ik)^{-1}(jk)^{-1}\}(ji)^{-1}(jl)^{\varepsilon}\\
&\equiv&(jk)(ik)(jk)^{-1}\{(jl)^{\varepsilon},\{(kl)^{-1},(il)^{-1}\}\{(il)^{-1},(kl)^{-1}(il)^{-1}\}
(jl)^{-1}\}(ik)^{-1}(jk)^{-1}(ji)^{-1}\\
&\equiv&(jk)(ik)\{\{(jl)^{\varepsilon},(kl)^{-1}(jl)^{-1}\},\{(kl)^{-1},(jl)^{-1}\}\\
&&\{(il)^{-1},(kl)(jl)(kl)^{-1}(jl)^{-1}\}\{(jl)^{-1},(kl)^{-1}(jl)^{-1}\}\}(jk)^{-1}(ik)^{-1}(jk)(ji)^{-1}\\
&\equiv&(jk)\{(jl)^{\varepsilon},\{(kl)^{-1},(il)(kl)\}(jl)^{-1}\{(kl)^{-1},(il)(kl)\}\\
&&\{(il)^{-1},(kl)\}(jl)^{-1}\}(ik)(jk)^{-1}(ik)^{-1}(jk)(ji)^{-1}\\
&\equiv&\{\{(jl)^{\varepsilon},(kl)\},\{\{(kl)^{-1},(jl)(kl),\{(il),(jl)^{-1}(kl)^{-1}(jl)(kl)\}\\
&&\{(kl),(jl)(kl)\}\}\{(jl)^{-1},(kl)\}\{(kl)^{-1},(jl)(kl)\}\{(il)^{-1},(jl)^{-1}(kl)^{-1}(jl)(kl)\}\\
&&\{(jl)^{-1},(kl)\}\}(jk)(ik)(jk)^{-1}(ik)^{-1}(jk)(ji)^{-1}\\
&\equiv&\{(jl)^{\varepsilon},(il)^{-1}(jl)^{-1}(kl)^{-1}(jl)(il)(jl)^{-1}(il)^{-1}(jl)^{-1}\}
\{(jk)^{-1},(ik)^{-1}(jk)^{-1}\}(ji)^{-1}.
\end{eqnarray*}

$(9)\wedge(10)$\\
Let $f=(ji)^{-1}(jk)-\{(jk),(ik)^{-1}(jk)^{-1}\}(ji)^{-1},\
g=(jk)(jl)^{\varepsilon}-\{(jl)^{\varepsilon},(kl)\}(jk),j<i<k<l$.
Then $w=(ji)^{-1}(jk)(jl)^{\varepsilon}$ and
$$
(f,g)_w=(ji)^{-1}\{(jl)^{\varepsilon},(kl)\}(jk)-\{(jk),(ik)^{-1}(jk)^{-1}\}(ji)^{-1}(jl)^{\varepsilon}\equiv
0
$$
since
\begin{eqnarray*}
&&(ji)^{-1}\{(jl)^{\varepsilon},(kl)\}(jk)\\
&\equiv&\{(jl)^{\varepsilon},(il)^{-1}(jl)^{-1}(kl)\}\{(jk),(ik)^{-1}(jk)^{-1}\}(ji)^{-1} \ \ \ \ \ \ \ \ and \\
&&\\
&&\{(jk),(ik)^{-1}(jk)^{-1}\}(ji)^{-1}(jl)^{\varepsilon}\\
&\equiv&(jk)(ik)(jk)\{(jl)^{\varepsilon},(kl)^{-1}(il)^{-1}(jl)^{-1}\}(ik)^{-1}(jk)^{-1}(ji)^{-1}\\
&\equiv&(jk)(ik)\{\{(jl)^{\varepsilon},(kl)\},\{(kl)^{-1},(jl)(kl)\}\{(il)^{-1},(jl)^{-1}(kl)^{-1}(jl)(kl)\}\\
&&\{(jl)^{-1},(kl)\}\}(jk)(ik)^{-1}(jk)^{-1}(ji)^{-1}\\
&\equiv&(jk)\{(jl)^{\varepsilon},\{(il)^{-1},(kl)\}(jl)^{-1}\}(ik)(jk)(ik)^{-1}(jk)^{-1}(ji)^{-1}\\
&\equiv&\{\{(jl)^{\varepsilon},(kl)\},\{\{(il)^{-1},(jl)^{-1}(kl)^{-1}(jl)(kl)\},\{(kl),(jl)(kl)\}\}\\
&&\{(jl)^{-1},(kl)\}\}(jk)(ik)(jk)(ik)^{-1}(jk)^{-1}(ji)^{-1}\\
&\equiv&\{(jl)^{\varepsilon},(il)^{-1}(jl)^{-1}(kl)\}\{(jk),(ik)^{-1}(jk)^{-1}\}(ji)^{-1}.
\end{eqnarray*}

$(9)\wedge(11)$\\
Let
$f=(ip)^{-1}(ik)^{-1}-\{(ik)^{-1},(pk)^{-1}(ik)^{-1}\}(ip)^{-1},i<p<k,\
g=(ik)^{-1}(jl)^{\varepsilon}-\{(jl)^{\varepsilon},(kl)(il)(kl)^{-1}(il)^{-1}\}(ik)^{-1},i<j<k<l.$
 \ Then \ $w=(ip)^{-1}(ik)^{-1}(jl)^{\varepsilon}$ and
$$
(f,g)_w=(ip)^{-1}\{(jl)^{\varepsilon},(kl)(il)(kl)^{-1}(il)^{-1}\}(ik)^{-1}-
\{ik)^{-1},(pk)^{-1}(ik)^{-1}\}(ip)^{-1}(jl)^{\varepsilon}.
$$
There are three cases to consider.
\begin{enumerate}
\item[1)]\ $p<j,i<p<j<k<l.$
\begin{eqnarray*}
&&(ip)^{-1}\{(jl)^{\varepsilon},(kl)(il)(kl)^{-1}(il)^{-1}\}(ik)^{-1}\\
&\equiv&\{(jl)^{\varepsilon},(kl)\{(il),(pl)^{-1}(il)^{-1}\}(kl)^{-1}\{(il)^{-1},(pl)^{-1}(il)^{-1}\}\}\\
&&\{(ik)^{-1},(pk)^{-1}(ik)^{-1}\}(ip)^{-1}\ \ \ \ \ \ \ \ \ \ \ \ \  and \\
&&\\
&&\{(ik)^{-1},(pk)^{-1}(ik)^{-1}\}(ip)^{-1}(jl)^{\varepsilon}\\
&\equiv&(ik)(pk)(ik)^{-1}(pk)^{-1}\{(jl)^{\varepsilon},(kl)(il)(kl)^{-1}(il)^{-1}\}(ik)^{-1}(ip)^{-1}\\
&\equiv&(ik)(pk)(ik)^{-1}\{\{(jl)^{\varepsilon},(kl)(pl)(kl)^{-1}(pl)^{-1}\},\{(kl),(pl)^{-1}\}(il)\\
&&\{(kl)^{-1},(pl)^{-1}\}(il)^{-1}\}(pk)^{-1}(ik)^{-1}(ip)^{-1}\\
&\equiv&(ik)(pk)\{\{(jl)^{\varepsilon},(kl)(il)(kl)^{-1}(il)^{-1}\},\{(kl),(il)^{-1}\}\{\{(kl)^{-1},(il)^{-1}\}\\
&&,\{(pl)^{-1},(kl)(il)(kl)^{-1}(il)^{-1}\}\{(il)^{-1},(kl)^{-1}(il)^{-1}\}\}\}(ik)^{-1}(pk)^{-1}(ik)^{-1}(ip)^{-1}\\
&\equiv&(ik)\{\{(jl)^{\varepsilon},(pl)^{-1}(kl)^{-1}(pl)(kl)\},\{(kl),(pl)(kl)\}\{\{(kl)^{-1},(pl)(kl)\},(il)^{-1}\\
&&\{(kl)^{-1},(pl)(kl)\}\{(pl)^{-1},(kl)\}(il)^{-1}\}\}\}(pk)(ik)^{-1}(pk)^{-1}(ik)^{-1}(ip)^{-1}\\
&\equiv&\{\{(jl)^{\varepsilon},(il)^{-1}(kl)^{-1}(il)(kl)\},\{(kl),(il)(kl)\}\{\{(kl)^{-1},(il)(kl)\},\\
&&\{\{(il)^{-1},(kl)\},\{(kl)^{-1},(il)(kl)\}\{(pl)^{-1},(il)^{-1}(kl)^{-1}(il)(kl)\}\}\{(il)^{-1},(kl)\}\}\}\\
&&(ik)(pk)(ik)^{-1}(pk)^{-1}(ik)^{-1}(ip)^{-1}\\
&\equiv&\{(jl)^{\varepsilon},(kl)\{(il),(pl)^{-1}(il)^{-1}\}(kl)^{-1}\{(il),(pl)^{-1}(il)^{-1}\}\}\\
&&\{(ik)^{-1},(pk)^{-1}(ik)^{-1}\}(ip)^{-1}.
\end{eqnarray*}
\item[2)]\ $p=j,i<p=j<k<l$.
\begin{eqnarray*}
&&(ip)^{-1}\{(jl)^{\varepsilon},(kl)(il)(kl)^{-1}(il)^{-1}\}(ik)^{-1}\\
&=&(ij)^{-1}\{(jl)^{\varepsilon},(kl)(il)(kl)^{-1}(il)^{-1}\}(ik)^{-1}\\
&\equiv&\{(jl)^{\varepsilon},(il)^{-1}(kl)\{(il),(jl)^{-1}(il)^{-1}(kl)^{-1}\}\{(il)^{-1},(jl)^{-1}(il)^{-1}\}\}\\
&&\{(ik)^{-1},(jk)^{-1}(ik)^{-1}\}(ij)^{-1}\ \ \ \ \ \ \ \ \ \ \ \ \ \ \ and \\
&&\\
&&\{(ik)^{-1},(pk)^{-1}(ik)^{-1}\}(ip)^{-1}(jl)^{\varepsilon}\\
&=&\{(ik)^{-1},(jk)^{-1}(ik)^{-1}\}(ij)^{-1}(jl)^{\varepsilon}\\
&\equiv&(ik)(jk)(ik)^{-1}(jk)^{-1}(ik)^{-1}\{(jl)^{\varepsilon},(il)^{-1}\}(ij)^{-1}\\
&\equiv&(ik)(jk)(ik)^{-1}(jk)^{-1}\{\{(jl)^{\varepsilon},(kl)(il)(kl)^{-1}(il)^{-1}\},\\
&&\{(il)^{-1},(kl)^{-1}(il)^{-1}\}\}(ik)^{-1}(ij)^{-1}\\
&\equiv&(ik)(jk)(ik)^{-1}\{\{(jl)^{\varepsilon},(kl)^{-1}(jl)^{-1}\},(il)^{-1}\}(jk)^{-1}(ik)^{-1}(ij)^{-1}\\
&\equiv&(ik)(jk)\{\{(jl)^{\varepsilon},(kl)(il)(kl)^{-1}(il)^{-1}\},\{(kl)^{-1},(il)^{-1}\}\{(jl)^{-1},(kl)(il)\\
&&(kl)^{-1}(il)^{-1}\}\{(il)^{-1},(kl)^{-1}(il)^{-1}\}\}(ik)^{-1}(jk)^{-1}(ik)^{-1}(ij)^{-1}\\
&\equiv&(ik)\{\{(jl)^{\varepsilon},(kl)\},\{\{(kl)^{-1},(jl)(kl)\},(il)^{-1}\{(kl)^{-1},(jl)(kl)\}\}\\
&&\{(jl)^{-1},(kl)\}(il)^{-1}\}(jk)(ik)^{-1}(jk)^{-1}(ik)^{-1}(ij)^{-1}\\
&\equiv&\{\{(jl)^{\varepsilon},(il)^{-1}(kl)^{-1}(il)(kl)\},\{(kl),(il)(kl)\}\{\{(kl)^{-1},(il)(kl)\},\\
&&\{(jl),(il)^{-1}(kl)^{-1}(il)(kl)\}\{(kl),(il)(kl)\}\},\{(il)^{-1},(kl)\}\{(kl)^{-1},(il)(kl)\}\\
&&\{(jl)^{-1},(il)^{-1}(kl)^{-1}(il)(kl)\}\}\{(il)^{-1},(kl)\}\}(ik)(jk)(ik)^{-1}(jk)^{-1}(ik)^{-1}(ij)^{-1}\\
&\equiv&\{(jl)^{\varepsilon},(il)^{-1}(kl)\{(il),(jl)^{-1}(il)^{-1}\}(kl)^{-1}\{(il)^{-1},(jl)^{-1}(il)^{-1}\}\}\\
&&\{(ik)^{-1},(jk)^{-1}(ik)^{-1}\}(ij)^{-1}.
\end{eqnarray*}
\item[3)]\ $p>j,i<j<p<k<l$.
\begin{eqnarray*}
&&(ip)^{-1}\{(jl)^{\varepsilon},(kl)(il)(kl)^{-1}(il)^{-1}\}(ik)^{-1}\\
&\equiv&\{(jl)^{\varepsilon},(pl)(il)(pl)^{-1}(il)^{-1}(kl)\{(il),(pl)^{-1}(il)^{-1}\}(kl)^{-1}\\
&&\{(il)^{-1},(pl)^{-1}(il)^{-1}\}\}\{(ik)^{-1},(pk)^{-1}(ik)^{-1}\}(ip)^{-1}\ \ \ \ \ and \\
&&\\
&&\{(ik)^{-1},(pk)^{-1}(ik)^{-1}\}(ip)^{-1}(jl)^{\varepsilon}\\
&\equiv&(ik)(pk)(ik)^{-1}(pk)^{-1}(ik)^{-1}\{(jl)^{\varepsilon},(pl)(il)(pl)^{-1}(il)^{-1}\}(ip)^{-1}\\
&\equiv&(ik)(pk)(ik)^{-1}(pk)^{-1}\{\{(jl)^{\varepsilon},(kl)(il)(kl)^{-1}(il)^{-1}\},\{(pl),(kl)(il)(kl)^{-1}(il)^{-1}\}\\
&&\{(il),(kl)^{-1}(il)^{-1}\}\{(pl)^{-1},(kl)(il)(kl)^{-1}(il)^{-1}\}\{(il)^{-1},(kl)^{-1}(il)^{-1}\}\}(ik)^{-1}(ip)^{-1}\\
&\equiv&(ik)(pk)(ik)^{-1}\{(jl)^{\varepsilon},\{(pl),(kl)^{-1}(pl)^{-1}\}\{(kl),(pl)^{-1}\}(il)\{(kl)^{-1},(pl)^{-1}\}\\
&&\{(pl)^{-1},(kl)^{-1}(pl)^{-1}\}(il)^{-1}\}(pk)^{-1}(ik)^{-1}(ip)^{-1}\\
&\equiv&(ik)(pk)\{\{(jl)^{\varepsilon},(kl)(il)(kl)^{-1}(il)^{-1}\},\{\{(il),(kl)^{-1}(il)^{-1}\},\{(kl)^{-1},(il)^{-1}\}\\
&&\{(pl)^{-1},(kl)(il)(kl)^{-1}(il)^{-1}\}\}\{(il)^{-1},(kl)^{-1}(il)^{-1}\}\}(ik)^{-1}(pk)^{-1}(ik)^{-1}(ip)^{-1}\\
&\equiv&(ik)\{(jl)^{\varepsilon},\{(il),\{(kl)^{-1},(pl)(kl)\}(il)^{-1}\{(kl)^{-1},(pl)(kl)\}\{(pl)^{-1},(kl)\}\}\\
&&(il)^{-1}\}(pk)(ik)^{-1}(pk)^{-1}(ik)^{-1}(ip)^{-1}\\
&\equiv&\{\{(jl)^{\varepsilon},(il)^{-1}(kl)^{-1}(il)(kl)\},\{\{(il),(kl)\},\{\{(kl)^{-1},(il)(kl)\},\\
&&\{(pl),(il)^{-1}(kl)^{-1}(il)(kl)\}\{(kl),(il)(kl)\}\}\{(il)^{-1},(kl)\}\{(kl)^{-1},(il)(kl)\}\\
&&\{(pl)^{-1},(il)^{-1}(kl)^{-1}(il)(kl)\}\}\{(il),(kl)\}\}\{(ik)^{-1},(pk)^{-1}(ik)^{-1}\}(ip)^{-1}\\
&\equiv&\{(jl)^{\varepsilon},(pl)(il)(pl)^{-1}(il)^{-1}(kl)\{(il),(pl)^{-1}(il)^{-1}\}(kl)^{-1}\\
&&\{(il)^{-1},(pl)^{-1}(il)^{-1}\}\}\{(ik)^{-1},(pk)^{-1}(ik)^{-1}\}(ip)^{-1}.
\end{eqnarray*}
\end{enumerate}

$(9)\wedge(12)$\\
Let $f=(ip)^{-1}(ik)-\{(ik),(pk)^{-1}(ik)^{-1}\}(ip)^{-1},i<p<k,\ \
 g=(ik)(jl)^{\varepsilon}-\{(jl)^{\varepsilon},(il)^{-1}(kl)^{-1}$\\
 $(il)(kl)\}(ik),
 \ i<j<k<l$. \ Then $w=(ip)^{-1}(ik)(jl)^{\varepsilon}$ and
$$
(f,g)_w=(ip)^{-1}\{(jl)^{\varepsilon},(il)^{-1}(kl)^{-1}(il)(kl)\}(ik)-
\{(ik),(pk)^{-1}(ik)^{-1}\}(ip)^{-1}(jl)^{\varepsilon}.
$$
There are three cases to consider.
\begin{enumerate}
\item[1)]\  $p<j,i<p<j<k<l$.
\begin{eqnarray*}
&&(ip)^{-1}\{(jl)^{\varepsilon},(il)^{-1}(kl)^{-1}(il)(kl)\}(ik)\\
&\equiv&\{(jl)^{\varepsilon},\{(il)^{-1},(pl)^{-1}(il)^{-1}\}(kl)^{-1}\{(il),(pl)^{-1}(il)^{-1}\}(kl)\}\\
&&\{(ik),(pk)^{-1}(ik)^{-1}\}(ip)^{-1}\ \ \ \ \ \ \ \ \ \ \  \ \ \ \ \ and \\
&&\\
&&\{(ik),(pk)^{-1}(ik)^{-1}\}(ip)^{-1}(jl)^{\varepsilon}\\
&\equiv&(ik)(pk)(ik)\{(jl)^{\varepsilon},(kl)\{(kl)^{-1},(pl)^{-1}(il)^{-1}\}\}(pk)^{-1}(ik)^{-1}(ip)^{-1}\\
&\equiv&(ik)(pk)\{\{(jl)^{\varepsilon},(il)^{-1}(kl)^{-1}(il)(kl)\},\{(kl),(il)(kl)\}\{\{(kl)^{-1},(il)(kl)\},\\
&&\{(pl)^{-1},(il)^{-1}(kl)^{-1}(il)(kl)\}\{(il)^{-1},(kl)\}\}\}(ik)(pk)^{-1}(ik)^{-1}(ip)^{-1}\\
&\equiv&(ik)\{\{(jl)^{\varepsilon},(pl)^{-1}(kl)^{-1}(pl)(kl)\},\{(kl),(pl)(kl)\}\{\{(kl)^{-1},(pl)(kl)\},\\
&&(il)\{(pl)^{-1},(kl)\}(il)^{-1}\}\}(pk)(ik)(pk)^{-1}(ik)^{-1}(ip)^{-1}\\
&\equiv&\{\{(jl)^{\varepsilon},(il)^{-1}(kl)^{-1}(il)(kl)\},\{(kl),(il)(kl)\}\{\{(kl)^{-1},(il)(kl)\},\{\{(il),(kl)\},\\
&&\{(kl)^{-1},(il)(kl)\}\{(pl)^{-1},(il)^{-1}(kl)^{-1}(il)(kl)\}\}\{(kl),(il)(kl)\}\\
&&\{(il)^{-1},(kl)\}\}\}\{(ik),(pk)^{-1}(ik)^{-1}\}(ip)^{-1}\\
&\equiv&\{(jl)^{\varepsilon},\{(il)^{-1},(pl)^{-1}(il)^{-1}\}(kl)^{-1}\{(il),(pl)^{-1}(il)^{-1}\}(kl)\}\\
&&\{(ik),(pk)^{-1}(ik)^{-1}\}(ip)^{-1}.
\end{eqnarray*}
\item[2)]\  $p=j,i<p=j<k<l$.
\begin{eqnarray*}
&&(ip)^{-1}\{(jl)^{\varepsilon},(il)^{-1}(kl)^{-1}(il)(kl)\}(ik)\\
&=&(ij)^{-1}\{(jl)^{\varepsilon},(il)^{-1}(kl)^{-1}(il)(kl)\}(ik)\\
&\equiv&\{(jl)^{\varepsilon},(il)^{-1}(jl)^{-1}(il)^{-1}(kl)^{-1}\{(il),(jl)^{-1}(il)^{-1}\}(kl)\}\\
&&\{(ik),(jk)^{-1}(ik)^{-1}\}(ij)^{-1}\ \ \ \ \ \ \ \ \ \ \  \ \ \ \ and \\
&&\\
&&\{(ik),(pk)^{-1}(ik)^{-1}\}(ip)^{-1}(jl)^{\varepsilon}\\
&\equiv&(ik)(jk)(ik)\{(jl)^{\varepsilon},(kl)^{-1}(jl)^{-1}(il)^{-1}\}(jk)^{-1}(ik)^{-1}(ij)^{-1}\\
&\equiv&(ik)(jk)\{\{(jl)^{\varepsilon},(il)^{-1}(kl)^{-1}(il)(kl)\},\{(kl)^{-1},(il)(kl)\}\{(jl)^{-1},\\
&&(il)^{-1}(kl)^{-1}(il)(kl)\}\{(il)^{-1},(kl)\}\}\times(ik)(jk)^{-1}(ik)^{-1}(ij)^{-1}\\
&\equiv&(ik)\{\{(jl)^{\varepsilon},(kl)\},(il)^{-1}\{(kl)^{-1},(jl)(kl)\}(il)\{(jl)^{-1},(kl)\}(il)^{-1}\}\\
&&(jk)(ik)(jk)^{-1}(ik)^{-1}(ij)^{-1}\\
&\equiv&\{\{\{(jl)^{\varepsilon},(il)^{-1}(kl)^{-1}(il)(kl)\},\{(kl),(il)(kl)\}\},\{(il)^{-1},(kl)\}
\{\{(kl)^{-1},(il)(kl)\},\\
&&\{(jl),(il)^{-1}(kl)^{-1}(il)(kl)\}\{(kl),(il)(kl)\}\}\{(il),(kl)\}\{\{(jl)^{-1},(il)^{-1}(kl)^{-1}\\
&&(il)(kl)\},\{(kl),(il)(kl)\}\}\{(il)^{-1},(kl)\}\}(ik)(jk)(ik)(jk)^{-1}(ik)^{-1}(ji)^{-1}\\
&\equiv&\{(jl)^{\varepsilon},(il)^{-1}(jl)^{-1}(il)^{-1}(kl)^{-1}\{(il),(jl)^{-1}(il)^{-1}\}(kl)\}\\
&&\{(ik),(jk)^{-1}(ik)^{-1}\}(ij)^{-1}.
\end{eqnarray*}
\item[3)] \  $p>j,i<j<p<k<l$.
\begin{eqnarray*}
&&(ip)^{-1}\{(jl)^{\varepsilon},(il)^{-1}(kl)^{-1}(il)(kl)\}(ik)\\
&\equiv&\{\{(jl)^{\varepsilon},(pl)(il)(pl)^{-1}(il)^{-1}\},\{(il)^{-1},(pl)^{-1}(il)^{-1}\}(kl)^{-1}\\
&&\{(il),(pl)^{-1}(il)^{-1}\}(kl)\}\{(ik),(pk)^{-1}(jk)^{-1}\}(ip)^{-1}\\
&\equiv&\{(jl)^{\varepsilon},(il)^{-1}(kl)^{-1}(il)(pl)(il)(pl)^{-1}(il)^{-1}(kl)\}
\{(ik),(pk)^{-1}(ik)^{-1}\}(ip)^{-1}\  \ and \\
&&\\
&&\{(ik),(pk)^{-1}(ik)^{-1}\}(ip)^{-1}(jl)^{\varepsilon}\\
&=&(ik)(pk)(ik)\{(jl)^{\varepsilon},\{(il),(kl)^{-1}(pl)^{-1}\}(il)^{-1}\}(pk)^{-1}(ik)^{-1}(ip)^{-1}\\
&\equiv&(ik)(pk)\{\{(jl)^{\varepsilon},(il)^{-1}(kl)^{-1}(il)(kl)\},\{\{(il),(kl)\},\{(kl)^{-1},(il)(kl)\}\\
&&\{(pl)^{-1},(il)^{-1}(kl)^{-1}(il)(kl)\}\{(il)^{-1},(kl)\}\}(ik)(pk)^{-1}(ik)^{-1}(ip)^{-1}\\
&\equiv&(ik)\{(jl)^{\varepsilon},\{(pl),(kl)\}(il)\{(pl)^{-1},(kl)\}(il)^{-1}\}(pk)(ik)(pk)^{-1}(ik)^{-1}(ip)^{-1}\\
&\equiv&\{\{(jl)^{\varepsilon},(il)^{-1}(kl)^{-1}(il)(kl)\},\{\{(pl),(il)^{-1}(kl)^{-1}(il)(kl)\},\{(kl),(il)(kl)\}\}\\
&&\{(il),(kl)\}\{\{(pl)^{-1},(il)^{-1}(kl)^{-1}(il)(kl)\},\{(kl),(il)(kl)\}\}\{(il)^{-1},(kl)\}\}\\
&&(ik)(pk)(ik)(pk)^{-1}(ik)^{-1}(ip)^{-1}\\
&\equiv&\{(jl)^{\varepsilon},(il)^{-1}(kl)^{-1}(il)(pl)(il)(pl)^{-1}(il)^{-1}(kl)\}\{(ik),(pk)^{-1}(ik)^{-1}\}(ip)^{-1}.
\end{eqnarray*}
\end{enumerate}

$(9)\wedge(13)$\\
Let
$f=(ip)^{-1}(ik)^\delta-\{(ik)^\delta,(pk)^{-1}(ik)^{-1}\}(ip)^{-1},\
g=(ik)^\delta(jl)^{\varepsilon}-(jl)^{\varepsilon}(ik)^\delta,\
j<i<p<k<l\ \ or\ \  i<p<k<j<l$. Then
$w=(ip)^{-1}(ik)^\delta(jl)^{\varepsilon}$ and
\begin{eqnarray*}
(f,g)_w&=&(ip)^{-1}(jl)^{\varepsilon}(ik)^\delta-\{(ik)^\delta,(pk)^{-1}(ik)^{-1}\}(ip)^{-1}(jl)^{\varepsilon}.\\
&\equiv&(jl)^{\varepsilon}\{(ik)^\delta,(pk)^{-1}(ik)^{-1}\}(ip)^{-1}-
(jl)^{\varepsilon}\{(ik)^\delta,(pk)^{-1}(ik)^{-1}\}(ip)^{-1}\\
&\equiv&0.
\end{eqnarray*}

$(10)\wedge(7)$\\
Let $f=(ji)(jk)^{-1}-\{(jk)^{-1},(ik)\}(ji),\
g=(jk)^{-1}(kl)^{\varepsilon}-\{(kl)^{\varepsilon},(jl)^{-1}\}(jk)^{-1},\
j<i<k<l$. Then $w=(ji)(jk)^{-1}(kl)^{\varepsilon}$ and
$$
(f,g)_w=(ji)\{(kl)^{\varepsilon},(jl)^{-1}\}(jk)^{-1}-\{(jk)^{-1},(ik)\}(ji)(kl)^{\varepsilon}\equiv
0
$$
since
\begin{eqnarray*}
&&(ji)\{(kl)^{\varepsilon},(jl)^{-1}\}(jk)^{-1}
\equiv\{(kl)^{\varepsilon},(il)^{-1}(jl)^{-1}(il)\}\{(jk)^{-1},(ik)\}(ji) \ \ \ \ \ \ \ and \\
&&\\
&&\{(jk)^{-1},(ik)\}(ji)(kl)^{\varepsilon}\\
&\equiv&(ik)^{-1}(jk)^{-1}(ik)(kl)^{\varepsilon}(ji)\\
&\equiv&(ik)^{-1}(jk)^{-1}\{(kl)^{\varepsilon},(il)(kl)\}(ik)(ji)\\
&\equiv&(ik)^{-1}\{\{(kl)^{\varepsilon},(jl)^{-1}\},\{(il),(kl)(jl)(kl)^{-1}(jl)^{-1}\}\{(kl),(jl)^{-1}\}\}
(jk)^{-1}(ik)(ji)\\
&\equiv&\{\{(kl)^{\varepsilon},(il)^{-1}\},(jl)^{-1}\{(kl)^{-1},(il)^{-1}\}\{(il),(kl)^{-1}(il)^{-1}\}\\
&&\{(kl),(il)^{-1}\}\}(ik)^{-1}(jk)^{-1}(ik)(ji)\\
&\equiv&\{(kl)^{\varepsilon},(il)^{-1}(jl)^{-1}(il)\}\{(jk)^{-1},(ik)\}(ji).
\end{eqnarray*}

$(10)\wedge(8)$\\
Let $f=(ji)(jk)-\{(jk),(ik)\}(ji),\
g=(jk)(kl)^{\varepsilon}-\{(kl)^{\varepsilon},(jl)(kl)\}(jk),j<i<k<l$.
Then $w=(ji)(jk)(kl)^{\varepsilon}$ and
$$
(f,g)_w=(ji)\{(kl)^{\varepsilon},(jl)(kl)\}(jk)-\{(jk),(ik)\}(ji)(kl)^{\varepsilon}\equiv
0
$$
since
\begin{eqnarray*}
&&(ji)\{(kl)^{\varepsilon},(jl)(kl)\}(jk)\equiv\{(kl)^{\varepsilon},(il)^{-1}(jl)(il)(kl)\}\{(jk),(ik)\}(ji) \ \ \ \ \ and\\
&&\\
&&\{(jk),(ik)\}(ji)(kl)^{\varepsilon}\\
&\equiv&(ik)^{-1}(jk)\{(kl)^{\varepsilon},(il)(kl)\}(ik)(ji)\\
&\equiv&(ik)^{-1}\{\{(kl)^{\varepsilon},(jl)(kl)\},\{(il),(jl)^{-1}(kl)^{-1}(jl)(kl)\}\{(kl),(jl)(kl)\}\}(jk)(ik)(ji)\\
&\equiv&\{\{(kl)^{\varepsilon},(il)^{-1}\},(jl)\{(il),(kl)^{-1}(il)^{-1}\}\{(kl),(il)^{-1}\}\}(ik)^{-1}(jk)(ik)(ji)\\
&\equiv&\{(kl)^{\varepsilon},(il)^{-1}(jl)(il)(kl)\}\{(jk),(ik)\}(ji).
\end{eqnarray*}

$(10)\wedge(9)$\\
Let $f=(ji)(jk)^{-1}-\{(jk)^{-1},(ik)\}(ji),\
g=(jk)^{-1}(jl)^{\varepsilon}-\{(jl)^{\varepsilon},(kl)^{-1}(jl)^{-1}\}(jk)^{-1},\
j<i<k<l$. Then $w=(ji)(jk)^{-1}(jl)^{\varepsilon}$ and
$$
(f,g)_w=(ji)\{(jl)^{\varepsilon},(kl)^{-1}(jl)^{-1}\}(jk)^{-1}-\{(jk)^{-1},(ik)\}(ji)(jl)^{\varepsilon}\equiv
0
$$
since
\begin{eqnarray*}
&&(ji)\{(jl)^{\varepsilon},(kl)^{-1}(jl)^{-1}\}(jk)^{-1}\\
&&\equiv\{(jl)^{\varepsilon},(il)(kl)^{-1}(il)^{-1}(jl)^{-1}(il)\}\{(jk)^{-1},(ik)\}(ji)\ \ \ \ \ and \\
&&\\
&&\{(jk)^{-1},(ik)\}(ji)(jl)^{\varepsilon}\\
&\equiv&(ik)^{-1}(jk)^{-1}(ik)\{(jl)^{\varepsilon},(il)\}(ji)\\
&\equiv&(ik)^{-1}(jk)^{-1}\{(jl)^{\varepsilon},\{(il),(kl)\}\}(ik)(ji)\\
&\equiv&(ik)^{-1}\{\{(jl)^{\varepsilon},(kl)^{-1}(jl)^{-1}\},\{(kl)^{-1},(jl)^{-1}\}\\
&&\{(il),(kl)(jl)(kl)^{-1}(jl)^{-1}\}\{(kl),(jl)^{-1}\}\}(jk)^{-1}(ik)(ji)\\
&\equiv&\{(jl)^{\varepsilon},\{(kl)^{-1},(il)^{-1}\}(jl)^{-1}\{(kl)^{-1},(il)^{-1}\}\\
&&\{(il),(kl)^{-1}(il)^{-1}\}\{(kl),(il)^{-1}\}\}(ik)^{-1}(jk)^{-1}(ik)(ji)\\
&\equiv&\{(jl)^{\varepsilon},(il)(kl)^{-1}(il)^{-1}(jl)^{-1}(il)\}\{(jk)^{-1},(ik)\}(ji).
\end{eqnarray*}

$(10)\wedge(10)$\\
Let $f=(ji)(jk)-\{(jk),(ik)\}(ji),
g=(jk)(jl)^{\varepsilon}-\{(jl)^{\varepsilon},(kl)\}(jk), j<i<k<l.$
Then $w=(ji)(jk)(jl)^{\varepsilon}$ and
$$
(f,g)_w=(ji)\{(jl)^{\varepsilon},(kl)\}(jk)-\{(jk),(ik)\}(ji)(jl)^{\varepsilon}\equiv
0
$$
since
\begin{eqnarray*}
&&(ji)\{(jl)^{\varepsilon},(kl)\}(jk)\\
&\equiv&\{(jl)^{\varepsilon},(il)(kl)\}\{(jk),(ik)\}(ji)\ \ \ \ \ \ and \\
&&\\
&&\{(jk),(ik)\}(ji)(jl)^{\varepsilon}\\
&\equiv&(ik)^{-1}(jk)\{(jl)^{\varepsilon},\{(il),(kl)\}\}(ik)(ji)\\
&\equiv&(ik)^{-1}\{\{(jl)^{\varepsilon},(kl)\},\{\{(il),(jl)^{-1}(kl)^{-1}(jl)(kl)\},\{(kl),(jl)(kl)\}\}\}(jk)(ik)(ji)\\
&\equiv&\{(jl)^{\varepsilon},\{(il),(kl)^{-1}(il)^{-1}\}\{(kl),(il)^{-1}\}\}(ik)^{-1}(jk)(ik)(ji)\\
&\equiv&\{(jl)^{\varepsilon},(il)(kl)\}\{(jk),(ik)\}(ji).
\end{eqnarray*}

$(10)\wedge(11)$\\
Let $f=(ip)(ik)^{-1}-\{(ik)^{-1},(pk)\}(ip),
g=(ik)^{-1}(jl)^{\varepsilon}-\{(jl)^{\varepsilon},(kl)(il)(kl)^{-1}$\\
$(il)^{-1}\}(ik)^{-1}, i<p<k,i<j<k<l.$ Then
 $w=(ip)(ik)^{-1}(jl)^{\varepsilon}$ and
$$
(f,g)_w=(ip)\{(jl)^{\varepsilon},(kl)(il)(kl)^{-1}(il)^{-1}\}(ik)^{-1}-\{(ik)^{-1},(pk)\}(ip)(jl)^{\varepsilon}.
$$
There are three cases to consider.
\begin{enumerate}
\item[1)]\ $p<j,i<p<j<k<l$.
\begin{eqnarray*}
&&(ip)\{(jl)^{\varepsilon},(kl)(il)(kl)^{-1}(il)^{-1}\}(ik)^{-1}\\
&\equiv&\{(jl)^{\varepsilon},(kl)\{(il),(pl)\}(kl)^{-1}\{(il)^{-1},(pl)\}\}\{(ik)^{-1},(pk)\}(ip)\ \ \ \ \ and \\
&&\\
&&\{(ik)^{-1},(pk)\}(ip)(jl)^{\varepsilon}\\
&\equiv&(pk)^{-1}(ik)^{-1}\{(jl)^{\varepsilon},(pl)^{-1}(kl)^{-1}(pl)(kl)\}(pk)(ip)\\
&\equiv&(pk)^{-1}\{\{(jl)^{\varepsilon},(kl)(il)(kl)^{-1}(il)^{-1}\},\{(pl)^{-1},(kl)(il)(kl)^{-1}(il)^{-1}\}\\
&&\{(kl)^{-1},(il)^{-1}\}\{(pl),(kl)(il)(kl)^{-1}(il)^{-1}\}\{(kl),(il)^{-1}\}(ik)^{-1}(pk)(ip)\\
&\equiv&\{\{(jl)^{\varepsilon},(kl)(pl)(kl)^{-1}(pl)^{-1}\},\{\{(kl)^{-1},(pl)^{-1}\}(il)^{-1}\{(kl)^{-1},(pl)^{-1}\}\\
&&\{(pl),(kl)^{-1}(pl)^{-1}\}\}\{(kl),(pl)^{-1}\}\}(pk)^{-1}(ik)^{-1}(pk)(ip)\\
&\equiv&\{(jl)^{\varepsilon},(kl)\{(il),(pl)\}(kl)^{-1}\{(il)^{-1},(pl)\}\}\{(ik)^{-1},(pk)\}(ip).
\end{eqnarray*}
\item[2)]\ $p=j,i<p=j<k<l$.
\begin{eqnarray*}
&&(ip)\{(jl)^{\varepsilon},(kl)(il)(kl)^{-1}(il)^{-1}\}(ik)^{-1}\\
&=&(ij)\{(jl)^{\varepsilon},(kl)(il)(kl)^{-1}(il)^{-1}\}(ik)^{-1}\\
&\equiv&\{(jl)^{\varepsilon},(il)(jl)(kl)\{(il),(jl)\}(kl)^{-1}\{(il)^{-1},(jl)\}\}\{(ik)^{-1},(jk)\}(ij)\ \ \ \ \ and \\
&&\\
&&\{(ik)^{-1},(pk)\}(ip)(jl)^{\varepsilon}=(jk)^{-1}(ik)^{-1}(jk)(ij)(jl)^{\varepsilon}\\
&\equiv&(jk)^{-1}(ik)^{-1}(jk)\{(jl)^{\varepsilon},(il)(jl)\}(ij)\\
&\equiv&(jk)^{-1}(ik)^{-1}\{\{(jl)^{\varepsilon},(kl)\},(il)\{(jl),(kl)\}\}(jk)(ij)\\
&\equiv&(jk)^{-1}\{\{(jl)^{\varepsilon},(kl)(il)(kl)^{-1}(il)^{-1}\},\{(kl),(il)^{-1}\}\{(il),(kl)^{-1}(il)^{-1}\}\\
&&\{\{(jl),(kl)(il)(kl)^{-1}(il)^{-1}\},\{(kl),(il)^{-1}\}\}\}(ik)^{-1}(jk)(ij)\\
&\equiv&\{\{(jl)^{\varepsilon},(kl)^{-1}(jl)^{-1}\},\{(il),\{(kl)^{-1},(jl)^{-1}\}(il)^{-1}\{(kl)^{-1},(jl)^{-1}\}\}\\
&&\{(jl),(kl)^{-1}(jl)^{-1}\}\{(kl),(jl)^{-1}\}\}(jk)^{-1}(ik)^{-1}(ik)(ij)\\
&\equiv&\{(jl)^{\varepsilon},(il)(jl)(kl)\{(il),(jl)\}(kl)^{-1}\{(il)^{-1},(jl)\}\}\{(ik)^{-1},(jk)\}(ij).
\end{eqnarray*}
\item[3)]\ $p>j,i<j<p<k<l$.
\begin{eqnarray*}
&&(ip)\{(jl)^{\varepsilon},(kl)(il)(kl)^{-1}(il)^{-1}\}(ik)^{-1}\\
&\equiv&\{(jl)^{\varepsilon},(il)^{-1}(pl)^{-1}(il)(pl)(kl)\{(il),(pl)\}(kl)^{-1}\\
&&\{(il)^{-1},(pl)\}\}\{(ik)^{-1},(pk)\}(ip)\ \ \ and \\
&&\\
&& \{(ik)^{-1},(pk)\}(ip)(jl)^{\varepsilon}\\
&\equiv&(pk)^{-1}(ik)^{-1}(pk)\{(jl)^{\varepsilon},(il)^{-1}(pl)^{-1}(il)(pl)\}(ip)\\
&\equiv&(pk)^{-1}(ik)^{-1}\{(jl)^{\varepsilon},(il)^{-1}\{(pl)^{-1},(kl)\}(il)\{(pl),(kl)\}\}(pk)(ip)\\
&\equiv&(pk)^{-1}\{\{(jl)^{\varepsilon},(kl)(il)(kl)^{-1}(il)^{-1}\},\{(il)^{-1},(kl)^{-1}(il)^{-1}\}\{\{(pl)^{-1},\\
&&(kl)(il)(kl)^{-1}(il)^{-1}\},\{(kl),(il)^{-1}\}\}\{(il),(kl)^{-1}(jl)^{-1}\}\{\{(pl),(kl)(il)\\
&&(kl)^{-1}(il)^{-1}\},\{(kl),(il)^{-1}\}\}\}(ik)^{-1}(pk)(ip)\\
&\equiv&\{(jl)^{\varepsilon},(il)^{-1}\{\{(pl)^{-1},(kl)^{-1}(pl)^{-1}\},\{(kl),(pl)^{-1}\}\}
\{(il),\{(kl)^{-1},(pl)^{-1}\}\\
&&(il)^{-1}\}\{\{(pl),(kl)^{-1}(pl)^{-1}\},\{(kl),(pl)^{-1}\}\}\}(pk)^{-1}(ik)^{-1}(pk)(ip)\\
&\equiv&\{(jl)^{\varepsilon},(il)^{-1}(pl)^{-1}(il)(pl)(kl)\{(il),(pl)\}(kl)^{-1}\{(il)^{-1},(pl)\}\}
\{(ik)^{-1},(pk)\}(ip).
\end{eqnarray*}
\end{enumerate}

$(10)\wedge(12)$\\
Let $f=(ip)(ik)-\{(ik),(pk)\}(ip),\
g=(ik)(jl)^{\varepsilon}-\{(jl)^{\varepsilon},(il)^{-1}(kl)^{-1}(il)(kl)\}(ik),i<p<k,i<j<k<l$.
Then $w=(ip)(ik)(jl)^{\varepsilon}$ and
$$
(f,g)_w=(ip)\{(jl)^{\varepsilon},(il)^{-1}(kl)^{-1}(il)(kl)\}(ik)-\{(ik),(pk)\}(ip)(jl)^{\varepsilon}.
$$
There are three cases to consider.
\begin{enumerate}
\item[1)]\ $p<j,i<p<j<k<l$.
\begin{eqnarray*}
&&(ip)\{(jl)^{\varepsilon},(il)^{-1}(kl)^{-1}(il)(kl)\}(ik)\\
&\equiv&\{(jl)^{\varepsilon},\{(il)^{-1},(pl)\}(kl)^{-1}\{(il),(pl)\}(kl)\}\{(ik),(pk)\}(ip)  \ \ \ \ \ \ and \\
&&\\
&&\{(ik),(pk)\}(ip)(jl)^{\varepsilon}\\
&\equiv&(pk)^{-1}(ik)\{(jl)^{\varepsilon},(pl)^{-1}(kl)^{-1}(pl)(kl)\}(pk)(ip)\\
&\equiv&(pk)^{-1}\{\{(jl)^{\varepsilon},(il)^{-1}(kl)^{-1}(il)(kl)\},\{(pl)^{-1},(il)^{-1}(kl)^{-1}(il)(kl)\}\\
&&\{(kl)^{-1},(il)(kl)\}\{(pl),(il)^{-1}(kl)^{-1}(il)(kl)\}\{(kl),(il)(kl)\}\}(ik)(pk)(ip)\\
&\equiv&\{\{(jl)^{\varepsilon},(kl)(pl)(kl)^{-1}(pl)^{-1}\},\{\{(kl)^{-1},(pl)^{-1}\},(il)\{(pl),(kl)^{-1}(pl)^{-1}\}\}\\
&&\{(kl),(pl)^{-1}\}\}(jk)^{-1}(ik)(pk)(ip)\\
&\equiv&\{(jl)^{\varepsilon},\{(il)^{-1},(pl)\}(kl)^{-1}\{(il),(pl)\}(kl)\}\{(ik),(pk)\}(ip).
\end{eqnarray*}
\item[2)]\ $p=j,i<p=j<k<l$.
\begin{eqnarray*}
&&(ip)\{(jl)^{\varepsilon},(il)^{-1}(kl)^{-1}(il)(kl)\}(ik)\\
&=&(ij)\{(jl)^{\varepsilon},(il)^{-1}(kl)^{-1}(il)(kl)\}(ik)\\
&\equiv&\{\{(jl)^{\varepsilon},(il)(jl)\},\{(il)^{-1},(jl)\}(kl)^{-1}\{(il),(jl)\}(kl)\}\{(ik),(jk)\}(ij)\\
&\equiv&\{(jl)^{\varepsilon},(kl)^{-1}(jl)^{-1}(il)(jl)(kl)\}(jk)^{-1}(ik)(jk)(ij) \ \ \ \ \ \ and\\
&&\\
&&\{(ik),(pk)\}(ip)(jl)^{\varepsilon}\\
&\equiv&(jk)^{-1}(ik)(jk)(ij)(jl)^{\varepsilon}\\
&\equiv&(jk)^{-1}(ik)\{(jl)^{\varepsilon},(kl)(il)\{(jl),(kl)\}\}(jk)(ij)\\
&\equiv&(jk)^{-1}\{\{(jl)^{\varepsilon},(il)^{-1}(kl)^{-1}(il)(kl)\},\{(kl),(il)(kl)\}\{(il),(kl)\}\\
&&\{\{(jl),(il)^{-1}(kl)^{-1}(il)(kl)\},\{(kl),(il)(kl)\}\}\}(ik)(jk)(ij)\\
&\equiv&\{\{(jl)^{\varepsilon},(kl)^{-1}(jl)^{-1}\},(il)\{(jl),(kl)^{-1}(jl)^{-1}\}\{(kl),(jl)^{-1}\}\}
(jk)^{-1}(ik)(jk)(ij)\\
&\equiv&\{(jl)^{\varepsilon},(kl)^{-1}(jl)^{-1}(il)(jl)(kl)\}(jk)^{-1}(ik)^{-1}(jk)(ij).
\end{eqnarray*}
\item[3)]\ $p>j,i<j<p<k<l$.
\begin{eqnarray*}
&&(ip)\{(jl)^{\varepsilon},(il)^{-1}(kl)^{-1}(il)(kl)\}(ik)\\
&\equiv&\{\{(jl)^{\varepsilon},(il)^{-1}(pl)^{-1}(il)(pl)\},\{(il)^{-1},(pl)\}(kl)^{-1}\{(il),(pl)\}(kl)\}
\{(ik),(pk)\}(ip)\\
&\equiv&\{(jl)^{\varepsilon},(il)^{-1}(kl)^{-1}(pl)^{-1}(il)(pl)(kl)\}(pk)^{-1}(ik)(pk)(ip)\ \ \ \ \ \ and \\
&&\\
&&\{(ik),(pk)\}(ip)(jl)^{\varepsilon}\\
&\equiv&(pk)^{-1}(ik)\{(jl)^{\varepsilon},(il)^{-1}\{(pl)^{-1},(kl)\}(il)\{(pl),(kl)\}\}(pk)(ip)\\
&\equiv&(pk)^{-1}\{\{(jl)^{\varepsilon},(il)^{-1}(kl)^{-1}(il)(kl)\},\{(il)^{-1},(kl)\}\\
&&\{\{(pl)^{-1},(il)^{-1}(kl)^{-1}(il)(kl)\},\{(kl),(il)(kl)\}\}\{(il),(kl)\}\\
&&\{\{(pl),(il)^{-1}(kl)^{-1}(il)(kl)\},\{(kl),(il)(kl)\}\}\}(ik)(pk)(ip)\\
&\equiv&\{(jl)^{\varepsilon},(il)^{-1}\{(kl)^{-1},(pl)^{-1}\}\{(pl)^{-1},(kl)^{-1}(pl)^{-1}\}(il)\\
&&\{(pl),(kl)^{-1}(pl)^{-1}\}\{(kl),(pl)^{-1}\}\}(pk)^{-1}(ik)(pk)(ip)\\
&\equiv&\{(jl)^{\varepsilon},(il)^{-1}(kl)^{-1}(pl)^{-1}(il)(pl)(kl)\}(pk)^{-1}(ik)(pk)(ip).
\end{eqnarray*}
\end{enumerate}

$(10)\wedge(13)$\\
Let $f=(ip)(ik)^\delta-\{(ik)^\delta,(pk)\}(ip),\
g=(ik)^\delta(jl)^{\varepsilon}-(jl)^{\varepsilon}(ik)^\delta,\
j<i<p<k<l,or,i<p<k<j<l$. Then $w=(ip)(ik)^\delta(jl)^{\varepsilon}$
and
\begin{eqnarray*}
(f,g)_w&=&(ip)(jl)^{\varepsilon}(ik)^\delta-\{(ik)^\delta,(pk)\}(ip)(jl)^{\varepsilon}\\
&\equiv&(jl)^{\varepsilon}(ip)(ik)^\delta-(jl)^{\varepsilon}\{(ik)^\delta,(pk)\}(ip)\\
&\equiv&0.
\end{eqnarray*}

$(11)\wedge(7)$\\
Let
$f=(pq)^{-1}(jk)^{-1}-\{(jk)^{-1},(qk)(pk)(qk)^{-1}(pk)^{-1}\}(pq)^{-1},\
\
g=(jk)^{-1}(kl)^{\varepsilon}-\{(kl)^{\varepsilon},$\\
$(jl)^{-1}\}(jk)^{-1},\ p<j<q<k<l$. Then
$w=(pq)^{-1}(jk)^{-1}(kl)^{\varepsilon}$ and
$$
(f,g)_w=(pq)^{-1}\{(kl)^{\varepsilon},(jl)^{-1}\}(jk)^{-1}-
\{(jk)^{-1},(qk)(pk)(qk)^{-1}(pk)^{-1}\}(pq)^{-1}(kl)^{\varepsilon}\equiv
0
$$
since
\begin{eqnarray*}
&&(pq)^{-1}\{(kl)^{\varepsilon},(jl)^{-1}\}(jk)^{-1}\\
&\equiv&\{(kl)^{\varepsilon},\{(jl)^{-1},(ql)(pl)(ql)^{-1}(pl)^{-1}\}\}
\{(jk)^{-1},(qk)(pk)(qk)^{-1}(pk)^{-1}\}(pq)^{-1} \ \ \ \ and\\
&&\\
&&\{(jk)^{-1},(qk)(pk)(qk)^{-1}(pk)^{-1}\}(pq)^{-1}(kl)^{\varepsilon}\\
&\equiv&(pk)(qk)(pk)^{-1}(qk)^{-1}(jk)^{-1}(qk)(pk)(qk)^{-1}(pk)^{-1}(kl)^{\varepsilon}(pq)^{-1}\\
&\equiv&(pk)(qk)(pk)^{-1}(qk)^{-1}(jk)^{-1}(qk)(pk)(qk)^{-1}\{(kl)^{\varepsilon},(pl)^{-1}\}(pk)^{-1}(pq)^{-1}\\
&\equiv&(pk)(qk)(pk)^{-1}(qk)^{-1}(jk)^{-1}(qk)(pk)\{\{(kl)^{\varepsilon},(ql)^{-1}\},(pl)^{-1}\}
(qk)^{-1}(pk)^{-1}(pq)^{-1}\\
&\equiv&(pk)(qk)(pk)^{-1}(qk)^{-1}(jk)^{-1}(qk)\{\{(kl)^{\varepsilon},(pl)(kl)\},\{(ql)^{-1},(pl)^{-1}(kl)^{-1}(pl)(kl)\}\\
&&\{(pl)^{-1},(kl)\}\}(pk)(qk)^{-1}(pk)^{-1}(pq)^{-1}\\
&\equiv&(pk)(qk)(pk)^{-1}(qk)^{-1}(jk)^{-1}\{\{(kl)^{\varepsilon},(ql)(kl)\},(pl)\{(ql)^{-1},(kl)\}(pl)^{-1}\}\\
&&(qk)(pk)(qk)^{-1}(pk)^{-1}(pq)^{-1}\\
&\equiv&(pk)(qk)(pk)^{-1}(qk)^{-1}\{\{\{(kl)^{\varepsilon},(jl)^{-1}\},\{(ql),(kl)(jl)(kl)^{-1}(jl)^{-1}\}\\
&&\{(kl),(jl)^{-1}\}\},(pl)\{\{(ql)^{-1},(kl)(jl)(kl)^{-1}(jl)^{-1}\},\{(kl),(jl)^{-1}\}\}(pl)^{-1}\}\\
&&(jk)^{-1}(qk)(pk)(qk)^{-1}(pk)^{-1}(pq)^{-1}\\
&\equiv&(pk)(qk)(pk)^{-1}\{\{(kl)^{\varepsilon},(ql)^{-1}\},(jl)^{-1}\{(pl),\{(kl)^{-1},(ql)^{-1}\}\\
&&\{(ql)^{-1},(kl)^{-1}(ql)^{-1}\}\{(kl),(ql)^{-1}\}\}(pl)^{-1}\}(qk)^{-1}(jk)^{-1}(qk)(pk)(qk)^{-1}(pk)^{-1}(pq)^{-1}\\
&\equiv&(pk)(qk)\{\{(kl)^{\varepsilon},(pl)^{-1}\},\{(ql)^{-1},(kl)(pl)(kl)^{-1}(pl)^{-1}\}\{(jl)^{-1},(kl)(pl)\\
&&(kl)^{-1}(pl)^{-1}\}\{(ql),(kl)(pl)(kl)^{-1}(pl)^{-1}\}\{(pl),(kl)^{-1}(pl)^{-1}\}\{(ql)^{-1},(kl)(pl)\\
&&(kl)^{-1}(pl)^{-1}\}\{pl()^{-1},(kl)^{-1}(pl)^{-1}\}\}(pk)^{-1}(qk)^{-1}(jk)^{-1}(qk)(pk)(qk)^{-1}(pk)^{-1}(pq)^{-1}\\
&\equiv&(pk)\{\{(kl)^{\varepsilon},(ql)(kl)\},\{(jl)^{-1},\{(ql),(kl)\}\{(kl),(ql)(kl)\}(pl)\}\{(kl)^{-1},(ql)(kl)\}\\
&&\{(ql)^{-1},(kl)\}(pl)^{-1}\}(qk)(pk)^{-1}(qk)^{-1}(jk)^{-1}(qk)(pk)(qk)^{-1}(pk)^{-1}(pq)^{-1}\\
&\equiv&\{\{(kl)^{\varepsilon},(pl)(kl)\},\{\{(jl)^{-1},(pl)^{-1}(kl)^{-1}(pl)(kl)\},\{(ql),(pl)^{-1}(kl)^{-1}(pl)(kl)\}\\
&&\{(kl),(pl)(kl)\}\{(pl),(kl)\}\{(kl)^{-1},(pl)(kl)\}\{(ql)^{-1},(pl)^{-1}(kl)^{-1}(pl)(kl)\}\\
&&\{(pl)^{-1},(kl)\}\{(jk)^{-1},(qk)(pk)(qk)^{-1}(pk)^{-1}\}(pq)^{-1}\\
&\equiv&\{(kl)^{\varepsilon},\{(jl)^{-1},(ql)(pl)(ql)^{-1}(pl)^{-1}\}\}
\{(jk)^{-1},(qk)(pk)(qk)^{-1}(pk)^{-1}\}(pq)^{-1}.
\end{eqnarray*}

$(11)\wedge(8)$\\
Let $f=(pq)^{-1}(jk)-\{(jk),(qk)(pk)(qk)^{-1}(pk)^{-1}\}(pq)^{-1},\
g=(jk)(kl)^{\varepsilon}-\{(kl)^{\varepsilon},(jl)(kl)\}(jk),\
p<j<q<k<l$. Then $w=(pq)^{-1}(jk)(kl)^{\varepsilon}$ and
$$
(f,g)_w=(pq)^{-1}\{(kl)^{\varepsilon},(jl)(kl)\}(jk)-\{(jk),(qk)(pk)(qk)^{-1}(pk)^{-1}\}(pq)^{-1}(kl)^{\varepsilon}\equiv
0
$$
since
\begin{eqnarray*}
&&(pq)^{-1}\{(kl)^{\varepsilon},(jl)(kl)\}(jk)\\
&\equiv&\{(kl)^{\varepsilon},\{(jl),(ql)(pl)(ql)^{-1}(pl)^{-1}\}(kl)\}
\{(jk),(qk)(pk)(qk)^{-1}(pk)^{-1}\}(pq)^{-1} \ \ \ \ and \\
&&\\
&&\{(jk),(qk)(pk)(qk)^{-1}(pk)^{-1}\}(pq)^{-1}(kl)^{\varepsilon}\\
&\equiv&(pk)(qk)(pk)^{-1}(qk)^{-1}(jk)\{(kl)^{\varepsilon},(ql)(kl)(pl)\{(ql)^{-1},(kl)\}(pl)^{-1}\}\\
&&(qk)(pk)(qk)^{-1}(pk)^{-1}(pq)^{-1}\\
&\equiv&(pk)(qk)(pk)^{-1}(qk)^{-1}\{\{(kl)^{\varepsilon},(jl)(kl)\},\{(ql),(jl)^{-1}(kl)^{-1}(jl)(kl)\}\\
&&\{(kl),(jl)(kl)\}(pl)\{\{(ql)^{-1},(jl)^{-1}(kl)^{-1}(jl)(kl)\}\{(kl),(jl)(kl)\}\}(pl)^{-1}\}\\
&&(jk)(qk)(pk)(qk)^{-1}(pk)^{-1}(pq)^{-1}\\
&\equiv&(pk)(qk)(pk)^{-1}\{\{(kl)^{\varepsilon},(ql)^{-1}\},(jl)\{(pl),\{(kl)^{-1},(ql)^{-1}\}
\{(ql)^{-1},(kl)^{-1}\\
&&(ql)^{-1}\}\}\{(kl),(ql)^{-1}\}(pl)^{-1}\}(qk)^{-1}(jk)(qk)(pk)(qk)^{-1}(pk)^{-1}(pq)^{-1}\\
&\equiv&(pk)(qk)\{\{(kl)^{\varepsilon},(pl)^{-1}\},\{(ql)^{-1},(kl)(pl)(kl)^{-1}(pl)^{-1}\}\{(jl),(kl)(pl)(kl)^{-1}\\
&&(pl)^{-1}\}\{(ql),(kl)(pl)(kl)^{-1}(pl)^{-1}\}\{(kl),(pl)^{-1}\}\{(pl),(kl)^{-1}(pl)^{-1}\}\{(ql)^{-1},(kl)(pl)\\
&&(kl)^{-1}(pl)^{-1}\}\{(pl)^{-1},(kl)^{-1}(pl)^{-1}\}\}(pk)^{-1}(qk)^{-1}(jk)(qk)(pk)(qk)^{-1}(pk)^{-1}(pq)^{-1}\\
&\equiv&(pk)\{\{(kl)^{\varepsilon},(ql)(kl)\},\{(jl),\{(ql),(kl)\}\{(kl),(ql)(kl)\}(pl)\}\{(ql)^{-1},(kl)\}(pl)^{-1}\}\\
&&(qk)(pk)^{-1}(qk)^{-1}(jk)(qk)(pk)(qk)^{-1}(pk)^{-1}(pq)^{-1}\\
&\equiv&\{\{(kl)^{\varepsilon},(pl)(kl)\},\{\{(jl),(pl)^{-1}(kl)^{-1}(pl)(kl)\},\{(ql),(pl)^{-1}(kl)^{-1}(pl)(kl)\}\\
&&\{(kl),(pl)(kl)\}\{(pl),(kl)\}\{(kl)^{-1},(pl)(kl)\}\{(ql)^{-1},(pl)^{-1}(kl)^{-1}(pl)(kl)\}\}\\
&&\{(kl),(pl)(kl)\}\{(pl)^{-1},(kl)\}\}\{(jk),(qk)(pk)(qk)^{-1}(pk)^{-1}\}(pq)^{-1}\\
&\equiv&\{(kl)^{\varepsilon},\{(jl),(ql)(pl)(ql)^{-1}(pl)^{-1}\}(kl)\}\{(jk),(qk)(pk)(qk)^{-1}(pk)^{-1}\}(pq)^{-1}.
\end{eqnarray*}

$(11)\wedge(9)$\\
Let
$f=(pq)^{-1}(jk)^{-1}-\{(jk)^{-1},(qk)(pk)(qk)^{-1}(pk)^{-1}\}(pq)^{-1},\
g=(jk)^{-1}(jl)^{\varepsilon}-\{(jl)^{\varepsilon},(kl)^{-1}$\\
$(jl)^{-1}\}(jk)^{-1},\ p<j<q<k<l$. Then
$w=(pq)^{-1}(jk)^{-1}(jl)^{\varepsilon}$ and
$$
(f,g)_w=(pq)^{-1}\{(jl)^{\varepsilon},(kl)^{-1}(jl)^{-1}\}(jk)^{-1}-
\{(jk)^{-1},(qk)(pk)(qk)^{-1}(pk)^{-1}\}(pq)^{-1}(jl)^{\varepsilon}\equiv
0
$$
since
\begin{eqnarray*}
&&(pq)^{-1}\{(jl)^{\varepsilon},(kl)^{-1}(jl)^{-1}\}(jk)^{-1}\\
&\equiv&\{(jl)^{\varepsilon},(ql)(pl)(ql)^{-1}(pl)^{-1}(kl)^{-1}\{(jl)^{-1},(ql)(pl)(ql)^{-1}(pl)^{-1}\}\}\\
&&\{(jk)^{-1},(qk)(pk)(qk)^{-1}(pk)^{-1}\}(pq)^{-1} \ \ \ \ \ \ \ \ \ and \\
&&\\
&&\{(jk)^{-1},(qk)(pk)(qk)^{-1}(pk)^{-1}\}(pq)^{-1}(jl)^{\varepsilon}\\
&\equiv&(pk)(qk)(pk)^{-1}(qk)^{-1}(jk)^{-1}(qk)(pk)(qk)^{-1}(pk)^{-1}\{(jl)^{\varepsilon},(ql)(pl)(ql)^{-1}(pl)^{-1}\}
(pq)^{-1}\\
&\equiv&(pk)(qk)(pk)^{-1}(qk)^{-1}(jk)^{-1}(qk)(pk)(qk)^{-1}\{\{(jl)^{\varepsilon},(kl)(pl)(kl)^{-1}(pl)^{-1}\},\\
&&\{(ql),(kl)(pl)(kl)^{-1}(pl)^{-1}\}\{(pl),(kl)^{-1}(pl)^{-1}\}\{(ql)^{-1},(kl)(pl)(kl)^{-1}(pl)^{-1}\}\\
&&\{(pl)^{-1},(kl)^{-1}(pl)^{-1}\}\}(pk)^{-1}(pq)^{-1}\\
&\equiv&(pk)(qk)(pk)^{-1}(qk)^{-1}(jk)^{-1}(qk)(pk)\{(jl)^{\varepsilon},\{(pl),\{(kl)^{-1},(ql)^{-1}\}\\
&&\{(ql)^{-1},(kl)^{-1}(ql)^{-1}\}\}(pl)^{-1}\}(qk)^{-1}(pk)^{-1}(pq)^{-1}\\
&\equiv&(pk)(qk)(pk)^{-1}(qk)^{-1}(jk)^{-1}(qk)\{\{(jl)^{\varepsilon},(pl)^{-1}(kl)^{-1}(pl)(kl)\},\{\{(pl),(kl)\},\\
&&\{(kl)^{-1},(pl)(kl)\}\{(ql)^{-1},(pl)^{-1}(kl)^{-1}(pl)(kl)\}\}\{(pl)^{-1},(kl)\}\}(pk)(qk)^{-1}(pk)^{-1}(pq)^{-1}\\
&\equiv&(pk)(qk)(pk)^{-1}(qk)^{-1}(jk)^{-1}\{(jl)^{\varepsilon},\{(ql),(kl)\}(pl)\{(ql)^{-1},(kl)\}(pl)^{-1}\}\\
&&(qk)(pk)(qk)^{-1}(pk)^{-1}(pq)^{-1}\\
&\equiv&(pk)(qk)(pk)^{-1}(qk)^{-1}\{\{(jl)^{\varepsilon},(kl)^{-1}(jl)^{-1}\},\{\{(ql),(kl)(jl)(kl)^{-1}(jl)^{-1}\},\\
&&\{(kl),(jl)^{-1}\}\}(pl)\{\{(ql)^{-1},(kl)(jl)(kl)^{-1}(jl)^{-1}\},\{(kl),(jl)^{-1}\}\}(pl)^{-1}\}\\
&&(jk)^{-1}(qk)(pk)(qk)^{-1}(pk)^{-1}(pq)^{-1}\\
&\equiv&(pk)(qk)(pk)^{-1}\{(jl)^{\varepsilon},\{(kl)^{-1},(ql)^{-1}\}(jl)^{-1}\{\{(ql),(kl)^{-1}(ql)^{-1}\},\\
&&\{(kl),(ql)^{-1}\}\}(pl)\{\{(ql)^{-1},(kl)^{-1}(ql)^{-1}\},\{(kl),(ql)^{-1}\}\}(pl)^{-1}\}\\
&&(qk)^{-1}(jk)^{-1}(qk)(pk)(qk)^{-1}(pk)^{-1}(pq)^{-1}\\
&\equiv&(pk)(qk)\{\{(jl)^{\varepsilon},(kl)(pl)(kl)^{-1}(pl)^{-1}\},\{(ql),(kl)(pl)(kl)^{-1}(pl)^{-1}\}\\
&&\{(kl)^{-1},(pl)^{-1}\}\{(ql)^{-1},(kl)(pl)(kl)^{-1}(pl)^{-1}\}\{(jl)^{-1},(kl)(pl)(kl)^{-1}(pl)^{-1}\}\\
&&\{(ql),(kl)(pl)(kl)^{-1}(pl)^{-1}\}\{(pl),(kl)^{-1}(pl)^{-1}\}\{(ql)^{-1},(kl)(pl)(kl)^{-1}(pl)^{-1}\}\\
&&\{(pl)^{-1},(kl)^{-1}(pl)^{-1}\}\},(pk)^{-1}(qk)^{-1}(jk)^{-1}(qk)(pk)(qk)^{-1}(pk)^{-1}(pq)^{-1}\\
&\equiv&(pk)\{(jl)^{\varepsilon},\{(ql),(kl)\}\{(kl),(ql)(kl)\}(pl)\{(kl),(ql)(kl)\}\{(jl)^{-1},\{(ql),(kl)\}\\
&&\{(kl),(ql)(kl)\}(pl)\}\{(kl)^{-1},(ql)(kl)\}\{(ql)^{-1},(kl)\}(pl)^{-1}\}\\
&&(qk)(pk)^{-1}(qk)^{-1}(jk)^{-1}(qk)(pk)(qk)^{-1}(pk)^{-1}(pq)^{-1}\\
&\equiv&\{\{(jl)^{\varepsilon},(pl)^{-1}(kl)^{-1}(pl)(kl)\},
\{\{(pl),(kl)\},\{(kl)^{-1},(pl)(kl)\}\{(ql)^{-1},(pl)^{-1}(kl)^{-1}\\
&&(pl)(kl)\}\}\{(kl),(pl)(kl)\}\{\{(jl)^{-1},(pl)^{-1}(kl)^{-1}(pl)(kl)\},\{\{(pl),(kl)\},\{(kl)^{-1},(pl)(kl)\}\\
&&\{(ql)^{-1},(pl)^{-1}(kl)^{-1}(pl)(kl)\}\}\}\{(pl)^{-1},(kl)\}\{(jk)^{-1},(qk)(pk)(qk)^{-1}(pk)^{-1}\}(pq)^{-1}\\
&\equiv&\{\{(jl)^{\varepsilon},(ql)(pl)(ql)^{-1}(pl)^{-1}(kl)^{-1}\{(jl)^{-1},(ql)(pl)(ql)^{-1}(pl)^{-1}\}\}\\
&&\{(jk)^{-1},(qk)(pk)(qk)^{-1}(pk)^{-1}\}(pq)^{-1}.
\end{eqnarray*}

$(11)\wedge(10)$\\
Let $f=(pq)^{-1}(jk)-\{(jk),(qk)(pk)(qk)^{-1}(pk)^{-1}\}(pq)^{-1},\
g=(jk)(jl)^{\varepsilon}-\{(jl)^{\varepsilon},(kl)\}(jk),\
p<j<q<k<l$. Then $w=(pq)^{-1}(jk)(jl)^{\varepsilon}$ and
$$
(f,g)_w=(pq)^{-1}\{(jl)^{\varepsilon},(kl)\}(jk)-\{(jk),(qk)(pk)(qk)^{-1}(pk)^{-1}\}(pq)^{-1}(jl)^{\varepsilon}\equiv
0
$$
since
\begin{eqnarray*}
&&(pq)^{-1}\{(jl)^{\varepsilon},(kl)\}(jk)\\
&\equiv&\{(jl)^{\varepsilon},(ql)(pl)(ql)^{-1}(pl)^{-1}(kl)\}\{(jk),(qk)(pk)(qk)^{-1}(pk)^{-1}\}(pq)^{-1} \ \ \ \ \ and\\
&&\\
&&\{(jk),(qk)(pk)(qk)^{-1}(pk)^{-1}\}(pq)^{-1}(jl)^{\varepsilon}\\
&\equiv&(pk)(qk)(pk)^{-1}(qk)^{-1}(jk)\{(jl)^{\varepsilon},\{(ql),(kl)\}(pl)\{(ql)^{-1},(kl)\}(pl)^{-1}\}\\
&&(qk)(pk)(qk)^{-1}(pk)^{-1}(pq)^{-1}\\
&\equiv&(pk)(qk)(pk)^{-1}(qk)^{-1}\{\{(jl)^{\varepsilon},(kl)\},\{\{(ql),(jl)^{-1}(kl)^{-1}(jl)(kl)\},
\{(kl),(jl)(kl)\}\}(pl)\\
&&\{\{(ql)^{-1},(jl)^{-1}(kl)^{-1}(jl)(kl)\},\{(kl),(jl)(kl)\}(pl)^{-1}\}(jk)(qk)(pk)(qk)^{-1}(pk)^{-1}(pq)^{-1},\\
&\equiv&(pk)(qk)(pk)^{-1}\{(jl)^{\varepsilon},\{(kl),(ql)^{-1}\}\{\{(ql),(kl)^{-1}(ql)^{-1}\},\{(kl),(ql)^{-1}\}\}(pl)\\
&&\{\{(ql)^{-1},(kl)^{-1}(ql)^{-1}\}\{(kl),(ql)^{-1}\}\}(pl)^{-1}\}(qk)^{-1}(jk)(qk)(pk)(qk)^{-1}(pk)^{-1}(pq)^{-1}\\
&\equiv&(pk)(qk)\{\{(jl)^{\varepsilon},(kl)(pl)(kl)^{-1}(pl)^{-1}\},\{(ql),(kl)(pl)(kl)^{-1}(pl)^{-1}\}\{(kl),(pl)^{-1}\}\\
&&\{(pl),(kl)^{-1}(pl)^{-1}\}\{(ql)^{-1},(kl)(pl)(kl)^{-1}(pl)^{-1}\}\{(pl)^{-1},(kl)^{-1}(pl)^{-1}\}\}\\
&&(pk)^{-1}(qk)^{-1}(jk)(qk)(pk)(qk)^{-1}(pk)^{-1}(pq)^{-1}\\
&\equiv&(pk)\{(jl)^{\varepsilon},\{(ql),(kl)\}\{(kl),(ql)(kl)\}(pl)\{(ql)^{-1},(kl)\}(pl)^{-1}\}\\
&&(qk)(pk)^{-1}(qk)^{-1}(jk)(qk)(pk)(qk)^{-1}(pk)^{-1}(pq)^{-1}\\
&\equiv&\{\{(jl)^{\varepsilon},(pl)^{-1}(kl)^{-1}(pl)(kl)\},\{(ql),(pl)^{-1}(kl)^{-1}(pl)(kl)\}\{(kl),(pl)(kl)\}\\
&&\{(pl),(kl)\}\{(kl)^{-1},(pl)(kl)\}\{(ql)^{-1},(pl)^{-1}(kl)^{-1}(pl)(kl)\}\{(kl),(pl)(kl)\}\\
&&\{(pl)^{-1},(kl)\}\}\{(jk),(qk)(pk)(qk)^{-1}(pk)^{-1}\}(pq)^{-1}\\
&\equiv&\{(jl)^{\varepsilon},(ql)(pl)(ql)^{-1}(pl)^{-1}(kl)\}\{(jk),(qk)(pk)(qk)^{-1}(pk)^{-1}\}(pq)^{-1}.
\end{eqnarray*}

$(11)\wedge(11)$\\
Let
$f=(pq)^{-1}(ik)^{-1}-\{(ik)^{-1},(qk)(pk)(qk)^{-1}(pk)^{-1}\}(pq)^{-1},p<i<q<k,
g=(ik)^{-1}(jl)^{\varepsilon}-\{(jl)^{\varepsilon},(kl)(il)(kl)^{-1}(il)^{-1}\}(ik)^{-1},
i<j<k<l.$ Then $w=(pq)^{-1}(ik)^{-1}(jl)^{\varepsilon}$ \ and
$$
(f,g)_w=(pq)^{-1}\{(jl)^{\varepsilon},(kl)(il)(kl)^{-1}(il)^{-1}\}(ik)^{-1}
-\{(ik)^{-1},(qk)(pk)(qk)^{-1}(pk)^{-1}\}(pq)^{-1}(jl)^{\varepsilon}.
$$
There are three cases to consider.
\begin{enumerate}
\item[1)]\ $q<j,p<i<q<j<k<l$.
\begin{eqnarray*}
&&(pq)^{-1}\{(jl)^{\varepsilon},(kl)(il)(kl)^{-1}(il)^{-1}\}(ik)^{-1}\\
&\equiv&\{(jl)^{\varepsilon},(kl)\{(kl)^{-1},\{(il)^{-1},(ql)(pl)(ql)^{-1}(pl)^{-1}\}\}\}\\
&&\{(ik)^{-1},(qk)(pk)(qk)^{-1}(pk)^{-1}\}(pq)^{-1}\ \ \ \ \ \ \ and \\
&&\\
&&\{(ik)^{-1},(qk)(pk)(qk)^{-1}(pk)^{-1}\}(pq)^{-1}(jl)^{\varepsilon}\\
&\equiv&(pk)(qk)(pk)^{-1}(qk)^{-1}(ik)^{-1}(qk)(pk)(qk)^{-1}\{(jl)^{\varepsilon},(kl)(pl)(kl)^{-1}(pl)^{-1}\}\\
&&(pk)^{-1}(pq)^{-1}\\
&\equiv&(pk)(qk)(pk)^{-1}(qk)^{-1}(ik)^{-1}(qk)(pk)\{\{(jl)^{\varepsilon},(kl)(ql)(kl)^{-1}(ql)^{-1}\},\{(kl),(ql)^{-1}\}\\
&&(pl)\{(kl)^{-1},(ql)^{-1}\}(pl)^{-1}\}(qk)^{-1}(pk)^{-1}(pq)^{-1}\\
&\equiv&(pk)(qk)(pk)^{-1}(qk)^{-1}(ik)^{-1}(qk)\{\{(jl)^{\varepsilon},(pl)^{-1}(kl)^{-1}(pl)(kl)\},\{(kl),(pl)(kl)\}\\
&&\{\{(kl)^{-1},(pl)(kl)\},\{(ql)^{-1},(pl)^{-1}(kl)^{-1}(pl)(kl)\}\{(pl)^{-1},(kl)\}\}\}\\
&&(pk)(qk)^{-1}(pk)^{-1}(pq)^{-1}\\
&\equiv&(pk)(qk)(pk)^{-1}(qk)^{-1}(ik)^{-1}\{\{(jl)^{\varepsilon},(ql)^{-1}(kl)^{-1}(ql)(kl)\},\{(kl),(ql)(kl)\}\\
&&\{\{(kl)^{-1},(ql)(kl)\},(pl)\{(ql)^{-1},(kl)\}(pl)^{-1}\}\}(qk)(pk)(qk)^{-1}(pk)^{-1}(pq)^{-1}\\
&\equiv&(pk)(qk)(pk)^{-1}(qk)^{-1}\{\{(jl)^{\varepsilon},(kl)(il)(kl)^{-1}(il)^{-1}\},\{(kl),(il)^{-1}\}\{\{(kl)^{-1},\\
&&(il)^{-1}\},\{(ql),(kl)(il)(kl)^{-1}(il)^{-1}\}\{(kl),(il)^{-1}\}(pl)\{(kl)^{-1},(il)^{-1}\}\{(ql)^{-1},(kl)(il)\\
&&(kl)^{-1}(il)^{-1}\}\{(kl),(il)^{-1}\}(pl)^{-1}\}\}(ik)^{-1}(qk)(pk)(qk)^{-1}(pk)^{-1}(pq)^{-1}\\
&\equiv&(pk)(qk)(pk)^{-1}\{\{(jl)^{\varepsilon},(kl)(ql)(kl)^{-1}(ql)^{-1}\},\{(kl),(ql)^{-1}\}\{\{(kl)^{-1},(ql)^{-1}\},\\
&&(il)^{-1}\{(kl)^{-1},(ql)^{-1}\}\{(ql),(kl)^{-1}(ql)^{-1}\}\{(kl),(ql)^{-1}\}(pl)\{(kl)^{-1},(ql)^{-1}\}\{(ql)^{-1},\\
&&(kl)^{-1}(ql)^{-1}\}\{(kl),(ql)^{-1}\}(pl)^{-1}\}(qk)^{-1}(ik)^{-1}(qk)(pk)(qk)^{-1}(pk)^{-1}(pq)^{-1}\\
&\equiv&(pk)(qk)\{\{(jl)^{\varepsilon},(kl)(pl)(kl)^{-1}(pl)^{-1}\},\{(kl),(pl)^{-1}\}\{\{(kl)^{-1},(pl)^{-1}\},\\
&&\{(ql)^{-1},(kl)(pl)(kl)^{-1}(pl)^{-1}\}\{(il)^{-1},(kl)(pl)(kl)^{-1}(pl)^{-1}\}\{(ql),(kl)(pl)\\
&&(kl)^{-1}(pl)^{-1}\}\{(pl),(kl)^{-1}(pl)^{-1}\}\{(ql)^{-1},(kl)(pl)(kl)^{-1}(pl)^{-1}\}\\
&&\{(pl)^{-1},(kl)^{-1}(pl)^{-1}\}\}\}(pk)^{-1}(qk)^{-1}(ik)^{-1}(qk)(pk)(qk)^{-1}(pk)^{-1}(pq)^{-1}\\
&\equiv&(pk)\{\{(jl)^{\varepsilon},(ql)^{-1}(kl)^{-1}(ql)(kl)\},\{(kl),(ql)(kl)\}\{\{(kl)^{-1},(ql)(kl)\},\{(il)^{-1},\\
&&\{(ql),(kl)\}\{(kl),(ql)(kl)\}(pl)\}\{(kl)^{-1},(ql)(kl)\}\{(ql)^{-1},(kl)\}(pl)^{-1}\}\}\\
&&(qk)(pk)^{-1}(qk)^{-1}(ik)^{-1}(qk)(pk)(qk)^{-1}(pk)^{-1}(pq)^{-1}\\
&\equiv&\{\{(jl)^{\varepsilon},(pl)^{-1}(kl)^{-1}(pl)(kl)\},\{(kl),(pl)(kl)\}\{\{(kl)^{-1},(pl)(kl)\},\{\{(il)^{-1},\\
&&(pl)^{-1}(kl)^{-1}(pl)(kl)\},\{(ql),(pl)^{-1}(kl)^{-1}(pl)(kl)\}\{(kl),(pl)(kl)\}\{(pl),(kl)\}\\
&&\{(kl)^{-1},(pl)(kl)\}\{(ql)^{-1},(pl)^{-1}(kl)^{-1}(pl)(kl)\}\}\{(pl)^{-1},(kl)\}\}\}\\
&&\{(ik)^{-1},(qk)(pk)(qk)^{-1}(pk)^{-1}(pq)^{-1}\}\\
&\equiv&\{(jl)^{\varepsilon},(kl)\{(kl)^{-1},\{(il)^{-1},(ql)(pl)(ql)^{-1}(pl)^{-1}\}\}\}\\
&&\{(ik)^{-1},(qk)(pk)(qk)^{-1}(pk)^{-1}\}(pq)^{-1}.
\end{eqnarray*}
\item[2)]\ $q=j,p<i<q=<j<k<l$.
\begin{eqnarray*}
&&(pq)^{-1}\{(jl)^{\varepsilon},(kl)(il)(kl)^{-1}(il)^{-1}\}(ik)^{-1}\\
&=&(pj)^{-1}\{(jl)^{\varepsilon},(kl)(il)(kl)^{-1}(il)^{-1}\}(ik)^{-1}\\
&\equiv&\{\{(jl)^{\varepsilon},(pl)^{-1}\},(kl)\{(kl)^{-1},\{(il)^{-1},(jl)(pl)(jl)^{-1}(pl)^{-1}\}\}\}\\
&&\{(ik)^{-1},(jk)(pk)(jk)^{-1}(pk)^{-1}\}(pj)^{-1}\ \ \ \ \ \ \ \ \ and \\
&&\\
&&\{(ik)^{-1},(qk)(pk)(qk)^{-1}(pk)^{-1}\}(pq)^{-1}(jl)^{\varepsilon}\\
&=&\{(ik)^{-1},(jk)(pk)(jk)^{-1}(pk)^{-1}\}(pj)^{-1}(jl)^{\varepsilon}\\
&\equiv&(pk)(jk)(pk)^{-1}(jk)^{-1}(ik)^{-1}(jk)(pk)(jk)^{-1}(pk)^{-1}\{(jl)^{\varepsilon},(pl)^{-1}\}(pj)^{-1}\\
&\equiv&(pk)(jk)(pk)^{-1}(jk)^{-1}(ik)^{-1}(jk)(pk)(jk)^{-1}\{\{(jl)^{\varepsilon},(kl)(pl)(kl)^{-1}(pl)^{-1}\},\\
&&\{(pl)^{-1},(kl)^{-1}(pl)^{-1}\}\}(pk)^{-1}(pj)^{-1}\\
&\equiv&(pk)(jk)(pk)^{-1}(jk)^{-1}(ik)^{-1}(jk)(pk)\{\{(jl)^{\varepsilon},(kl)^{-1}(jl)^{-1}\},(pl)^{-1}\}\\
&&(jk)^{-1}(pk)^{-1}(pj)^{-1}\\
&\equiv&(pk)(jk)(pk)^{-1}(jk)^{-1}(ik)^{-1}(jk)\{\{(jl)^{\varepsilon},(pl)^{-1}(kl)^{-1}(pl)(kl)\},\{(kl)^{-1},(pl)\\
&&(kl)\}\{(jl)^{-1},(pl)^{-1}(kl)^{-1}(pl)(kl)\}\{(pl)^{-1},(kl)\}\}(pk)(jk)^{-1}(pk)^{-1}(pj)^{-1}\\
&\equiv&(pk)(jk)(pk)^{-1}(jk)^{-1}(ik)^{-1}\{\{(jl)^{\varepsilon},(kl)\},(pl)^{-1}\{(kl)^{-1},(jl)(kl)\}(pl)\\
&&\{(jl)^{-1},(kl)\}(pl)^{-1}\}(jk)(pk)(jk)^{-1}(pk)^{-1}(pj)^{-1}\\
&\equiv&(pk)(jk)(pk)^{-1}(jk)^{-1}\{\{(jl)^{\varepsilon},(kl)(il)(kl)^{-1}(il)^{-1}\},\{(kl),(il)^{-1}\}(pl)^{-1}\\
&&\{\{(kl)^{-1},(il)^{-1}\},\{(jl),(kl)(il)(kl)^{-1}(il)^{-1}\}\{(kl),(il)^{-1}\}\}(pl)\{\{(jl)^{-1},(kl)\\
&&(il)(kl)^{-1}(il)^{-1}\},\{(kl),(il)^{-1}\}\}(pl)^{-1}\}(ik)^{-1}(jk)(pk)(jk)^{-1}(pk)^{-1}(pj)^{-1}\\
&\equiv&(pk)(jk)(pk)^{-1}\{\{(jl)^{\varepsilon},(kl)^{-1}(jl)^{-1}\},\{(kl),(jl)^{-1}\}\{\{(kl)^{-1},(jl)^{-1}\},(il)^{-1}\\
&&\{(kl)^{-1},(jl)^{-1}\}\{(jl),(kl)^{-1}(jl)^{-1}\}\{(kl),(jl)^{-1}\}(pl)\}\{(kl)^{-1},(jl)^{-1}\}\{(jl)^{-1},\\
&&(kl)^{-1}(jl)^{-1}\}\{(kl),(jl)^{-1}\}(pl)^{-1}\}(jk)^{-1}(ik)^{-1}(jk)(pk)(jk)^{-1}(pk)^{-1}(pj)^{-1}\\
&\equiv&(pk)(jk)\{\{(jl)^{\varepsilon},(kl)(pl)(kl)^{-1}(pl)^{-1}\},\{\{(kl)^{-1},(pl)^{-1}\},\{(jl)^{-1},(kl)(pl)(kl)^{-1}\\
&&(pl)^{-1}\}\{(il)^{-1},(kl)(pl)(kl)^{-1}(pl)^{-1}\}\{(jl),(kl)(pl)(kl)^{-1}(pl)^{-1}\}\\
&&\{(pl),(kl)^{-1}(pl)^{-1}\}\}\{(jl)^{-1},(kl)(pl)(kl)^{-1}(pl)^{-1}\}\{(pl)^{-1},(kl)^{-1}(pl)^{-1}\}\}\\
&&(pk)^{-1}(jk)^{-1}(ik)^{-1}(jk)(pk)(jk)^{-1}(pk)^{-1}(pj)^{-1}\\
&\equiv&(pk)\{\{(jl)^{\varepsilon},(kl)\},\{\{(kl)^{-1},(jl)(kl)\},\{(il)^{-1},\{(jl),(kl)\}\{(kl),(jl)(kl)\}(pl)\}\\
&&\{(kl)^{-1},(jl)(kl)\}\}\{(jl)^{-1},(kl)\}(pl)^{-1}\}(jk)(pk)^{-1}(jk)^{-1}(ik)^{-1}\\
&&(jk)(pk)(jk)^{-1}(pk)^{-1}(pj)^{-1}\\
&\equiv&\{\{(jl)^{\varepsilon},(pl)^{-1}(kl)^{-1}(pl)(kl)\},\{(kl),(pl)(kl)\}\{\{(kl)^{-1},(pl)(kl)\},\{(jl),(pl)^{-1}(kl)^{-1}\\
&&(pl)(kl)\}\{(kl),(pl)(kl)\}\{\{(il)^{-1},(pl)(kl)(pl)^{-1}(kl)^{-1}\},\{(jl),(pl)^{-1}(kl)^{-1}(pl)(kl)\}\\
&&\{(kl),(pl)(kl)\}\{(pl),(kl)\}\}\{\{(kl)^{-1},(pl)(kl)\},\{(jl),(pl)^{-1}(kl)^{-1}(pl)(kl)\}\\
&&\{(kl),(pl)(kl)\}\}\}\{\{(jl),(pl)^{-1}(kl)^{-1}(pl)(kl)\},\{(kl),(pl)(kl)\}\}\\
&&\{(pl)^{-1},(kl)\}\}\{(ik)^{-1},(jk)(pk)(jk)^{-1}(pk)^{-1}\}(pj)^{-1}\\
&\equiv&\{\{(jl)^{\varepsilon},(pl)^{-1}\},(kl)\{(kl)^{-1},\{(il)^{-1},(jl)(pl)(jl)^{-1}(pl)^{-1}\}\}\}\\
&&\{(ik)^{-1},(jk)(pk)(jk)^{-1}(pk)^{-1}\}(pj)^{-1}.
\end{eqnarray*}
\item[3)]\ $q>j,p<i<j<q<k<l$.
\begin{eqnarray*}
&&(pq)^{-1}\{(jl)^{\varepsilon},(kl)(il)(kl)^{-1}(il)^{-1}\}(ik)^{-1}\\
&\equiv&\{\{(jl)^{\varepsilon},(ql)(pl)(pl)^{-1}(ql)^{-1}\},(kl)\{(il),(ql)(pl)(ql)^{-1}(pl)^{-1}\}(kl)^{-1}\\
&&\{(il)^{-1},(ql)(pl)(ql)^{-1}(pl)^{-1}\}\}\{(ik)^{-1},(qk)(pk)(qk)^{-1}(pk)^{-1}\}(pq)^{-1}\ \ \ \ \ \ and \\
&&\\
&&\{(ik)^{-1},(qk)(pk)(qk)^{-1}(pk)^{-1}\}(pq)^{-1}(jl)^{\varepsilon}\\
&\equiv&(pk)(qk)(pk)^{-1}(qk)^{-1}(ik)^{-1}(qk)(pk)(qk)^{-1}(pk)^{-1}\{(jl)^{\varepsilon},(ql)(pl)(ql)^{-1}(pl)^{-1}\}(pq)^{-1}\\
&\equiv&(pk)(qk)(pk)^{-1}(qk)^{-1}(ik)^{-1}(qk)(pk)(qk)^{-1}\{\{(jl)^{\varepsilon},(kl)(pl)(kl)^{-1}(pl)^{-1}\},\\
&&\{(ql),(kl)(pl)(kl)^{-1}(pl)^{-1}\}\{(pl),(kl)^{-1}(pl)^{-1}\}\{(ql)^{-1},(kl)(pl)(kl)^{-1}(pl)^{-1}\}\\
&&\{(pl)^{-1},(kl)^{-1}(pl)^{-1}\}\}(pk)^{-1}(pq)^{-1}\\
&\equiv&(pk)(qk)(pk)^{-1}(qk)^{-1}(ik)^{-1}(qk)(pk)\{(jl)^{\varepsilon},\{(ql),(kl)^{-1}(ql)^{-1}\}\{(kl),(ql)^{-1}\}\\
&&(pl)\{(kl)^{-1},(ql)^{-1}\}\{(ql)^{-1},(kl)^{-1}(ql)^{-1}\}(pl)^{-1}\}(qk)^{-1}(pk)^{-1}(pq)^{-1}\\
&\equiv&(pk)(qk)(pk)^{-1}(qk)^{-1}(ik)^{-1}(qk)\{\{(jl)^{\varepsilon},(pl)^{-1}(kl)^{-1}(pl)(kl)\},\\
&&\{(ql),(pl)^{-1}(kl)^{-1}(pl)(kl)\}\{(kl),(pl)(kl)\}\{(pl),(kl)\}\{(kl)^{-1},(pl)(kl)\}\\
&&\{(ql)^{-1},(pl)^{-1}(kl)^{-1}(pl)(kl)\}\{(pl)^{-1},(kl)\}\}(pk)(qk)^{-1}(pk)^{-1}(pq)^{-1}\\
&\equiv&(pk)(qk)(pk)^{-1}(qk)^{-1}(ik)^{-1}\{(jl)^{\varepsilon},\{(ql),(kl)\}(pl)\{(ql)^{-1},(kl)\}(pl)^{-1}\}\\
&&(qk)(pk)(qk)^{-1}(pk)^{-1}(pq)^{-1}\\
&\equiv&(pk)(qk)(pk)^{-1}(qk)^{-1}\{\{(jl)^{\varepsilon},(kl)(il)(kl)^{-1}(il)^{-1}\},\{\{(ql),(kl)(il)(kl)^{-1}(il)^{-1}\},\\
&&\{(kl),(il)^{-1}\}\}(pl)\{\{(ql)^{-1},(kl)(il)(kl)^{-1}(il)^{-1}\},\{(kl),(il)^{-1}\}\}(pl)^{-1}\}\\
&&(ik)^{-1}(qk)(pk)(qk)^{-1}(pk)^{-1}(pq)^{-1}\\
&\equiv&(pk)(qk)(pk)^{-1}\{(jl)^{\varepsilon},\{(kl),(ql)^{-1}\}(il)\{(kl)^{-1},(ql)^{-1}\}(il)^{-1}\{(kl)^{-1},(ql)^{-1}\}\\
&&\{(ql),(kl)^{-1}(ql)^{-1}\}\{(kl),(ql)^{-1}\}(pl)\{(kl)^{-1},(ql)^{-1}\}\{(ql)^{-1},(kl)^{-1}(ql)^{-1}\}\\
&&\{(kl),(ql)^{-1}\}(pl)^{-1}\}(qk)^{-1}(ik)(qk)(pk)(qk)^{-1}(pk)^{-1}(pq)^{-1}\\
&\equiv&(pk)(qk)\{\{(jl)^{\varepsilon},(kl)(pl)(kl)^{-1}(pl)^{-1}\},\{\{(il),(kl)(pl)(kl)^{-1}(pl)^{-1}\},\\
&&\{(ql),(kl)(pl)(kl)^{-1}(pl)^{-1}\}\{(kl)^{-1},(pl)^{-1}\}\{(ql)^{-1},(kl)(pl)(kl)^{-1}(pl)^{-1}\}\}\\
&&\{(il)^{-1},(kl)(pl)(kl)^{-1}(pl)^{-1}\}\{(ql),(kl)(pl)(kl)^{-1}(pl)^{-1}\}\{(pl),(kl)^{-1}(pl)^{-1}\}\\
&&\{(ql)^{-1},(kl)(pl)(kl)^{-1}(pl)^{-1}\}\{(pl)^{-1},(kl)^{-1}(pl)^{-1}\}\}\\
&&(pk)^{-1}(qk)^{-1}(ik)(qk)(pk)(qk)^{-1}(pk)^{-1}(pq)^{-1}\\
&\equiv&(pk)\{(jl)^{\varepsilon},\{(il),\{\{(kl)^{-1},(ql)(kl)\},(pl)^{-1}\{(kl)^{-1},(ql)(kl)\}\{(ql)^{-1},(kl)\}\}\}\\
&&(il)^{-1}\{(ql),(kl)\}\{(pl),\{(kl)^{-1},(ql)(kl)\}\}\{(ql)^{-1},(kl)\}(pl)^{-1}\}\\
&&(qk)(pk)(qk)^{-1}(ik)(qk)(pk)(qk)^{-1}(pk)^{-1}(pq)^{-1}\\
&\equiv&\{\{(jl)^{\varepsilon},(pl)^{-1}(kl)^{-1}(pl)(kl)\},\{\{(il),(pl)^{-1}(kl)^{-1}(pl)(kl)\},\{\{(kl)^{-1},(pl)(kl)\},\\
&&\{(ql),(pl)^{-1}(kl)^{-1}(pl)(kl)\}\{(kl),(pl)(kl)\}\{(pl)^{-1},(kl)\}\{(kl)^{-1},(pl)(kl)\}\\
&&\{(ql)^{-1},(pl)^{-1}(kl)^{-1}(pl)(kl)\}\}\}\{(il)^{-1},(pl)^{-1}(kl)^{-1}(pl)(kl)\}\{(ql),(pl)^{-1}\\
&&(kl)^{-1}(pl)(kl)\}\{(kl),(pl)(kl)\}\{(pl),(kl)\}\{(kl)^{-1},(pl)(kl)\}\{(ql)^{-1},(pl)^{-1}\\
&&(kl)^{-1}(pl)(kl)\}\{(pl)^{-1},(kl)\}\}\{(ik)^{-1},(qk)(pk)(qk)^{-1}(pk)^{-1}\}(pq)^{-1}\\
&\equiv&\{\{(jl)^{\varepsilon},(ql)(pl)(ql)^{-1}(pl)^{-1}\},(kl)\{(il),(ql)(pl)(ql)^{-1}(pl)^{-1}\}(kl)^{-1}\\
&&\{(il)^{-1},(ql)(pl)(ql)^{-1}(pl)^{-1}\}\}\{(ik)^{-1},(qk)(pk)(qk)^{-1}(pk)^{-1}\}(pq)^{-1}.
\end{eqnarray*}
\end{enumerate}

$(11)\wedge(12)$\\
Let$f=(pq)^{-1}(ik)-\{(ik),(qk)(pk)(qk)^{-1}(pk)^{-1}\}(pq)^{-1},p<i<q<k,
g=(ik)(jl)^{\varepsilon}-\{(jl)^{\varepsilon},(il)^{-1}(kl)^{-1}(il)(kl)\}(ik),
i<j<k<l.$ Then $w=(pq)^{-1}(ik)(jl)^{\varepsilon}$ and
$$
(f,g)_w=(pq)^{-1}\{(jl)^{\varepsilon},(il)^{-1}(kl)^{-1}(il)(kl)\}(ik)
-\{(ik),(qk)(pk)(qk)^{-1}(pk)^{-1}\}(pq)^{-1}(jl)^{\varepsilon}.
$$
There are three cases to consider.
\begin{enumerate}
\item[1)]\ $q<j,p<i<q<j<k<l.$
\begin{eqnarray*}
&&(pq)^{-1}\{(jl)^{\varepsilon},(il)^{-1}(kl)^{-1}(il)(kl)\}(ik)\\
&\equiv&\{(jl)^{\varepsilon},\{(kl)^{-1},\{(il),(ql)(pl)(ql)^{-1}(pl)^{-1}\}\}(kl)\}\\
&&\{(ik),(qk)(pk)(qk)^{-1}(pk)^{-1}\}(pq)^{-1}\ \ \ \ \ \ \ and\\
&&\\
&&\{(ik),(qk)(pk)(qk)^{-1}(pk)^{-1}\}(pq)^{-1}(jl)^{\varepsilon}\\
&\equiv&(pk)(qk)(pk)^{-1}(qk)^{-1}(ik)\{(jl)^{\varepsilon},(kl)\{(kl)^{-1},(ql)(kl)(pl)(kl)^{-1}(ql)^{-1}(kl)(pl)^{-1}\}\}\\
&&(qk)(pk)(qk)^{-1}(pk)^{-1}(pq)^{-1}\\
&\equiv&(pk)(qk)(pk)^{-1}(qk)^{-1}\{\{(jl)^{\varepsilon},(il)^{-1}(kl)^{-1}(il)(kl)\},\{(kl),(il)(kl)\}\\
&&\{\{(kl)^{-1},(il)(kl)\},\{(ql),(il)^{-1}(kl)^{-1}(il)(kl)\}\{(kl),(il)(kl)\}(pl)\\
&&\{(kl)^{-1},(il)(kl)\}\{(ql)^{-1},(il)^{-1}(kl)^{-1}(il)(kl)\}\{(kl),(il)(kl)\}(pl)^{-1}\}\}\\
&&(ik)(qk)(pk)(qk)^{-1}(pk)^{-1}(pq)^{-1}\\
&\equiv&(pk)(qk)(pk)^{-1}\{\{(jl)^{\varepsilon},(kl)(ql)(kl)^{-1}(ql)^{-1}\},\{(kl),(ql)^{-1}\}\{\{(kl)^{-1},(ql)^{-1}\},\\
&&(il)\{(ql),(kl)^{-1}(ql)^{-1}\}\{(kl),(ql)^{-1}\}(pl)\{(kl)^{-1},(ql)^{-1}\}\{(ql)^{-1},(kl)^{-1}(ql)^{-1}\}\\
&&\{(kl),(ql)^{-1}\}(pl)^{-1}\}\}(qk)^{-1}(ik)(qk)(pk)(qk)^{-1}(pk)^{-1}(pq)^{-1}\\
&\equiv&(pk)(qk)\{\{(jl)^{\varepsilon},(kl)(pl)(kl)^{-1}(pl)^{-1}\},\{(kl),(pl)^{-1}\}\{\{(kl)^{-1},(pl)^{-1}\},\\
&&\{(ql)^{-1},(kl)(pl)(kl)^{-1}(pl)^{-1}\}\{(il),(kl)(pl)(kl)^{-1}(pl)^{-1}\}\{(ql),(kl)(pl)\\
&&(kl)^{-1}(pl)^{-1}\}\{(kl),(pl)^{-1}\}\{(pl),(kl)^{-1}(pl)^{-1}\}\{(ql)^{-1},(kl)(pl)(kl)^{-1}(pl)^{-1}\}\\
&&\{(pl)^{-1},(kl)^{-1}(pl)^{-1}\}\}\}(pk)^{-1}(qk)^{-1}(ik)(qk)(pk)(qk)^{-1}(pk)^{-1}(pq)^{-1}\\
&\equiv&(pk)\{\{(jl)^{\varepsilon},(ql)^{-1}(kl)^{-1}(ql)(kl)\},\{(kl),(ql)(kl)\}\{\{(kl)^{-1},(ql)(kl)\},\{(il),\\
&&\{(ql),(kl)\}\{(kl),(ql)(kl)\}(pl)\}\{(ql)^{-1},(kl)\}(pl)^{-1}\}\}\\
&&(qk)(pk)^{-1}(qk)^{-1}(ik)(qk)(pk)(qk)^{-1}(pk)^{-1}(pq)^{-1}\\
&\equiv&\{\{(jl)^{\varepsilon},(pl)^{-1}(kl)^{-1}(pl)(kl)\},\{(kl),(pl)(kl)\}\{\{(kl)^{-1},(pl)(kl)\},\{\{(il),(pl)^{-1}\\
&&(kl)^{-1}(pl)(kl)\},\{(ql),(pl)^{-1}(kl)^{-1}(pl)(kl)\}\{(kl),(pl)(kl)\}\{(pl),(kl)\}\\
&&\{(kl)^{-1},(pl)(kl)\}\{(ql)^{-1},(pl)^{-1}(kl)^{-1}(pl)(kl)\}\}\{(kl),(pl)(kl)\}\\
&&\{(pl)^{-1},(kl)\}\}\}\{(ik),(qk)(pk)(qk)^{-1}(pk)^{-1}\}(pq)^{-1}\\
&\equiv&\{(jl)^{\varepsilon},\{(kl)^{-1},\{(il),(ql)(pl)(ql)^{-1}(pl)^{-1}\}\}(kl)\}\\
&&\{(ik),(qk)(pk)(qk)^{-1}(pk)^{-1}\}(pq)^{-1}.
\end{eqnarray*}
\item[2)]\ $q=j,p<i<q=j<k<l$.
\begin{eqnarray*}
&&(pq)^{-1}\{(jl)^{\varepsilon},(il)^{-1}(kl)^{-1}(il)(kl)\}(ik)\\
&=&(pj)^{-1}\{(jl)^{\varepsilon},(il)^{-1}(kl)^{-1}(il)(kl)\}(ik)\\
&\equiv&\{\{(jl)^{\varepsilon},(pl)^{-1}\},\{(il)^{-1},(jl)(pl)(jl)^{-1}(pl)^{-1}\}(kl)^{-1}\{(il),(jl)(pl)\\
&&(jl)^{-1}(pl)^{-1}\}(kl)\}\{(ik),(jk)(pk)(jk)^{-1}(pk)^{-1}\}(pj)^{-1}\ \ \ \ \ \ and\\
&&\\
&&\{(ik),(qk)(pk)(qk)^{-1}(pk)^{-1}\}(pq)^{-1}(jl)^{\varepsilon}\\
&=&(pk)(jk)(pk)^{-1}(jk)^{-1}(ik)(jk)(pk)(jk)^{-1}(pk)^{-1}(pj)^{-1}(jl)^{\varepsilon}\\
&\equiv&(pk)(jk)(pk)^{-1}(jk)^{-1}(ik)\{(jl)^{\varepsilon},(kl)(pl)^{-1}\{(kl)^{-1},(jl)(kl)\}(pl)\{(jl)^{-1},(kl)\}\\
&&(pl)^{-1}\}(jk)(pk)(jk)^{-1}(pk)^{-1}(pj)^{-1}\\
&\equiv&(pk)(jk)(pk)^{-1}(jk)^{-1}\{\{(jl)^{\varepsilon},(il)^{-1}(kl)^{-1}(il)(kl)\},\{(kl),(il)(kl)\}(pl)^{-1}\\
&&\{\{(kl)^{-1},(il)(kl)\},\{(jl),(il)^{-1}(kl)^{-1}(il)(kl)\}\{(kl),(il)(kl)\}\}(pl)\{\{(jl)^{-1},(il)^{-1}\\
&&(kl)^{-1}(il)(kl)\},\{(kl),(il)(kl)\}\}(pl)^{-1}\}(ik)(jk)(pk)(jk)^{-1}(pk)^{-1}(pj)^{-1}\\
&\equiv&(pk)(jk)(pk)^{-1}\{\{(jl)^{\varepsilon},(kl)^{-1}(jl)^{-1}\},\{(kl),(jl)^{-1}\}\{\{(kl)^{-1},(jl)^{-1}\},(il)\\
&&\{(jl),(kl)^{-1}(jl)^{-1}\}\{(kl),(jl)^{-1}\}(pl)\}\{\{(jl)^{-1},(kl)^{-1}(jl)^{-1}\}\{(kl),(jl)^{-1}\}\}\\
&&(pl)^{-1}\}(jk)^{-1}(ik)(jk)(pk)(jk)^{-1}(pk)^{-1}(pj)^{-1}\\
&\equiv&(pk)(jk)\{\{(jl)^{\varepsilon},(kl)(pl)(kl)^{-1}(pl)^{-1}\},\{\{(kl)^{-1},(pl)^{-1}\},\{(jl)^{-1},(kl)(pl)\\
&&(kl)^{-1}(pl)^{-1}\}\{(il),(kl)(pl)(kl)^{-1}(pl)^{-1}\}\{(jl),(kl)(pl)(kl)^{-1}(pl)^{-1}\}\\
&&\{(kl),(pl)^{-1}\}\{(pl),(kl)^{-1}(pl)^{-1}\}\}\{(jl)^{-1},(kl)(pl)(kl)^{-1}(pl)^{-1}\}\{(pl)^{-1},\\
&&(kl)^{-1}(pl)^{-1}\}\}(pk)^{-1}(jk)^{-1}(ik)(jk)(pk)(jk)^{-1}(pk)^{-1}(pj)^{-1}\\
&\equiv&(pk)\{\{(jl)^{\varepsilon},(kl)\},\{\{(kl)^{-1},(jl)(kl)\},\{(il),\{(jl),(kl)\}\{(kl),(jl)(kl)\}(pl)\}\}\\
&&\{(jl)^{-1},(kl)\}(pl)^{-1}\}(jk)(pk)^{-1}(jk)^{-1}(ik)(jk)(pk)(jk)^{-1}(pk)^{-1}(pj)^{-1}\\
&\equiv&\{\{(jl)^{\varepsilon},(pl)^{-1}(kl)^{-1}(pl)(kl)\},\{\{(kl)^{-1},(pl)(kl)\},\{(jl),(pl)^{-1}(kl)^{-1}(pl)(kl)\}\\
&&\{\{(il),(pl)^{-1}(kl)^{-1}(pl)(kl)\},\{(jl),(pl)^{-1}(kl)^{-1}(pl)(kl)\}\{(kl),(pl)(kl)\}\\
&&\{(pl),(kl)\}\{(kl)^{-1},(pl)(kl)\}\}\}\{(jl)^{-1},(pl)^{-1}(kl)^{-1}(pl)(kl)\}\{(kl),(pl)(kl)\}\\
&&\{(pl)^{-1},(kl)\}\}\{(ik),(jk)(pk)(jk)^{-1}(pk)^{-1}\}(pj)^{-1}\\
&\equiv&\{\{(jl)^{\varepsilon},(pl)^{-1}\},\{(il)^{-1},(jl)(pl)(jl)^{-1}(pl)^{-1}\}(kl)^{-1}\\
&&\{(il),(jl)(pl)(jl)^{-1}(pl)^{-1}\}(kl)\}\{(ik),(jk)(pk)(jk)^{-1}(pk)^{-1}\}(pj)^{-1}.
\end{eqnarray*}
\item[3)]\ $q>j,p<i<j<q<k<l$.
\begin{eqnarray*}
&&(pq)^{-1}\{(jl)^{\varepsilon},(il)^{-1}(kl)^{-1}(il)(kl)\}(ik)\\
&\equiv&\{(jl)^{\varepsilon},(il)^{-1}(ql)(pl)(ql)^{-1}(pl)^{-1}\{\{(il),(ql)(pl)(ql)^{-1}(pl)^{-1}\},(kl)\}\}\\
&&\{(ik),(qk)(pk)(qk)^{-1}(pk)^{-1}\}(pq)^{-1}\ \ \ \ \ \ and \\
&&\\
&&\{(ik),(qk)(pk)(qk)^{-1}(pk)^{-1}\}(pq)^{-1}(jl)^{\varepsilon}\\
&\equiv&(pk)(qk)(pk)^{-1}(qk)^{-1}(ik)\{(jl)^{\varepsilon}\{(ql),(kl)\}(pl)\\
&&\{(ql)^{-1},(kl)\}(pl)^{-1}\}(qk)(pk)(qk)^{-1}(pk)^{-1}(pq)^{-1}\\
&\equiv&(pk)(qk)(pk)^{-1}(qk)^{-1}\{\{(jl)^{\varepsilon},(il)^{-1}(kl)^{-1}(il)(kl)\},\{\{(ql),(il)^{-1}(kl)^{-1}\\
&&(il)(kl)\},\{(kl),(il)(kl)\}\}(pl)\{\{(ql)^{-1},(il)^{-1}(kl)^{-1}(il)(kl)\},\{(kl),(il)(kl)\}\}\\
&&(pl)^{-1}\}(ik)(qk)(pk)(qk)^{-1}(pk)^{-1}(pq)^{-1}\\
&\equiv&(pk)(qk)(pk)^{-1}\{(jl)^{\varepsilon},(il)^{-1}\{(kl)^{-1},(ql)^{-1}\}(il)\{(ql),(kl)^{-1}(ql)^{-1}\}\\
&&\{(kl),(ql)^{-1}\}\{\{(ql)^{-1},(kl)^{-1}(ql)^{-1}\},\{(kl),(ql)^{-1}\}(pl)^{-1}\}\}\\
&&(qk)^{-1}(ik)(qk)(pk)(qk)^{-1}(pk)^{-1}(pq)^{-1}\\
&\equiv&(pk)(qk)\{\{(jl)^{\varepsilon},(kl)(pl)(kl)^{-1}(pl)^{-1}\},\{(il)^{-1},(kl)(pl)(kl)^{-1}(pl)^{-1}\}\\
&&\{\{(kl)^{-1},(pl)^{-1}\},\{(ql)^{-1},(kl)(pl)(kl)^{-1}(pl)^{-1}\}\}\{(il),(kl)(pl)(kl)^{-1}(pl)^{-1}\}\\
&&\{(ql),(kl)(pl)(kl)^{-1}(pl)^{-1}\}\{(kl),(pl)^{-1}\}\{\{(ql)^{-1},(kl)(pl)(kl)^{-1}(pl)^{-1}\},\\
&&\{(pl)^{-1},(kl)^{-1}(pl)^{-1}\}\}\}(pk)^{-1}(qk)^{-1}(ik)(qk)(pk)(qk)^{-1}(pk)^{-1}(pq)^{-1}\\
&\equiv&(pk)\{(jl)^{\varepsilon},\{\{(kl)^{-1},(ql)(kl)\},(pl)^{-1}\{(kl)^{-1},(ql)(kl)\}\{(ql)^{-1},(kl)\}(il)\}\{(ql),\\
&&(kl)\}\{(kl),(ql)(kl)\}\{\{(ql)^{-1},(kl)\},(pl)^{-1}\}\}(qk)(pk)^{-1}(qk)^{-1}(ik)\\
&&(qk)(pk)(qk)^{-1}(pk)^{-1}(pq)^{-1}\\
&\equiv&\{\{(jl)^{\varepsilon},(pl)^{-1}(kl)^{-1}(pl)(kl)\},\{\{(kl)^{-1},(pl)(kl)\},\{\{(pl)^{-1},(kl)\},\\
&&\{(kl)^{-1},(pl)(kl)\}\{(ql)^{-1},(pl)^{-1}(kl)^{-1}(pl)(kl)\}\}\{(il),(pl)^{-1}(kl)^{-1}(pl)(kl)\}\}\\
&&\{(ql),(pl)^{-1}(kl)^{-1}(pl)(kl)\}\{(kl),(pl)(kl)\}\{\{(ql)^{-1},(pl)^{-1}(kl)^{-1}(pl)(kl)\},\\
&&\{(kl),(pl)(kl)\}\{(pl)^{-1},(kl)\}\}\}\{(ik),(qk)(pk)(qk)^{-1}(pk)^{-1}\}(pq)^{-1}\\
&\equiv&\{(jl)^{\varepsilon},(il)^{-1}(ql)(pl)(ql)^{-1}(pl)^{-1}(kl)^{-1}\{(il),(ql)(pl)(ql)^{-1}(pl)^{-1}\}(kl)\}\\
&&\{(ik),(qk)(pk)(qk)^{-1}(pk)^{-1}\}(pq)^{-1}.
\end{eqnarray*}
\end{enumerate}

$(11)\wedge(13)$\\
Let
$f=(pq)^{-1}(ik)^\delta-\{(ik)^\delta,(qk)(pk)(qk)^{-1}(pk)^{-1}\}(pq)^{-1},\
g=(ik)^\delta(jl)^{\varepsilon}-(jl)^{\varepsilon}(ik)^\delta$. Then
$w=(pq)^{-1}(ik)^\delta(jl)^{\varepsilon}$ and
$$
(f,g)_w=(pq)^{-1}(jl)^{\varepsilon}(ik)^\delta-\{(ik)^\delta,(qk)(pk)(qk)^{-1}(pk)^{-1}\}(pq)^{-1}(jl)^{\varepsilon}.
$$
There are two cases to consider.
\begin{enumerate}
\item[1)] \  $p<i<q<k,\ j<i<k<l $. In this case, there are three subcases to consider.
\begin{enumerate}
\item[a)] \  $p<j,p<j<i<q<k<l$.
\begin{eqnarray*}
&&(pq)^{-1}(jl)^{\varepsilon}(ik)^\delta\\
&\equiv&\{(jl)^{\varepsilon},(ql)(pl)(ql)^{-1}(pl)^{-1}\}\{(ik)^\delta,(qk)(pk)(qk)^{-1}(pk)^{-1}\}(pq)^{-1} \ \ \ and\\
&&\\
&&\{(ik)^\delta,(qk)(pk)(qk)^{-1}(pk)^{-1}\}(pq)^{-1}(jl)^{\varepsilon}\\
&\equiv&(pk)(qk)(pk)^{-1}(qk)^{-1}(ik)^\delta(qk)(pk)(qk)^{-1}(pk)^{-1}\{(jl)^{\varepsilon},(ql)(pl)\\
&&(ql)^{-1}(pl)^{-1}\}(pq)^{-1}\\
&\equiv&(pk)(qk)(pk)^{-1}(qk)^{-1}(ik)^\delta\{(jl)^{\varepsilon},\{(ql),(kl)\}(pl)\{(ql)^{-1},(kl)\}(pl)^{-1}\}\\
&&(qk)(pk)(qk)^{-1}(pk)^{-1}.\\
\end{eqnarray*}
If $\delta=1$, then
\begin{eqnarray*}
&&\{(ik)^\delta,(qk)(pk)(qk)^{-1}(pk)^{-1}\}(pq)^{-1}(jl)^{\varepsilon}\\
&\equiv&(pk)(qk)(pk)^{-1}(qk)^{-1}\{(jl)^{\varepsilon},\{\{(ql),(il)^{-1}(kl)^{-1}(il)(kl)\},\{(kl),(il)(kl)\}\}\\
&&\{\{(ql)^{-1},(il)^{-1}(kl)^{-1}(il)(kl)\},\{(kl),(il)(kl)\}(pl)^{-1}\}\}(ik)^\delta(qk)(pk)(qk)^{-1}(pk)^{-1}\\
&\equiv&(pk)(qk)(pk)^{-1}\{(jl)^{\varepsilon},\{\{(ql),(kl)^{-1}(ql)^{-1}\},\{(kl),(ql)^{-1}\}\}\{(ql)^{-1},(kl)^{-1}\\
&&(ql)^{-1}\},\{(kl),(ql)^{-1}\}(pl)^{-1}\}\}(qk)^{-1}(ik)^\delta(qk)(pk)(qk)^{-1}(pk)^{-1}(pq)^{-1}\\
&\equiv&(pk)(qk)\{\{(jl)^{\varepsilon},(kl)(pl)(kl)^{-1}(pl)^{-1}\},\{(ql),(kl)(pl)(kl)^{-1}(pl)^{-1}\}\\
&&\{(pl),(kl)^{-1}(pl)^{-1}\}\{(ql)^{-1},(kl)(pl)(kl)^{-1}(pl)^{-1}\}\{(pl)^{-1},(kl)^{-1}(pl)^{-1}\}\\
&&(pk)^{-1}(qk)^{-1}(ik)^\delta(qk)(pk)(qk)^{-1}(pk)^{-1}(pq)^{-1}\\
&\equiv&(pk)\{(jl)^{\varepsilon},\{(pl),\{(kl)^{-1},(ql)(kl)\}\{(ql)^{-1},(kl)\}\}(pl)^{-1}\}\\
&&\{(ik)^\delta,(qk)(pk)(qk)^{-1}\}(pk)^{-1}(pq)^{-1}\\
&\equiv&\{\{(jl)^{\varepsilon},(pl)^{-1}(kl)^{-1}(pl)(kl)\},\{\{(pl),(kl)\},\{(kl)^{-1},(pl)(kl)\}\{(ql)^{-1},\\
&&(pl)^{-1}(kl)^{-1}(pl)(kl)\}\}\{(pl)^{-1},(kl)\}\}\{(ik)^\delta,(qk)(pk)(qk)^{-1}(pk)^{-1}\}(pq)^{-1}\\
&\equiv&\{(jl)^{\varepsilon},(ql)(pl)(ql)^{-1}(pl)^{-1}\}\{(ik)^\delta,(qk)(pk)(qk)^{-1}(pk)^{-1}\}(pq)^{-1}.
\end{eqnarray*}
If $\delta=-1$, then
\begin{eqnarray*}
&&\{(ik)^\delta,(qk)(pk)(qk)^{-1}(pk)^{-1}\}(pq)^{-1}(jl)^{\varepsilon}\\
&\equiv&(pk)(qk)(pk)^{-1}(qk)^{-1}\{(jl)^{\varepsilon},\{\{(ql),(kl)(il)(kl)^{-1}(il)^{-1}\},\{(kl),(il)^{-1}\}\}\\
&&\{(ql)^{-1},(kl)(il)(kl)^{-1}(il)^{-1}\},\{(kl),(il)^{-1}\}(pl)^{-1}\}\}(ik)^\delta(qk)(pk)(qk)^{-1}(pk)^{-1}\\
&\equiv&(pk)(qk)(pk)^{-1}(qk)^{-1}\{(jl)^{\varepsilon},\{(ql),(kl)\}\{(ql)^{-1},(kl)(pl)^{-1}\}\}\\
&&(ik)^\delta(qk)(pk)(qk)^{-1}(pk)^{-1}(pq)^{-1}\\
&\equiv&\{(jl)^{\varepsilon},(ql)(pl)(ql)^{-1}(pl)^{-1}\}\{(ik)^\delta,(qk)(pk)(qk)^{-1}(pk)^{-1}\}(pq)^{-1}.
\end{eqnarray*}
\item[b)]\  $p=j,p=j<i<q<k<l$.
\begin{eqnarray*}
&&(pq)^{-1}(jl)^{\varepsilon}(ik)^\delta=(jq)^{-1}(jl)^{\varepsilon}(ik)^\delta\\
&\equiv&\{(jl)^{\varepsilon},(ql)^{-1}(jl)^{-1}\}\{(ik)^\delta,(qk)(jk)(qk)^{-1}(jk)^{-1}\}(jq)^{-1} \ \ \ \ \ \ and\\
&&\\
&&\{(ik)^\delta,(qk)(pk)(qk)^{-1}(pk)^{-1}\}(pq)^{-1}(jl)^{\varepsilon}\\
&=&\{(ik)^\delta,(qk)(jk)(qk)^{-1}(jk)^{-1}\}(jq)^{-1}(jl)^{\varepsilon}\\
&\equiv&(jk)(qk)(jk)^{-1}(qk)^{-1}(ik)^\delta(qk)(jk)(qk)^{-1}(jk)^{-1}\{(jl)^{\varepsilon},(ql)^{-1}(jl)^{-1}\}(jq)^{-1}\\
&\equiv&(jk)(qk)(jk)^{-1}(qk)^{-1}(ik)^\delta(qk)(jk)(qk)^{-1}\{\{(jl)^{\varepsilon},(kl)^{-1}(jl)^{-1}\},\\
&&\{(ql)^{-1},(kl)(jl)(kl)^{-1}(jl)^{-1}\}\{(jl)^{-1},(kl)^{-1}(jl)^{-1}\}\}(jk)^{-1}(jq)^{-1}\\
&\equiv&(jk)(qk)(jk)^{-1}(qk)^{-1}(ik)^\delta(qk)(jk)\{(jl)^{\varepsilon},\{(kl)^{-1},(ql)^{-1}\}\\
&&\{(ql)^{-1},(kl)^{-1}(ql)^{-1}\}(jl)^{-1}\}(qk)^{-1}(jk)^{-1}(jq)^{-1}\\
&\equiv&(jk)(qk)(jk)^{-1}(qk)^{-1}(ik)^\delta(qk)\{\{(jl)^{\varepsilon},(kl)\},\{(kl)^{-1},(jl)(kl)\}\\
&&\{(ql)^{-1},(jl)^{-1}(kl)^{-1}(jl)(kl)\}\{(jl)^{-1},(kl)\}\}(jk)(qk)^{-1}(jk)^{-1}(jq)^{-1}\\
&\equiv&(jk)(qk)(jk)^{-1}(qk)^{-1}(ik)^\delta\{(jl)^{\varepsilon},\{(ql)^{-1},(kl)\}(jl)^{-1}\}\\
&&(qk)(jk)(qk)^{-1}(jk)^{-1}(jq)^{-1}.
\end{eqnarray*}
If $\delta=1$, then
\begin{eqnarray*}
&&\{(ik)^\delta,(qk)(pk)(qk)^{-1}(pk)^{-1}\}(pq)^{-1}(jl)^{\varepsilon}\\
&\equiv&(jk)(qk)(jk)^{-1}(qk)^{-1}\{(jl)^{\varepsilon},\{(kl)^{-1},(il)(kl)\}\{(ql)^{-1},(il)^{-1}(kl)^{-1}(il)(kl)\}\\
&&\{(kl),(il)(kl)\}(jl)^{-1}\}(ik)^\delta(qk)(jk)(qk)^{-1}(jk)^{-1}(jq)^{-1}\\
&\equiv&(jk)(qk)(jk)^{-1}\{(jl)^{\varepsilon},\{(kl)^{-1},(ql)^{-1}\}\{(ql)^{-1},(kl)^{-1}(ql)^{-1}\}\\
&&\{(kl),(ql)^{-1}\}(jl)^{-1}\}(qk)^{-1}(ik)^\delta(qk)(jk)(qk)^{-1}(jk)^{-1}(jq)^{-1}\\
&\equiv&(jk)(qk)\{\{(jl)^{\varepsilon},(kl)^{-1}(jl)^{-1}\},\{(ql)^{-1},(kl)(jl)(kl)^{-1}(jl)^{-1}\}\\
&&\{(jl)^{-1},(kl)^{-1}(jl)^{-1}\}\}(jk)^{-1}(qk)^{-1}(ik)^\delta(qk)(jk)(qk)^{-1}(jk)^{-1}(jq)^{-1}\\
&\equiv&(jk)\{(jl)^{\varepsilon},\{(kl)^{-1},(ql)(kl)\}\{(ql)^{-1},(kl)\}(jl)^{-1}\}\\
&&(qk)(jk)^{-1}(qk)^{-1}(ik)^\delta(qk)(jk)(qk)^{-1}(jk)^{-1}(jq)^{-1}\\
&\equiv&\{\{(jl)^{\varepsilon},(kl)\},\{(kl)^{-1},(jl)(kl)\}\{(ql)^{-1},(jl)^{-1}(kl)^{-1}(jl)(kl)\}\\
&&\{(jl)^{-1},(kl)\}\}\{(ik)^\delta,(qk)(jk)(qk)^{-1}(jk)^{-1}\}(jq)^{-1}\\
&\equiv&\{(jl)^{\varepsilon},(ql)^{-1}(jl)^{-1}\}\{(ik)^\delta,(qk)(jk)(qk)^{-1}(jk)^{-1}\}(jq)^{-1}.
\end{eqnarray*}
If $\delta=-1$, then
\begin{eqnarray*}
&&\{(ik)^\delta,(qk)(pk)(qk)^{-1}(pk)^{-1}\}(pq)^{-1}(jl)^{\varepsilon}\\
&\equiv&(jk)(qk)(jk)^{-1}(qk)^{-1}\{(jl)^{\varepsilon},\{(kl)^{-1},(il)^{-1}\}\{(ql)^{-1},(kl)(il)(kl)^{-1}(il)^{-1}\}\\
&&\{(kl),(il)^{-1}\}(jl)^{-1}\}(ik)^\delta(qk)(jk)(qk)^{-1}(jk)^{-1}(jq)^{-1}\\
&\equiv&(jk)(qk)(jk)^{-1}(qk)^{-1}\{(jl)^{\varepsilon},(kl)^{-1}(ql)^{-1}(kl)(jl)^{-1}\}\\
&&(ik)^\delta(qk)(jk)(qk)^{-1}(jk)^{-1}(jq)^{-1}\\
&\equiv&\{(jl)^{\varepsilon},(ql)^{-1}(jl)^{-1}\}\{(ik)^\delta,(qk)(jk)(qk)^{-1}(jk)^{-1}\}(jq)^{-1}.
\end{eqnarray*}
\item[c)]\  $p>j,j<p<i<q<k<l$.
\begin{eqnarray*}
(f,g)_{w}&\equiv&(jl)^{\varepsilon}\{(ik)^\delta,(qk)(pk)(qk)^{-1}(pk)^{-1}\}(pq)^{-1}\\
&&-(jl)^{\varepsilon}\{(ik)^\delta,(qk)(pk)(qk)^{-1}(pk)^{-1}\}(pq)^{-1}\\
&\equiv&0.
\end{eqnarray*}
\end{enumerate}
\item[2)]\ $ i<k<j<l,p<i<q<k$, then \ $p<i<q<k<j<l$.
\begin{eqnarray*}
(f,g)_w&=&(pq)^{-1}(jl)^{\varepsilon}(ik)^\delta-\{(ik)^\delta,(qk)(pk)(qk)^{-1}(pk)^{-1}\}(pq)^{-1}(jl)^{\varepsilon}\\
&\equiv&(jl)^{\varepsilon}(pq)^{-1}(ik)^\delta-(jl)^{\varepsilon}\{(ik)^\delta,(qk)(pk)(qk)^{-1}(pk)^{-1}\}(pq)^{-1}\\
&\equiv&0.
\end{eqnarray*}
\end{enumerate}

$(12)\wedge(7)$\\
Let $f=(pq)(jk)^{-1}-\{(jk)^{-1},(pk)^{-1}(qk)^{-1}(pk)(qk)\}(pq),\
p<j<q<k,\
g=(jk)^{-1}(kl)^{\varepsilon}-\{(kl)^{\varepsilon},(jl)^{-1}\}(jk)^{-1},\
 j<k<l$. Then $p<j<q<k<l,\ w=(pq)(jk)^{-1}(kl)^{\varepsilon}$ and
$$
(f,g)_w=(pq)\{(kl)^{\varepsilon},(jl)^{-1}\}(jk)^{-1}-\{(jk)^{-1},(pk)^{-1}(qk)^{-1}(pk)(qk)\}(pq)(kl)^{\varepsilon}\equiv
0
$$
since
\begin{eqnarray*}
&&(pq)\{(kl)^{\varepsilon},(jl)^{-1}\}(jk)^{-1}\\
&\equiv&\{(kl)^{\varepsilon},\{(jl)^{-1},(pl)^{-1}(ql)^{-1}(pl)(ql)\}\}
\{(jk)^{-1},(pk)^{-1}(qk)^{-1}(pk)(qk)\}(pq) \ \ \ and\\
&&\\
&&\{(jk)^{-1},(pk)^{-1}(qk)^{-1}(pk)(qk)\}(pq)(kl)^{\varepsilon}\\
&\equiv&(qk)^{-1}(pk)^{-1}(qk)(pk)(jk)^{-1}(pk)^{-1}(qk)^{-1}(pk)(qk)(kl)^{\varepsilon}(pq)\\
&\equiv&(qk)^{-1}(pk)^{-1}(qk)(pk)(jk)^{-1}(pk)^{-1}(qk)^{-1}(pk)\{(kl)^{\varepsilon},(ql)(kl)\}(qk)(pq)\\
&\equiv&(qk)^{-1}(pk)^{-1}(qk)(pk)(jk)^{-1}(pk)^{-1}(qk)^{-1}\{\{(kl)^{\varepsilon},(pl)(kl)\},\{(ql),(pl)^{-1}\\
&&(kl)^{-1}(pl)(kl)\}\{(kl),(pl)(kl)\}\}(pk)(qk)(pq)\\
&\equiv&(qk)^{-1}(pk)^{-1}(qk)(pk)(jk)^{-1}(pk)^{-1}\{\{(kl)^{\varepsilon},(ql)^{-1}\},(pl)\{(ql),(kl)^{-1}(ql)^{-1}\}\\
&&\{(kl),(ql)^{-1}\}\}(qk)^{-1}(pk)(qk)(pq)\\
&\equiv&(qk)^{-1}(pk)^{-1}(qk)(pk)(jk)^{-1}\{\{(kl)^{\varepsilon},(pl)^{-1}\},\{\{(pl),(kl)^{-1}(pl)^{-1}\},\\
&&\{(ql),(kl)(pl)(kl)^{-1}(pl)^{-1}\}\}\{(kl),(ql)^{-1}\}\}(pk)^{-1}(qk)^{-1}(pk)(qk)(pq)\\
&\equiv&(qk)^{-1}(pk)^{-1}(qk)(pk)\{\{(kl)^{\varepsilon},(jl)^{-1}\}(pl)^{-1}\{(pl),\{(kl)^{-1}(jl)^{-1}\}\\
&&\{(ql),(kl)(jl)(kl)^{-1}(jl)^{-1}\}\{(kl),(jl)^{-1}\}\}\}(jk)^{-1}(pk)^{-1}(qk)^{-1}(pk)(qk)(pq)\\
&\equiv&(qk)^{-1}(pk)^{-1}(qk)\{\{(kl)^{\varepsilon},(pl)(kl)\},\{(jl)^{-1},(pl)^{-1}(kl)^{-1}(pl)(kl)\}\\
&&\{(pl)^{-1},(kl)\}\{\{(pl),(kl)\}\{\{(ql),(pl)^{-1}(kl)^{-1}(pl)(kl)\},\{(kl),(pl)(kl)\}\}\}\}\\
&&(pk)(jk)^{-1}(pk)^{-1}(qk)^{-1}(pk)(qk)(pq)\\
&\equiv&(qk)^{-1}(pk)^{-1}\{\{(kl)^{\varepsilon},(ql)(kl)\},(pl)(jl)^{-1}(pl)^{-1}\{(pl),\{(ql),(kl)\}\\
&&\{(kl),(ql)(kl)\}\}\}(qk)(pk)(jk)^{-1}(pk)^{-1}(qk)^{-1}(qk)(qk)(pq)\\
&\equiv&(qk)^{-1}\{\{(kl)^{\varepsilon},(pl)^{-1}\},\{(ql),(kl)(pl)(kl)^{-1}(pl)^{-1}\}\{(kl),(pl)^{-1}\}\\
&&\{(pl),(kl)^{-1}(pl)^{-1}\}\{(jl)^{-1},(kl)(pl)(kl)^{-1}(pl)^{-1}\}\{(pl)^{-1},(kl)^{-1}(pl)^{-1}\}\\
&&\{\{(pl),(kl)^{-1}(pl)^{-1}\},\{(ql),(kl)(pl)(kl)^{-1}(pl)^{-1}\}\{(kl),(pl)^{-1}\}\}\}\\
&&(pk)^{-1}(qk)(pk)(jk)^{-1}(pk)^{-1}(qk)^{-1}(pk)(qk)(pq)\\
&\equiv&\{\{(kl)^{\varepsilon},(pl)^{-1}\},\{(jl)^{-1},\{\{(ql)^{-1},(kl)^{-1}(ql)^{-1}\},\{(kl),(ql)^{-1}\}(pl)\}\}\\
&&\{(kl)^{-1},(ql)^{-1}\}\{(ql),(kl)^{-1}(ql)^{-1}\}\{(kl),(ql)^{-1}\}\}\{(jk)^{-1},(pk)^{-1}(qk)^{-1}(pk)(qk)\}(pq)\\
&\equiv&\{(kl)^{\varepsilon},\{(jl)^{-1},(pl)^{-1}(ql)^{-1}(pl)(ql)\}\}\{(jk)^{-1},(pk)^{-1}(qk)^{-1}(pk)(qk)\}(pq).
\end{eqnarray*}

$(12)\wedge(8)$\\
Let $f=(pq)(jk)-\{(jk)^{-1},(pk)^{-1}(qk)^{-1}(pk)(qk)\}(pq),\ \
g=(jk)(kl)^{\varepsilon}-\{(kl)^{\varepsilon},(jl)(kl)\}(jk),\
 p<j<q<k<l$. Then $w=(pq)(jk)(kl)^{\varepsilon}$ and
$$
(f,g)_w=(pq)\{(kl)^{\varepsilon},(jl)(kl)\}(jk)-\{(jk),(pk)^{-1}(qk)^{-1}(pk)(qk)\}(pq)(kl)^{\varepsilon}\equiv
0
$$
since
\begin{eqnarray*}
&&(pq)\{(kl)^{\varepsilon},(jl)(kl)\}(jk)\\
&\equiv&\{(kl)^{\varepsilon},\{jl),(pl)^{-1}(ql)^{-1}(pl)(ql)\}(kl)\}\{(jk),(pk)^{-1}(qk)^{-1}(pk)(qk)\}(pq) \ \ \ \ \ and\\
&&\\
&&\{(jk),(pk)^{-1}(qk)^{-1}(pk)(qk)\}(pq)(kl)^{\varepsilon}\\
&\equiv&(qk)^{-1}(pk)^{-1}(qk)(pk)(jk)\{(kl)^{\varepsilon},(pl)^{-1}\{(pl),(kl)^{-1}(ql)(kl)\}\}
(pk)^{-1}(qk)^{-1}(pk)(qk)(pq)\\
&\equiv&(qk)^{-1}(pk)^{-1}(qk)(pk)\{\{(kl)^{\varepsilon},(jl)(kl)\},(pl)^{-1}\{(pl),\{(kl)^{-1},(jl)(kl)\}\\
&&\{(ql),(jl)^{-1}(kl)^{-1}(jl)(kl)\}\{(kl),(jl)(kl)\}\}\}(jk)(pk)^{-1}(qk)^{-1}(pk)(qk)(pq)\\
&\equiv&(qk)^{-1}(pk)^{-1}(qk)\{\{(kl)^{\varepsilon},(pl)(kl)\},\{(jl),(pl)^{-1}(kl)^{-1}(pl)(kl)\}\\
&&\{(kl),(pl)(kl)\}\{(pl)^{-1},(kl)\}\{\{(pl),(kl)\},\{(kl)^{-1},(pl)(kl)\}\{(ql),(pl)^{-1}(kl)^{-1}(pl)(kl)\}\\
&&\{(kl),(pl)(kl)\}\}\}(pk)(kl)(pk)^{-1}(qk)^{-1}(pk)(qk)(pq)\\
&\equiv&(qk)^{-1}(pk)^{-1}\{\{(kl)^{\varepsilon},(pl)(kl)\},(pl)(jl)(pl)^{-1}\{(ql)^{-1},(kl)\}(pl)\{(ql),(kl)\}\\
&&\{(kl),(ql)(kl)\}\}(qk)(pk)(jk)(pk)^{-1}(qk)^{-1}(pk)(qk)(pq)\\
&\equiv&(qk)^{-1}\{\{(kl)^{\varepsilon},(pl)^{-1}\},\{\{(jl),(kl)(pl)(kl)^{-1}(pl)^{-1}\},\{(pl)^{-1},(kl)^{-1}(pl)^{-1}\}\\
&&\{(kl)^{-1},(pl)^{-1}\}\{(ql)^{-1},(kl)(pl)(kl)^{-1}(pl)^{-1}\}\}\{(kl),(pl)^{-1}\}\{(pl),(kl)^{-1}(pl)^{-1}\}\\
&&\{(ql),(kl)(pl)(kl)^{-1}(pl)^{-1}\}\{(kl),(pl)^{-1}\}\}(pk)^{-1}(qk)(pk)(jk)(pk)^{-1}(qk)^{-1}(pk)(qk)(pq)\\
&\equiv&\{\{(kl)^{\varepsilon},(ql)^{-1}\}\{(jl),(pl)^{-1}\{(kl)^{-1},(ql)^{-1}\}\{(ql)^{-1},(kl)^{-1}(ql)^{-1}\}\\
&&\{(kl),(ql)^{-1}\}(pl)\}\{(ql),(kl)^{-1}(ql)^{-1}\}\{(kl),(ql)^{-1}\}\}\{(jk),(pk)^{-1}(qk)^{-1}(pk)(qk)\}(pq)\\
&\equiv&\{(kl)^{\varepsilon},\{(jl),(pl)^{-1}(ql)^{-1}(pl)(ql)\}(kl)\}\{(jk),(pk)^{-1}(qk)^{-1}(pk)(qk)\}(pq).
\end{eqnarray*}

$(12)\wedge(9)$\\
Let $f=(pq)(jk)^{-1}-\{(jk)^{-1},(pk)^{-1}(qk)^{-1}(pk)(qk)\}(pq),\
g=(jk)^{-1}(jl)^{\varepsilon}-\{(jl)^{\varepsilon},(kl)^{-1}$\\
$(jl)^{-1}\}(jk)^{-1},\ p<j<q<k<l$. Then
$w=(pq)(jk)^{-1}(jl)^{\varepsilon}$ and
$$
(f,g)_w=(pq)\{(jl)^{\varepsilon},(kl)^{-1}(jl)^{-1}\}(jk)^{-1}-
\{(jk)^{-1},(pk)^{-1}(qk)^{-1}(pk)(qk)\}(pq)(jl)^{\varepsilon}\equiv
0
$$
since
\begin{eqnarray*}
&&(pq)\{(jl)^{\varepsilon},(kl)^{-1}(jl)^{-1}\}(jk)^{-1}\\
&\equiv&\{(jl)^{\varepsilon},(pl)^{-1}(ql)^{-1}(pl)(ql)(kl)^{-1}\{(jl)^{-1},(pl)^{-1}(ql)^{-1}(pl)(ql)\}\}\\
&&\{(jk)^{-1},(pk)^{-1}(qk)^{-1}(pk)(qk)\}(pq) \ \ \ \ \ \  \ and\\
&&\\
&&\{(jk)^{-1},(pk)^{-1}(qk)^{-1}(pk)(qk)\}(pq)(jl)^{\varepsilon}\\
&\equiv&(qk)^{-1}(pk)^{-1}(qk)(pk)(jk)^{-1}(pk)^{-1}(qk)^{-1}(pk)(qk)\{(jl)^{\varepsilon},(pl)^{-1}(ql)^{-1}(pl)(ql)\}(pq)\\
&\equiv&(qk)^{-1}(pk)^{-1}(qk)(pk)(jk)^{-1}(pk)^{-1}(qk)^{-1}(pk)\\
&&\{(jl)^{\varepsilon},(pl)^{-1}\{(ql)^{-1},(kl)\}(pl)\{(ql),(kl)\}\}(qk)(pq)\\
&\equiv&(qk)^{-1}(pk)^{-1}(qk)(pk)(jk)^{-1}(pk)^{-1}(qk)^{-1}\{\{(jl)^{\varepsilon},(pl)^{-1}(kl)^{-1}(pl)(kl)\},\\
&&\{(pl)^{-1},(kl)\}\{\{(ql)^{-1},(pl)^{-1}(kl)^{-1}(pl)(kl)\},\{(kl),(pl)(kl)\}\}\{(pl),(kl)\}\\
&&\{\{(ql),(pl)^{-1}(kl)^{-1}(pl)(kl)\},\{(kl),(pl)(kl)\}\}\}(pk)(qk)(pq)\\
&\equiv&(qk)^{-1}(pk)^{-1}(qk)(pk)(jk)^{-1}(pk)^{-1}\{(jl)^{\varepsilon},(pl)^{-1}\{(pl),\{(ql),(kl)^{-1}(ql)^{-1}\}\\
&&\{(kl),(ql)^{-1}\}\}\}(qk)^{-1}(pk)(qk)(pq)\\
&\equiv&(qk)^{-1}(pk)^{-1}(qk)(pk)(jk)^{-1}\{\{(jl)^{\varepsilon},(kl)(pl)(kl)^{-1}(pl)^{-1}\},\{(pl)^{-1},(kl)^{-1}\\
&&(pl)^{-1}\}\{\{(pl),(kl)^{-1}(pl)^{-1}\}\{(ql),(kl)(pl)(kl)^{-1}(pl)^{-1}\}\{(kl),(pl)^{-1}\}\}\}\\
&&(pk)^{-1}(qk)^{-1}(pk)(qk)(pq)\\
&\equiv&(qk)^{-1}(pk)^{-1}(qk)(pk)\{\{(jl)^{\varepsilon},(kl)^{-1}(jl)^{-1}\},(pl)^{-1}\{(pl),\{(kl)^{-1},(jl)^{-1}\}\\
&&\{(ql),(kl)(jl)(kl)^{-1}(jl)^{-1}\}\{(kl),(jl)^{-1}\}\}\}(jk)^{-1}(pk)^{-1}(qk)^{-1}(pk)(qk)(pq)\\
&\equiv&(qk)^{-1}(pk)^{-1}(qk)\{\{(jl)^{\varepsilon},(pl)^{-1}(kl)^{-1}(pl)(kl)\},\{(kl)^{-1},(pl)(kl)\}\\
&&\{(jl)^{-1},(pl)^{-1}(kl)^{-1}(pl)(kl)\}\{(pl)^{-1}(kl)\}\{\{(pl),(kl)\},\{(kl)^{-1},(pl)(kl)\}\\
&&\{(ql),(pl)^{-1}(kl)^{-1}(pl)(kl)\}\{(kl),(pl)(kl)\}\}\}(pk)(jk)^{-1}(pk)^{-1}(qk)^{-1}(pk)(qk)(pq)\\
&\equiv&(qk)^{-1}(pk)^{-1}\{(jl)^{\varepsilon},(pl)^{-1}\{(kl)^{-1},(ql)(kl)\}(pl)(jl)^{-1}(pl)^{-1}\{(pl),\{(ql),(kl)\}\\
&&\{(kl),(ql)(kl)\}\}\}(qk)(pk)(jk)^{-1}(pk)^{-1}(qk)^{-1}(pk)(qk)(pq)\\
&\equiv&(qk)^{-1}\{\{(jl)^{\varepsilon},(kl)(pl)(kl)^{-1}(pl)^{-1}\},\{(pl)^{-1},(kl)^{-1}(pl)^{-1}\}\\
&&\{\{(kl)^{-1},(pl)^{-1}\},\{(ql),(kl)(pl)(kl)^{-1}(pl)^{-1}\}\{(kl),(pl)^{-1}\}\}\{(jl)^{-1},(kl)(pl)\\
&&(kl)^{-1}(pl)^{-1}\},\{(pl)^{-1},(kl)^{-1}(pl)^{-1}\}\}\{\{(pl),(kl)^{-1}(pl)^{-1}\}, \{(ql),(kl)(pl)\\
&&(kl)^{-1}(pl)^{-1}\}\{(kl),(pl)^{-1}\}\}\}(pk)^{-1}(qk)(pk)(jk)^{-1}(pk)^{-1}(qk)^{-1}(pk)(qk)(pq)\\
&\equiv&\{(jl)^{\varepsilon},(pl)^{-1}\{(kl)^{-1},(ql)^{-1}\},(pl)^{-1}\{(kl)^{-1},(ql)^{-1}\}\\
&&\{(ql),(kl)^{-1}(ql)^{-1}\}\{(kl),(ql)^{-1}\}\}\{(jl)^{-1},(pl)^{-1}\}\{(pl),\{(kl)^{-1},(ql)^{-1}\}\\
&&\{(ql),(kl)^{-1}(ql)^{-1}\}\{(kl),(ql)^{-1}\}\}\}\{(jk)^{-1},(pk)^{-1}(qk)^{-1}(pk)(qk)\}(pq)\\
&\equiv&\{(jl)^{\varepsilon},(pl)^{-1}(ql)^{-1}(pl)(ql)(kl)^{-1}\{(jl)^{-1},(pl)^{-1}(ql)^{-1}(pl)(ql)\}\}\\
&&\{(jk)^{-1},(pk)^{-1}(qk)^{-1}(pk)(qk)\}(pq).
\end{eqnarray*}

$(12)\wedge(10)$\\
Let $f=(pq)(jk)-\{(jk),(pk)^{-1}(qk)^{-1}(pk)(qk)\}(pq),\
g=(jk)(jl)^{\varepsilon}-\{(jl)^{\varepsilon},(kl)\}(jk),\ \
p<j<q<k<l$. Then $w=(pq)(jk)(jl)^{\varepsilon}$ and
$$
(f,g)_w=(pq)\{(jl)^{\varepsilon},(kl)\}(jk)-\{(jk),(pk)^{-1}(qk)^{-1}(pk)(qk)\}(pq)(jl)^{\varepsilon}\equiv
0
$$
since
\begin{eqnarray*}
&&(pq)\{(jl)^{\varepsilon},(kl)\}(jk)\\
&\equiv&\{(jl)^{\varepsilon},(pl)^{-1}(ql)^{-1}(pl)(ql)(kl)\}\{(jk),(pk)^{-1}(qk)^{-1}(pk)(qk)\}(pq) \ \ \ \  \ and\\
&&\\
&&\{(jk),(pk)^{-1}(qk)^{-1}(pk)(qk)\}(pq)(jl)^{\varepsilon}\\
&\equiv&(qk)^{-1}(pk)^{-1}(qk)(pk)(jk)\{(jl)^{\varepsilon},(pl)^{-1}\{(pl),(kl)^{-1}(ql)(kl)\}\}\\
&&(pk)^{-1}(qk)^{-1}(pk)(qk)(pq)\\
&\equiv&(qk)^{-1}(pk)^{-1}(qk)(pk)\{\{(jl)^{\varepsilon},(kl)\},(pl)^{-1}\{(pl),\{(kl)^{-1},(jl)(kl)\}\\
&&\{(ql),(jl)^{-1}(kl)^{-1}(jl)(kl)\} \{(kl),(jl)(kl)\}\}\}(jk)(pk)^{-1}(qk)^{-1}(pk)(qk)(pq)\\
&\equiv&(qk)^{-1}(pk)^{-1}(qk)\{\{(jl)^{\varepsilon},(pl)^{-1}(kl)^{-1}(ql)(kl)\},\{(kl),(pl)(kl)\}\\
&&\{(pl)^{-1},(kl)\}\{\{(pl),(kl)\},\{(kl)^{-1},(pl)(kl)\}\{(ql),(pl)^{-1}(kl)^{-1}(ql)(kl)\}\\
&&\{(kl),(pl)(kl)\}\}\}(pk)(jk)(pk)^{-1}(qk)^{-1}(pk)(qk)(pq)\\
&\equiv&(qk)^{-1}(pk)^{-1}\{(jl)^{\varepsilon},(pl)^{-1}\{(kl),(ql)(kl)\}\{(pl),\{(ql),(kl)\}\{(kl),(ql)(kl)\}\}\}\\
&&(qk)(pk)(jk)(pk)^{-1}(qk)^{-1}(pk)(qk)(pq)\\
&\equiv&(qk)^{-1}\{\{(jl)^{\varepsilon},(kl)(pl)(kl)^{-1}(pl)^{-1}\},\{\{(ql)^{-1},(kl)(pl)(kl)^{-1}(pl)^{-1}\},\\
&&\{(kl),(pl)^{-1}\}\{(pl),(kl)^{-1}(pl)^{-1}\}\}\{(ql),(kl)(pl)(kl)^{-1}(pl)^{-1}\}\{(kl),(pl)^{-1}\}\}\\
&&(pk)^{-1}(qk)(pk)(jk)(pk)^{-1}(qk)^{-1}(pk)(qk)(pq)\\
&\equiv&\{(jl)^{\varepsilon},\{\{(ql)^{-1},(kl)^{-1}(ql)^{-1}\},\{(kl),(ql)^{-1}\}(pl)\}\{(ql),(kl)^{-1}(ql)^{-1}\}\\
&&\{(kl),(ql)^{-1}\}\}\{ (jk),(pk)^{-1}(qk)^{-1}(pk)(qk)\}(pq)\\
&\equiv&\{(jl)^{\varepsilon},(pl)^{-1}(ql)^{-1}(pl)(ql)(kl)\}\{(jk),(pk)^{-1}(qk)^{-1}(pk)(qk)\}(pq).
\end{eqnarray*}

$(12)\wedge(11)$\\
Let $f=(pq)(ik)^{-1}-\{(ik)^{-1},(pk)^{-1}(qk)^{-1}(pk)(qk)\}(pq), \
p<i<q<k,
g=(ik)^{-1}(jl)^{\varepsilon}-\{(jl)^{\varepsilon},(kl)(il)(kl)^{-1}(il)^{-1}\}(ik)^{-1},i<j<k<l.$
\ Then $w=(pq)(ik)^{-1}(jl)^{\varepsilon}$ and
$$
(f,g)_w=(pq)\{(jl)^{\varepsilon},(kl)(il)(kl)^{-1}(il)^{-1}\}(ik)^{-1}-
\{(ik)^{-1},(pk)^{-1}(qk)^{-1}(pk)(qk)\}(pq)(jl)^{\varepsilon} .
$$
There are three cases to consider.
\begin{enumerate}
\item[1)]\ $\ \ q<j,\ \ p<i<q<j<k<l.$
\begin{eqnarray*}
&&(pq)\{(jl)^{\varepsilon},(kl)(il)(kl)^{-1}(il)^{-1}\}(ik)^{-1}\\
&\equiv&\{(jl)^{\varepsilon} ,(kl)\{(il),(pl)^{-1}(ql)^{-1}(pl)(ql)\}(kl)^{-1}\{(il)^{-1},(pl)^{-1}(ql)^{-1}(pl)(ql)\}\}\\
&&\{(ik)^{-1},(pk)^{-1}(qk)^{-1}(pk)(qk)\}(pq)\ \ \ \ \ \ \ \ and \\
&&\\
&&\{(ik)^{-1},(pk)^{-1}(qk)^{-1}(pk)(qk)\}(pq)(jl)^{\varepsilon}.\\
&\equiv&(qk)^{-1}(pk)^{-1}(qk)(pk)(ik)^{-1}(pk)^{-1}(qk)^{-1}(pk)(qk)(jl)^{\varepsilon}(pq)\\
&\equiv&(qk)^{-1}(pk)^{-1}(pk)(qk)(ik)^{-1}(pk)^{-1}(qk)^{-1}(pk)\{(jl)^{\varepsilon},(ql)^{-1}(kl)^{-1}(ql)(kl)\}\\
&&(qk)(pq)\\
&\equiv&(qk)^{-1}(pk)^{-1}(qk)(pk)(ik)^{-1}(pk)^{-1}(qk)^{-1}\{\{(jl)^{\varepsilon},(pl)^{-1}(kl)^{-1}(pl)(kl)\},\\
&&\{(ql)^{-1},(pl)^{-1}(kl)^{-1}(pl)(kl)\}\{(kl)^{-1},(pl)(kl)\}\{(ql),(pl)^{-1}(kl)^{-1}(pl)(kl)\}\\
&&\{(kl),(pl)(kl)\}\}(pk)(qk)(pq)\\
&\equiv&(qk)^{-1}(pk)^{-1}(qk)(pk)(ik)^{-1}(pk)^{-1}\{\{(jl)^{\varepsilon},(kl)(ql)(kl)^{-1}(ql)^{-1}\},\{\{(kl)^{-1},\\
&&(ql)^{-1}\},(pl)\{(ql),(kl)^{-1}(ql)^{-1}\}\}\{(kl),(ql)^{-1}\}\}(qk)^{-1}(pk)(qk)(pq)\\
&\equiv&(qk)^{-1}(pk)^{-1}(qk)(pk)(ik)^{-1}\{\{(jl)^{\varepsilon},(kl)(pl)(kl)^{-1}(pl)^{-1}\},\{\{(kl)^{-1},(pl)^{-1}\},\\
&&\{(ql)^{-1},(kl)(pl)(kl)^{-1}(pl)^{-1}\}\{(pl),(kl)^{-1}(pl)^{-1}\}\{(ql),(kl)(pl)(kl)^{-1}(pl)^{-1}\}\}\\
&&\{(kl),(pl)^{-1}\}\}(pk)^{-1}(qk)^{-1}(pk)(qk)(pq)\\
&\equiv&(qk)^{-1}(pk)^{-1}(qk)(pk)\{\{(jl)^{\varepsilon},(kl)(il)(kl)^{-1}(il)^{-1}\},\{\{(kl)^{-1},(il)^{-1}\},\\
&&\{\{(ql)^{-1},(kl)(il)(kl)^{-1}(il)^{-1}\},\{(kl),(il)^{-1}\}(pl)\}\{(kl)^{-1},(il)^{-1}\}\{(ql),(kl)(il)\\
&&(kl)^{-1}(il)^{-1}\}\}\{(kl),(il)^{-1}\}\}(ik)^{-1}(pk)^{-1}(qk)^{-1}(pk)(qk)(pq)\\
&\equiv&(qk)^{-1}(pk)^{-1}(qk)\{\{(jl)^{\varepsilon},(pl)^{-1}(kl)^{-1}(pl)(kl)\},\{\{(kl)^{-1},(pl)(kl)\},\{(il)^{-1},\\
&&(pl)^{-1}(kl)^{-1}(pl)(kl)\}\{\{(ql)^{-1},(pl)^{-1}(kl)^{-1}(pl)(kl)\},\{(kl),(pl)(kl)\}\\
&&\{(pl),(kl)\}\}\{(kl)^{-1},(pl)(kl)\}\{(ql),(pl)^{-1}(kl)^{-1}(pl)(kl)\}\}\{(kl),(pl)(kl)\}\}\\
&&(pk)(ik)^{-1}(pk)^{-1}(qk)^{-1}(pk)(qk)(pq)\\
&\equiv&(qk)^{-1}(pk)^{-1}\{\{(jl)^{\varepsilon},(ql)^{-1}(kl)^{-1}(ql)(kl)\},\{\{(kl)^{-1},(ql)(kl)\},(pl)(il)^{-1}\\
&&(pl)^{-1}\{(kl)^{-1},(ql)(kl)\}\{(ql)^{-1},(kl)\}(pl)\{(ql),(kl)\}\}\{(kl),(ql)(kl)\}\}\\
&&(qk)(pk)(ik)^{-1}(pk)^{-1}(qk)^{-1}(pk)(qk)(pq)\\
&\equiv&(qk)^{-1}\{\{(jl)^{\varepsilon},(kl)(pl)(kl)^{-1}(pl)^{-1}\},\{\{(kl)^{-1},(pl)^{-1}\},\{\{(il)^{-1},(kl)(pl)\\
&&(kl)^{-1}(pl)^{-1}\},\{(pl)^{-1},(kl)^{-1}(pl)^{-1}\}\{(kl)^{-1},(pl)^{-1}\}\{(ql)^{-1},(kl)(pl)\\
&&(kl)^{-1}(pl)^{-1}\}\}\{(pl),(kl)^{-1}(pl)^{-1}\}\{(ql),(kl)(pl)(kl)^{-1}(pl)^{-1}\}\}\\
&&\{(kl),(pl)^{-1}\}\}(pk)^{-1}(qk)(pk)(ik)^{-1}(pk)^{-1}(qk)^{-1}(pk)(qk)(pq)\\
&\equiv&\{\{(jl)^{\varepsilon},(kl)(ql)(kl)^{-1}(ql)^{-1}\},\{\{(kl)^{-1},(ql)^{-1}\},\{(il)^{-1},(pl)^{-1}\\
&&\{(kl)^{-1},(ql)^{-1}\}\{(ql)^{-1},(kl)^{-1}(ql)^{-1}\}\{(kl),(ql)^{-1}\}(pl)\{(kl)^{-1},(ql)^{-1}\}\}\\
&&\{(ql),(kl)^{-1}(ql)^{-1}\}\}\{(kl),(ql)^{-1}\}\}\{(ik)^{-1},(pk)^{-1}(qk)^{-1}(pk)(qk)\}(pq)\\
&\equiv&\{(jl)^{\varepsilon} ,(kl)\{(il),(pl)^{-1}(ql)^{-1}(pl)(ql)\}(kl)^{-1}\{(il)^{-1},(pl)^{-1}(ql)^{-1}(pl)(ql)\}\}\\
&&\{(ik)^{-1},(pk)^{-1}(qk)^{-1}(pk)(qk)\}(pq).
\end{eqnarray*}
\item[2)]\  $q=j, \ p<i<q=j<k<l$.
\begin{eqnarray*}
&&(pq)\{(jl)^{\varepsilon},(kl)(il)(kl)^{-1}(il)^{-1}\}(ik)^{-1}\\
&=&(pj)\{(jl)^{\varepsilon},(kl)(il)(kl)^{-1}(il)^{-1}\}(ik)^{-1}\\
&\equiv&\{\{(jl)^{\varepsilon},(pl)(jl)\},(kl)\{(il),(pl)^{-1}(jl)^{-1}(pl)(jl)\}(kl)^{-1}\{(il)^{-1},(pl)^{-1}\\
&&(jl)^{-1}(pl)(jl)\}\}\{(ik)^{-1},(pk)^{-1}(jk)^{-1}(pk)(jk)\}(pj) \ \ \ \ \ \ \ \ and \\
&&\\
&&\{(ik)^{-1},(pk)^{-1}(qk)^{-1}(pk)(qk)\}(pq)(jl)^{\varepsilon}\\
&=&\{(ik)^{-1},(pk)^{-1}(jk)^{-1}(pk)(jk)\}(pj)(jl)^{\varepsilon}\\
&\equiv&(jk)^{-1}(pk)^{-1}(jk)(pk)(ik)^{-1}(pk)^{-1}(jk)^{-1}(pk)(jk)\{(jl)^{\varepsilon},(pl)(jl)\}(pj)\\
&\equiv&(jk)^{-1}(pk)^{-1}(jk)(pk)(ik)^{-1}(pk)^{-1}(jk)^{-1}(pk)\{\{(jl)^{\varepsilon},(kl)\},(pl)\{(jl),(kl)\}\}\\
&&(jk)(pj)\\
&\equiv&(jk)^{-1}(pk)^{-1}(jk)(pk)(ik)^{-1}(pk)^{-1}(jk)^{-1}\{\{(jl)^{\varepsilon},(pl)^{-1}(kl)^{-1}(pl)(kl)\},\{(kl),\\
&&(pl)(kl)\}\{(pl),(kl)\}\{(kl)^{-1},(pl)(kl)\}\{(jl),(pl)^{-1}(kl)^{-1}(pl)(kl)\}\\
&&\{(kl),(pl)(kl)\}\}(pk)(jk)(pj)\\
&\equiv&(jk)^{-1}(pk)^{-1}(jk)(pk)(ik)^{-1}(pk)^{-1}\{\{(jl)^{\varepsilon},(kl)^{-1}(jl)^{-1}\},\\
&&(pl)\{(jl),(kl)^{-1}(jl)^{-1}\}\}(jk)^{-1}(pk)(jk)(pj)\\
&\equiv&(jk)^{-1}(pk)^{-1}(jk)(pk)(ik)^{-1}\{\{(jl)^{\varepsilon},(kl)(pl)(kl)^{-1}(pl)^{-1}\},\{\{(pl),(kl)^{-1}\\
&&(pl)^{-1}\},\{(jl),(kl)(pl)(kl)^{-1}(pl)^{-1}\}\{(kl),(pl)^{-1}\}\}\}(pk)^{-1}(jk)^{-1}(pk)(jk)(pj)\\
&\equiv&(jk)^{-1}(pk)^{-1}(jk)(pk)\{\{(jl)^{\varepsilon},(kl)(il)(kl)^{-1}(il)^{-1}\},\{(kl),(il)^{-1}\}(pl)\\
&&\{(kl)^{-1},(il)^{-1}\}(pl)^{-1}\{(pl),\{(kl)^{-1}(il)^{-1}\}\{(jl),(kl)(il)(kl)^{-1}(il)^{-1}\}\\
&&\{(kl),(il)^{-1}\}\}\}(ik)^{-1}(pk)^{-1}(jk)^{-1}(pk)(jk)(pj)\\
&\equiv&(jk)^{-1}(pk)^{-1}(jk)\{\{(jl)^{\varepsilon},(pl)^{-1}(kl)^{-1}(pl)(kl)\},\{(kl),(pl)(kl)\}\{\{(kl)^{-1},(pl)(kl)\},\\
&&\{(il)^{-1},(pl)^{-1}(kl)^{-1}(pl)(kl)\}\{(pl)^{-1},(kl)\}\}\{\{(pl),(kl)\},\{(kl)^{-1},(pl)(kl)\}\\
&&\{(jl),(pl)^{-1}(kl)^{-1}(pl)(kl)\}\{(kl),(pl)(kl)\}\}\}(pk)(ik)^{-1}(pk)^{-1}(jk)^{-1}(pk)(jk)(pj)\\
&\equiv&(jk)^{-1}(pk)^{-1}\{\{(jl)^{\varepsilon},(kl)\},\{(kl),(jl)(kl)\}\{\{(kl)^{-1},(jl)(kl)\},(pl)(il)^{-1}\\
&&(pl)^{-1}\}\{(pl),(jl),(kl)\}\{(kl),(jl)(kl)\}\}\}(jk)(pk)(ik)^{-1}(pk)^{-1}(jk)^{-1}(pk)(jp)(pj)\\
&\equiv&(jk)^{-1}\{\{(jl)^{\varepsilon},(kl)(pl)(kl)^{-1}(pl)^{-1}\},\{(kl),(pl)^{-1}\}\{(jl),(kl)(pl)(kl)^{-1}(pl)^{-1}\}\\
&&\{(kl),(pl)^{-1}\}\{\{(kl)^{-1},(pl)^{-1}\},\{(jl),(kl)(pl)(kl)^{-1}(pl)^{-1}\}\{(kl),(pl)^{-1}\}\\
&&\{\{(il)^{-1},(kl)(pl)(kl)^{-1}(pl)^{-1}\},\{(pl)^{-1},(kl)^{-1}(pl)^{-1}\}\}\}\{\{(pl),(kl)^{-1}(pl)^{-1}\},\\
&&\{(jl),(kl)(pl)(kl)^{-1}(pl)^{-1}\}\{(kl),(pl)^{-1}\}\}\}(pk)^{-1}(jk)(pk)(ik)^{-1}\\
&&(pk)^{-1}(jk)^{-1}(pk)(jk)(pj)\\
&\equiv&\{\{(jl)^{\varepsilon},(kl)^{-1}(jl)^{-1}\},\{(il),(pl)^{-1}\{(kl)^{-1},(jl)^{-1}\}\{(jl)^{-1},(kl)^{-1}(jl)^{-1}\}\\
&&\{(kl),(jl)^{-1}\}(pl)\{(kl)^{-1},(jl)^{-1}\}(pl)^{-1}\{(kl)^{-1},(jl)^{-1}\}\}\{\{(il)^{-1},(pl)^{-1}\\
&&\{(kl)^{-1},(jl)^{-1}\}\{(jl)^{-1},(kl)^{-1}(jl)^{-1}\}\}\{(kl),(jl)^{-1}\}(pl)\{(kl)^{-1},(jl)^{-1}\}\\
&&\{(jl),(kl)^{-1}(jl)^{-1}\}\{(kl),(jl)^{-1}\}\}\{(ik)^{-1},(pk)^{-1}(jk)^{-1}(pk)(jk)\}(pj)\\
&\equiv&\{\{(jl)^{\varepsilon},(pl)(jl)\},(kl)\{(il),(pl)^{-1}(jl)^{-1}(pl)(jl)\}(kl)^{-1}\{(il)^{-1},(pl)^{-1}(jl)^{-1}\\
&&(pl)(jl)\}\}\{(ik)^{-1},(pk)^{-1}(jk)^{-1}(pk)(jk)\}(pj).
\end{eqnarray*}
\item[3)]\  \ $  q>j, \ p<i<j<q<k<l$.
\begin{eqnarray*}
&&(pq)\{(jl)^{\varepsilon},(kl)(il)(kl)^{-1}(il)^{-1}\}(ik)^{-1}\\
&\equiv&\{\{(jl)^{\varepsilon},(pl)^{-1}(ql)^{-1}(pl)(ql)\},(kl)\{(il),(pl)^{-1}(ql)^{-1}(pl)(ql)\}(kl)^{-1}\\
&&\{(il)^{-1},(pl)^{-1}(ql)^{-1}(pl)(ql)\}\}\{(ik)^{-1},(pk)^{-1}(qk)^{-1}(pk)(qk)\}(pq)\ \ \ \ \ \ and \\
&&\\
&&\{(ik)^{-1},(pk)^{-1}(qk)^{-1}(pk)(qk)\}(pq)(jl)^{\varepsilon}\\
&\equiv&(qk)^{-1}(pk)^{-1}(qk)(pk)(ik)^{-1}(pk)^{-1}(qk)^{-1}(pk)(qk)\{(jl)^{\varepsilon},(pl)^{-1}(ql)^{-1}(pl)(ql)\}(pq)\\
&\equiv&(qk)^{-1}(pk)^{-1}(qk)(pk)(ik)^{-1}(pk)^{-1}(qk)^{-1}(pk)\{(jl)^{\varepsilon},(pl)^{-1}\{(ql)^{-1},(kl)\}(pl)\\
&&\{(ql),(kl)\}\}(qk)(pq)\\
&\equiv&(qk)^{-1}(pk)^{-1}(qk)(pk)(ik)^{-1}(pk)^{-1}(qk)^{-1}\{\{(jl)^{\varepsilon},(pl)^{-1}(kl)^{-1}(pl)(kl)\},\\
&&\{(pl)^{-1},(kl)\}\{\{(ql)^{-1},(pl)^{-1}(kl)^{-1}(pl)(kl)\},\{(kl),(pl)(kl)\}\{(pl),(kl)\}\\
&&\{\{(ql),(pl)^{-1}(kl)^{-1}(ql)(kl)\},\{(kl),(pl)(kl)\}\}\}(pk)(qk)(pq)\\
&\equiv&(qk)^{-1}(pk)^{-1}(qk)(pk)(ik)^{-1}(pk)^{-1}\{(jl)^{\varepsilon},(pl)^{-1}\{(pl),\{(ql),(kl)^{-1}(ql)^{-1}\}\\
&&\{(kl),(ql)^{-1}\}\}\}(qk)^{-1}(pk)(qk)(pq)\\
&\equiv&(qk)^{-1}(pk)^{-1}(qk)(pk)(ik)^{-1}\{\{(jl)^{\varepsilon},(kl)(pl)(kl)^{-1}(pl)^{-1}\},\{(pl)^{-1},(kl)^{-1}\\
&&(pl)^{-1}\}\{\{(pl),(kl)^{-1}(pl)^{-1}\},\{(ql),(kl)(pl)(kl)^{-1}(pl)^{-1}\}\{(kl),(pl)^{-1}\}\}\}\\
&&(pk)^{-1}(qk)^{-1}(pk)(qk)(pq)\\
&\equiv&(qk)^{-1}(pk)^{-1}(qk)(pk)\{\{(jl)^{\varepsilon},(kl)(il)(kl)^{-1}(il)^{-1}\},(pl)^{-1}\{(pl),\{(kl)^{-1}(il)^{-1}\}\\
&&\{(ql),(kl)(il)(kl)^{-1}(il)^{-1}\{(kl),(il)^{-1}\}\}\}(ik)^{-1}(pk)^{-1}(qk)^{-1}(pk)(qk)(pq)\\
&\equiv&(qk)^{-1}(pk)^{-1}(qk)\{\{(jl)^{\varepsilon},(pl)^{-1}(kl)^{-1}(pl)(kl)\},\{(kl),(pl)(kl)\}\{(il),(pl)^{-1}\\
&&(kl)^{-1}(pl)(kl)\}\{(kl)^{-1},(pl)(kl)\}\{(il)^{-1},(pl)^{-1}(kl)^{-1}(pl)(kl)\}\{\{(pl),(kl)\},\\
&&\{(kl)^{-1},(pl)(kl)\}\{(ql),(pl)^{-1}(kl)^{-1}(pl)(kl)\}\{(kl),(pl)(kl)\}\}\}\\
&&(pk)(ik)^{-1}(pk)^{-1}(qk)^{-1}(pk)(qk)(pq)\\
&\equiv&(qk)^{-1}(pk)^{-1}\{(jl)^{\varepsilon},(pl)^{-1}\{\{(kl)^{-1},(ql)(kl)\},(pl)(il)^{-1}(pl)^{-1}\{(kl)^{-1},(ql)(kl)\}\}\\
&&\{(ql)^{-1},(kl)\}(pl)\{(ql),(kl)\}\{(kl),(ql)(kl)\}\}(qk)(pk)(ik)^{-1}\\
&&(pk)^{-1}(qk)^{-1}(pk)(qk)(pq)\\
&\equiv&(qk)^{-1}\{\{(jl)^{\varepsilon},(kl)(pl)(kl)^{-1}(pl)^{-1}\},\{\{(il),(kl)(pl)(kl)^{-1}(pl)^{-1}\},\{\{(kl)^{-1},\\
&&(pl)^{-1}\},\{(ql),(kl)(pl)(kl)^{-1}(pl)^{-1}\}\{(kl),(pl)^{-1}\}\{(pl),(kl)^{-1}(pl)^{-1}\}\}\}\\
&&\{(il)^{-1},(kl)(pl)(kl)^{-1}(pl)^{-1}\}\{(pl)^{-1},(kl)^{-1}(pl)^{-1}\}\{(kl)^{-1},(pl)^{-1}\}\\
&&\{(ql)^{-1},(kl)(pl)(kl)^{-1}(pl)^{-1}\}\{(pl),(kl)^{-1}(pl)^{-1}\}\{(ql),(kl)(pl)(kl)^{-1}(pl)^{-1}\}\\
&&\{(kl),(pl)^{-1}\}\}(pk)^{-1}(qk)(pk)(ik)^{-1}(pk)^{-1}(qk)^{-1}(pk)(qk)(pq)\\
&\equiv&\{(jl)^{\varepsilon},\{(il),\{\{(kl)^{-1},(ql)^{-1}\},(pl)^{-1}\{\{(ql),(kl)^{-1}(ql)^{-1}\},\{(kl),(ql)^{-1}\}\}\\
&&(pl)\}\}(il)^{-1}\{\{(ql)^{-1},(kl)^{-1}(ql)^{-1}\},\{(kl),(ql)^{-1}\}(pl)\}\\
&&\{\{(ql),(kl)^{-1}(ql)^{-1},\{(kl),(ql)^{-1}\}\}\}\{(ik)^{-1}n,(pk)^{-1}(qk)^{-1}(pk)(qk)\}(pq)\\
&\equiv&\{\{(jl)^{\varepsilon},(pl)^{-1}(ql)^{-1}(pl)(ql)\},(kl)\{(il),(pl)^{-1}(ql)^{-1}(pl)(ql)\}(kl)^{-1}\\
&&\{(il)^{-1},(pl)^{-1}(ql)^{-1}(pl)(ql)\}\}\{(ik)^{-1},(pk)^{-1}(qk)^{-1}(pk)(qk)\}(pq).
\end{eqnarray*}
\end{enumerate}

$(12)\wedge(12)$ \\
Let $f=(pq)(ik)-\{(ik),(pk)^{-1}(qk)^{-1}(pk)(qk)\}(pq),p<i<q<k, \
g=(ik)(jl)^{\varepsilon}-\{(jl)^{\varepsilon},(il)^{-1}(kl)^{-1}(il)(kl)\}(ik),i<j<k<l$.
Then $w=(pq)(ik)(jl)^{\varepsilon}$ and
$$
(f,g)_{w}=(pq)\{(jl)^{\varepsilon},(il)^{-1}(kl)^{-1}(il)(kl)\}(ik)-\{(ik),
(pk)^{-1}(qk)^{-1}(pk)(qk)\}(pq)(jl)^{\varepsilon}.
$$
There are three cases to consider.
\begin{enumerate}
\item[1)]\ $q<j,p<i<q<j<k<l$.
\begin{eqnarray*}
&&(pq)\{(jl)^{\varepsilon},(il)^{-1}(kl)^{-1}(il)(kl)\}(ik)\\
&\equiv&\{(jl)^{\varepsilon},\{(il)^{-1},(pl)^{-1}(ql)^{-1}(pl)(ql)\}(kl)^{-1}\{(il),
(pl)^{-1}(ql)^{-1}(pl)(ql)\}(kl)\}\\
&& \ \ \{(ik),(pk)^{-1}(qk)^{-1}(pk)(qk)\}(pq) \ \ \ \  \ \ and \\
&&\\
&&\{(ik),(pk)^{-1}(qk)^{-1}(pk)(qk)\}(pq)(jl)^{\varepsilon}\\
&\equiv&(qk)^{-1}(pk)^{-1}(qk)(pk)(ik)\{(jl)^{\varepsilon},\{(kl)^{-1},\{(ql)^{-1},
(kl)(pl)\}(kl)^{-1}(ql)\}(kl)\}\\
&&(pk)^{-1}(qk)^{-1}(pk)(qk)(pq)\\
&\equiv&(qk)^{-1}(pk)^{-1}(qk)(pk)\{\{(jl)^{\varepsilon},(il)^{-1}(kl)^{-1}
(il)(kl)\},\{\{(kl)^{-1},(il)(kl)\},\{\{(ql)^{-1},\\
&&(il)^{-1}(kl)^{-1}(il)(kl)\},
\{(kl),(il)(kl)\}(pl)\}\{(kl)^{-1},(il)(kl)\}\{(ql),
(il)^{-1}(kl)^{-1}\\
&&(il)(kl)\}\}\{(kl),(il)(kl)\}\}
(ik)(pk)^{-1}(qk)^{-1}(pk)(qk)(pq)\\
&\equiv&(qk)^{-1}(qk)(pk)\{\{(jl)^{\varepsilon},(pl)^{-1}(kl)^{-1}(pl)(kl)\},\{\{(kl)^{-1},(pl)(kl)\},\{(il),(pl)^{-1}\\
&&(kl)^{-1}(pl)(kl)\}\{\{(ql)^{-1},(pl)^{-1}(kl)^{-1}(pl)(kl)\},\{(kl),(pl)(kl)\}\{(pl),(kl)\}\\
&&\{(kl)^{-1},(pl)(kl)\}\}\{(ql),(pl)^{-1}(kl)^{-1}(pl)(kl)\}\}\\
&&\{(kl),(pl)(kl)\}\}(qk)(ik)(pk)^{-1}(qk)^{-1}(pk)(qk)(pq)\\
&\equiv&(qk)^{-1}(pk)^{-1}\{\{(jl)^{\varepsilon},(ql)^{-1}(kl)^{-1}(ql)(kl)\},
\{\{(kl)^{-1},(ql)(kl)\},(pl)(il)(pl)^{-1}\\
&&\{(ql)^{-1},(kl)\}(pl)\{(ql),(kl)\}\}\{(kl),(ql)(kl)\}\}(qk)(pk)(ik)(pk)^{-1}(qk)^{-1}(pk)(qk)(pq)\\
&\equiv&(qk)^{-1}\{\{(jl)^{\varepsilon},(kl)(pl)(kl)^{-1}(pl)^{-1}\},\{\{(kl)^{-1},\\
&&(pl)^{-1}\},\{\{(il),(kl)(pl)(kl)^{-1}(pl)^{-1}\},\{(pl)^{-1},(kl)^{-1}(pl)^{-1}\}\{(kl)^{-1},(pl)^{-1}\}\\
&&\{(ql)^{-1},(kl)(pl)(kl)^{-1}(pl)^{-1}\}\}\{(kl),(pl)^{-1}\}\{(pl),(kl)^{-1}(pl)^{-1}\}\{(ql),(kl)(pl)\\
&&(kl)^{-1}(pl)^{-1}\}\}\{(kl),(pl)^{-1}\}\}(pk)^{-1}(qk)(pk)(ik)(pk)^{-1}(qk)^{-1}(pk)(qk)(pq)\\
&\equiv&\{\{(jl)^{\varepsilon},(kl)(ql)(kl)^{-1}(ql)^{-1}\},\{\{(kl)^{-1},(ql)^{-1}\},\{(il)^{-1},(pl)^{-1}\\
&&\{(kl)^{-1},(ql)^{-1}\}\{(ql)^{-1},(kl)^{-1}(ql)^{-1}\}\{(kl),(ql)^{-1}\}(pl)\}\{(kl)^{-1},(ql)^{-1}\}\{(ql),\\
&&(kl)^{-1}(ql)^{-1}\}\}\{(kl),(ql)^{-1}\}\}\{(ik),(pk)^{-1}(qk)^{-1}(pk)(qk)\}(pq)\\
&\equiv&\{(jl)^{\varepsilon},\{(il),(pl)^{-1}(ql)^{-1}(pl)(ql)\}(kl)^{-1}\{(il)^{-1},(pl)^{-1}(ql)^{-1}(pl)(ql)\}(kl)\}\\
&&\{(ik),(pk)^{-1}(qk)^{-1}(pk)(qk)\}(pq).
\end{eqnarray*}
 \item[2)]\ $q=j,p<i<q=j<k<l$.\\
\begin{eqnarray*}
&&(pq)\{(jl)^{\varepsilon},(il)^{-1}(kl)^{-1}(il)(kl)\}(ik)=(pj)\{(jl)^{\varepsilon},(il)^{-1}(kl)^{-1}(il)(kl)\}(ik)\\
&\equiv&\{\{(jl)^{\varepsilon},(pl)(jl)\},\{(il)^{-1},(pl)^{-1}(jl)^{-1}(pl)(jl)\}(kl)^{-1}\{(il),(pl)^{-1}\\
&&(jl)^{-1}(pl)(jl)\}(kl)\}\{(ik),(pk)^{-1}(jk)^{-1}(pk)(jk)\}(pj) \ \ \ \ \ \ \ and \\
&&\\
&&\{(ik),(pk)^{-1}(qk)^{-1}(pk)(qk)\}(pq)(jl)^{\varepsilon}=\{(ik),(pk)^{-1}(jk)^{-1}(pk)(jk)\}(pj)(jl)^{\varepsilon}\\
&\equiv&(jk)^{-1}(pk)^{-1}(jk)(pk)(ik)\{(jl)^{\varepsilon},(kl)(pl)(kl)^{-1}(pl)^{-1}\{(pl),(jl)(kl)\}\}\\
&&(pk)^{-1}(jk)^{-1}(pk)(jk)(pj)\\
&\equiv&(jk)^{-1}(pk)^{-1}(jk)(pk)\{\{(jl)^{\varepsilon},(il)^{-1}(kl)^{-1}(il)(kl)\},\{(kl),(il)(kl)\}\\
&&(pl)\{(kl)^{-1}(il)(kl)\}(pl)^{-1}\{(pl),\{(kl)^{-1},(il)(kl)\}\{(jl),(il)^{-1}(kl)^{-1}(il)(kl)\}\\
&&\{(kl),(il)(kl)\}\}\}(ik)(pk)^{-1}(jk)^{-1}(pk)(jk)(pj)\\
&\equiv&(jk)^{-1}(pk)^{-1}(jk)\{\{(jl)^{\varepsilon},(pl)^{-1}(kl)^{-1}(pl)(kl)\},\{(kl),(pl)(kl)\}
\{\{(kl)^{-1},(pl)(kl)\},\\
&&\{(il),(pl)^{-1}(kl)^{-1}(pl)(kl)\}\{(kl),(pl)(kl)\}\{(pl)^{-1},(kl)\}\}\{\{(pl),(kl)\},\{(kl)^{-1},(pl)(kl)\}\\
&&\{(jl),(pl)^{-1}(kl)^{-1}(pl)(kl)\}\{(kl),(pl)(kl)\}\}\}(ik)(pk)^{-1}(jk)^{-1}(pk)(jp)(pj)\\
&\equiv&(jk)^{-1}(pk)^{-1}\{\{(jl)^{\varepsilon},(kl)\},\{\{(kl)^{-1},(jl)(kl)\},(pl)(il)(pl)^{-1}\}\\
&&\{(pl),\{(jl),(kl)\}\}\{(kl),(jl)(kl)\}\}(jk)(pk)(ik)(pk)^{-1}(jk)^{-1}(pk)(jk)(pj)\\
&\equiv&(jk)^{-1}\{\{(jl)^{\varepsilon},(kl)(pl)(kl)^{-1}(pl)^{-1}\},\{(kl),(pl)^{-1}\}\{\{(kl)^{-1},(pl)^{-1}\},\\
&&\{(jl),(kl)(pl)(kl)^{-1}(pl)^{-1}\}\{(kl),(pl)^{-1}\}\{\{(il),(kl)(pl)(kl)^{-1}(pl)^{-1}\},\{(pl)^{-1},(kl)^{-1}\\
&&(pl)^{-1}\}\}\}\{(kl)^{-1},(pl)^{-1}\}\{(jl)^{-1},(kl)(pl)(kl)^{-1}(pl)^{-1}\}\{(kl),(pl)^{-1}\}\\
&&\{(pl),(kl)^{-1}(pl)^{-1}\}\{(jl),(kl)(pl)(kl)^{-1}(pl)^{-1}\}\{(kl),(pl)^{-1}\}\}\\
&&(pk)^{-1}(jk)(pk)(ik)(pk)^{-1}(jk)^{-1}(pk)(jk)(pj)\\
&\equiv&\{\{(jl)^{\varepsilon},(kl)^{-1}(jl)^{-1}\},\{(kl),(jl)^{-1}\}\{\{(kl),(jl)^{-1}\},(pl)^{-1}\\
&&\{(kl)^{-1},(jl)^{-1}\}\{(jl),(kl)^{-1}(jl)^{-1}\}\{(kl),(jl)^{-1}\}(pl)(il)(pl)^{-1}\}\\
&&\{(kl)^{-1},(jl)^{-1}\}\{(jl)^{-1},(kl)^{-1}(jl)^{-1}\}\{(kl),(jl)^{-1}\}(pl)\{(jl),(kl)^{-1}(jl)^{-1}\}\\
&&\{(kl),(jl)^{-1}\}\{(ik),(pk)^{-1}(jk)^{-1}(pk)(jk)\}(pj)\\
&\equiv&\{\{(jl)^{\varepsilon},(pl)(jl)\},\{(il)^{-1},(pl)^{-1}(jl)^{-1}(pl)(jl)\}(kl)^{-1}\{(il),(pl)^{-1}\\
&&(jl)^{-1}(pl)(jl)\}(kl)\}\{(ik),(pk)^{-1}(jk)^{-1}(pk)(jk)\}(pj).
\end{eqnarray*}
\item[3)]\  $q>j,\ p<i<j<q<k<l$.
\begin{eqnarray*}
&&(pq)\{(jl)^{\varepsilon},(il)^{-1}(kl)^{-1}(il)(kl)\}(ik)\\
&\equiv&\{\{(jl)^{\varepsilon},(pl)^{-1}(ql)^{-1}(pl)(ql)\},\{(il)^{-1},(pl)^{-1}(ql)^{-1}(pl)(ql)\}\\
&&(kl)^{-1}\{(il),(pl)^{-1}(ql)^{-1}(pl)(ql)\}(kl)\}\{(ik),(pk)^{-1}(qk)^{-1}(pk)(qk)\}(pq)\ \ \ \ \ and \\
&&\\
&&\{(ik),(pk)^{-1}(qk)^{-1}(pk)(qk)\}(pq)(jl)^{\varepsilon}\\
&\equiv&(qk)^{-1}(pk)^{-1}(qk)(pk)(ik)\{(jl)^{\varepsilon},(pl)^{-1}\{(pl),(kl)^{-1}(ql)(kl)\}\}
(pk)^{-1}(qk)^{-1}(pk)(qk)(pq)\\
&\equiv&(qk)^{-1}(pk)^{-1}(qk)(pk)\{\{(jl)^{\varepsilon},(il)^{-1}(kl)^{-1}(il)(kl)\},(pl)^{-1}\{(pl),\{(kl)^{-1},(il)\\
&&(kl)\}\{(ql),(il)^{-1}(kl)^{-1}(il)(kl)\}\{(kl),(il)(kl)\}\}\}(ik)(pk)^{-1}(qk)^{-1}(pk)(qk)(pq)\\
&\equiv&(qk)^{-1}(pk)^{-1}(qk)\{\{(jl)^{\varepsilon},(pl)^{-1}(kl)^{-1}(pl)(kl)\},\{\{(kl)^{-1},(pl)(kl)\},\{(il),(pl)^{-1}\\
&&(kl)^{-1}(pl)(kl)\}\}\{(kl),(pl)(kl)\}\{(pl)^{-1},(kl)\}\{\{(pl),(kl)\},\{(kl)^{-1},(pl)(kl)\}\\
&&\{(ql),(pl)^{-1}(kl)^{-1}(pl)(kl)\}\{(kl),(pl)(kl)\}\}\}(pk)(ik)(pk)^{-1}(qk)^{-1}(pk)(qk)(pq)\\
&\equiv&(qk)^{-1}(pk)^{-1}\{(jl)^{\varepsilon},\{\{(kl)^{-1},(ql)(kl)\},(pl)(il)\}(pl)^{-1}\{(pl),\{(ql),(kl)\}\}\\
&&\{(kl),(ql)(kl)\}\}(qk)(pk)(ik)(pk)^{-1}(qk)^{-1}(pk)(qk)(pq)\\
&\equiv&(qk)^{-1}\{\{(jl)^{\varepsilon},(kl)(pl)(kl)^{-1}(pl)^{-1}\},\{\{(kl)^{-1},(pl)^{-1}\},\{(ql),(kl)(pl)\\
&&(kl)^{-1}(pl)^{-1}\}\{(kl),(pl)^{-1}\}\{(pl),(kl)^{-1}(pl)^{-1}\}\{(il),(kl)(pl)(kl)^{-1}(pl)^{-1}\}\}\\
&&\{\{(ql)^{-1},(kl)(pl)(kl)^{-1}(pl)^{-1}\},\{(kl),(pl)^{-1}\}\{(pl),(kl)^{-1}(pl)^{-1}\}\}\{(ql),(kl)(pl)\\
&&(kl)^{-1}(pl)^{-1}\}\{(kl),(pl)^{-1}\}\}(pk)^{-1}(qk)(pk)(ik)(pk)^{-1}(qk)^{-1}(pk)(qk)(pq)\\
&\equiv&\{(jl)^{\varepsilon},\{\{(kl)^{-1},(ql)^{-1}\},(pl)^{-1}\{(kl)^{-1},(ql)^{-1}\}\{(ql),(kl)^{-1}(ql)^{-1}\}\\
&&\{(kl),(ql)^{-1}\}(pl)(il)\}\{\{(ql)^{-1},(kl)^{-1}(ql)^{-1}\}\{(kl),(ql)^{-1}\}(pl)\{(ql),(kl)^{-1}(ql)^{-1}\}\\
&&\{(kl),(ql)^{-1}\}\}\{(ik),(pk)^{-1}(qk)^{-1}(pk)(qk)\}(pq)\\
&\equiv&\{\{(jl)^{\varepsilon},(pl)^{-1}(ql)^{-1}(pl)(ql)\},\{(il)^{-1},(pl)^{-1}(ql)^{-1}(pl)(ql)\}\\
&&(kl)^{-1}\{(il),(pl)^{-1}(ql)^{-1}(pl)(ql)\}(kl)\}\{(ik),(pk)^{-1}(qk)^{-1}(pk)(qk)\}(pq).
\end{eqnarray*}
\end{enumerate}

$(12)\wedge(13)$.\\
Let
$f=(pq)(ik)^{\delta}-\{(ik)^{\delta},(pk)^{-1}(qk)^{-1}(pk)(qk)\}(pq),
g=(ik)^{\delta}(jl)^{\varepsilon}-(jl)^{\varepsilon}(ik)^{\delta}$.\
Then $w=(pq)(ik)^{\delta}(jl)^{\varepsilon}$\  and
$$
(f,g)_{w}=(pq)(jl)^{\varepsilon}(ik)^{\delta}-\{(ik)^{\delta},(pk)^{-1}(qk)^{-1}(pk)(qk)\}(pq)(jl)^{\varepsilon}.
$$
There are two cases to consider.
\begin{enumerate}
\item[1)] \ $p<i<q<k$, $j<i<k<l$.\ In this case, there are three subcases
to consider.
\begin{enumerate}
\item[a)] \ $p<j,p<j<i<q<k<l$.
\begin{eqnarray*}
&&(pq)(jl)^{\varepsilon}(ik)^{\delta}\\
&\equiv&\{(jl)^{\varepsilon},(pl)^{-1}(ql)^{-1}(pl)(ql)\}\{(ik)^{\delta},(pk)^{-1}(qk)^{-1}(pk)(qk)\}(pq)\ \ \ \ \ \ \ \ and \\
&&\\
&&\{(ik)^{\delta},(pk)^{-1}(qk)^{-1}(pk)(qk)\}(pq)(jl)^{\varepsilon}\\
&\equiv&(qk)^{-1}(pk)^{-1}(qk)(pk)(ik)^{\delta}\{(jl)^{\varepsilon},(pl)^{-1}\{(pl),(kl)^{-1}(ql)(kl)\}\}\\
&&(pk)^{-1}(qk)^{-1}(pk)(qk)(pq).
\end{eqnarray*}
If $\delta =1$, then
\begin{eqnarray*}
&&\{(ik)^{\delta},(pk)^{-1}(qk)^{-1}(pk)(qk)\}(pq)(jl)^{\varepsilon}\\
&\equiv&(qk)^{-1}(pk)^{-1}(qk)(pk)\{(jl)^{\varepsilon},(pl)^{-1}\{(pl),\{(kl)^{-1},(il)(kl)\}\{(ql),(il)^{-1}(kl)^{-1}\\
&&(il)(kl)\}\{(kl),(il)(kl)\}\}\}(ik)^{\delta}(pk)^{-1}(qk)^{-1}(pk)(qk)(pq)\\
&\equiv&(qk)^{-1}(pk)^{-1}(qk)(pk)\{(jl)^{\varepsilon},(pl)^{-1}\{(pl),(kl)^{-1}(ql)(kl)\}\}\\
&&\{(ik)^{\delta}(pk)^{-1}(qk)^{-1}(pk)(qk)(pq)\\
&\equiv&(qk)^{-1}(pk)^{-1}(qk)\{\{(jl)^{\varepsilon},(pl)^{-1}(kl)^{-1}(pl)(kl)\},\{(pl)^{-1},(kl)\}\{\{(pl),(kl)\},\\
&&\{(kl)^{-1},(pl)(kl)\}\{(ql),(pl)^{-1}(kl)^{-1}(pl)(kl)\}\{(kl),(pl)(kl)\}\}\}\\
&&(pk)(ik)^{}(pk)^{-1}(qk)^{-1}(pk)(qk)(pq)\\
&\equiv&(qk)^{-1}(pk)^{-1}\{(jl)^{\varepsilon},(pl)^{-1}\{(pl),\{(ql),(kl)\}\{(kl),(ql)(kl)\}\}\}\\
&&(qk)(pk)(ik)^{\delta}(pk)^{-1}(qk)^{-1}(pk)(qk)(pq)\\
&\equiv&(qk)^{-1}\{\{jl^{\varepsilon},(kl)(pl)(kl)^{-1}(pl)^{-1}\},\{(pl)^{-1},(kl)^{-1}(pl)^{-1}\}\\
&&\{\{(pl),(kl)^{-1}(pl)^{-1}\},\{(ql),(kl)(pl)(kl)^{-1}(pl)^{-1}\}\{(kl),(pl)^{-1}\}\}\}\\
&&(pk)^{-1}(qk)(pk)(ik)^{}(pk)^{-1}(qk)^{-1}(pk)(qk)(pq)\\
&\equiv&\{(jl)^{\varepsilon},(pl)^{-1}\{(pl),\{(kl)^{-1},(ql)^{-1}\}\{(ql),(kl)^{-1}(ql)^{-1}\}\{(kl),(ql)^{-1}\}\}\}\\
&&\{(ik)^{\delta},(pk)^{-1}(qk)^{-1}(pk)(qk)\}(pq)\\
&\equiv&\{(jl)^{\varepsilon},(pl)^{-1}(ql)^{-1}(pl)(ql)\}\{(ik)^{\delta},(pk)^{-1}(qk)^{-1}(pk)(qk)\}(pq).
\end{eqnarray*}
If $\delta =-1$, then by using the result of the case of $\delta
=1$, we have
\begin{eqnarray*}
&&\{(ik)^{\delta},(pk)^{-1}(qk)^{-1}(pk)(qk)\}(pq)(jl)^{\varepsilon}\\
&\equiv&(qk)^{-1}(pk)^{-1}(qk)(pk)\{(jl)^{\varepsilon},(pl)^{-1}\{(pl),\{(kl)^{-1},(il)^{-1}\}\{(ql),(kl)(il)\\
&&(kl)^{-1}(il)^{-1}\}\{(kl),(il)^{-1}\}\}\}(ik)^{\delta}(pk)^{-1}(qk)^{-1}(pk)(qk)(pq)\\
&\equiv&(qk)^{-1}(pk)^{-1}(qk)(pk)\{(jl)^{\varepsilon},(pl)^{-1}\{(pl),(kl)^{-1}(ql)(kl)\}\}\\
&&\{(ik)^{\delta},(pk)^{-1}(qk)^{-1}(pk)(qk)\}(pq)\\
&\equiv&\{(jl)^{\varepsilon},(pl)^{-1}(ql)^{-1}(pl)(ql)\}
\{(ik)^{\delta},(pk)^{-1}(qk)^{-1}(pk)(qk)\}(pq).
\end{eqnarray*}
\item[b)] \ $p=j,p=j<i<q<k<l$.
\begin{eqnarray*}
&&(pq)(jl)^{\varepsilon}(ik)^{\delta}=(jq)(jl)^{-1}(ik)^{\delta}\\
&\equiv&\{(jl)^{\varepsilon},(ql)\}\{(ik)^{\delta},(jk)^{-1}(qk)^{-1}(jk)(qk)\}(jq)\ \ \ \ \ \ \ and \\
&&\\
&&\{(ik)^{\delta},(pk)^{-1}(qk)^{-1}(pk)(qk)\}(pq)(jl)^{\varepsilon}\\
&=&\{(ik)^{\delta},(jk)^{-1}(qk)^{-1}(jk)(qk)\}(jq)(jl)^{\varepsilon}\\
&\equiv&(qk)^{-1}(jk)^{-1}(qk)(jk)(ik)^{\delta}(jk)^{-1}(qk)^{-1}(jk)(qk)\{(jl)^{\varepsilon},(ql)\}(jq)\\
&\equiv&(qk)^{-1}(jk)^{-1}(qk)(jk)(ik)^{\delta}(jk)^{-1}(qk)^{-1}(jk)\{(jl)^{\varepsilon},\{(ql),(kl)\}\}(qk)(jq)\\
&\equiv&(qk)^{-1}(jk)^{-1}(qk)(jk)(ik)^{\delta}(jk)^{-1}(qk)^{-1}\{\{(jl)^{\varepsilon},(kl)\},\\
&&\{\{(ql),(jl)^{-1}(kl)^{-1}(jl)(kl)\},\{(kl),(jl)(kl)\}\}\}(jk)(qk)(jq)\\
&\equiv&(qk)^{-1}(jk)^{-1}(qk)(jk)(ik)^{\delta}(jk)^{-1}\{(jl)^{\varepsilon},\{(ql),(kl)^{-1}(ql)^{-1}\}\\
&&\{(kl),(ql)^{-1}\}\}(qk)^{-1}(jk)(qk)(jq)\\
&\equiv&(qk)^{-1}(jk)^{-1}(qk)(jk)(ik)^{\delta}\{\{(jl)^{\varepsilon},(kl)^{-1}(jl)^{-1}\},\{(ql),(kl)(jl)\\
&&(kl)^{-1}(jl)^{-1}\}\{(kl),(jl)^{-1}\}\}(jk)^{-1}(qk)^{-1}(jk)(qk)(jq)\\
&\equiv&(qk)^{-1}(jk)^{-1}(qk)(jk)(ik)^{\delta}\{(jl)^{\varepsilon},(kl)^{-1}(ql)(kl)\}(jk)^{-1}(qk)^{-1}(jk)(qk)(jq).
\end{eqnarray*}
If $\delta =1$, then
\begin{eqnarray*}
&&\{(ik)^{\delta},(pk)^{-1}(qk)^{-1}(pk)(qk)\}(pq)(jl)^{\varepsilon}\\
&\equiv&(qk)^{-1}(jk)^{-1}(qk)(jk)\{(jl)^{\varepsilon},\{(kl)^{-1},(il)(kl)\}\{(ql),(il)^{-1}(kl)^{-1}(il)(kl)\}\\
&&\{(kl),(il)(kl)\}\}(ik)^{\delta}(jk)^{-1}(qk)^{-1}(jk)(qk)(jq)\\
&\equiv&(qk)^{-1}(jk)^{-1}(qk)(jk)\{(jl)^{\varepsilon},(kl)^{-1}(ql)(kl)\}(ik)^{\delta}jk)^{-1}(qk)^{-1}(jk)(qk)(jq)\\
&\equiv&(qk)^{-1}(jk)^{-1}(qk)\{\{(jl)^{\varepsilon},(kl)\},\{(kl)^{-1},(jl)(kl)\}\{(ql),(jl)^{-1}(kl)^{-1}(jl)(kl)\}\\
&&\{(kl),(jl)(kl)\}\}(jk)(ik)^{\delta}(jk)^{-1}(qk)^{-1}(jk)(qk)(jq)\\
&\equiv&(qk)^{-1}(jk)^{-1}\{(jl)^{\varepsilon},\{(ql),(kl)\}\{(kl),(ql)(kl)\}\}\\
&&(qk)(jk)(ik)^{\delta}(jk)^{-1}(qk)^{-1}(jk)(qk)(jq)\\
&\equiv&(qk)^{-1}\{\{(jl)^{\varepsilon},(kl)^{-1}(jl)^{-1}\},\{(ql),(kl)(jl)(kl)^{-1}(jl)^{-1}\}\{(kl),(jl)^{-1}\}\}\\
&&(jk)^{-1}(qk)(jk)(ik)^{\delta}(jk)^{-1}(qk)^{-1}(jk)(qk)(jq)\\
&\equiv&\{(jl)^{\varepsilon},\{(kl)^{-1},(ql)^{-1}\}\{(ql),(kl)^{-1}(ql)^{-1}\}\{(kl),(ql)^{-1}\}\}\\
&&\{(ik)^{\delta},(jk)^{-1}(qk)^{-1}(jk)(qk)\}(jq)\\
&\equiv&\{(jl)^{\varepsilon},(ql)\}\{(ik)^{\delta},(jk)^{-1}(qk)^{-1}(jk)(qk)\}(jq).
\end{eqnarray*}
If $ \delta =-1$, then by using the result of the case of $\delta
=1$, we have
\begin{eqnarray*}
&&\{(ik)^{\delta},(pk)^{-1}(qk)^{-1}(pk)(qk)\}(pq)(jl)^{\varepsilon}\\
&\equiv&(qk)^{-1}(jk)^{-1}(qk)(jk)\{(jl)^{\varepsilon},\{(kl)^{-1},(il)^{-1}\}\{(ql),(kl)(il)(kl)^{-1}(il)^{-1}\}\\
&&\{(kl),(il)^{-1}\}\}(ik)^{\delta}(jk)^{-1}(qk)^{-1}(jk)(qk)(jq)\\
&\equiv&(qk)^{-1}(jk)^{-1}(qk)(jk)\{(jl)^{\varepsilon},(kl)^{-1}(ql)(kl)\}(ik)^{\delta}(jk)^{-1}(qk)^{-1}(jk)(qk)(jq)\\
&\equiv&\{(jl)^{\varepsilon},(ql)\}\{(ik)^{\delta},(jl)^{-1}(qk)^{-1}(jk)(qk)\}(jq).
\end{eqnarray*}
\item[c)] \ $p>j,j<p<i<q<k<l$.
\begin{eqnarray*}
(f,g)_{w}&\equiv&(jl)^{\varepsilon}\{(ik)^{\delta},(pk)^{-1}(qk)^{-1}(pk)(qk)\}(pq)\\
&&-(jl)^{\varepsilon}\{(ik)^{\delta},(pk)^{-1}(qk)^{-1}(pk)(qk)\}(pq)\\
&\equiv&0.
\end{eqnarray*}
\end{enumerate}
\item[2)] \ $i<k<j<l$.

In this case, we have $p<i<q<k<j<l$ and
\begin{eqnarray*}
(f,g)_{w}&\equiv&(jl)^{\varepsilon}\{(ik)^{\delta},(pk)^{-1}(qk)^{-1}(pk)(qk)\}(pq)\\
&&-(jl)^{\varepsilon}\{(ik)^{\delta},(pk)^{-1}(qk)^{-1}(pk)(qk)\}(pq)\\
&\equiv&0.
\end{eqnarray*}
\end{enumerate}

$(13)\wedge(7)$\\
Let $f=(ip)^\delta(jk)^{-1}-(jk)^{-1}(ip)^\delta,\
g=(jk)^{-1}(kl)^{\varepsilon}-\{(kl)^{\varepsilon},(jl)^{-1}\}(jk)^{-1},\
j<i<p<k<l\ or\  i<p<j<k<l$. Then
$w=(ip)^\delta(jk)^{-1}(kl)^{\varepsilon}$ and
\begin{eqnarray*}
(f,g)_w&=&(ip)^\delta\{(kl)^{\varepsilon},(jl)^{-1}\}(jk)^{-1}-(jk)^{-1}(ip)^\delta(kl)^{\varepsilon}\\
&\equiv&\{(kl)^{\varepsilon},(jl)^{-1}\}(jk)^{-1}(ip)^\delta-(jk)^{-1}(kl)^{\varepsilon}(ip)^\delta\\
&\equiv&0.
\end{eqnarray*}

$(13)\wedge(8)$\\
Let $f=(ip)^\delta(jk)-(jk)(ip)^\delta,\
g=(jk)(kl)^{\varepsilon}-\{(kl)^{\varepsilon},(jl)(kl)\}(jk), \
j<i<p<k<l\ or\ i<p<j<k<l$. Then
$w=(ip)^\delta(jk)(kl)^{\varepsilon}$ and
\begin{eqnarray*}
(f,g)_w&=&(ip)^\delta\{(kl)^{\varepsilon},(jl)(kl)\}(jk)-(jk)(ip)^\delta(kl)^{\varepsilon}\\
&\equiv&\{(kl)^{\varepsilon},(jl)(kl)\}(jk)(ip)^\delta-(jk)(kl)^{\varepsilon}(ip)^\delta\\
&\equiv&0.
\end{eqnarray*}

$(13)\wedge(9)$\\
Let $f=(ip)^\delta(jk)^{-1}-(jk)^{-1}(ip)^\delta,\
g=(jk)^{-1}(jl)^{\varepsilon}-\{(jl)^{\varepsilon},(kl)^{-1}(jl)^{-1}\}(jk)^{-1},\
j<i<p<k<l\ or\ i<p<j<k<l$. Then
$w=(ip)^\delta(jk)^{-1}(jl)^{\varepsilon}$ and
\begin{eqnarray*}
(f,g)_w&=&(ip)^\delta\{(jl)^{\varepsilon},(kl)^{-1}(jl)^{-1}\}(jk)^{-1}-(jk)^{-1}(ip)^\delta(jl)^{\varepsilon}\\
&\equiv&\{(jl)^{\varepsilon},(kl)^{-1}(jl)^{-1}\}(jk)^{-1}(ip)^\delta-(jk)^{-1}(jl)^{\varepsilon}(ip)^\delta\\
&\equiv&0.
\end{eqnarray*}

$(13)\wedge(10)$\\
Let $f=(ip)^\delta(jk)-(jk)(ip)^\delta,\
g=(jk)(jl)^{\varepsilon}-\{(jl)^{\varepsilon},(kl)\}(jk),\
j<i<p<k<l\ or\ i<p<j<k<l$. Then
$w=(ip)^\delta(jk)(jl)^{\varepsilon}$ and
\begin{eqnarray*}
(f,g)_w&=&(ip)^\delta\{(jl)^{\varepsilon},(kl)\}(jk)-(jk)(ip)^\delta(jl)^{\varepsilon}\\
&\equiv&\{(jl)^{\varepsilon},(kl)\}(jk)(ip)^\delta-(jk)(jl)^{\varepsilon}(ip)^\delta\\
&\equiv&0.
\end{eqnarray*}

$(13)\wedge(11)$\\
Let $f=(pq)^\delta(ik)^{-1}-(ik)^{-1}(pq)^\delta,\
g=(ik)^{-1}(jl)^{\varepsilon}-\{(jl)^{\varepsilon},(kl)(il)(kl)^{-1}(il)^{-1}\}(ik)^{-1}$.
Then $w=(pq)^\delta(ik)^{-1}(jl)^{\varepsilon}$ and
$$
(f,g)_w=(pq)^\delta\{(jl)^{\varepsilon},(kl)(il)(kl)^{-1}(il)^{-1}\}(ik)^{-1}-(ik)^{-1}(pq)^\delta(jl)^{\varepsilon}.
$$
There are two cases to consider.
\begin{enumerate}
\item[1)]\  $i<p<q<k,i<j<k<l$. In this case, there are five subcases
to consider.
\begin{enumerate}
\item[a)]\  $j<p,i<j<p<q<k<l$.
$$
(f,g)_w\equiv\{(jl)^{\varepsilon},(kl)(il)(kl)^{-1}(il)^{-1}\}(ik)^{-1}(pq)^\delta-
(ik)^{-1}(jl)^{\varepsilon}(pq)^\delta\equiv0.
$$
\item[b)]\  $j=p,i<j=p<q<k<l$.
$$
(f,g)_w=(jq)^\delta\{(jl)^{\varepsilon},(kl)(il)(kl)^{-1}(il)^{-1}\}(ik)^{-1}-(ik)^{-1}(jq)^\delta(jl)^{\varepsilon}.
$$
If $\delta=1$, then
\begin{eqnarray*}
(f,g)_w&\equiv&\{\{(jl)^{\varepsilon},(ql)\},(kl)(il)(kl)^{-1}(il)^{-1}\}(ik)^{-1}(jq)^\delta\\
&&-(ik)^{-1}\{(jl)^{\varepsilon},(ql)\}(jq)^\delta\\
&\equiv&\{(jl)^{\varepsilon},(ql)(kl)(il)(kl)^{-1}(il)^{-1}\}(jq)^\delta(ik)^{-1}\\
&&-\{\{(jl)^{\varepsilon},(kl)(il)(kl)^{-1}(il)^{-1}\},\{(ql),(kl)(il)(kl)^{-1}(il)^{-1}\}\}(ik)^{-1}(jq)^\delta\\
&\equiv&0.
\end{eqnarray*}
If $\delta=-1$, then
\begin{eqnarray*}
(f,g)_w&\equiv&\{\{(jl)^{\varepsilon},(ql)^{-1}(jl)^{-1}\},(kl)(il)(kl)^{-1}(il)^{-1}\}(ik)^{-1}(jq)^\delta\\
&&-(ik)^{-1}\{(jl)^{\varepsilon},(ql)^{-1}(jl)^{-1}\}(jq)^\delta\\
&\equiv&\{(jl)^{\varepsilon},(ql)^{-1}(jl)^{-1}(kl)(il)(kl)^{-1}(il)^{-1}\}(ik)^{-1}(jq)^\delta\\
&&-\{\{(jl)^{\varepsilon},(kl)(il)(kl)^{-1}(il)^{-1}\},\{(ql)^{-1},(kl)(il)(kl)^{-1}(il)^{-1}\}\\
&&\{(jl)^{-1},(kl)(il)(kl)^{-1}(il)^{-1}\}\}(ik)^{-1}(jq)^\delta\\
&\equiv&0.
\end{eqnarray*}
\item[c)]\  $p<j<q,i<p<j<q<k<l$.\\
If $\delta=1$, then
\begin{eqnarray*}
(f,g)_w&\equiv&\{\{(jl)^{\varepsilon},(pl)^{-1}(ql)^{-1}(pl)(ql)\},(kl)(il)(kl)^{-1}(il)^{-1}\}(ik)^{-1}(pq)^\delta\\
&&-(ik)^{-1}\{(jl)^{\varepsilon},(pl)^{-1}(ql)^{-1}(pl)(ql)\}(pq)^\delta\\
&\equiv&\{(jl)^{\varepsilon},(pl)^{-1}(ql)^{-1}(pl)(ql)(kl)(il)(kl)^{-1}(il)^{-1}\}(ik)^{-1}(pq)^\delta\\
&&-\{\{(jl)^{\varepsilon},(il)^{-1}(kl)^{-1}(il)(kl)\},\{\{(ql)^{-1},(il)^{-1}(kl)^{-1}(il)(kl)\},\\
&&\{(pl),(il)^{-1}(kl)^{-1}(il)(kl)\}\}\{(ql),(il)^{-1}(kl)^{-1}(il)(kl)\}\}(ik)^{-1}(pq)^\delta\\
&\equiv&0.
\end{eqnarray*}
If $\delta=-1$, then
\begin{eqnarray*}
(f,g)_w&\equiv&\{(jl)^{\varepsilon},(ql)(pl)(ql)^{-1}(pl)^{-1}(kl)(il)(kl)^{-1}(il)^{-1}\}(ik)^{-1}(pq)^\delta\\
&&-(ik)^{-1}\{(jl)^{\varepsilon},(ql)(pl)(ql)^{-1}(pl)^{-1}\}(pq)^\delta\\
&\equiv&\{(jl)^{\varepsilon},(ql)(pl)(ql)^{-1}(pl)^{-1}(kl)(il)(kl)^{-1}(il)^{-1}\}(ik)^{-1}(pq)^\delta-\\
&&\{\{(jl)^{\varepsilon},(kl)(il)(kl)^{-1}(il)^{-1}\},\{\{(pl),(kl)(il)(kl)^{-1}(il)^{-1}\},\\
&&\{(ql)^{-1},(kl)(il)(kl)^{-1}(il)^{-1}\}\}\{(pl)^{-1},(kl)(il)(kl)^{-1}(il)^{-1}\}\}(ik)^{-1}(pq)^\delta\\
&\equiv&0.
\end{eqnarray*}
\item[d)]\  $j=q,i<p<q=j<k<l$.
$$
(f,g)_w=(pj)^\delta\{(jl)^{\varepsilon},(kl)(il)(kl)^{-1}(il)^{-1}\}(ik)^{-1}-(ik)^{-1}(pj)^\delta(jl)^{\varepsilon}.
$$
If $\delta=1$, then
\begin{eqnarray*}
(f,g)_w&\equiv&\{\{(jl)^{\varepsilon},(pl)(jl)\},(kl)(il)(kl)^{-1}(il)^{-1}\}(ik)^{-1}(pj)^\delta\\
&&-(ik)^{-1}\{(jl)^{\varepsilon},(pl)(jl)\}(pj)^\delta\\
&\equiv&(jl)^{\varepsilon},(pl)(jl)(kl)(il)(kl)^{-1}(il)^{-1}\}(ik)^{-1}(pj)^\delta-\\
&&\{\{(jl)^{\varepsilon},(kl)(il)(kl)^{-1}(il)^{-1}\},\{(pl),(kl)(il)(kl)^{-1}(il)^{-1}\}\\
&&\{(jl),(kl)(il)(kl)^{-1}(il)^{-1}\}\}(ik)^{-1}(pj)^\delta\\
&\equiv&0.
\end{eqnarray*}
If $\delta=-1$, then
\begin{eqnarray*}
(f,g)_w&\equiv&\{\{(jl)^{\varepsilon},(pl)^{-1}\},(kl)(il)(kl)^{-1}(il)^{-1}\}(ik)^{-1}(pj)^\delta-\\
&&(ik)^{-1}\{(jl)^{\varepsilon},(pl)^{-1}\}(pj)^\delta\\
&\equiv&\{(jl)^{\varepsilon},(pl)^{-1}(kl)(il)(kl)^{-1}(il)^{-1}\}(ik)^{-1}(pj)^\delta-\{\{(jl)^{\varepsilon},(kl)(il)\\
&&(kl)^{-1}(il)^{-1}\},\{(pl)^{-1},(kl)(il)(kl)^{-1}(il)^{-1}\}\}(ik)^{-1}(pj)^\delta\\
&\equiv&0.
\end{eqnarray*}
\item[e)]\  $j>q,i<p<q<j<k<l$.
$$
(f,g)_w\equiv\{(jl)^{\varepsilon},(kl)(il)(kl)^{-1}(il)^{-1}\}(ik)^{-1}(pq)^\delta-
(ik)^{-1}(jl)^{\varepsilon}(pq)^\delta\equiv0.
$$
\end{enumerate}
\item[2)]\   $p<q<i<k,i<j<k<l$, that is, \ $p<q<i<j<k<l$.
$$
(f,g)_w\equiv\{(jl)^{\varepsilon},(kl)(il)(kl)^{-1}(il)^{-1}\}(ik)^{-1}(pq)^\delta-
(ik)^{-1}(jl)^{\varepsilon}(pq)^\delta\equiv0.
$$
\end{enumerate}

$(13)\wedge(12)$\\
Let $f=(pq)^\delta(ik)-(ik)(pq)^\delta,\
g=(ik)(jl)^{\varepsilon}-\{(jl)^{\varepsilon},(il)^{-1}(kl)^{-1}(il)(kl)\}(ik)$.
Then $w=(pq)^\delta(ik)(jl)^{\varepsilon}$ and
$$
(f,g)_w=(pq)^\delta\{(jl)^{\varepsilon},(il)^{-1}(kl)^{-1}(il)(kl)\}(ik)-(ik)(pq)^\delta(jl)^{\varepsilon}.
$$
There are two cases to consider.
\begin{enumerate}
\item[1)]\ $i<p<q<k,i<j<k<l$. In this case, there are five subcases
to consider.
\begin{enumerate}
\item[a)]\ $j<p,i<j<p<q<k<l$.
$$
(f,g)_w\equiv\{(jl)^{\varepsilon},(il)^{-1}(kl)^{-1}(il)(kl)\}(ik)(pq)^\delta-(ik)(jl)^{\varepsilon}(pq)^\delta\equiv0.
$$
\item[b)]\ $p=j,i<j=p<q<k<l$.
$$
(f,g)_w\equiv(jq)^\delta\{(jl)^{\varepsilon},(il)^{-1}(kl)^{-1}(il)(kl)\}(ik)-(ik)(jq)^\delta(jl)^{\varepsilon}.
$$
If $\delta=1$, then
\begin{eqnarray*}
(f,g)_w&\equiv&\{\{(jl)^{\varepsilon},(ql)\},(il)^{-1}(kl)^{-1}(il)(kl)\}(ik)(jq)^\delta
-(ik)\{(jl)^{\varepsilon},(ql)\}(jq)^\delta\\
&\equiv&\{(jl)^{\varepsilon},(ql)(il)^{-1}(kl)^{-1}(il)(kl)\}(ik)(jq)^\delta-\\
&&\{\{(jl)^{\varepsilon},(il)^{-1}(kl)^{-1}(il)(kl)\},\{(ql),(il)^{-1}(kl)^{-1}(il)(kl)\}\}(ik)(jq)^\delta\\
&\equiv&0.
\end{eqnarray*}
If $\delta=-1$, then
\begin{eqnarray*}
(f,g)_w&\equiv&\{\{(jl)^{\varepsilon},(ql)^{-1}(jl)^{-1}\},(il)^{-1}(kl)^{-1}(il)(kl)\}(ik)(jq)^\delta\\
&&-(ik)\{(jl)^{\varepsilon},(ql)^{-1}(jl)^{-1}\}(jq)^\delta\\
&\equiv&\{(jl)^{\varepsilon},(ql)^{-1}(jl)^{-1}(il)^{-1}(kl)^{-1}(il)(kl)\}(ik)(jq)^\delta\\
&&-\{\{(jl)^{\varepsilon},(il)^{-1}(kl)^{-1}(il)(kl)\},\{(ql)^{-1},(il)^{-1}(kl)^{-1}(il)(kl)\}\\
&&\{(jl)^{-1},(il)^{-1}(kl)^{-1}(il)(kl)\}\}(ik)(jq)^\delta\\
&\equiv&0.
\end{eqnarray*}
\item[c)]\ $p<j<q,i<p<j<q<k<l$.\\
If $\delta=1$, then
\begin{eqnarray*}
(f,g)_w&\equiv&\{\{(jl)^{\varepsilon},(pl)^{-1}(ql)^{-1}(pl)(ql)\},(il)^{-1}(kl)^{-1}(il)(kl)\}(ik)(pq)^\delta\\
&&-(ik)\{(jl)^{\varepsilon},(pl)^{-1}(ql)^{-1}(pl)(ql)\}(pq)^\delta\\
&\equiv&\{(jl)^{\varepsilon},(pl)^{-1}(ql)^{-1}(pl)(ql)(il)^{-1}(kl)^{-1}(il)(kl)\}(ik)(pq)^\delta-\\
&&\{\{(jl)^{\varepsilon},(il)^{-1}(kl)^{-1}(il)(kl)\},\{(pl)^{-1},(il)^{-1}(kl)^{-1}(il)(kl)\}\\
&&\{\{(pl),(il)^{-1}(kl)^{-1}(il)(kl)\},\{(ql),(il)^{-1}(kl)^{-1}(il)(kl)\}\}\}(ik)(pq)^\delta\\
&\equiv&0.
\end{eqnarray*}
If $ \delta=-1$, then
\begin{eqnarray*}
(f,g)_w&\equiv&\{\{(jl)^{\varepsilon},(ql)(pl)(ql)^{-1}(pl)^{-1}\},(il)^{-1}(kl)^{-1}(il)(kl)\}(ik)(pq)^\delta\\
&&-(ik)\{(jl)^{\varepsilon},(ql)(pl)(ql)^{-1}(pl)^{-1}\}(pq)^\delta\\
&\equiv&\{(jl)^{\varepsilon},(ql)(pl)(ql)^{-1}(pl)^{-1}(il)^{-1}(kl)^{-1}(il)(kl)\}(ik)(pq)^\delta\\
&&-\{\{(jl)^{\varepsilon},(il)^{-1}(kl)^{-1}(il)(kl)\},\{(ql)^{-1},(il)^{-1}(kl)^{-1}(il)(kl)\}\\
&&\{\{(ql)^{-1},(il)^{-1}(kl)^{-1}(il)(kl)\},\{(pl)^{-1},(il)^{-1}(kl)^{-1}(il)(kl)\}\}\}(ik)(pq)^\delta\\
&\equiv&0.
\end{eqnarray*}
\item[d)]\ $j=q,i<p<q=j<k<l$.
$$
(f,g)_w=(pj)^\delta\{(jl)^{\varepsilon},(il)^{-1}(kl)^{-1}(il)(kl)\}(ik)-(ik)(pj)^\delta(jl)^{\varepsilon}.
$$
If $\delta=1$, then
\begin{eqnarray*}
(f,g)_w&\equiv&\{\{(jl)^{\varepsilon},(pl)\},(il)^{-1}(kl)^{-1}(il)(kl)\}(ik)(pj)^\delta\\
&&-(ik)\{(jl)^{\varepsilon},(pl)\}(pj)^\delta\\
&\equiv&\{(jl)^{\varepsilon},(pl)(il)^{-1}(kl)^{-1}(il)(kl)\}(ik)(pj)^\delta-\\
&&\{\{(jl)^{\varepsilon},(il)^{-1}(kl)^{-1}(il)(kl)\}\{(pl),(il)^{-1}(kl)^{-1}(il)(kl)\}\}(kl)(pj)^\delta\\
&\equiv&0.
\end{eqnarray*}
If $\delta=-1$, then
\begin{eqnarray*}
(f,g)_w&\equiv&\{(jl)^{\varepsilon},(pl)^{-1}(jl)^{-1}(il)^{-1}(kl)^{-1}(il)(kl)\}(ik)(pj)^\delta\\
&&-(ik)\{(jl)^{\varepsilon},(pl)^{-1}(jl)^{-1}\}(pj)^\delta\\
&\equiv&\{(jl)^{\varepsilon},(pl)^{-1}(jl)^{-1}(il)^{-1}(kl)^{-1}(il)(kl)\}(ik)(pj)^\delta\\
&&-\{\{(jl)^{\varepsilon},(il)^{-1}(kl)^{-1}(il)(kl)\},\{(pl)^{-1},(il)^{-1}(kl)^{-1}(il)(kl)\}\\
&&\{(jl)^{-1},(il)^{-1}(kl)^{-1}(il)(kl)\}\}(ik)(pj)^\delta\\
&\equiv&0.
\end{eqnarray*}
\item[e)]\ $q>j,i<p<q<j<k<l$.
$$
(f,g)_w\equiv\{(jl)^{\varepsilon},(il)^{-1}(kl)^{-1}(il)(kl)\}(ik)(pq)^\delta-(ik)(jl)^{\varepsilon}(pq)^\delta\equiv0.
$$
\end{enumerate}
\item[2)]\ $p<q<i<k, \ i<j<k<l$,\ that is, \ $p<q<i<j<k<l$.
$$
(f,g)_w\equiv\{(jl)^{\varepsilon},(il)^{-1}(kl)^{-1}(il)(kl)\}(ik)(pq)^\delta-(ik)(jl)^{\varepsilon}(pq)^\delta\equiv0.
$$
\end{enumerate}

$(13)\wedge(13)$\\
Let
$f=(pq)^\gamma(ik)^\delta-(ik)^\delta(pq)^\gamma,g=(ik)^\delta(jl)^{\varepsilon}-(jl)^{\varepsilon}(ik)^\delta$.\
Then $w=(pq)^\gamma(ik)^\delta(jl)^{\varepsilon}$ and
$$
(f,g)_w=(pq)^\gamma(jl)^{\varepsilon}(ik)^\delta-(ik)^\delta(pq)^\gamma(jl)^{\varepsilon}.
$$
There are four cases to consider.
\begin{enumerate}
\item[1)]\ $i<p<q<k,j<i<k<l$,\ that  is, \ $j<i<p<q<k<l$.
\begin{eqnarray*}
(f,g)_w&\equiv&(jl)^{\varepsilon}(ik)^\delta(pq)^\gamma-(ik)^\delta(jl)^{\varepsilon}(pq)^\gamma\\
&\equiv&(jl)^{\varepsilon}(ik)^\delta(pq)^\gamma-(jl)^{\varepsilon}(ik)^\delta(pq)^\gamma\\
&\equiv&0.
\end{eqnarray*}
\item[2)]\ $i<p<q<k,i<k<j<l$,\ that is, \ $i<p<q<k<j<l$.\\
$$
(f,g)_w\equiv(jl)^{\varepsilon}(ik)^\delta(pq)^\gamma-(jl)^{\varepsilon}(ik)^\delta(jl)^{\varepsilon}\equiv0.
$$
\item[3)]\ $p<q<i<k,j<i<k<l$.\ In this case,\ there are five subcases to
consider.
\begin{enumerate}
\item[a)]\ $j<p$,\ that is, \ $j<p<q<i<k<l$.\\
$$
(f,g)_w\equiv(jl)^{\varepsilon}(ik)^\delta(pq)^\gamma-(jl)^{\varepsilon}(ik)^\delta(pq)^\gamma\equiv0.
$$
\item[b)]\ $j=p$,\ that is, \ $j=p<q<i<k<l$.
\begin{eqnarray*}
(f,g)_w&=&(jq)^\gamma(jl)^{\varepsilon}(ik)^\delta-(ik)^\delta(jq)^\gamma(jl)^{\varepsilon}.
\end{eqnarray*}
If $\gamma=1$, then
\begin{eqnarray*}
(f,g)_w&\equiv&\{(jl)^{\varepsilon},(ql)\}(ik)^\delta(jq)^\gamma-(ik)^\delta\{(jl)^{\varepsilon},(ql)\}(jq)^\gamma\\
&\equiv&\{(jl)^{\varepsilon},(ql)\}(ik)^\delta(jq)^\gamma-\{(jl)^{\varepsilon},(ql)\}(ik)^\delta(jq)^\gamma\\
&\equiv&0.
\end{eqnarray*}
If $\gamma=-1$, then
\begin{eqnarray*}
(f,g)_w&\equiv&\{(jl)^{\varepsilon},(ql)^{-1}(jl)^{-1}\}(ik)^\delta(jq)^\gamma
-(ik)^\delta\{(jl)^{\varepsilon},(ql)^{-1}(jl)^{-1}\}(jq)^\gamma\\
&\equiv&\{(jl)^{\varepsilon},(ql)^{-1}(jl)^{-1}\}(ik)^\delta(jq)^\gamma
-\{(jl)^{\varepsilon},(ql)^{-1}(jl)^{-1}\}(ik)^\delta(jq)^\gamma\\
&\equiv&0.
\end{eqnarray*}
\item[c)]\ $p<j<q,p<j<q<i<k<l$.\\
If $\gamma=1$, then
\begin{eqnarray*}
(f,g)_w&\equiv&\{(jl)^{\varepsilon},(pl)^{-1}(ql)^{-1}(pl)(ql)\}(ik)^\delta(pq)^\gamma\\
&&-(ik)^\delta\{(jl)^{\varepsilon},(pl)^{-1}(ql)^{-1}(pl)(ql)\}(pq)^\gamma\\
&\equiv&\{(jl)^{\varepsilon},(pl)^{-1}(ql)^{-1}(pl)(ql)\}(ik)^\delta(pq)^\gamma\\
&&-\{(jl)^{\varepsilon},(pl)^{-1}(ql)^{-1}(pl)(ql)\}(ik)^\delta(pq)^\gamma\\
&\equiv&0.
\end{eqnarray*}
If$ \gamma=-1$, then
\begin{eqnarray*}
(f,g)_w&\equiv&\{(jl)^{\varepsilon},(ql)(pl)(ql)^{-1}(pl)^{-1}\}(ik)^\delta(pq)^\gamma\\
&&-(ik)^\delta\{(jl)^{\varepsilon},(ql)(pl)(ql)^{-1}(pl)^{-1}\}(pq)^\gamma\\
&\equiv&\{(jl)^{\varepsilon},(ql)(pl)(ql)^{-1}(pl)^{-1}\}(ik)^\delta(pq)^\gamma\\
&&-\{(jl)^{\varepsilon},(ql)(pl)(ql)^{-1}(pl)^{-1}\}(ik)^\delta(pq)^\gamma\\
&\equiv&0.
\end{eqnarray*}
\item[d)]\ $j=q$,\ that is, \ $p<q=j<i<k<l$.
$$
(f,g)_w=(pj)^\gamma(jl)^{\varepsilon}(ik)^\delta-(ik)^\delta(pj)^\gamma(jl)^{\varepsilon}.
$$
If $\gamma=1$, then
\begin{eqnarray*}
(f,g)_w&\equiv&\{(jl)^{\varepsilon},(pl)\}(ik)^\delta(pj)^\gamma-(ik)^\delta\{(jl)^{\varepsilon},(pl)\}(pj)^\gamma\\
&\equiv&\{(jl)^{\varepsilon},(pl)\}(ik)^\delta(pj)^\gamma-\{(jl)^{\varepsilon},(pl)\}(ik)^\delta(pj)^\gamma\\
&\equiv&0.
\end{eqnarray*}
If $\gamma=-1$, then
\begin{eqnarray*}
(f,g)_w&\equiv&\{(jl)^{\varepsilon},(pl)^{-1}(jl)^{-1}\}(ik)^\delta(pj)^\gamma
-(ik)^\delta\{(jl)^{\varepsilon},(pl)^{-1}(jl)^{-1}\}(pj)^\gamma\\
&\equiv&\{(jl)^{\varepsilon},(pl)^{-1}(jl)^{-1}\}(ik)^\delta(pj)^\gamma
-\{(jl)^{\varepsilon},(pl)^{-1}(jl)^{-1}\}(ik)^\delta(pj)^\gamma\\
&\equiv&0.
\end{eqnarray*}
\item[e)]\ $j>q$,\ that is, \ $p<q<j<i<k<l$.\\
$$
(f,g)_w\equiv(jl)^{\varepsilon}(ik)^\delta(pq)^\gamma-(jl)^{\varepsilon}(ik)^\delta(pq)^\gamma\equiv0.
$$
\end{enumerate}
\item[4)]\ $p<q<i<k$,\ that is, \ $i<k<j<l,p<q<i<k<j<l$.\\
$$
(f,g)_w\equiv(jl)^{\varepsilon}(ik)^\delta(pq)^\gamma-(jl)^{\varepsilon}(ik)^\delta(pq)^\gamma\equiv0.
$$
\end{enumerate}

$(14)\wedge(1)$\\
Let
$f=\sigma_{j}^{-1}\sigma_{k}^{-1}-\sigma_{k}^{-1}\sigma_{j}^{-1},j<k-1,\
\ g=\sigma_{k}^{-1}(pq)^\delta-(pq)^\delta\sigma_{k}^{-1},k\neq
p,p-1,q,q-1$. Then $w=\sigma_{j}^{-1}\sigma_{k}^{-1}(pq)^\delta $
and
$$
(f,g)_w=\sigma_{j}^{-1}(pq)^\delta\sigma_{k}^{-1}-\sigma_{k}^{-1}\sigma_{j}^{-1}(pq)^\delta
\equiv(pq)^\delta\sigma_{k}^{-1}\sigma_{j}^{-1}-(pq)^\delta\sigma_{k}^{-1}\sigma_{j}^{-1}\equiv0.
$$

$(14)\wedge(2)$\\
Let
$f=\sigma_{j}^{-1}\sigma_{k}^{-1}-\sigma_{k}^{-1}\sigma_{j}^{-1},j<k-1,\
g=\sigma_{k}^{-1}(k,k+1)^\delta-(k,k+1)^\delta\sigma_{k}^{-1}$. Then
$w=\sigma_{j}^{-1}\sigma_{k}^{-1}(k,k+1)^\delta$ and
$$
(f,g)_w=\sigma_{j}^{-1}(k,k+1)^\delta\sigma_{k}^{-1}-\sigma_{k}^{-1}\sigma_{j}^{-1}(k,k+1)^\delta
\equiv(k,k+1)^\delta\sigma_{k}^{-1}\sigma_{j}^{-1}-(k,k+1)^\delta\sigma_{k}^{-1}\sigma_{j}^{-1}\equiv0.
$$

$(14)\wedge(3)$\\
Let
$f=\sigma_{j}^{-1}\sigma_{k}^{-1}-\sigma_{k}^{-1}\sigma_{j}^{-1},j<k-1,\
g=\sigma_{k}^{-1}(k+1,q)^\delta-(k,q)^\delta\sigma_{k}^{-1},k+1<q$.
Then $w=\sigma_{j}^{-1}\sigma_{k}^{-1}(k+1,q)^\delta$ and
$$
(f,g)_w=\sigma_{j}^{-1}(k,q)^\delta\sigma_{k}^{-1}-\sigma_{k}^{-1}\sigma_{j}^{-1}(k+1,q)^\delta
\equiv(k,q)^\delta\sigma_{k}^{-1}\sigma_{j}^{-1}-(k,q)^\delta\sigma_{k}^{-1}\sigma_{j}^{-1}\equiv0.
$$

$(14)\wedge(4)$\\
Let
$f=\sigma_{j}^{-1}\sigma_{k}^{-1}-\sigma_{k}^{-1}\sigma_{j}^{-1},j<k-1,\
g=\sigma_{k}^{-1}(kq)^\delta-\{(k+1,q)^\delta,(k,k+1)\}\sigma_{k}^{-1},i+1<q$.
Then $w=\sigma_{j}^{-1}\sigma_{k}^{-1}(kq)^\delta$ and
\begin{eqnarray*}
(f,g)_w&=&\sigma_{j}^{-1}\{(k+1,q)^\delta,(k,k+1)\}\sigma_{k}^{-1}-\sigma_{k}^{-1}\sigma_{j}^{-1}(kq)^\delta\\
&\equiv&\{(k+1,q)^\delta,(k,k+1)\}\sigma_{k}^{-1}\sigma_{j}^{-1}
-\{(k+1,q)^\delta,(k,k+1)\}\sigma_{k}^{-1}\sigma_{j}^{-1}\\
&\equiv&0.
\end{eqnarray*}

$(14)\wedge(5)$\\
Let
$f=\sigma_{j}^{-1}\sigma_{k}^{-1}-\sigma_{k}^{-1}\sigma_{j}^{-1},j<k-1,\
g=\sigma_{k}^{-1}(q,k+1)^\delta\sigma_{k}^{-1}-(q,k)^\delta\sigma_{k}^{-1},q<k$.
Then $w=\sigma_{j}^{-1}\sigma_{k}^{-1}(q,k+1)^\delta$ and
$$
(f,g)_w=\sigma_{j}^{-1}(q,k)^\delta\sigma_{k}^{-1}-\sigma_{k}^{-1}\sigma_{j}^{-1}(q,k+1)^\delta.
$$
There are four cases to consider.
\begin{enumerate}
\item[1)]\ $q=k-1,j<q$.
$$
(f,g)_w=\sigma_{j}^{-1}(k-1,k)^\delta\sigma_{k}^{-1}-\sigma_{k}^{-1}\sigma_{j}^{-1}(k-1,k+1)^\delta.
$$
If  $ j=k-2$, then
\begin{eqnarray*}
(f,g)_w&=&\sigma_{k-2}^{-1}(k-1,k)^\delta\sigma_{k}^{-1}-\sigma_{k}^{-1}\sigma_{k-2}^{-1}(k-1,k+1)^\delta\\
&\equiv&(k-2,k)^\delta\sigma_{k}^{-1}\sigma_{k-2}^{-1}-\sigma_{k}^{-1}(k-2,k+1)^\delta\sigma_{k-2}^{-1}\\
&\equiv&(k-2,k)^\delta\sigma_{k}^{-1}\sigma_{k-2}^{-1}-(k-2,k)^\delta\sigma_{k}^{-1}\sigma_{k-2}^{-1}\\
&\equiv&0.
\end{eqnarray*}
If $j<k-2$, then
\begin{eqnarray*}
(f,g)_w&=&(k-1,k)^\delta\sigma_{k}^{-1}\sigma_{j}^{-1}-\sigma_{k}^{-1}(k-1,k+1)^\delta\sigma_{j}^{-1}\\
&\equiv&(k-1,k)^\delta\sigma_{k}^{-1}\sigma_{j}^{-1}-(k-1,k)^\delta\sigma_{k}^{-1}\sigma_{j}^{-1}\\
&\equiv&0.
\end{eqnarray*}
\item[2)]\ $j<q<k-1$.\\
If $ j=q-1$ , then
\begin{eqnarray*}
(f,g)_w&=&\sigma_{q-1}^{-1}(q,k)^\delta\sigma_{k}^{-1}-\sigma_{k}^{-1}\sigma_{q-1}^{-1}(q,k+1)^\delta\\
&\equiv&(q-1,k)^\delta\sigma_{k}^{-1}\sigma_{q-1}^{-1}-\sigma_{k}^{-1}(q-1,k+1)^\delta\sigma_{q-1}^{-1}\\
&\equiv&(q-1,k)^\delta\sigma_{k}^{-1}\sigma_{q-1}^{-1}-(q-1,k)^\delta\sigma_{k}^{-1}\sigma_{q-1}^{-1}\\
&\equiv&0.
\end{eqnarray*}
If $j<q-1$, then
$$
(f,g)_w\equiv(q,k)^\delta\sigma_{k}^{-1}\sigma_{j}^{-1}-\sigma_{k}^{-1}(q,k+1)^\delta\sigma_{j}^{-1}
\equiv(q,k)^\delta\sigma_{k}^{-1}\sigma_{j}^{-1}-(q,k)^\delta\sigma_{k}^{-1}\sigma_{j}^{-1}\equiv0.
$$
\item[3)]\ $q=j<k-1$.
\begin{eqnarray*}
(f,g)_w&=&\sigma_{q}^{-1}(q,k)^\delta\sigma_{k}^{-1}-\sigma_{k}^{-1}\sigma_{q}^{-1}(q,k+1)^\delta,\\
&\equiv&\{(q+1,k)^\delta,(q,q+1)\}\sigma_{k}^{-1}\sigma_{q}^{-1}-
\sigma_{k}^{-1}\{(q+1,k+1)^\delta,(q,q+1)\}\sigma_{q}^{-1}\\
&\equiv&\{(q+1,k)^\delta,(q,q+1)\}\sigma_{k}^{-1}\sigma_{q}^{-1}-
\{(q+1,k)^\delta,(q,q+1)\}\sigma_{k}^{-1}\sigma_{q}^{-1}\\
&\equiv&0.
\end{eqnarray*}
\item[4)]\ $q<j<k-1$.
$$
(f,g)_w\equiv(q,k)^\delta\sigma_{k}^{-1}\sigma_{j}^{-1}-\sigma_{k}^{-1}(q,k+1)^\delta\sigma_{j}^{-1}
\equiv(q,k)^\delta\sigma_{k}^{-1}\sigma_{j}^{-1}-(q,k)^\delta\sigma_{k}^{-1}\sigma_{j}^{-1}\equiv0.
$$
\end{enumerate}

$(14)\wedge(6)$\\
Let
$f=\sigma_{j}^{-1}\sigma_{k}^{-1}-\sigma_{k}^{-1}\sigma_{j}^{-1},j<k-1,\
g=\sigma_{k}^{-1}(q,k)^\delta-\{(q,k+1)^\delta,(k,k+1)\}\sigma_{k}^{-1},q<k$.
Then $w=\sigma_{j}^{-1}\sigma_{k}^{-1}(q,k)^\delta$ and
$$
(f,g)_w=\sigma_{j}^{-1}\{(q,k+1)^\delta,(k,k+1)\}\sigma_{k}^{-1}-\sigma_{k}^{-1}\sigma_{j}^{-1}(q,k)^\delta.
$$
There are four cases to consider.
\begin{enumerate}
\item[1)]\ $j=q-1<k-1$.
\begin{eqnarray*}
(f,g)_w&=&\sigma_{q-1}^{-1}\{(q,k+1)^\delta,(k,k+1)\}\sigma_{k}^{-1}-\sigma_{k}^{-1}\sigma_{q-1}^{-1}(q,k)^\delta\\
&\equiv&\{(q-1,k+1)^\delta,(k,k+1)\}\sigma_{k}^{-1}\sigma_{q-1}^{-1}-\sigma_{k}^{-1}(q-1,k)^\delta\sigma_{q-1}^{-1}\\
&\equiv&\{(q-1,k+1)^\delta,(k,k+1)\}\sigma_{k}^{-1}\sigma_{q-1}^{-1}-
\{(q-1,k+1)^\delta,(k,k+1)\}\sigma_{k}^{-1}\sigma_{q-1}^{-1}\\
&\equiv&0.
\end{eqnarray*}
\item[2)]\ $j<q-1,q<k$.
\begin{eqnarray*}
(f,g)_w&=&\{(q,k+1)^\delta,(k,k+1)\}\sigma_{k}^{-1}\sigma_{j}^{-1}-\sigma_{k}^{-1}(qk)^\delta\sigma_{j}^{-1}\\
&\equiv&\{(q,k+1)^\delta,(k,k+1)\}\sigma_{k}^{-1}\sigma_{j}^{-1}-
\{(q,k+1)^\delta,(k,k+1)\}\sigma_{k}^{-1}\sigma_{j}^{-1}\\
&\equiv&0.
\end{eqnarray*}
\item[3)]\ $j=q<k-1$.
\begin{eqnarray*}
(f,g)_w&=&\sigma_{q}^{-1}\{(q,k+1)^\delta,(k,k+1)\}\sigma_{k}^{-1}-\sigma_{k}^{-1}\sigma_{q}^{-1}(q,k)^\delta\\
&\equiv&\{\{(q+1,k+1)^\delta,(q,q+1)\},(k,k+1)\}\sigma_{k}^{-1}\sigma_{q}^{-1}\\
&&-\sigma_{k}^{-1}\{(q+1,k)^\delta,(q,q+1)\}\sigma_{q}^{-1}\\
&\equiv&\{(q+1,k+1)^\delta,(q,q+1)(k,k+1)\}\sigma_{k}^{-1}\sigma_{q}^{-1}\\
&&-\{\{(q+1,k+1)^\delta,(k,k+1)\},(q,q+1)\}\sigma_{k}^{-1}\sigma_{q}^{-1}\\
&\equiv&\{(q+1,k+1)^\delta,(k,k+1)(q,q+1)\}\sigma_{k}^{-1}\sigma_{q}^{-1}\\
&&-\{(q+1,k+1)^\delta,(k,k+1)\},(q,q+1)\}\sigma_{k}^{-1}\sigma_{q}^{-1}\\
&\equiv&0.
\end{eqnarray*}
\item[4)]\ $q<j<k-1$.
$$
(f,g)_w\equiv\{(q,k+1)^\delta,(k,k+1)\}\sigma_{k}^{-1}\sigma_{j}^{-1}-
\sigma_{k}^{-1}(q,k)^\delta\sigma_{j}^{-1}\equiv0.
$$
\end{enumerate}

$(14)\wedge(14)$\\
Let
$f=\sigma_{j}^{-1}\sigma_{k}^{-1}-\sigma_{k}^{-1}\sigma_{j}^{-1},\
g=\sigma_{k}^{-1}\sigma_{l}^{-1}-\sigma_{l}^{-1}\sigma_{k}^{-1},j<k-1<k<l-1<l$.
Then $ w=\sigma_{j}^{-1}\sigma_{k}^{-1}\sigma_{l}^{-1}$ and
$$
(f,g)_w=\sigma_{j}^{-1}\sigma_{l}^{-1}\sigma_{k}^{-1}-\sigma_{k}^{-1}\sigma_{j}^{-1}\sigma_{l}^{-1}
\equiv\sigma_{l}^{-1}\sigma_{k}^{-1}\sigma_{j}^{-1}-\sigma_{l}^{-1}\sigma_{k}^{-1}\sigma_{j}^{-1}\equiv0.
$$

$(14)\wedge(15)$\\
Let
$f=\sigma_{q}^{-1}\sigma_{j}^{-1}-\sigma_{j}^{-1}\sigma_{q}^{-1},q<j-1,\
\
g=\sigma_{j}^{-1}\sigma_{k,j+1}-\sigma_{k,j+1}\sigma_{j-1}^{-1},k<j$.
Then $w=\sigma_{q}^{-1}\sigma_{j}^{-1}\sigma_{k,j+1}$ and
$$
(f,g)_w=\sigma_{q}^{-1}\sigma_{k,j+1}\sigma_{j-1}^{-1}-\sigma_{j}^{-1}\sigma_{q}^{-1}\sigma_{k,j+1}.
$$
There are four cases to consider.
\begin{enumerate}
\item[1)]\ $k<q<j-1$.
\begin{eqnarray*}
(f,g)_w&\equiv&\sigma_{k}^{-1}\cdots\sigma_{j}^{-1}\sigma_{q-1}^{-1}\sigma_{j-1}^{-1}-
\sigma_{j}^{-1}\sigma_{k}^{-1}\cdots\sigma_{j}^{-1}\sigma_{q-1}^{-1}\\
&\equiv&\sigma_{k}^{-1}\cdots\sigma_{j}^{-1}\sigma_{j-1}^{-1}\sigma_{q-1}^{-1}-
\sigma_{k}^{-1}\cdots\sigma_{j}^{-1}\sigma_{j-1}^{-1}\sigma_{q-1}^{-1}\\
&\equiv&0.
\end{eqnarray*}
\item[2)]\ $q=k<j-1$.
\begin{eqnarray*}
(f,g)_w&=&\sigma_{k}^{-1}\sigma_{k}^{-1}\cdots\sigma_{j}^{-1}\sigma_{j-1}^{-1}-
\sigma_{j}^{-1}\sigma_{k}^{-1}\sigma_{k}^{-1}\cdots\sigma_{j}^{-1}\\
&\equiv&(k,k+1)^{-1}\sigma_{k+1}^{-1}\cdots\sigma_{j}^{-1}\sigma_{j-1}^{-1}-
\sigma_{j}^{-1}(k,k+1)^{-1}\sigma_{k+1}^{-1}\cdots\sigma_{j}^{-1}\\
&\equiv&(k,k+1)^{-1}\sigma_{k+1}^{-1}\cdots\sigma_{j}^{-1}\sigma_{j-1}^{-1}-
(k,k+1)^{-1}\sigma_{j}^{-1}\sigma_{k+1}^{-1}\cdots\sigma_{j}^{-1}\\
&\equiv&(k,k+1)^{-1}\sigma_{k+1}^{-1}\cdots\sigma_{j}^{-1}\sigma_{j-1}^{-1}-
(k,k+1)^{-1}\sigma_{k+1}^{-1}\cdots\sigma_{j}^{-1}\sigma_{j-1}^{-1}\\
&\equiv&0.
\end{eqnarray*}
\item[3)]\ $q=k-1<j-1$.
\begin{eqnarray*}
(f,g)_w&=&\sigma_{k-1}^{-1}\sigma_{k}^{-1}\cdots\sigma_{j}^{-1}\sigma_{j-1}^{-1}
-\sigma_{j}^{-1}\sigma_{k-1}^{-1}\sigma_{k}^{-1}\cdots\sigma_{j}^{-1}\\
&\equiv&\sigma_{k-1}^{-1}\sigma_{k}^{-1}\cdots\sigma_{j}^{-1}\sigma_{j-1}^{-1}
-\sigma_{k-1}^{-1}\sigma_{k}^{-1}\cdots\sigma_{j}^{-1}\sigma_{j-1}^{-1}\\
&\equiv&0.
\end{eqnarray*}
\item[4)]\ $q<k-1<j-1$.
\begin{eqnarray*}
(f,g)_w&\equiv&\sigma_{k}^{-1}\cdots\sigma_{j}^{-1}\sigma_{j-1}^{-1}\sigma_{q}^{-1}
-\sigma_{j}^{-1}\sigma_{k}^{-1}\cdots\sigma_{j}^{-1}\sigma_{q}^{-1}\\
&\equiv&\sigma_{k}^{-1}\cdots\sigma_{j}^{-1}\sigma_{j-1}^{-1}\sigma_{q}^{-1}
-\sigma_{k}^{-1}\cdots\sigma_{j}^{-1}\sigma_{j-1}^{-1}\sigma_{q}^{-1}\\
&\equiv&0.
\end{eqnarray*}
\end{enumerate}

$(14)\wedge(16)$\\
Let
$f=\sigma_{j}^{-1}\sigma_{k}^{-1}-\sigma_{k}^{-1}\sigma_{j}^{-1},j<k-1,\
 g=\sigma_{k}^{-1}\sigma_{k}^{-1}-(k,k+1)^{-1}$. Then $w=\sigma_{j}^{-1}\sigma_{k}^{-1}\sigma_{k}^{-1}$
 and
\begin{eqnarray*}
(f,g)_w&=&\sigma_{j}^{-1}(k,k+1)^{-1}-\sigma_{k}^{-1}\sigma_{j}^{-1}\sigma_{k}^{-1}\\
&\equiv&(k,k+1)^{-1}\sigma_{j}^{-1}-\sigma_{k}^{-1}\sigma_{k}^{-1}\sigma_{j}^{-1}\\
&\equiv&(k,k+1)^{-1}\sigma_{j}^{-1}-(k,k+1)^{-1}\sigma_{j}^{-1}\\
&\equiv&0.
\end{eqnarray*}

$(15)\wedge(1)$\\
Let
$f=\sigma_{j}^{-1}\sigma_{k}^{-1}\cdots\sigma_{j}^{-1}-\sigma_{k}^{-1}\cdots\sigma_{j}^{-1}\sigma_{j-1}^{-1},
g=\sigma_{j}^{-1}(pq)^\delta-(pq)^\delta\sigma_{j}^{-1},k<j,j\neq
p,p-1,q,q-1$.\ Then
$w=\sigma_{j}^{-1}\sigma_{k}^{-1}\cdots\sigma_{j}^{-1}(pq)^\delta$\
\ and\\
$$
(f,g)_w=\sigma_{j}^{-1}\sigma_{k}^{-1}\cdots\sigma_{j-1}^{-1}(pq)^\delta\sigma_{j}^{-1}
-\sigma_{k}^{-1}\cdots\sigma_{j}^{-1}\sigma_{j-1}^{-1}(pq)^\delta.
$$
There are four cases to consider.
\begin{enumerate}
\item[1)]\ $p>j+1$.
\begin{eqnarray*}
(f,g)_w&\equiv&(pq)^\delta\sigma_{j}^{-1}\sigma_{k}^{-1}\cdots\sigma_{j-1}^{-1}\sigma_{j}^{-1}-(pq)^\delta
\sigma_{k}^{-1}\cdots\sigma_{j}^{-1}\sigma_{j-1}^{-1}\\
&\equiv&(pq)^\delta\sigma_{k}^{-1}\cdots\sigma_{j}^{-1}\sigma_{j-1}^{-1}-(pq)^\delta\sigma_{k}^{-1}
\cdots\sigma_{j}^{-1}\sigma_{j-1}^{-1}\\
&\equiv&0.
\end{eqnarray*}
\item[2)]\  $p=j-1$.

In this case, we have $q>j+1$ and
\begin{eqnarray*}
&&\sigma_{j}^{-1}\sigma_{k}^{-1}\cdots\sigma_{j-1}^{-1}(j-1,q)^\delta\sigma_{j}^{-1}\\
&\equiv&\sigma_{j}^{-1}\sigma_{k}^{-1}\cdots\sigma_{j-2}^{-1}\{(j,q)^\delta,(j-1,j)\}
\sigma_{j-1}^{-1}\sigma_{j}^{-1}\\
&\equiv&\sigma_{j}^{-1}\sigma_{k}^{-1}\cdots\sigma_{j-2}^{-1}(j-1,q)(jq)^\delta(j-1,q)^{-1}
\sigma_{j-1}^{-1}\sigma_{j}^{-1}\\
&\equiv&\sigma_{j}^{-1}(kq)(jq)^\delta(kq)^{-1}\sigma_{k}^{-1}\cdots\sigma_{j}^{-1}\\
&\equiv&(kq)\{(j+1,q)^\delta,(j,j+1)\}(kq)^{-1}\sigma_{k}^{-1}\cdots\sigma_{j}^{-1}\sigma_{j-1}^{-1}\\
&\equiv&(kq)(jq)(j+1,q)^\delta(jq)^{-1}(kq)^{-1}\sigma_{k}^{-1}\cdots\sigma_{j}^{-1}\sigma_{j-1}^{-1},\ \ \ \ \ and \\
&&\\
&&\sigma_{k}^{-1}\cdots\sigma_{j}^{-1}\sigma_{j-1}^{-1}(j-1,q)^\delta\\
&\equiv&\sigma_{k}^{-1}\cdots\sigma_{j}^{-1}\{(jq)^\delta,(j-1,j)\}\sigma_{j-1}^{-1}\\
&\equiv&\sigma_{k}^{-1}\cdots\sigma_{j}^{-1}(j-1,q)(jq)^\delta(j-1,q)^{-1}\sigma_{j-1}^{-1}\\
&\equiv&\sigma_{k}^{-1}\cdots\sigma_{j-1}^{-1}(j-1,q)(jq)(j+1,q)^\delta(j-1,q)^{-1}\sigma_{j}^{-1}\sigma_{j-1}^{-1}\\
&\equiv&(kq)(jq)(kq)^{-1}(kq)(j+i,q)^\delta(kq)^{-1}(kq)(jq)^{-1}(kq)^{-1}\sigma_{k}^{-1}\cdots\sigma_{j}^{-1}
\sigma_{j-1}^{-1}\\
&\equiv&(kq)(jq)(j+1,q)^\delta(jq)^{-1}(kq)^{-1}\sigma_{k}^{-1}\cdots\sigma_{j}^{-1}\sigma_{j-1}^{-1}.
\end{eqnarray*}
\item[3)]\ $k\leq p<j-1$. In this case, there are three subcases to
consider.
\begin{enumerate}
\item[a)]\ $ q=j-1$.
\begin{eqnarray*}
&&\sigma_{j}^{-1}\sigma_{k}^{-1}\cdots\sigma_{j-1}^{-1}(pq)^\delta\sigma_{j}^{-1}\\
&\equiv&\sigma_{j}^{-1}\sigma_{k}^{-1}\cdots\sigma_{j-2}^{-1}\{(pj)^\delta,(j-1,j)\}
\sigma_{j-1}^{-1}\sigma_{j}^{-1}\\
&\equiv&\sigma_{j}^{-1}\sigma_{k}^{-1}\cdots\sigma_{p}^{-1}\{(pj)^\delta,(p+1,j)\}
\sigma_{p+1}^{-1}\cdots\sigma_{j}^{-1}\\
&\equiv&\sigma_{j}^{-1}\sigma_{k}^{-1}\cdots\sigma_{p-1}^{-1}\{(p+1,j)^\delta,(p,p+1)(pj)\}
\sigma_{p}^{-1}\cdots\sigma_{j}^{-1}\\
&\equiv&\sigma_{j}^{-1}\sigma_{k}^{-1}\cdots\sigma_{p-1}^{-1}(p+1,j)^\delta\sigma_{p}^{-1}\cdots\sigma_{j}^{-1}\\
&\equiv&\{(p+1,j+1)^\delta,(j,j+1)\}\sigma_{k}^{-1}\cdots\sigma_{j}^{-1}\sigma_{j-1}^{-1}\ \ \ \ \ \ \ and \\
&&\\
&&\sigma_{k}^{-1}\cdots\sigma_{j}^{-1}\sigma_{j-1}^{-1}(pq)^\delta\\
&\equiv&\sigma_{k}^{-1}\cdots\sigma_{j}^{-1}\{(pj)^\delta,(j-1,j)\}\sigma_{j-1}^{-1}\\
&\equiv&\sigma_{k}^{-1}\cdots\sigma_{j-1}^{-1}\{\{(p,j+1)^\delta,(j,j+1)\},\{(j-1,j+1),(j,j+1)\}\}
\sigma_{j}^{-1}\sigma_{j-1}^{-1}\\
&\equiv&\sigma_{k}^{-1}\cdots\sigma_{j-1}^{-1}\{(p,j+1)^\delta,(j-1,j+1)(j,j+1)\}\sigma_{j}^{-1}
\sigma_{j-1}^{-1}\\
&\equiv&\sigma_{k}^{-1}\cdots\sigma_{j-2}^{-1}\{(p,j+1)^\delta,\{(j,j+1),(j-1,j+1)\}(j-1,j+1)\}
\sigma_{j-1}^{-1}\sigma_{j}^{-1}\sigma_{j-1}^{-1}\\
&\equiv&\sigma_{k}^{-1}\cdots\sigma_{j-2}^{-1}\{(p,j+1)^\delta,(j-1,j+1)^{-1}(j,j+1)\}
\sigma_{j-1}^{-1}\sigma_{j}^{-1}\sigma_{j-1}^{-1}\\
&\equiv&\sigma_{k}^{-1}\cdots\sigma_{p}^{-1}\{(p,j+1)^\delta,(p+1,j+1)^{-1}(j,j+1)\}
\sigma_{p+1}^{-1}\cdots\sigma_{j}^{-1}\sigma_{j-1}^{-1}\\
&\equiv&\sigma_{k}^{-1}\cdots\sigma_{p-1}^{-1}\{(p+1,j+1)^\delta,(p,p+1)(p,j+1)^{-1}(j,j+1)\}
\sigma_{p}^{-1}\cdots\sigma_{j}^{-1}\sigma_{j-1}^{-1}\\
&\equiv&\sigma_{k}^{-1}\cdots\sigma_{p-1}^{-1}(j,j+1)^{-1}(p+1,j+1)^\delta(j,j+1)
\sigma_{p}^{-1}\cdots\sigma_{j}^{-1}\sigma_{j-1}^{-1}\\
&\equiv&\{(p+1,j+1)^\delta,(j,j+1)\}\sigma_{k}^{-1}\cdots\sigma_{j}^{-1}\sigma_{j-1}^{-1}.
\end{eqnarray*}
\item[b)]\ $q>j+1$.
\begin{eqnarray*}
&&\sigma_{j}^{-1}\sigma_{k}^{-1}\cdots\sigma_{j-1}^{-1}(pq)^\delta\sigma_{j}^{-1}\\
&\equiv&\sigma_{j}^{-1}\sigma_{k}^{-1}\cdots\sigma_{p}^{-1}(pq)^\delta\sigma_{p+1}^{-1}\cdots\sigma_{j}^{-1}\\
&\equiv&\sigma_{j}^{-1}\sigma_{k}^{-1}\cdots\sigma_{p-1}^{-1}\{(p+1,q)^\delta,(p,p+1)\}
\sigma_{p}^{-1}\cdots\sigma_{j}^{-1}\\
&\equiv&\sigma_{j}^{-1}\sigma_{k}^{-1}\cdots\sigma_{p-1}^{-1}(pq)^{-1}(p+1,q)^\delta(pq)
\sigma_{p}^{-1}\cdots\sigma_{j}^{-1}\\
&\equiv&\sigma_{j}^{-1}(kq)^{-1}(p+1,q)^\delta(kq)\sigma_{k}^{-1}\cdots\sigma_{j}^{-1}\\
&\equiv&(kq)^{-1}(p+1,q)^\delta(kq)\sigma_{k}^{-1}\cdots\sigma_{j}^{-1}\sigma_{j-1}^{-1}\ \ \ \ \ \ \ and \\
&&\\
&&\sigma_{k}^{-1}\cdots\sigma_{j}^{-1}\sigma_{j-1}^{-1}(pq)^\delta\\
&\equiv&\sigma_{k}^{-1}\cdots\sigma_{j}^{-1}(pq)^\delta\sigma_{j-1}^{-1}\\
&\equiv&\sigma_{k}^{-1}\cdots\sigma_{p}^{-1}(pq)^\delta\sigma_{j}^{-1}\sigma_{j-1}^{-1}\\
&\equiv&\sigma_{kq}^{-1}(p+1,q)^\delta(kq)\sigma_{k}^{-1}\cdots\sigma_{j}^{-1}\sigma_{j-1}^{-1}.
\end{eqnarray*}
\item[c)]\ $q<j-1$.

In this case, we have $k\leqslant p<q<j-1$ and
\begin{eqnarray*}
&&\sigma_{j}^{-1}\sigma_{k}^{-1}\cdots\sigma_{j-1}^{-1}(pq)^\delta\sigma_{j}^{-1}\\
&\equiv&\sigma_{j}^{-1}\sigma_{k}^{-1}\cdots\sigma_{q}^{-1}(pq)^\delta\sigma_{q+1}^{-1}\cdots\sigma_{j}^{-1}\\
&\equiv&\sigma_{j}^{-1}\sigma_{k}^{-1}\cdots\sigma_{q-1}^{-1}\{(p,q+1)^\delta,(q,q+1)\}
\sigma_{q}^{-1}\cdots\sigma_{j}^{-1}\\
&\equiv&\sigma_{j}^{-1}\sigma_{k}^{-1}\cdots\sigma_{p}^{-1}\{(p,q+1)^\delta,(p+1,q+1)\}
\sigma_{p+1}^{-1}\cdots\sigma_{j}^{-1}\\
&\equiv&\sigma_{j}^{-1}\sigma_{k}^{-1}\cdots\sigma_{p-1}^{-1}\{(p+1,q+1)^\delta,(p,p+1)(q,q+1)\}
\sigma_{p}^{-1}\cdots\sigma_{j}^{-1}\\
&\equiv&\sigma_{j}^{-1}\sigma_{k}^{-1}\cdots\sigma_{p-1}^{-1}(p+1,q+1)^\delta\sigma_{p}^{-1}\cdots\sigma_{j}^{-1}\\
&\equiv&\sigma_{j}^{-1}(p+1,q+1)^\delta\sigma_{k}^{-1}\cdots\sigma_{j}^{-1}\\
&\equiv&(p+1,q+1)^\delta\sigma_{k}^{-1}\cdots\sigma_{j}^{-1}\sigma_{j-1}^{-1},\ \ \ \ \ \ \ \ \ and\\
&&\\
&&\sigma_{k}^{-1}\cdots\sigma_{j}^{-1}\sigma_{j-1}^{-1}(pq)^\delta\\
&\equiv&\sigma_{k}^{-1}\cdots\sigma_{j}^{-1}(pq)^\delta\sigma_{j-1}^{-1}\\
&\equiv&\sigma_{k}^{-1}\cdots\sigma_{q}^{-1}(pq)^\delta\sigma_{q+1}^{-1}\cdots\sigma_{j}^{-1}\sigma_{j-1}^{-1}\\
&\equiv&(p+1,q+1)^\delta\sigma_{k}^{-1}\cdots\sigma_{j}^{-1}\sigma_{j-1}^{-1}.
\end{eqnarray*}
\end{enumerate}
\item[4)]\ $p<k$. In this case, there are three subcases to
consider.
\begin{enumerate}
\item[a)]\ $\ q>j+1,p<k<j<j+1<q$.\\
\begin{eqnarray*}
&&\sigma_{j}^{-1}\sigma_{k}^{-1}\cdots\sigma_{j-1}^{-1}(pq)^\delta\sigma_{j}^{-1}
\equiv\sigma_{j}^{-1}(pq)^\delta\sigma_{k}^{-1}\cdots\sigma_{j}^{-1}
\equiv(pq)^\delta\sigma_{k}^{-1}\cdots\sigma_{j}^{-1}\sigma_{j-1}^{-1}\\
&&\\
&&  \mbox{ and } \ \ \ \ \
\sigma_{k}^{-1}\cdots\sigma_{j}^{-1}\sigma_{j-1}^{-1}(pq)^\delta
\equiv(pq)^\delta\sigma_{k}^{-1}\cdots\sigma_{j}^{-1}\sigma_{j-1}^{-1}.
\end{eqnarray*}
\item[b)]\ $ \ k\leqslant q\leqslant j-1<j$.
\begin{eqnarray*}
&&\sigma_{j}^{-1}\sigma_{k}^{-1}\cdots\sigma_{j-1}^{-1}(pq)^\delta\sigma_{j}^{-1}\\
&\equiv&\sigma_{j}^{-1}\sigma_{k}^{-1}\cdots\sigma_{q}^{-1}(pq)^\delta\sigma_{q+1}^{-1}\cdots\sigma_{j}^{-1}\\
&\equiv&\sigma_{j}^{-1}\sigma_{k}^{-1}\cdots\sigma_{q-1}^{-1}\{(p,q+1)^\delta,(q,q+1)\}\sigma_{q}^{-1}\cdots
\sigma_{j}^{-1}\\
&\equiv&\sigma_{j}^{-1}\{(p,q+1)^\delta,(k,q+1)\}\sigma_{k}^{-1}\cdots\sigma_{j}^{-1}\\
&\equiv& \left\{
\begin{array}{ll} \{(p,j+1)^\delta,(k,j+1)(j,j+1)\}\sigma_{k}^{-1}
\cdots\sigma_{j}^{-1}\sigma_{j-1}^{-1} \ \ \ \ \ \ \ \  if\ \ q+1=j,\\
\{(p,q+1)^\delta,(k,q+1)\}\sigma_{k}^{-1}\cdots\sigma_{j}^{-1}\sigma_{j-1}^{-1}
\ \ \ \ \ \ \ \ \ \ \ \ \ \ \ \ \ \ \ \ if\ \ q+1<j.
\end{array}\right.
\end{eqnarray*}
If \ $q+1=j$, \  then
\begin{eqnarray*}
&&\sigma_{k}^{-1}\cdots\sigma_{j}^{-1}\sigma_{j-1}^{-1}(pq)^\delta\\
&\equiv&\sigma_{k}^{-1}\cdots\sigma_{j}^{-1}\{(p,j)^\delta,(j-1,j)\}\sigma_{j-1}^{-1}\\
&\equiv&\sigma_{k}^{-1}\cdots\sigma_{j-1}^{-1}\{\{(p,j+1)^\delta,(j,j+1)\},\{(j-1,j+1),(j,j+1)\}\}
\sigma_{j}^{-1}\sigma_{j-1}^{-1}\\
&\equiv&\sigma_{k}^{-1}\cdots\sigma_{j-1}^{-1}\{(p,j+1)^\delta,(j-1,j+1)(j,j+1)\}\sigma_{j}^{-1}\sigma_{j-1}^{-1}\\
&\equiv&\sigma_{k}^{-1}\cdots\sigma_{j-2}^{-1}\{(p,j+1)^\delta,\{(j,j+1),(j-1,j)\}(j-1,j+1)\}\sigma_{j-1}^{-1}
\sigma_{j}^{-1}\sigma_{j-1}^{-1}\\
&\equiv&\sigma_{k}^{-1}\cdots\sigma_{j-2}^{-1}\{(p,j+1)^\delta,(j-1,j+1)(j,j+1)\}\sigma_{j-1}^{-1}\sigma_{j}^{-1}
\sigma_{j-1}^{-1}\\
&\equiv&\{(p,j+1)^\delta,(k,j+1)(j,j+1)\}\sigma_{k}^{-1}\cdots\sigma_{j}^{-1}\sigma_{j-1}^{-1}.
\end{eqnarray*}
If\ $q+1<j$, \ then
\begin{eqnarray*}
&&\sigma_{k}^{-1}\cdots\sigma_{j}^{-1}\sigma_{j-1}^{-1}(pq)^\delta
\equiv\sigma_{k}^{-1}\cdots\sigma_{j}^{-1}(pq)^\delta\sigma_{j-1}^{-1}\\
&\equiv&\sigma_{k}^{-1}\cdots\sigma_{q}^{-1}(pq)^\delta\sigma_{q+1}^{-1}\cdots\sigma_{j}^{-1}\sigma_{j-1}^{-1}\\
&\equiv&\{(p,q+1)^\delta,(k,q+1)\}\sigma_{k}^{-1}\cdots\sigma_{j}^{-1}\sigma_{j-1}^{-1}.
\end{eqnarray*}
\item[c)]\ $q<k,p<q<k<j$.\\
\begin{eqnarray*}
&&\sigma_{j}^{-1}\sigma_{k}^{-1}\cdots\sigma_{j-1}^{-1}(pq)^\delta\sigma_{j}^{-1}
\equiv\sigma_{j}^{-1}(pq)^\delta\sigma_{k}^{-1}\cdots\sigma_{j}^{-1}
\equiv(pq)^\delta\sigma_{k}^{-1}\cdots\sigma_{j}^{-1}\sigma_{j-1}^{-1}\\
&&\\
&\mbox{ and }& \ \ \ \ \ \
\sigma_{k}^{-1}\cdots\sigma_{j}^{-1}\sigma_{j-1}^{-1}(pq)^\delta
\equiv(pq)^\delta\sigma_{k}^{-1}\cdots\sigma_{j}^{-1}\sigma_{j-1}^{-1}.
\end{eqnarray*}
\end{enumerate}
\end{enumerate}

$(15)\wedge(2)$\\
Let
$f=\sigma_{j}^{-1}\sigma_{k}^{-1}\cdots\sigma_{j}^{-1}-\sigma_{k}^{-1}\cdots\sigma_{j}^{-1}\sigma_{j-1}^{-1},k<j,\
g=\sigma_{j}^{-1}(j,j+1)^\delta-(j,j+1)^\delta\sigma_{j}^{-1}$. Then
$w=\sigma_{j}^{-1}\sigma_{k}^{-1}\cdots\sigma_{j}^{-1}(j,j+1)^\delta$
and
\begin{eqnarray*}
(f,g)_w&=&\sigma_{j}^{-1}\sigma_{k}^{-1}\cdots\sigma_{j-1}^{-1}(j,j+1)^\delta\sigma_{j}^{-1}
-\sigma_{k}^{-1}\cdots\sigma_{j}^{-1}\sigma_{j-1}^{-1}(j,j+1)^\delta\\
&\equiv&\sigma_{j}^{-1}(k,k+1)^\delta\sigma_{k}^{-1}\cdots\sigma_{j-1}^{-1}\sigma_{j}^{-1}
-\sigma_{k}^{-1}\cdots\sigma_{j}^{-1}(j-1,j+1)^\delta\sigma_{j-1}^{-1}\\
&\equiv&(k,j)^\delta\sigma_{k}^{-1}\cdots\sigma_{j}^{-1}\sigma_{j-1}^{-1}
-(k,j)^\delta\sigma_{k}^{-1}\cdots\sigma_{j}^{-1}\sigma_{j-1}^{-1}\\
&\equiv&0.
\end{eqnarray*}

$(15)\wedge(3)$\\
Let
$f=\sigma_{j}^{-1}\sigma_{k}^{-1}\cdots\sigma_{j}^{-1}-\sigma_{k}^{-1}\cdots\sigma_{j}^{-1}\sigma_{j-1}^{-1},k<j,\
g=\sigma_{j}^{-1}(j+1,q)^\delta-(jq)^\delta\sigma_{j}^{-1},j+1<q$.
Then
$w=\sigma_{j}^{-1}\sigma_{k}^{-1}\cdots\sigma_{j}^{-1}(j+1,q)^\delta$
and
\begin{eqnarray*}
(f,g)_w&=&\sigma_{j}^{-1}\sigma_{k}^{-1}\cdots\sigma_{j-1}^{-1}(j,q)^\delta\sigma_{j}^{-1}
-\sigma_{k}^{-1}\cdots\sigma_{j}^{-1}\sigma_{j-1}^{-1}(j+1,q)^\delta\\
&\equiv&\sigma_{j}^{-1}(kq)^\delta\sigma_{k}^{-1}\cdots\sigma_{j-1}^{-1}\sigma_{j}^{-1}
-\sigma_{k}^{-1}\cdots\sigma_{j}^{-1}(j+1,q)^\delta\sigma_{j-1}^{-1}\\
&\equiv&(kq)^\delta\sigma_{k}^{-1}\cdots\sigma_{j}^{-1}\sigma_{j-1}^{-1}
-(kq)^\delta\sigma_{k}^{-1}\cdots\sigma_{j}^{-1}\sigma_{j-1}^{-1}\\
&\equiv&0.
\end{eqnarray*}

$(15)\wedge(4)$\\
Let
$f=\sigma_{j}^{-1}\sigma_{k}^{-1}\cdots\sigma_{j}^{-1}-\sigma_{k}^{-1}\cdots\sigma_{j}^{-1}\sigma_{j-1}^{-1},k<j,\
g=\sigma_{j}^{-1}(jq)^\delta-\{(j+1,q)^\delta,(j,j+1)\}\sigma_{j}^{-1},j+1<q$.
Then
$w=\sigma_{j}^{-1}\sigma_{k}^{-1}\cdots\sigma_{j}^{-1}(jq)^\delta$
and
\begin{eqnarray*}
(f,g)_w&=&\sigma_{j}^{-1}\sigma_{k}^{-1}\cdots\sigma_{j-1}^{-1}\{(j+1,q)^\delta,(j,j+1)\}\sigma_{j}^{-1}
-\sigma_{k}^{-1}\cdots\sigma_{j}^{-1}\sigma_{j-1}^{-1}(jq)^\delta\\
&\equiv&\sigma_{j}^{-1}\{(j+1,q)^\delta,(k,j+1)\}\sigma_{k}^{-1}\cdots\sigma_{j}^{-1}
-\sigma_{k}^{-1}\cdots\sigma_{j}^{-1}(j-1,q)^\delta\sigma_{j-1}^{-1}\\
&\equiv&\{(jq)^\delta,(kj)\}\sigma_{k}^{-1}\cdots\sigma_{j}^{-1}\sigma_{j-1}^{-1}
-\sigma_{k}^{-1}\cdots\sigma_{j-1}^{-1}(j-1,q)^\delta\sigma_{j}^{-1}\sigma_{j-1}^{-1}\\
&\equiv&\{(jq)^\delta,(kj)\}\sigma_{k}^{-1}\cdots\sigma_{j}^{-1}\sigma_{j-1}^{-1}
-\{(jq)^\delta,(kj)\}\sigma_{k}^{-1}\cdots\sigma_{j}^{-1}\sigma_{j-1}^{-1}\\
&\equiv&0.
\end{eqnarray*}

$(15)\wedge(5)$\\
Let
$f=\sigma_{j}^{-1}\sigma_{k}^{-1}\cdots\sigma_{j}^{-1}-\sigma_{k}^{-1}\cdots\sigma_{j}^{-1}\sigma_{j-1}^{-1},k<j,\
g=\sigma_{j}^{-1}(i,j+1)^\delta-(ij)^\delta\sigma_{j}^{-1},i<j$.
Then $w=\sigma_{j}^{-1}\sigma_{k}^{-1}\cdots\sigma_{j}^{-1}(i,j+1)^
\delta $ and
$$
(f,g)_w=\sigma_{j}^{-1}\sigma_{k}^{-1}\cdots\sigma_{j-1}^{-1}(ij)^{\delta}\sigma_{j}^{-1}
-\sigma_{k}^{-1}\cdots\sigma_{j}^{-1}\sigma_{j-1}^{-1}(i,j+1)^\delta.
$$
There are two cases to consider.
\begin{enumerate}
\item[1)]\ $k\leqslant i<j$.\\
If $i=j-1$, then
\begin{eqnarray*}
&&\sigma_{j}^{-1}\sigma_{k}^{-1}\cdots\sigma_{j-1}^{-1}(ij)^\delta\sigma_{j}^{-1}
=\sigma_{j}^{-1}\sigma_{k}^{-1}\cdots\sigma_{j-1}^{-1}(j-1,j)^\delta\sigma_{j}^{-1}\\
&\equiv&\sigma_{j}^{-1}\sigma_{k}^{-1}\cdots\sigma_{j-2}^{-1}(j-1,j)^\delta\sigma_{j-1}^{-1}\sigma_{j}^{-1}\\
&\equiv&\sigma_{j}^{-1}(kj)^\delta\sigma_{k}^{-1}\cdots\sigma_{j}^{-1}\\
&\equiv&\{(k,j+1)^\delta,(j,j+1)\}\sigma_{k}^{-1}\cdots\sigma_{j}^{-1}\sigma_{j-1}^{-1}\ \ \ \ \ \ \  and\\
&&\\
&&\sigma_{k}^{-1}\cdots\sigma_{j}^{-1}\sigma_{j-1}^{-1}(j-1,j+1)^\delta\\
&\equiv&\sigma_{k}^{-1}\cdots\sigma_{j}^{-1}\{(j,j+1)^\delta,(j-1,j)\}\sigma_{j-1}^{-1}\\
&\equiv&\sigma_{k}^{-1}\cdots\sigma_{j}^{-1}(j-1,j+1)(j,j+1)^\delta(j-1,j+1)^{-1}\sigma_{j-1}^{-1}\\
&\equiv&\sigma_{k}^{-1}\cdots\sigma_{j-1}^{-1}(j-1,j)(j,j+1)^\delta(j-1,j)^{-1}\sigma_{j}^{-1}\sigma_{j-1}^{-1}\\
&\equiv&(kj)(k,j+1)^\delta(kj)^{-1}\sigma_{k}^{-1}\cdots\sigma_{j}^{-1}\sigma_{j-1}^{-1}\\
&\equiv&\{(k,j+1)^\delta,(j,j+1)\}\sigma_{k}^{-1}\cdots\sigma_{j}^{-1}\sigma_{j-1}^{-1}.
\end{eqnarray*}
If $i<j-1$, then
\begin{eqnarray*}
&&\sigma_{j}^{-1}\sigma_{k}^{-1}\cdots\sigma_{j-1}^{-1}(ij)^\delta\sigma_{j}^{-1}\\
&\equiv&\sigma_{j}^{-1}\sigma_{k}^{-1}\cdots\sigma_{i}^{-1}(i,i+1)^\delta\sigma_{i+1}^{-1}\cdots\sigma_{j}^{-1}\\
&\equiv&\sigma_{j}^{-1}(k,i+1)^\delta\sigma_{k}^{-1}\cdots\sigma_{j}^{-1}\\
&\equiv&(k,i+1)^\delta\sigma_{k}^{-1}\cdots\sigma_{j}^{-1}\sigma_{j-1}^{-1} \ \ \ \ \ \ and \\
&&\\
&&\sigma_{k}^{-1}\cdots\sigma_{j}^{-1}\sigma_{j-1}^{-1}(i,j+1)^\delta\\
&\equiv&\sigma_{k}^{-1}\cdots\sigma_{j}^{-1}(i,j+1)^\delta\sigma_{j-1}^{-1}\\
&\equiv&\sigma_{k}^{-1}\cdots\sigma_{i}^{-1}(i,i+1)^\delta\sigma_{i+1}^{-1}\cdots\sigma_{j}^{-1}\sigma_{j-1}^{-1}\\
&\equiv&(k.i+1)^\delta\sigma_{k}^{-1}\cdots\sigma_{j}^{-1}\sigma_{j-1}^{-1}.
\end{eqnarray*}
\item[2)]\ $i<k<j$. We have
$$
\sigma_{j}^{-1}\sigma_{k}^{-1}\cdots\sigma_{j-1}^{-1}(ij)^\delta\sigma_{j}^{-1}
\equiv\sigma_{j}^{-1}(ik)^\delta\sigma_{k}^{-1}\cdots\sigma_{j}^{-1}
\equiv(ik)^\delta\sigma_{k}^{-1}\cdots\sigma_{j}^{-1}\sigma_{j-1}^{-1}
$$
and
$$
\sigma_{k}^{-1}\cdots\sigma_{j}^{-1}\sigma_{j-1}^{-1}(i,j+1)^\delta
\equiv\sigma_{k}^{-1}\cdots\sigma_{j}^{-1}(i,j+1)^\delta\sigma_{j-1}^{-1}
\equiv(ik)^\delta\sigma_{k}^{-1}\cdots\sigma_{j}^{-1}\sigma_{j-1}^{-1}.
$$
\end{enumerate}

$(15)\wedge(6)$\\
Let
$f=\sigma_{j}^{-1}\sigma_{k}^{-1}\cdots\sigma_{j}^{-1}-\sigma_{k}^{-1}\cdots\sigma_{j}^{-1}\sigma_{j-1}^{-1},k<j,\
g=\sigma_{j}^{-1}(ij)^\delta-\{(i,j+1)^\delta,(j,j+1)\}\sigma_{j}^{-1},i<j$.
Then
$w=\sigma_{j}^{-1}\sigma_{k}^{-1}\cdots\sigma_{j}^{-1}(ij)^\delta$
and
$$
(f,g)_w=\sigma_{j}^{-1}\sigma_{k}^{-1}\cdots\sigma_{j-1}^{-1}\{(i,j+1)^\delta,(j,j+1)\}\sigma_{j}^{-1}
-\sigma_{k}^{-1}\cdots\sigma_{j}^{-1}\sigma_{j-1}^{-1}(ij)^\delta.
$$
There are two cases to consider.
\begin{enumerate}
\item[1)]\  $k\leqslant i<j$.\\
If $i=j-1$, then
\begin{eqnarray*}
&&\sigma_{j}^{-1}\sigma_{k}^{-1}\cdots\sigma_{j-1}^{-1}\{(i,j+1)^\delta,(j,j+1)\}\sigma_{j}^{-1}\\
&\equiv&\sigma_{j}^{-1}\sigma_{k}^{-1}\cdots\sigma_{j-2}^{-1}\{\{(j,j+1)^\delta,(j-1,j)\},(j-1,j+1)\}
\sigma_{j-1}^{-1}\sigma_{j}^{-1}\\
&\equiv&\sigma_{j}^{-1}\sigma_{k}^{-1}\cdots\sigma_{j-2}^{-1}(j-1,j+1)^{-1}(j-1,j)^{-1}(j,j+1)^{\delta}\\
&&(j-1,j)(j-1,j+1)\sigma_{j-1}^{-1}\sigma_{j}^{-1}\\
&\equiv&\sigma_{j}^{-1}\sigma_{k}^{-1}\cdots\sigma_{j-2}^{-1}(j-1,j+1)^{-1}\{(j,j+1)^{\delta},(j-1,j+1^{-1})\}\\
&&(j-1,j+1)\sigma_{j-1}^{-1}\sigma_{j}^{-1}\\
&\equiv&\sigma_{j}^{-1}\sigma_{k}^{-1}\cdots\sigma_{j-2}^{-1}(j,j+1)^{\delta}\sigma_{j-1}^{-1}\sigma_{j}^{-1}\\
&\equiv&\sigma_{j}^{-1}(j,j+1)^{\delta}\sigma_{k}^{-1}\cdots\sigma_{j}^{-1}\\
&\equiv&(j,j+1)^{\delta}\sigma_{k}^{-1}\cdots\sigma_{j}^{-1}\sigma_{j-1}^{-1}\ \ \ \ \ \  and \\
&&\\
&&\sigma_{k}^{-1}\cdots\sigma_{j}^{-1}\sigma_{j-1}^{-1}(j-1,j)^{\delta}\\
&\equiv&\sigma_{k}^{-1}\cdots\sigma_{j}^{-1}(j-1,j)^{\delta}\sigma_{j-1}^{-1}\\
&\equiv&\sigma_{k}^{-1}\cdots\sigma_{j-1}^{-1}\{(j-1,j+1)^{\delta},(j,j+1)\}\sigma_{j}^{-1}\sigma_{j-1}^{-1}\\
&\equiv&\sigma_{k}^{-1}\cdots\sigma_{j-2}^{-1}\{\{(j,j+1)^{\delta},(j-1,j)\},(j-1,j+1)\}\sigma_{j-1}^{-1}
\sigma_{j}^{-1}\sigma_{j-1}^{-1}\\
&\equiv&\sigma_{k}^{-1}\cdots\sigma_{j-2}^{-1}(j-1,j+1)^{-1}(j-1,j)^{-1}(j,j+1)^{\delta}(j-1,j)(j-1,j+1)\\
&&\sigma_{j-1}^{-1}\sigma_{j}^{-1}\sigma_{j-1}^{-1}\\
&\equiv&(j,j+1)^{\delta}\sigma_{k}^{-1}\cdots\sigma_{j}^{-1}\sigma_{j-1}^{-1}.
\end{eqnarray*}
If $i<j-1$, then
\begin{eqnarray*}
&&\sigma_{j}^{-1}\sigma_{k}^{-1}\cdots\sigma_{j-1}^{-1}\{(i,j+1)^{\delta},(j,j+1)\}\sigma_{j}^{-1}\\
&\equiv&\sigma_{j}^{-1}\sigma_{k}^{-1}\cdots\sigma_{i}^{-1}\{(i,j+1)^{\delta},(i+1,j+1)\}
\sigma_{i+1}^{-1}\cdots\sigma_{j}^{-1}\\
&\equiv&\sigma_{j}^{-1}\sigma_{k}^{-1}\cdots\sigma_{i-1}^{-1}\{\{(i+1,j+1)^{\delta},(i,i+1)\},(i,j+1)\}
\sigma_{i}^{-1}\cdots\sigma_{j}^{-1}\\
&\equiv&\sigma_{j}^{-1}\sigma_{k}^{-1}\cdots\sigma_{i-1}^{-1}(i+1,j+1)^{\delta}\sigma_{i}^{-1}\cdots\sigma_{j}^{-1}\\
&\equiv&\sigma_{j}^{-1}(i+1,j+1)^{\delta}\sigma_{k}^{-1}\cdots\sigma_{j}^{-1}\\
&\equiv&(i+1,j)^{\delta}\sigma_{k}^{-1}\cdots\sigma_{j}^{-1}\sigma_{j-1}^{-1} \ \ \ \ \ \ and \\
&&\\
&&\sigma_{k}^{-1}\cdots\sigma_{j}^{-1}\sigma_{j-1}^{-1}(ij)^{\delta}\\
&\equiv&\sigma_{k}^{-1}\cdots\sigma_{j}^{-1}(i,j-1)^{\delta}\sigma_{j-1}^{-1}\\
&\equiv&\sigma_{k}^{-1}\cdots\sigma_{j-1}^{-1}(i,j-1)^{\delta}\sigma_{j}^{-1}\sigma_{j-1}^{-1}\\
&\equiv&\sigma_{k}^{-1}\cdots\sigma_{j-2}^{-1}\{(ij)^{\delta},(j-1,j)\}\sigma_{j-1}^{-1}\sigma_{j}^{-1}\sigma_{j-1}^{-1}.
\end{eqnarray*}
If \ \ $ i=j-2 $ ,\ then
\begin{eqnarray*}
&&\sigma_{k}^{-1}\cdots\sigma_{j-2}^{-1}\{(ij)^{\delta},(j-1,j)\}\sigma_{j-1}^{-1}\sigma_{j}^{-1}\sigma_{j-1}^{-1}\\
&\equiv&\sigma_{k}^{-1}\cdots\sigma_{j-3}^{-1}\{\{(j-1,j)^{\delta},(j-2,j-1)\},(j-2,j)\}\sigma_{j-2}^{-1}
\sigma_{j-1}^{-1}\sigma_{j}^{-1}\sigma_{j-1}^{-1}\\
&\equiv&\sigma_{k}^{-1}\cdots\sigma_{j-3}^{-1}(j-2,j)^{-1}(j-2,j-1)^{-1}(j-1,j)^{\delta}(j-2,j-1)(j-2,j)\\
&&\sigma_{j-2}^{-1}\sigma_{j-1}^{-1}\sigma_{j}^{-1}\sigma_{j-1}^{-1}\\
&\equiv&\sigma_{k}^{-1}\cdots\sigma_{j-3}^{-1}(j-1,j)^{\delta}\sigma_{j-2}^{-1}\sigma_{j-1}^{-1}
\sigma_{j}^{-1}\sigma_{j-1}^{-1}\\
&\equiv&(j-1,j)^{\delta}\sigma_{k}^{-1}\cdots\sigma_{j}^{-1}\sigma_{j-1}^{-1}\\
&\equiv&(i+1,j)^{\delta}\sigma_{k}^{-1}\cdots\sigma_{j}^{-1}\sigma_{j-1}^{-1}.
\end{eqnarray*}
If \ \ $i<j-2$ , \ then\
\begin{eqnarray*}
&&\sigma_{k}^{-1}\cdots\sigma_{j-2}^{-1}\{(ij)^{\delta},(j-1,j)\}\sigma_{j-1}^{-1}\sigma_{j}^{-1}\sigma_{j-1}^{-1}\\
&\equiv&\sigma_{k}^{-1}\cdots\sigma_{i}^{-1}\{(ij)^{\delta},(i+1,j)\}\sigma_{i+1}^{-1}\cdots
\sigma_{j}^{-1}\sigma_{j-1}^{-1}\\
&\equiv&\sigma_{k}^{-1}\cdots\sigma_{i-1}^{-1}\{\{(i+1,j)^{\delta},(i,i+1)\},(ij)\}
\sigma_{i}^{-1}\cdots\sigma_{j}^{-1}\sigma_{j-1}^{-1}\\
&\equiv&\sigma_{k}^{-1}\cdots\sigma_{i-1}^{-1}(i+1,j)^{\delta}\sigma_{i}^{-1}\cdots
\sigma_{j}^{-1}\sigma_{j-1}^{-1}\\
&\equiv&(i+1,j)^{\delta}\sigma_{k}^{-1}\cdots\sigma_{j}^{-1}\sigma_{j-1}^{-1}.
\end{eqnarray*}
\item[2)]\  \ $i<k<j$.
\begin{eqnarray*}
&&\sigma_{j}^{-1}\sigma_{k}^{-1}\cdots\sigma_{j-1}^{-1}\{(i,j+1)^{\delta},(j,j+1)\}\sigma_{j}^{-1}\\
&\equiv&\sigma_{j}^{-1}\{(i,j+1)^{\delta},(k,j+1)\}\sigma_{k}^{-1}\cdots\sigma_{j}^{-1}\\
&\equiv&\{(ij)^{\delta},(kj)\}\sigma_{k}^{-1}\cdots\sigma_{j}^{-1}\sigma_{j-1}^{-1} \ \ \ \ \ \ and \\
&&\\
&&\sigma_{k}^{-1}\cdots\sigma_{j}^{-1}\sigma_{j-1}^{-1}(ij)^{\delta}\\
&\equiv&\sigma_{k}^{-1}\cdots\sigma_{j}^{-1}(i,j-1)^{\delta}\sigma_{j-1}^{-1}\\
&\equiv&\sigma_{k}^{-1}\cdots\sigma_{j-1}^{-1}(i,j-1)^{\delta}\sigma_{j}^{-1}\sigma_{j-1}^{-1}\\
&\equiv&\sigma_{k}^{-1}\cdots\sigma_{j-2}^{-1}\{(ij)^{\delta},(j-1,j)\}\sigma_{j-1}^{-1}\sigma_{j}^{-1}
\sigma_{j-1}^{-1}\\
&\equiv&\{(ij)^{\delta},(kj)\}\sigma_{k}^{-1}\cdots\sigma_{j}^{-1}\sigma_{j-1}^{-1}.
\end{eqnarray*}
\end{enumerate}

$(15)\wedge(14)$\\
Let
$f=\sigma_{j}^{-1}\sigma_{k}^{-1}\cdots\sigma_{j}^{-1}-\sigma_{k}^{-1}\cdots\sigma_{j}^{-1}\sigma_{j-1}^{-1},k<j,\
 g=\sigma_{j}^{-1}\sigma_{q}^{-1}-\sigma_{q}^{-1}\sigma_{j}^{-1},j<q-1$.
 Then
$w=\sigma_{j}^{-1}\sigma_{k}^{-1}\cdots\sigma_{j}^{-1}\sigma_{q}^{-1}$
and
\begin{eqnarray*}
(f,g)_w&=&\sigma_{j}^{-1}\sigma_{k}^{-1}\cdots\sigma_{j-1}^{-1}\sigma_{q}^{-1}\sigma_{j}^{-1}
-\sigma_{k}^{-1}\cdots\sigma_{j}^{-1}\sigma_{j-1}^{-1}\sigma_{q}^{-1}\\
&\equiv&\sigma_{q}^{-1}\sigma_{j}^{-1}\sigma_{k}^{-1}\cdots\sigma_{j-1}^{-1}\sigma_{j}^{-1}
-\sigma_{q}^{-1}\sigma_{k}^{-1}\cdots\sigma_{j}^{-1}\sigma_{j-1}^{-1}\\
&\equiv&\sigma_{q}^{-1}\sigma_{k}^{-1}\cdots\sigma_{j-1}^{-1}\sigma_{j}^{-1}\sigma_{j-1}^{-1}
-\sigma_{q}^{-1}\sigma_{k}^{-1}\cdots\sigma_{j}^{-1}\sigma_{j-1}^{-1}\\
&\equiv&0.
\end{eqnarray*}

$(15)\wedge(15)$\\
Let
$f=\sigma_{j}^{-1}\sigma_{k}^{-1}\cdots\sigma_{j}^{-1}-\sigma_{k}^{-1}\cdots\sigma_{j}^{-1}\sigma_{j-1}^{-1},\
 \ k<j,\ g=\sigma_{j}^{-1}\sigma_{l}^{-1}\cdots\sigma_{j}^{-1}-\sigma_{l}^{-1}\cdots\sigma_{j}^{-1}\sigma_{j-1}^{-1},\
 \ l<j$. Then
$w=\sigma_{j}^{-1}\sigma_{k}^{-1}\cdots\sigma_{j}^{-1}\sigma_{l}^{-1}\cdots\sigma_{j}^{-1}$
and
$$
(f,g)_{w}=-\sigma_{k}^{-1}\cdots\sigma_{j}^{-1}\sigma_{j-1}^{-1}\sigma_{l}^{-1}\cdots\sigma_{j}^{-1}
+\sigma_{j}^{-1}\sigma_{k}^{-1}\cdots\sigma_{j-1}^{-1}\sigma_{l}^{-1}\cdots\sigma_{j}^{-1}\sigma_{j-1}^{-1}.
$$
There are two cases to consider.
\begin{enumerate}
\item[1)]\ \ $l=j-1$.
\begin{eqnarray*}
&&\sigma_{k}^{-1}\cdots\sigma_{j}^{-1}\sigma_{j-1}^{-1}\sigma_{l}^{-1}\cdots\sigma_{j}^{-1}\\
&\equiv&\sigma_{k}^{-1}\cdots\sigma_{j}^{-1}(j-1,j)^{-1}\sigma_{j}^{-1}\\
&\equiv&(j,j+1)^{-1}\sigma_{k}^{-1}\cdots\sigma_{j-1}^{-1}(j,j+1)^{-1}\\
&\equiv&(j,j+1)^{-1}(k,j+1)^{-1}\sigma_{k}^{-1}\cdots\sigma_{j-1}^{-1} \ \ \ \ \ and\\
&&\\
&&\sigma_{j}^{-1}\sigma_{k}^{-1}\cdots\sigma_{j-1}^{-1}\sigma_{l}^{-1}\cdots\sigma_{j}^{-1}\sigma_{j-1}^{-1}\\
&\equiv&\sigma_{j}^{-1}\sigma_{k}^{-1}\cdots\sigma_{j-2}^{-1}(j-1,j)^{-1}\sigma_{j}^{-1}\sigma_{j-1}^{-1}\\
&\equiv&\sigma_{j}^{-1}(kj)^{-1}\sigma_{k}^{-1}\cdots\sigma_{j-2}^{-1}\sigma_{j}^{-1}\sigma_{j-1}^{-1}\\
&\equiv&\{(k,j+1)^{-1},(j,j+1)\}(j,j+1)^{-1}\sigma_{k}^{-1}\cdots\sigma_{j-1}^{-1}\\
&\equiv&(j,j+1)^{-1}(k,j+1)^{-1}\sigma_{k}^{-1}\cdots\sigma_{j-1}^{-1}.
\end{eqnarray*}
\item[2)]\  \ $l<j-1$.
\begin{eqnarray*}
&&\sigma_{k}^{-1}\cdots\sigma_{j}^{-1}\sigma_{j-1}^{-1}\sigma_{l}^{-1}\cdots\sigma_{j}^{-1}\\
&\equiv&\sigma_{k}^{-1}\cdots\sigma_{j}^{-1}\sigma_{l}^{-1}\cdots\sigma_{j}^{-1}\sigma_{j-2}^{-1}\\
&\equiv&\sigma_{kj+1}\sigma_{lj+1}\sigma_{j-2}^{-1}\\
&\equiv& \left\{
\begin{array}{ll}\sigma_{lj+1}\sigma_{k-1j}\sigma_{j-2}^{-1} \ \ \ \ \ \ \ \ \ \ \
\ \ \ \ \ \ \ \ \ \ \ \ \ \ \ \ \ \ \ \ \ \ \  if\ \ l<k,\\
(k,l+1)^{-1}\sigma_{l+1j+1}\sigma_{kj}\sigma_{j-2}^{-1} \ \ \ \ \ \
\ \ \ \ \ \ \ \ \ \ \ \ \ \ if \ \ \ k \leq l.
\end{array}\right.
\end{eqnarray*}
If\ $l<k $,\ then
\begin{eqnarray*}
&&\sigma_{j}^{-1}\sigma_{k}^{-1}\cdots\sigma_{j-1}^{-1}\sigma_{l}^{-1}\cdots\sigma_{j}^{-1}\sigma_{j-1}^{-1} \\
&\equiv&\sigma_{j}^{-1}\sigma_{kj}\sigma_{lj}\sigma_{j}^{-1}\sigma_{j-1}^{-1}\\
&\equiv&\sigma_{j}^{-1}\sigma_{lj}\sigma_{k-1j-1}\sigma_{j}^{-1}\sigma_{j-1}^{-1}\\
&\equiv&\sigma_{j}^{-1}\sigma_{lj+1}\sigma_{k-1j}\\
&\equiv&\sigma_{lj+1}\sigma_{k-1j}\sigma_{j-2}^{-1}.
\end{eqnarray*}
If\ $k \leq l$, \ then $k \leq l<j-1$ and
\begin{eqnarray*}
&&\sigma_{j}^{-1}\sigma_{k}^{-1}\cdots\sigma_{j-1}^{-1}\sigma_{l}^{-1}\cdots\sigma_{j}^{-1}\sigma_{j-1}^{-1}\\
&\equiv&\sigma_{j}^{-1}\sigma_{kj}\sigma_{lj}\sigma_{j}^{-1}\sigma_{j-1}^{-1}\\
&\equiv&\sigma_{j}^{-1}(k,l+1)^{-1}\sigma_{l+1j}\sigma_{kj-1}\sigma_{j}^{-1}\sigma_{j-1}^{-1}\\
&\equiv&\sigma_{j}^{-1}(k,l+1)^{-1}\sigma_{l+1j+1}\sigma_{kj}\\
&\equiv&(k,l+1)^{-1}\sigma_{l+1j+1}\sigma_{kj}\sigma_{j-2}^{-1}.
\end{eqnarray*}
\end{enumerate}

$(15)\wedge(16)$ \\
Let
$f=\sigma_{j}^{-1}\sigma_{k}^{-1}\cdots\sigma_{j}^{-1}-\sigma_{k}^{-1}\cdots\sigma_{j}^{-1}\sigma_{j-1}^{-1},k<j,\
g=\sigma_{j}^{-1}\sigma_{j}^{-1}-(j,j+1)^{-1}$. Then
$w=\sigma_{j}^{-1}\sigma_{k}^{-1}\cdots\sigma_{j}^{-1}\sigma_{j}^{-1}$
and
\begin{eqnarray*}
(f,g)_w&=&\sigma_{j}^{-1}\sigma_{k}^{-1}\cdots\sigma_{j-1}^{-1}(j,j+1)^{-1}
-\sigma_{k}^{-1}\cdots\sigma_{j}^{-1}\sigma_{j-1}^{-1}\sigma_{j}^{-1}\\
&\equiv&\sigma_{j}^{-1}(k,j+1)^{-1}\sigma_{k}^{-1}\cdots\sigma_{j-1}^{-1}
-\sigma_{k}^{-1}\cdots\sigma_{j-1}^{-1}\sigma_{j-1}^{-1}\sigma_{j}^{-1}\sigma_{j-1}^{-1}.
\end{eqnarray*}
If  $\ \ k<j-1$, then
\begin{eqnarray*}
(f,g)_w&\equiv&(kj)^{-1}\sigma_{j}^{-1}\sigma_{k}^{-1}\cdots\sigma_{j-1}^{-1}
-\sigma_{k}^{-1}\cdots\sigma_{j-2}^{-1}(j-1,j)^{-1}\sigma_{j}^{-1}\sigma_{j-1}^{-1}\\
&\equiv&(kj)^{-1}\sigma_{k}^{-1}\cdots\sigma_{j-2}^{-1}\sigma_{j}^{-1}\sigma_{j-1}^{-1}
-(kj)^{-1}\sigma_{k}^{-1}\cdots\sigma_{j-2}^{-1}\sigma_{j}^{-1}\sigma_{j-1}^{-1}\\
&\equiv&0.
\end{eqnarray*}
If $\ \ k=j-1$, then
\begin{eqnarray*}
(f,g)_w&\equiv&(j-1,j)^{-1}\sigma_{j}^{-1}\sigma_{j-1}^{-1}-\sigma_{j-1}^{-1}\sigma_{j-1}^{-1}\sigma_{j}^{-1}
\sigma_{j-1}^{-1}\\
&\equiv&(j-1,j)^{-1}\sigma_{j}^{-1}\sigma_{j-1}^{-1}-(j-1,j)^{-1}\sigma_{j}^{-1}\sigma_{j-1}^{-1}\\
&\equiv&0.
\end{eqnarray*}

$(16)\wedge(1)$\\
Let $f=\sigma_{k}^{-1}\sigma_{k}^{-1}-(k,k+1)^{-1},\ \
g=\sigma_{k}^{-1}(ij)^{\delta}-(ij)^{\delta}\sigma_{k}^{-1},\ \
k\neq i,i-1,j,j-1$. Then
$w=\sigma_{k}^{-1}\sigma_{k}^{-1}(ij)^{\delta}$ and
$$
(f,g)_w=\sigma_{k}^{-1}(ij)^{\delta}\sigma_{k}^{-1}-(k,k+1)^{-1}(ij)^{\delta}.
$$
There are three cases to consider.
\begin{enumerate}
\item[1)]\ If $\ \ \ \ k\geq j+1,\ \ i<j<j+1\leq k<k+1$, then
$$
(f,g)_w\equiv(ij)^{\delta}(k,k+1)^{-1}-(k,k+1)^{-1}(ij)^{\delta}\equiv(k,k+1)^{-1}(ij)^{\delta}-
(k,k+1)^{-1}(ij)^{\delta}\equiv0.
$$
\item[2)] \ \ \ If $ \ i<k<j-1,\ \ i<k<k+1<j$,  then
$$
(f,g)_w\equiv(ij)^{\delta}(k,k+1)^{-1}-(ij)^{\delta}(k,k+1)^{-1}\equiv0.
$$
\item[3)] \  If $ k<i-1, \ k<k+1<i<j$, then
$$
(f,g)_w\equiv(ij)^{\delta}(k,k+1)^{-1}-(ij)^{\delta}(k,k+1)^{-1}\equiv0.
$$
\end{enumerate}

$(16)\wedge(2)$\\
Let $f=\sigma_{i}^{-1}\sigma_{i}^{-1}-(i,i+1)^{-1},\ \
 g=\sigma_{i}^{-1}(i,i+1)^{\delta}-(i,i+1)^{\delta}\sigma_{i}^{-1}$. Then
$w=\sigma_{i}^{-1}\sigma_{i}^{-1}(i,i+1)^{\delta}$ and
\begin{eqnarray*}
(f,g)_w&=&\sigma_{i}^{-1}(i,i+1)^{\delta}\sigma_{i}^{-1}-(i,i+1)^{-1}(i,i+1)^{\delta}\\
&\equiv&(i,i+1)^{\delta}\sigma_{i}^{-1}\sigma_{i}^{-1}-(i,i+1)^{-1}(i,i+1)^{\delta}\\
&\equiv&(i,i+1)^{\delta}(i,i+1)^{-1}-(i,i+1)^{-1}(i,i+1)^{\delta}\\
&\equiv&0.
\end{eqnarray*}

$(16)\wedge(3)$\\
Let $f=\sigma_{i-1}^{-1}\sigma_{i-1}^{-1}-(i-1,i)^{-1},\ \
 g=\sigma_{i-1}^{-1}(ij)^{\delta}-(i-1,j)^{\delta}\sigma_{i-1}^{-1},\ \ i<j$. Then
 $w=\sigma_{i-1}^{-1}\sigma_{i-1}^{-1}(ij)^{\delta}$ and
\begin{eqnarray*}
(f,g)_w&=&\sigma_{i-1}^{-1}(i-1,j)^{\delta}\sigma_{i-1}^{-1}-(i-1,i)^{-1}(ij)^{\delta}\\
&\equiv&\{(ij)^{\delta},(i-1,i)\}(i-1,i)^{-1}-(i-1,i)^{-1}(ij)^{\delta}\\
&\equiv&(i-1,i)^{-1}(ij)^{\delta}-(i-1,i)^{-1}(ij)^{\delta}\\
&\equiv&0.
\end{eqnarray*}

$(16)\wedge(4)$\\
Let $f=\sigma_{i}^{-1}\sigma_{i}^{-1}-(i,i+1)^{-1},\ \
g=\sigma_{i}^{-1}(ij)^{\delta}-\{(i+1,j)^{\delta},(i,i+1)\}\sigma_{i}^{-1},\
\ i+1<j$. Then $w=\sigma_{i}^{-1}\sigma_{i}^{-1}(ij)^{\delta}$ and
\begin{eqnarray*}
(f,g)_w&=&\sigma_{i}^{-1}\{(i+1,j)^{\delta},(i,i+1)\}\sigma_{i}^{-1}-(i,i+1)^{-1}(ij)^{\delta}\\
&\equiv&\{(ij)^{\delta},(i,i+1)\}(i,i+1)^{-1}-(i,i+1)^{-1}(ij)^{\delta}\\
&\equiv&(i,i+1)^{-1}(ij)^{\delta}-(i,i+1)^{-1}(ij)^{\delta}\\
&\equiv&0.
\end{eqnarray*}

$(16)\wedge(5)$\\
Let $f=\sigma_{j-1}^{-1}\sigma_{j-1}^{-1}-(j-1,j)^{-1},\
 \ g=\sigma_{j-1}^{-1}(ij)^{\delta}-(i,j-1)^{\delta}\sigma_{j-1}^{-1},\ \ i<j-1$. Then
$w=\sigma_{j-1}^{-1}\sigma_{j-1}^{-1}(ij)^{\delta}$ and
\begin{eqnarray*}
(f,g)_w&=&\sigma_{j-1}^{-1}(i,j-1)^{\delta}\sigma_{j-1}^{-1}-(j-1,j)^{-1}(ij)^{\delta}\\
&\equiv&\{(ij)^{\delta},(j-1,j)\}(j-1,j)^{-1}-(j-1,j)^{-1}(ij)^{\delta}\\
&\equiv&0.
\end{eqnarray*}

$(16)\wedge(6)$\\
Let $f=\sigma_{j}^{-1}\sigma_{j}^{-1}-(j,j+1)^{-1},\ \
 g=\sigma_{j}^{-1}(ij)^{\delta}-\{(i,j+1)^{\delta},(j,j+1)\}\sigma_{j}^{-1},\ \ i<j$. Then
$w=\sigma_{j}^{-1}\sigma_{j}^{-1}(ij)^{\delta}$ and
\begin{eqnarray*}
(f,g)_w&=&\sigma_{j}^{-1}\{(i,j+1)^{\delta},(j,j+1)\}\sigma_{j}^{-1}-(j,j+1)^{-1}(ij)^{\delta}\\
&\equiv&\{(ij)^{\delta},(j,j+1)\}(j,j+1)^{-1}-(j,j+1)^{-1}(ij)^{\delta}\\
&\equiv&0.
\end{eqnarray*}

$(16)\wedge(14)$\\
Let $f=\sigma_{j}^{-1}\sigma_{j}^{-1}-(j,j+1)^{-1},\
g=\sigma_{j}^{-1}\sigma_{k}^{-1}-\sigma_{k}^{-1}\sigma_{j}^{-1},j<k-1$.
Then $w=\sigma_{j}^{-1}\sigma_{j}^{-1}\sigma_{k}^{-1}$ and
\begin{eqnarray*}
(f,g)_w&=&\sigma_{j}^{-1}\sigma_{k}^{-1}\sigma_{j}^{-1}-(j,j+1)^{-1}\sigma_{k}^{-1}\\
&\equiv&\sigma_{k}^{-1}(j,j+1)^{-1}-(j,j+1)^{-1}\sigma_{k}^{-1}\\
&\equiv&(j,j+1)^{-1}\sigma_{k}^{-1}-(j,j+1)^{-1}\sigma_{k}^{-1}\\
&\equiv&0.
\end{eqnarray*}

$(16)\wedge(15)$\\
Let $f=\sigma_{j}^{-1}\sigma_{j}^{-1}-(j,j+1)^{-1},\
g=\sigma_{j}^{-1}\sigma_{k}^{-1}\cdots\sigma_{j}^{-1}-\sigma_{k}^{-1}\cdots\sigma_{j}^{-1}\sigma_{j-1}^{-1},k<j$.
Then
$w=\sigma_{j}^{-1}\sigma_{j}^{-1}\sigma_{k}^{-1}\cdots\sigma_{j}^{-1}$
and
\begin{eqnarray*}
(f,g)_w&=&\sigma_{j}^{-1}\sigma_{k}^{-1}\cdots\sigma_{j}^{-1}\sigma_{j-1}^{-1}
-(j,j+1)^{-1}\sigma_{k}^{-1}\cdots\sigma_{j}^{-1}\\
&\equiv&\sigma_{k}^{-1}\cdots\sigma_{j}^{-1}\sigma_{j-1}^{-1}\sigma_{j-1}^{-1}
-(j,j+1)^{-1}\sigma_{k}^{-1}\cdots\sigma_{j}^{-1}\\
&\equiv&\sigma_{k}^{-1}\cdots\sigma_{j}^{-1}(j-1,j)^{-1}-(j,j+1)^{-1}\sigma_{k}^{-1}\cdots\sigma_{j}^{-1}\\
&\equiv&\sigma_{k}^{-1}\cdots\sigma_{j-1}^{-1}\{(j-1,j+1)^{-1},(j,j+1)\}\sigma_{j}^{-1}
-(j,j+1)^{-1}\sigma_{k}^{-1}\cdots\sigma_{j}^{-1}\\
&\equiv&\sigma_{k}^{-1}\cdots\sigma_{j-2}^{-1}\{\{(j,j+1)^{-1},(j-1,j)\},(j-1,j+1)\}\sigma_{j-1}^{-1}\sigma_{j}^{-1}\\
&&-(j,j+1)^{-1}\sigma_{k}^{-1}\cdots\sigma_{j}^{-1}\\
&\equiv&\sigma_{k}^{-1}\cdots\sigma_{j-2}^{-1}(j-1,j+1)^{-1}(j-1,j)^{-1}(j,j+1)^{-1}\\
&&(j-1,j)(j-1,j+1)\sigma_{j-1}^{-1}\sigma_{j}^{-1}-(j,j+1)^{-1}\sigma_{k}^{-1}\cdots\sigma_{j}^{-1}\\
&\equiv&\sigma_{k}^{-1}\cdots\sigma_{j-2}^{-1}(j,j+1)^{-1}\sigma_{j-1}^{-1}\sigma_{j}^{-1}-
(j,j+1)^{-1}\sigma_{k}^{-1}\cdots\sigma_{j}^{-1}\\
&\equiv&(j,j+1)^{-1}\sigma_{k}^{-1}\cdots\sigma_{j}^{-1}-(j,j+1)^{-1}\sigma_{k}^{-1}\cdots\sigma_{j}^{-1}\\
&\equiv&0.
\end{eqnarray*}

$(16)\wedge(16)$\\
Let $f=\sigma_{i}^{-1}\sigma_{i}^{-1}-(i,i+1)^{-1},\ \ \
g=\sigma_{i}^{-1}\sigma_{i}^{-1}-(i,i+1)^{-1}$. Then
$w=\sigma_{i}^{-1}\sigma_{i}^{-1}\sigma_{i}^{-1}$ and
$$
(f,g)_w=\sigma_{i}^{-1}(i,i+1)^{-1}-(i,i+1)^{-1}\sigma_{i}^{-1}\equiv(i,i+1)^{-1}\sigma_{i}
-(i,i+1)^{-1}\sigma_{i}\equiv0.
$$
For the cases of $(17)\wedge((1)\sim(17))$ \ and \
$((1)\sim(17))\wedge(17)$, the possible compositions are the cases
of $(17)\wedge((7)\sim(13),(17))$ \ and \
$((7)\sim(13),(17))\wedge(17)$. For these cases, we can easily check
that the compositions are trivial.

This completes our proof.

\ \

\noindent{\bf Acknowledgement}: The authors would like to express
their deepest gratitude to Professor L. A. Bokut for his kind
guidance, useful discussions and enthusiastic encouragement.

\end{document}